\documentclass[11p,english]{smfbook} 
\usepackage[french,english]{babel} 
\usepackage{smfthm} 

\usepackage{graphicx}
\newtheorem{Thm}{Theorem}
\newtheorem{Def}{Definition}[section]
\newtheorem{Rem}[Def]{Remark}
\newtheorem{Prop}[Def]{Proposition}
\newtheorem{Cor}[Def]{Corollary}
\newtheorem{Lem}[Def]{Lemma}
\renewcommand{\theequation}{\thesection.\arabic{equation}}
\def\longformule#1#2{
\displaylines{
\qquad{#1}
\hfill\cr
\hfill {#2}
\qquad\cr
}
}

\def \N{{\mathbf N}}
\def \R{{\mathbf R}}
\def \T{{\mathbf T}}
\def \Z{{\mathbf Z}}

\def \DD{{\mathcal D}}

\def \Ker{\hbox{{\rm Ker}}}

\def \eps{{\varepsilon}}
\def \e{{\varepsilon}}

\def \s{{\sigma}}

\def \AA{{\mathcal A}}

\def \EE{{\mathcal E}}
\def \FF{{\mathcal F}}

\def \LL{{\mathcal L}}

\def \d{{\partial}}
\def \iint{\int \! \! \! \int }

\def \DIV{\nabla\! \! \cdot }

\def\eqdefa{\buildrel\hbox{\footnotesize def}\over =}
\newcommand{\Fr}{ \hbox{Fr}}

\newcommand{\ul}{ \underline u}
\newcommand{\el}{ \underline \eta}

\def\eqdefa{\buildrel\hbox{\footnotesize def}\over =}

\def\virgp{\raise 2pt\hbox{,}}
\def\sumetage#1#2{
\sum_{\scriptstyle {#1}\atop\scriptstyle {#2}}}
\def\supetage#1#2{
\sup_{\scriptstyle {#1}\atop\scriptstyle {#2}}}

\author[I. Gallagher] {Isabelle Gallagher} 
\address[I. Gallagher]%
{Institut de Math{\'e}matiques de Jussieu \\
     Universit{\'e} Paris VII\\
175, rue du Chevaleret\\ 75013 Paris
    \\
     FRANCE}
\email{Isabelle.Gallagher@math.jussieu.fr}
 \author[L. Saint-Raymond]{Laure Saint-Raymond}
\address[L. Saint-Raymond]%
{ Laboratoire J.-L. Lions UMR 7598\\ Universit{\'e} Paris
VI\\
175, rue du Chevaleret\\ 75013 Paris\\FRANCE }
\email{saintray@ann.jussieu.fr }
\title[Mathematical study of the betaplane model]{Mathematical study of the betaplane model: 
Equatorial waves and convergence results}

\begin{document}
\frontmatter
 \begin{abstract}
We are interested in a model of rotating fluids, describing the motion of the ocean in the equatorial zone. This model
is known as the Saint-Venant, or 
shallow-water type system, to which a rotation term is added  whose
amplitude is linear with respect to the latitude; in particular it vanishes at the equator. After a physical introduction to 
the model, we describe the various waves involved and study in detail the resonances associated with those waves. We
then exhibit the formal limit system (as the rotation becomes large), obtained as usual by filtering out the waves, and prove
 its wellposedness. Finally we
prove three types of convergence results: a weak convergence result towards a linear, geostrophic equation,
a  strong convergence result of the filtered solutions
towards the unique strong solution to the limit system, and finally
 a ``hybrid" strong convergence result
of the filtered solutions towards a weak solution to the limit system. In particular we obtain that 
  there are no
confined equatorial waves  in the mean motion as the
rotation becomes large.  
\end{abstract}
 
 \begin{altabstract}
 On s'int{\'e}resse {\`a} un mod{\`e}le de fluides en rotation rapide, d{\'e}crivant le mouvement de l'oc{\'e}an dans la
 zone {\'e}quatoriale. Ce mod{\`e}le est connu sous le nom de Saint-Venant, ou syst{\`e}me ``shallow water'', auquel on 
 ajoute un terme de rotation dont l'amplitude est lin{\'e}aire en la latitude ; en particulier il s'annule {\`a} l'{\'e}quateur. 
 Apr{\`e}s une introduction physique au mod{\`e}le, on d{\'e}crit les diff{\'e}rentes ondes en jeu et l'on {\'e}tudie en d{\'e}tail
 les r{\'e}sonances associ{\'e}es {\`a} ces ondes. On exhibe ensuite un syst{\`e}me limite formel (dans la limite d'une forte rotation), 
 obtenu comme d'habitude en filtrant les ondes, et l'on d{\'e}montre qu'il est bien pos{\'e}. Enfin on d{\'e}montre trois
 types de r{\'e}sultats de convergence : un th{\'e}or{\`e}me de convergence faible vers un syst{\`e}me g{\'e}ostrophique
 lin{\'e}aire,  un th{\'e}or{\`e}me de convergence forte des solutions filtr{\'e}es
 vers la solution unique du syst{\`e}me limite, et enfin  un r{\'e}sultat ``hybride'' de convergence forte des solutions filtr{\'e}es 
 vers une solution faible du syst{\`e}me limite. En particulier
 on d{\'e}montre l'absence d'ondes {\'e}quatoriales confin{\'e}es dans le mouvement moyen, quand la rotation
 augmente.  
 \end{altabstract}
 
 \thanks{The authors wish to thank  warmly J.-Y. Chemin and   F. Golse  for their patience in answering our
  many questions on a number of subjects addressed in this work. They also thank the anonymous referee for pointing out
  an important improvement in the study of resonances (namely the absence of resonances in Rossby modes). }
\maketitle

\tableofcontents

\mainmatter

\chapter{Introduction}\label{intro}
\setcounter{equation}{0}

The aim of this paper is to obtain a description of geophysical flows, especially oceanic flows, in the equatorial zone.
 For the scales considered, i.e., on domains extending over many thousands of kilometers, the forces with dominating
  influence are the gravity and the Coriolis force. The question is therefore to understand how they counterbalance
   eachother to impose the so-called geostrophic constraint on the mean motion, and to describe the  oscillations 
   which are generated around this geostrophic equilibrium.

At mid-latitudes, on ``small" geographical zones,  the variations of the Coriolis force due to the curvature of the
 Earth are usually neglected, which leads to a  singular perturbation problem    with constant coefficients. 
 The corresponding asymptotics, called asymptotics of rotating fluids, have been studied by a number of authors.
  We refer for instance to the review  by R. Temam and M. Ziane \cite{TZ}, or to the work  by J.-Y. Chemin, B. Desjardins,
   I. Gallagher and E. Grenier \cite{cdggbook}.

In order to a get a more realistic description, which allows for instance to exhibit the specificity of the equatorial 
zone, one has to study more intricate models, taking into account especially the interaction between the fluid and the 
atmosphere (free surface), and the geometry of the Earth (variations of the local vertical component of the Earth 
rotation). The mathematical modelling of these various phenomena, as well as their respective importance according
 to the scales considered, have been studied in a rather systematic way by A. Majda \cite{majda}, and R. Klein and
  A. Majda \cite{klein/majda}.

\bigskip
Here we will focus on {\it quasigeostrophic, oceanic} flows, meaning that we will consider horizontal length scales
 of order 1000km and vertical length scales of order 5km, so that the aspect ratio is very small and the shallow-water
  approximation is relevant (see for instance the works by D. Bresch, B. Desjardins  and
   C.K. Lin \cite{bresch-desjardinslin} or by J.-F. Gerbeau and B. Perthame \cite{GP}). In this framework, 
   the asymptotics of homogeneous rotating fluids have been studied by D. Bresch and 
   B. Desjardins \cite{breschdesjardins}.

For the description of {\it equatorial} flows, one has to take further into account the variations of the Coriolis
 force, and especially the fact that it cancels at equator. The inhomogeneity of the Coriolis force has already
  been studied by B. Desjardins and E. Grenier \cite{desjardinsgrenier}  and by the authors \cite{gallagher/sr}
   for an incompressible fluid with  rigid lid upper boundary (see also~\cite{dutrifoymajda} for a study
   of the wellposedness and weak asymptotics of a non viscous model). The question here is then to understand the combination
    of the effects due to the free surface, and of the effects due to the variations of the Coriolis force.

Note that, for the sake of simplicity, we will not discuss the effects of the interaction with the boundaries, 
describing neither the vertical boundary layers, known as Ekman layers (see for instance the paper by D. 
G\'erard-Varet  \cite{dgv}), nor the lateral boundary layers, known as Munk and Stommel layers 
(see for instance \cite{desjardinsgrenier}). We will indeed consider a purely horizontal model, assuming 
periodicity with respect to the longitude (omitting the stopping conditions on the continents) and and infinite 
domain for the latitude (using the exponential decay of the equatorial waves to neglect the boundary).

\section{ Physical phenomena observed in the equatorial zone of the earth}

The rotation of the earth has a dominating influence on the way the atmosphere and the 
ocean
respond to imposed changes. The dynamic effect is caused
(see~\cite{gill}, \cite{pedlosky}) by the Coriolis acceleration, which is
equal to the product of the Coriolis parameter $f$ and the horizontal
velocity.

An important feature of the response of a rotating fluid to gravity is
that it does not adjust to a state of rest, but rather to an equilibrium
which  contains more potential energy than does the rest state. The
steady equilibrium solution is a geostrophic balance, i.e. a balance
between the Coriolis acceleration and the pressure gradient divided by
density. The equation determining this steady solution contains a length
scale $a$, called the Rossby radius of deformation, which is equal to
$c/|f|$ where $c$ is the wave speed in the absence of rotation effects. 
If $f$ tends to zero, then~$a$ tends to infinity, indicating that for length
scales small compared with $a$, rotation effects are small, whereas for
scales comparable to or larger than~$a$, rotation effects are
important. 
Added to that mean, geostrophic motion, are time oscillations which correspond to the 
so-called   ageostrophic motion.
The use of a constant-$f$ approximation to describe motion on
the earth is adequate to handle the adjustment process at mid-latitudes~:
Kelvin \cite{thomson} stated that his wave solutions (also known as  Poincar{\'e} waves) 
are
applicable {\sl ``in any narrow lake or portion of the sea covering not
more than a few degrees of the earth's surface, if for $\frac12 f$ we
take the component of the earth's angular velocity round a vertical
through the locality, that is to say
$$ \frac12 f =\Omega \sin \phi,$$
where $\Omega$ denotes the earth's angular velocity and $\phi$ the
latitude."}

The
adjustment processes are somewhat special when the Coriolis acceleration vanishes~:
the equatorial zone is actually found to be a waveguide: as explained in~\cite{gill}, there is an 
equatorial
Kelvin wave, and there are equatorially trapped waves, which are the equivalent of the
Poincar{\'e} waves in a uniformly rotating system. There is also  an important new class 
of waves
with much slower frequencies, called planetary or quasi-geostrophic waves. These owe their existence to the
variations in the undisturbed potential vorticity and thus exist at all latitudes. 
However, the
ray paths along which they propagate bend, as do the paths of gravity waves, because of 
the
variation of Coriolis parameter with latitude, and it is this bending that tends to 
confine the
waves to the equatorial waveguide.

\section{A mathematical model for the ocean in the equatorial zone}

In order to explore the qualitative features of the equatorial flow, we
 restrict our attention here to a very simplified model of
oceanography. More precisely,  we consider the ocean as an
incompressible viscous fluid with free surface submitted to gravitation,
and further make the following classical assumptions~:
$$\hbox{ the density of the fluid is homogeneous,}\leqno(H1)$$
$$\hbox{ the pressure law is given by the hydrostatic
approximation,}\leqno (H2)$$
$$\hbox{ the motion is essentially horizontal and does not
depend on the vertical coordinate,}\leqno(H3)$$ leading to the so-called
shallow water approximation.

We therefore consider a so-called viscous Saint-Venant model, which describes
vertically averaged flows in three dimensional shallow domains in terms
of the horizontal mean velocity field $u$ and the depth variation $h$ due
to the free surface. Taking into account the Coriolis force, a particular
model reads as\label{phiperp}
\begin{equation}
\label{SW}
\begin{aligned}
\d_t h +\DIV (hu) =0\\
\d_t (hu) +\DIV (hu \otimes u) +  {f}(hu)^\perp + {1\over \Fr^2} h\nabla h -
h\nabla K(h) -A(h, u)=0
\end{aligned}
\end{equation}
where $f$ denotes the vertical component of the earth rotation, $\Fr
$ the Froude number,
 and
$K$ and $A$ are the  capillarity and viscosity operators. We have
written~$ u^\perp $ for the vector~$ (u_2 , -u_1)$.

Note that, from a theoretical  point of view,  it is not clear that the use of
the shallow water approximation is relevant  in this context since the
Coriolis force is known to generate vertical oscillations which are
completely neglected in such an approach. Nevertheless, this  very simplified model is 
commonly
used by physicists \cite{gill, philander} and we will see that its study  already  gives  a
description of the horizontal motion corresponding to experimental observations.

\bigskip
Of course, in order that the curvature of the earth can be neglected, and
that latitude and longitude can be considered as cartesian coordinates,
we should consider only  a thin strip around the equator. This means that
we should study (\ref{SW}) on a bounded domain, and  supplement it
with boundary conditions. Nevertheless, as we expect the Coriolis force to
confine equatorial waves, we will perform our study on $\R\times \T$
where~$ \T$ is the one-dimensional torus~$ \R_{/2\pi \Z}$,  and check
a posteriori that oscillating modes vanish far from the equator, so that
it is reasonable to conjecture that they should not be disturbed by
boundary conditions.

\section{Some orders of magnitude  in the equatorial zone}

For motions near the equator, the approximations
$$\sin \phi \sim \phi,\quad \cos \phi \sim 1$$
may be used, giving what is called the equatorial betaplane approximation. Half of the 
earth's
surface lies at latitudes of less than $30^{o}$ and the maximum percentage error in the 
above
approximation in that range of latitudes is only 14$\%$. In this
approximation, $f$ is given by
$$f=\beta x_1,$$
where $x_1$ is distance northward from the equator, taking values in the range
$$x_1 \in [ -3000  \, {\rm km} \: , \: 3000  \, {\rm km} ] ,$$
and~$ \beta$
 is a
constant given by
$$\beta ={2\Omega\over r} =2.3 \times 10^{-11} \,  {\rm m}^{-1} \, {\rm s}^{-1}.$$

A formal analysis of the linearized versions of the equations shows then that rotation 
effects do
not allow the motion in each plane $x_1=const$ to be independent because a geostrophic 
balance
between the eastward velocity and the north-south pressure gradient is required. 
Equatorial
waves  actually decay in a distance of order $a_e$, the so-called equatorial radius of
deformation,
$$a_e =\left({c\over 2\beta}\right)^{1/2}$$
where $c$ is the square root of $gH$, $H$ being interpreted as the equivalent depth. For
baroclinic ocean waves, appropriate values of $c$ are typically in the
range $0.5  {\rm m} {\rm s}^{-1}\hbox{ to
}3  {\rm m} {\rm s}^{-1}$, so the order of the equatorial Rossby radius is
$$a_e \sim 100 \, {\rm km},$$
which is effectively very small compared with the range of validity of the betaplane
approximation.

\section{The Cauchy problem for the betaplane model}
\setcounter{equation}{0}
\label{cauchy}
Before describing the equatorial waves and the asymptotic behaviour of the ocean in the 
fast rotation
limit, we need to give the mathematical framework for our study, and therefore to specify 
the
dissipative operators~$A$ and~$K$ occuring in (\ref{SW}).

From a physical point of view, it would be relevant to model the viscous effects
by the following operator
$$A(h, u)=\nu \nabla \cdot (h\nabla  u),$$
meaning in particular that the viscosity cancels when $h$ vanishes. Then, in order for 
the Cauchy
problem to be globally well-posed, it is necessary to get some control on the cavitation. 
Results by Bresch and
Desjardins \cite{breschdesjardins} show  that capillary or
friction effects  can
 prevent the formation of
singularities in the Saint-Venant system (without Coriolis force). On the other hand, in 
the absence of such
dissipative effects, Mellet and Vasseur \cite{mellet-vasseur} have  proved the weak 
stability of this same system
under a suitable integrability assumption on the initial velocity field. All these 
results are based on a new
entropy inequality  \cite{breschdesjardins} which controls in particular the first 
derivative of
$\sqrt{h}$. In particular, they cannot be easily extended to (\ref{SW}) since the 
betaplane approximation of the
Coriolis force prevents from deriving  such an entropy inequality.

For the sake of simplicity, as we are interested in some asymptotic regime where the 
depth $h$ is just a
fluctuation around a mean value~$H$, we will thus consider the following viscosity operator
$$A(h, u)=\nu \Delta u,$$
and we will neglect the capillarity
$$K(h)=0,$$
so that the usual theory of the  isentropic Navier-Stokes equations
can be applied (see for instance~\cite{pll}).

\begin{Thm}[Existence of weak solutions]
\label{existence} { 
Let $(h^0,u^0)$  be some  measurable nonnegative function and
vector-field on $\R\times \T$ such that
\begin{equation}
\label{init-energy2}
\EE ^0\eqdefa\int \left(\frac{
(h^0-H)^2}{2\Fr^2} +\frac{h^0}2 |u^0|^2\right)dx
<+\infty.
\end{equation}
Then there exists a
 global weak solution to (\ref{SW}) satisfying the initial
condition
$$ h_{|t=0} =h^0,\quad u_{|t=0}=u^0,$$
and which furthermore satisfies for almost every $t\geq0$ the energy estimate
\begin{equation}
\label{energy} \int \left(\frac{
(h-H)^2}{2\Fr^2} +\frac{h}2
|u|^2\right)(t,x)dx+ \nu \int_0^t
\int    |\nabla u|^2 (t',x)dxdt'\leq
\EE^0.
\end{equation}
}
\end{Thm}

\bigskip
In this paper we are interested in describing the behaviour of the ocean in the 
equatorial zone. We
thus expect the Froude number $\Fr$, which is the ratio of the fluid speed to a measure of
the internal wave speed, to be small. More precisely we will consider depth variations
$$h=H(1+\eps  \eta)$$
where
$\eps$ stands for the order of magnitude of the Froude number.

As seen in the introduction, in order for gravity waves to be notably modified by rotation
effects, the Rossby radius of deformation has to be comparable to the typical length
scales. In order to derive the quasi-geostrophic equations with free-surface term
used in oceanography, we will assume that $\eps$ is also  the order of magnitude of
the Rossby number.

In non-dimensional variables, the viscous Saint-Venant system (\ref{SW}) can therefore be
rewritten (normalizing~$ H$ to~$ H = 1$ for simplicity)
\begin{eqnarray}
\label{SW-eps}
\d_t \eta +\frac1\eps \DIV \Bigl((1+\eps \eta) u\Bigr) =0, \nonumber\\
\d_t \Bigl((1+\eps \eta)u\Bigr) + \nabla \cdot \Bigl((1+\eps \eta)u\otimes u\Bigr) +  
\frac{ \beta x_1}{ \eps} (1+\eps \eta)u^\perp+
{1\over \eps}(1+\eps \eta)
\nabla
\eta
-\nu \Delta u =0, \nonumber \\
\eta_{|t=0}=\eta^0,\quad u_{|t=0} =u^0.
\end{eqnarray}

In
such a framework, the energy  inequality~(\ref{energy})
provides uniform bounds on any family~$( \eta_\eps,u_\eps)_{\eps>0}$ of weak solutions of
(\ref{SW-eps}).

In all the sequel we will denote respectively $\dot
H^{s } $ and $
H^{s } $  the homogeneous and inhomogeneous Sobolev spaces\label{Hs-not} of
order~$ s$, defined by
$$
\longformule{ 
\dot H^s(\R\times \T) = \Bigl\{
f \in {\mathcal S}' (\R\times \T)\: \Big/ \: {\mathcal F}  f \in L^{1}_{loc} (\R\times \T) } {\hbox{ and }  \|f\|_{\dot H^s(\R\times \T)}^2 =
\sum_{k \in \Z} \int_{\R} |\xi^2 + k^2|^s | {\mathcal F}  f (\xi,k)|^2 \: d\xi
< \infty
\Bigr\},}
$$
and
$$
H^s (\R\times \T)= \left\{
f \in {\mathcal S}'(\R\times \T) \: \Big/ \:
 \|f\|_{H^s(\R\times \T)}^2 =
\sum_{k \in \Z} \int_{\R} |1 + \xi^2 + k^2|^s| {\mathcal F} f (\xi,k)|^2 \: d\xi
< \infty
\right\},
$$
where~$ {\mathcal F}$ denotes the Fourier transform\label{ff-not}
$$
\forall k \in \Z, \: \forall \xi \in \R \quad  {\mathcal F} f (\xi,k)
= \int e^{-ik x_2}  e^{-i\xi x_1} f(x_1,x_2) \: dx_2 dx_1.
$$
We will also denote, for all  subsets~$\Omega$ of~$\R\times \T$ and\label{hs0} for all~$s > 0$,  
by~$H^s_0(\Omega)$, the closure of~${\mathcal D}(\Omega)$ for the~$H^s $ norm, and by~$H^{-s}(\Omega)$ its  dual space.

The following result is a consequence of
Theorem~\ref{existence}.
\begin{Cor}\label{corexistenceeta}
{ 
Let $(\eta^0,u^0) \in L^2(\R\times \T)$ and $(\eta^0_\eps,u_\eps^0)$ such that
\begin{equation}
\label{initial-data}
\begin{aligned}
\frac12 \int \left(
|\eta^0_\eps|^2 +(1 + \e
 \eta^0_\eps)
|u^0_\eps|^2\right)dx\leq \EE^0,\\
(\eta^0_\eps, u^0_\eps) \to (\eta^0,u^0) \hbox{ in } L^2(\R\times \T).
\end{aligned}
\end{equation}
Then, for all
$\eps>0$, System (\ref{SW-eps}) has at least one weak solution  $(\eta_\eps,u_\eps)$\label{sol-not} with initial data $(\eta_\eps^0,u_\eps^0)$ satisfying the uniform bound
\begin{equation}
\label{uniform}
 \frac12 \int (\eta_\eps^2
+(1+\eps \eta_\eps )|u_\eps|^2)(t,x) dx+\int_0^t \int \nu| \nabla
u_\eps|^2(s,x) dxds
\leq \EE^0.
\end{equation}
Furthermore $u_\eps $ is uniformly bounded in  $L^2_{loc}(\R^+;L^2(\R\times \T)).
$

 In particular, there exist $\eta\in L^\infty(\R^+;L^2(\R\times \T))$  and  $u\in
L^\infty(\R^+;L^2(\R\times
\T)) \cap L^2(\R^+;\dot H^1(\R\times \T))$
such that, up to extraction of a subsequence,
$$(\eta_\eps,u_\eps) \rightharpoonup (\eta,u) \hbox{ in w-}L^2_{loc}(\R^+\times \R\times 
\T).$$
}
\end{Cor}

\begin{proof}
Replacing $ h$ by~$ 1 + \e \eta$ in the
energy inequality~(\ref{energy}), we get (\ref{uniform}) from which we deduce that there 
exist $\eta\in L^\infty(\R^+;L^2(\R\times \T))$  and  $u\in
 L^2(\R^+;\dot H^1(\R\times \T))$
such that, up to extraction of a subsequence,
$$(\eta_\eps,u_\eps) \rightharpoonup (\eta,u) \hbox{ in w-}L^2_{loc}(\R^+\times \R\times 
\T).$$
Furthermore we have the following inequality:
$$\int_\Omega   |u_\eps|^2(t,x) dx \leq \int_\Omega  (1+\eps \eta_\eps )|u_\eps|^2(t,x) 
dx+\eps\left(\int_\Omega
\eta_\eps^2(t,x) dx \right)^{1/2} \left( \int_\Omega  |u_\eps|^4(t,x)dx \right)^{1/2}
 $$
 from which 
we deduce that
$$u_\e \in L^2_{loc}(\R^+;L^2(\R\times \T)),$$
where we have used the interpolation inequality
$$
\int_\Omega  |u_\eps|^4(t,x)dx \leq C \int_\Omega  |u_\eps|^2(t,x)dx \int_\Omega  |\nabla u_\eps|^2(t,x)dx .
$$
By Fatou's lemma we get that~$u$ belongs to~$L^\infty(\R^+;L^2(\R\times \T))$. That concludes the proof.
\end{proof}


\chapter{Equatorial waves} \label{equatorialwaves}
\setcounter{equation}{0}
The aim of this chapter is to describe precisely the various waves induced by the  singular perturbation\label{L-not}
\begin{equation}
\label{Ldef}
L :(\eta,u)\in L^2(\R \times \T) \mapsto (\DIV u , \beta x_1 u^\perp +\nabla\eta).
\end{equation}
In the first paragraph we study the kernel of the operator, which describes the mean flow as we will see in 
Chapter~\ref{convergence}. In the second paragraph we describe all the other waves, using the Hermite 
functions in~$x_{1}$ and the Fourier transform in~$x_{2}$; this enables us to    recover results which are well-known from
physicists (see for instance~\cite{gill},\cite{pedlosky},\cite{ripa1}-\cite{ripa3}, 
as well as~\cite{dutrifoymajda} for a mathematical study). Finally in the last paragraph we study the possible
 resonances between all those waves; that result will be useful in
 Chapter~\ref{envelope} to prove   regularity estimates for the limit system introduced in Paragraph~\ref{sctschochet} below.


\section{The geostrophic constraint}
\label{sing}
\setcounter{equation}{0}
In this section we are going to study the kernel of the singular perturbation~$L$ defined in~(\ref{Ldef}).
\begin{Prop} \label{kernel}
{ 
Define the linear operator $L$ by (\ref{Ldef}). Then
$(\eta,u)
\in  L^2(\R \times \T)$ belongs to~$\Ker L$ if and only if $(\eta,u)$ belongs to~$ L^2(\R_{x_1})$
and 
\begin{equation}
\label{constraint}
u_1=0,\quad \beta  x_1 u_2+\d_1 \eta =0.
\end{equation}
}
\end{Prop}
\begin{proof}
If $(\eta,u)$ belongs to $L^2(\R\times \T) \cap \Ker L$, then we have
$$
   \DIV u=0,\quad \beta  x_1 u^\perp +\nabla
\eta=0,
$$ in the sense of distributions. Computing the vorticity in the
second identity leads to
$$
\nabla^\perp\cdot ( \beta x_1 u^\perp +\nabla
\eta)=  \beta (x_1 \nabla \cdot u +  u_1) =0,
$$
from which we deduce, since~$  \nabla \cdot  u=0 $,  that
$u_1=0.$
Plugging this identity respectively in the divergence-free condition and
in the second component of the vectorial condition gives
$$
\d_2 u_2=0,\quad \d_2 \eta=0,
$$
meaning that $\eta$ and $u$ depend only on the $x_1$ variable. The last
condition can then be rewritten
$$
 \beta x_1 u_2+\d_1 \eta=0.
$$
Conversely, it is easy to check that any $(\eta,u) \in L^2(\R)$ satisfying
(\ref{constraint}) belongs to $\Ker L$. 
\end{proof}
In the following we will denote  by $\Pi_0$ the orthogonal
projection of $L^2(\R\times \T)$ onto $\Ker L$. It is given by the
following formula.

\begin{Prop}\label{explicitPi0} { 
Define the linear operator $L$ by (\ref{Ldef}). Denote by $\Pi_0$\label{pi0-not} the orthogonal
projection of $L^2(\R\times \T)$ onto $\Ker L$. Then, for all $(\eta,u) \in L^2(\R\times 
\T)$
\begin{equation}
\label{pi-def}
\Pi_0(\eta,u)= \left(  \int  (DD^T +Id)^{-1} (\eta
+Du_2)dx_2, 0, \int D^T(DD^T +Id)^{-1} (\eta +Du_2)dx_2  \right),
\end{equation}
where $D$ is the differential operator defined by $\displaystyle D\cdot=\d_1\left({\cdot 
\over
\beta x_1}\right)$. }
\end{Prop}
\begin{proof}
By Proposition \ref{kernel}, for all $(\eta,u) \in L^2(\R\times \T)$,
$(  \eta^*,u^*)\eqdefa\Pi_0(\eta,u) $ belongs to $ L^2(\R)$ and satisfies
$$ u_1^*=0,\quad \beta x_1  u_2^* +\d_1 \eta^* =0.$$

Averaging with respect to the $x_2$-variable, one is reduced to the case when $(\eta,u)
\in L^2(\R)$.

By definition  $(\eta- \eta^*,u- u^* )$ is orthogonal in~$ L^2$ to any element $(\rho,v)$ 
of $\Ker
L$~: that implies that
$$
\int \left((\eta- \eta^*) \rho +(u_2- u_2^*) v_2 \right)dx_1=\int \left((\eta- \eta^*)
\rho -(u_2+{1\over \beta x_1} \d_1  \eta^*) {1\over \beta x_1} \d_1 \rho
\right)dx_1=0 .
$$
An integration by parts leads then to
$$\left( -\d_1 {1\over \beta ^2 x_1^2} \d_1 +Id\right) \eta^* = \eta + \d_1  {u_2\over 
\beta
x_1}\cdotp$$ Plugging this identity in the second constraint equation gives  the expected 
formula
for
$ u_2^*$.

That proves Proposition~\ref{explicitPi0}. 
\end{proof}

\section{Description of the waves}\label{descriptionofthewaves}
\setcounter{equation}{0}
In this section we are going to describe precisely the various waves created by~$L$. 
In the first paragraph of this section
(Paragraph~\ref{precisedescription}) we compute the eigenvalues of~$L$ and present its  eigenvectors, which
 constitute  a Hilbertian
basis of $L^2(\R\times \T)$ (that is proved in Paragraph \ref{diagonalization}). That basis
enables us  in the last paragraph to introduce the
filtering operator and formally derive the limit filtered system, in the spirit of S.
Schochet~\cite{schochet} (see also~\cite{grenier}).

  \subsection{Precise description of the oscillations}\label{precisedescription}
In this paragraph, we are going to  explain how to obtain the various 
eigenmodes of $L$.  The crucial point is that the description of
these eigenmodes can be achieved using the Fourier transform with 
respect to~$x_2$ and the decomposition on the Hermite
functions
$(\psi_n)_{n\in
\N}$ with respect to~$x_1$. Here the Hermite functions are conveniently rescaled so that
$$
\psi_n (x_1) = e^{-\frac{\beta x_1^2}{2}}P_n (\sqrt \beta x_1),
$$
where~$P_{n}$ is a polynomial of degree~$n$, and 
$(\psi_n)_{n\in \N}$\label{hermite-not} satisfy
$$-\psi_n''+\beta^2 x_1^2 \psi_n =\beta (2n+1) \psi_n.$$
We recall that $(\psi_n)_{n\in
\N}$ constitutes a Hermitian basis of $L^2(\R)$. 

Moreover we have the identities
\begin{equation}
\label{psi-identities}
\begin{aligned}
\psi_n'(x_1) +\beta x_1 \psi_n (x_1) &=\sqrt{2\beta n} \psi_{n-1}(x_1),\\
\psi_n'(x_1) -\beta x_1 \psi_n (x_1) &=-\sqrt{2\beta (n+1)} \psi_{n+1}(x_1).
\end{aligned}
\end{equation}
We have used the convention that~$\psi_n = 0$ if~$n <0$.

In the following we will then denote by~$\widehat f (n,k) $\label{widehat-not}, for~$ n \in
\N$ and~$ k \in \Z$, the components
of any function~$ f$ in the Hermite-Fourier basis~$(2\pi)^{-1/2} \psi_n (x_1)
e^{ik x_2} $. In other words we have
$$
\forall (n,k) \in
\N \times \Z, \quad \widehat f (n,k) = \frac1{\sqrt{2\pi}}\int_{\R \times \T}  \psi_n 
(x_1) e^{-ik x_2} f(x_1,x_2) \: dx_1 dx_2,
$$
along with the inversion formula
$$
\forall (x_1,x_2) \in
\R \times \T, \quad   f(x_1,x_2)=  {1\over \sqrt{2\pi}}\sumetage{n \in \N}{k \in \Z}
\psi_n (x_1) e^{ik x_2}\widehat f (n,k).
$$

In order to investigate the spectrum of $L$ (which is an unbounded 
skew-symmetric operator), we are interested in the non
trivial solutions
  to
\begin{equation}
\label{spectrum}
L(\eta,u)=i\tau (\eta,u).
\end{equation}
If one looks for the~$L^2$ solutions of (\ref{spectrum}) with $u_1$ non 
identically zero, one gets as a necessary condition that
the Fourier transform of $u_1$ with respect to $x_2$ (denoted by $\FF_2 u_1$\label{defFF2}) satisfies
$$(\FF_2 u_1)''+\left( \tau^2-k^2 +{\beta k\over \tau} -\beta^2 x_1^2 
\right) \FF_2 u_1=0,$$
from which we deduce that $\FF_2 u_1$ is  proportional to some $\psi_n$ and that
\begin{equation}
\label{tau-def1}
\tau^3 -(k^2+\beta(2n+1)) \tau+\beta k=0 ,
\end{equation}
for some $n\in \N$.

The following lemma is proved by elementary algebraic computations.  
\begin{Lem}\label{polynomial}
{ 
For any~$\beta > 0$ and any~$(n,k) \in \N^{*} \times \Z$, the
 polynomial
  $$
  P(\tau) = \tau^3 -(k^2+\beta(2n+1)) \tau+\beta k
  $$
   has three distinct roots in~$\R$, denoted
  in the following way:
  \begin{equation}
  \label{tau-not}
\tau(n,k,-1) <\tau(n,k,0)<\tau(n,k,1).
\end{equation}
 Moreover if~$\tau (n,k,j) = \tau (n',k,j') \neq 0 $ for some~$(n,n') \in \N^2$ with~$n  \neq 0$ and~$(j,j') \in \{-1,0,1\}^{2}$,
  then necessarily~$n = n'$
 and~$j = j'$.

Finally the following asymptotics hold if~$ n$ or~$ |k|$ goes to infinity:
$$
\tau(n,k,\pm 1) \sim  \pm\sqrt{k^2+\beta(2n+1)} , \quad \mbox{and}  \quad 
\tau(n,k,0) \sim  {\beta k\over k^2+\beta(2n+1)} \cdotp
$$
}
\end{Lem}
\begin{proof}
To prove that the polynomial has three distinct roots we simply analyze the function~$P(\tau) $.
Its derivative~$P'(\tau)$ vanishes in~$\pm \alpha$, where
$$
\alpha = \sqrt{\frac{k^{2}+ \beta (2n+1)}{3}} \cdotp
$$
It is then enough to prove that~$P(\alpha) < 0$ and~$P(-\alpha) > 0$. Let us write the argument for~$P(\alpha)$.
We have
$$
P(\alpha) = - 2 \alpha^{3} + \beta k.
$$
But for~$n \neq 0$, one checks easily that~$
2 \alpha^{3} > \beta |k|.
$
Indeed one has
$$
4 \alpha^{6} = \frac{4}{27} (k^{2}+ \beta (2n+1))^{3} > k^{2} \beta^{2}
$$
as soon as~$n \geq 1$.  So the first
result of the lemma is proved.

To prove the second result,   we notice that if
$\tau(n,k,j)=\tau(n',k ,j')= \tau \neq 0$,  then~$ 2(n-n') \tau =0$
from which we deduce that $n=n'$,    and therefore that
$j=j'$ since the polynomial (\ref{tau-def1}) admits three
separate   roots for~$n \neq 0$.

Finally the asymptotics of the eigenvalues is an easy computation. The
lemma is proved. 
\end{proof}

\begin{Rem}\label{racinesconfondues}
{ 
In the case when~$n = 0$, the three roots of~$P$ are
\begin{equation}
\label{tau-notagain}
\tau (0,k,-1) = -\frac k2 - \frac 12 \sqrt{k^{2} + 4 \beta} , \:  \: 
 \tau (0,k,0) = k, \:  \: \mbox{and} \:  \:  \tau (0,k,1) =   -\frac k2 +\frac 12 \sqrt{k^{2} + 4 \beta} .
\end{equation}
It follows that in the case when~$\beta = 2 k^{2}$, the roots become~$k$ (double) and~$-2k$.
}
\end{Rem}

Now let us study more precisely the waves generated by~$L$.

\medskip
$\bullet$ If $k\neq 0$ and $n\neq 0$, (\ref{tau-def1}) admits three solutions according to 
Lemma~\ref{polynomial},  
and one can check (see Paragraph \ref{diagonalization} below) that these solutions are 
eigenvalues of $L$ associated to the following unitary
eigenvectors\label{pagepoincarenknotzero}\label{pagerossby}
$$\Psi_{n,k,j} (x_1,x_2)=C_{n,k,j} e^{ikx_2} \left(
  \begin{array}{c}
  \displaystyle {i\tau(n,k,j)\over k^2-\tau(n,k,j)^2} \psi'_n(x_1)-{ik 
\over k^2-\tau(n,k,j)^2}\beta x_1
\psi_n(x_1)\\
\psi_n(x_1)
   \\
\displaystyle {ik\over k^2-\tau(n,k,j)^2} \psi'_n(x_1)-{i\tau(n,k,j) 
\over k^2-\tau(n,k,j)^2}\beta x_1 \psi_n(x_1)
\end{array}
\right)$$
which can be rewritten\label{psi1-not}
\begin{equation}
\label{psi-def1}
\Psi_{n,k,j} (x_1,x_2)=C_{n,k,j} e^{ikx_2} \left(
  \begin{array}{c}
  \displaystyle {-i\over \tau(n,k,j)+k} \sqrt{\beta n\over 2} 
\psi_{n-1}(x_1)+{i\over \tau(n,k,j)-k}
\sqrt{\beta (n+1)\over 2}\psi_{n+1}(x_1)\\
\psi_n(x_1)
   \\
\displaystyle {i\over \tau(n,k,j)+k} \sqrt{\beta n\over 2} 
\psi_{n-1}(x_1)+{i\over \tau(n,k,j)-k}
\sqrt{\beta (n+1)\over 2}\psi_{n+1}(x_1)
\end{array}
\right)
\end{equation}
because of the identities
(\ref{psi-identities}). The factor~$C_{n,k,j}$ ensures that~$\Psi_{n,k,j}$ is of norm~$1$
 in~$L^{2}(\R\times\T)$, 
its precise value is given in~(\ref{Cnkj}) below.

The modes corresponding to $\tau(n,k,-1)$ and $\tau(n,k,1)$ are 
called {\it Poincar{\'e} modes}, and satisfy
$$\tau(n,k,\pm 1) \sim \pm\sqrt{k^2+\beta(2n+1)} \hbox{ as } |k| \:
\mbox{or} \: n\to \infty,$$
which are the  frequencies of the gravity waves.

The modes corresponding to $\tau(n,k,0)$ are called {\it Rossby modes}, and satisfy
$$\tau(n,k,0)\sim {\beta k\over k^2+\beta(2n+1)} \hbox{ as } |k| \:
\mbox{or} \: n\to \infty,$$
meaning that the oscillation frequency is very small~: the planetary 
waves $\Psi_{n,k,0}$ satisfy indeed the
quasi-geostrophic approximation.

\medskip
$\bullet$ If $k= 0$ and $n\neq 0$, the three distinct solutions to 
(\ref{tau-def1}) are the two Poincar{\'e} modes $$\tau(n,0,\pm1)=\pm
\sqrt{\beta(2n+1)}$$ \label{pagepoincarekzero}and the {\it non-oscillating, or geostrophic, mode} \label{pagenonoscillating}
$\tau(n,0,0)=0$. The  corresponding eigenvectors of $L$ are
given by (\ref{psi-def1}) if $j\neq 0$ and by\label{psi2-not}
\begin{equation}
\label{psi-def2}
\Psi_{n,0,0} (x_1)= \frac{1}{\sqrt{2\pi(2n+1)}}\left(
  \begin{array}{c}
  \displaystyle - \sqrt{n+1\over 2} \psi_{n-1}(x_1)-
\sqrt{ n\over 2}\psi_{n+1}(x_1)\\
0
   \\
\displaystyle  \sqrt{ n+1 \over 2} \psi_{n-1}(x_1)-
\sqrt{n\over 2}\psi_{n+1}(x_1)
\end{array}
\right).
\end{equation}

\medskip
$\bullet$ If $n=0$, the three solutions to 
(\ref{tau-def1}) are the two Poincar{\'e} \label{pagepoincarenzero}and {\it mixed
Poincar{\'e}-Rossby modes} \label{pagemixed}
\begin{equation} \label{psi-def11}\tau(0,k,\pm1)=-{k\over 2} \pm \frac12
\sqrt{k^2+4\beta}
\end{equation}
with asymptotic behaviours given by
$$
\begin{aligned}
\tau (0,k,-\hbox{ sgn} (k))\sim -k \hbox{ as } |k|\to \infty,
\\
\tau (0,k,\hbox{ sgn} (k))\sim {\beta \over k} \hbox{ as } |k|\to \infty,
\end{aligned}$$
and the {\it  Kelvin mode} \label{pagekelvin}
$\tau(0,k,0)=k$. The  corresponding eigenvectors of $L$ are
given by (\ref{psi-def1}) if $j\neq 0$ and by\label{psi3-not}
\begin{equation}
\label{psi-def3}
\Psi_{0,k,0} (x_1,x_2)={1\over \sqrt{4\pi}} e^{ikx_2}\left(
  \begin{array}{c}
  \displaystyle \psi_0 (x_1)\\
0
   \\
\displaystyle \psi_0 (x_1)
\end{array}
\right).
\end{equation}
Note that in the case when the fluid studied is the atmosphere rather than the ocean, the mixed Poincar\'e-Rossby waves
are known as {\it Yanai waves}.

We recall that the functions~$\psi_n$ are defined by
$$
\psi_n (x_1) = e^{-\frac{\beta x_1^2}{2}}P_n (\sqrt\beta  x_1),
$$
where~$P_n$ is a polynomial of degree~$n$. We therefore have an
exponential decay far from the equator.

As mentioned in the introduction, the adjustment processes are
somewhat special in the vicinity of the equator (when 
the Coriolis
acceleration vanishes). A very important property of the equatorial zone is that
it acts as a {\sl waveguide}, i.e., disturbances are 
trapped in
the vicinity of the equator. The waveguide effect is due entirely to the
variation of Coriolis parameter with latitude.

Note that another important  effect of the
waveguide is the separation into a discrete set of modes
$n=0,1,2,...$ as occurs in a channel.
%

\subsection{Diagonalization of the singular perturbation}\label{diagonalization}
In this paragraph we are going to show that the previous study does provide   a Hilbertian basis of
eigenvectors.
\begin{Prop}\label{diag-prop} 
{ 
For all $(n,k,j)\in \N\times \Z\times 
\{-1,0,1\}$, denote by $\tau(n,k,j)$ the three roots   of~(\ref{tau-def1})  
and by $\Psi_{n,k,j}$ the unitary vector defined in Paragraph~\ref{precisedescription}.

Then $(\Psi_{n,k,j})_{(n,k,j)\in \N\times \Z\times \{-1,0,1\}}$ is a 
Hilbertian basis of $L^2(\R\times \T)$ constituted of
eigenvectors of $L$~:
\begin{equation}
\label{eigenvectors}
\forall (n,k,j)\in \N\times \Z\times \{-1,0,1\}, \quad L\Psi_{n,k,j} =i\tau(n,k,j) \Psi_{n,k,j}.
\end{equation}

Furthermore we have the following  estimates~: for all $s\geq 0$, there
exists  a nonnegative
constant~$C_s$ such that, for all
$(n,k,j)\in
\N\times
\Z
\times
\{-1,0,1\}$,
\begin{equation}
\label{eigenvectors-est}
\begin{aligned}
\| \Psi_{n,k,j} \|_{L^\infty(\R\times \T)} \leq C_0, \quad 
\| \Psi_{n,k,j} \|_{W^{s,\infty} (\R\times \T)}\leq C_s
(1+|k|^2+n)^{s/2}, \\
 C_s^{-1} (1+|k|^2+n)^{s/2} \leq 
\| \Psi_{n,k,j} \|_{H^s (\R\times \T)}\leq C_s (1+|k|^2+n)^{s/2},
\end{aligned}
\end{equation}
where~$W^{s,\infty}$ denotes the usual Sobolev space.
Moreover  the  eigenspace 
associated with any non zero eigenvalue is of finite dimension. 
}\end{Prop}

\begin{proof}
In order to establish the diagonalization result, the three points to be checked are  the
identity  (\ref{eigenvectors}), the orthonormality of the family
$(\Psi_{n,k,j})$, and the fact that it generates the whole space 
$L^2(\R\times \T)$.

$\bullet$ {$\Psi_{n,k,j}$ is an eigenvector of $L$}$ $

We start by establishing the identity (\ref{eigenvectors}), where 
  $\tau(n,k,j)$ is
defined by (\ref{tau-not}) and~(\ref{tau-notagain}) and
$\Psi_{n,k,j}$ is defined either by (\ref{psi-def1})  (for the 
Poincar{\'e} and Rossby modes)  or by (\ref{psi-def2}) (for 
the non-oscillating modes), or by
(\ref{psi-def3}) (for the Kelvin modes).

 For the Poincar{\'e}, Rossby  and mixed  Poincar{\'e}-Rossby modes, we start from formula 
(\ref{psi-def1})
$$\Psi_{n,k,j} =C_{n,k,j} e^{ikx_2} \left(
  \begin{array}{c}
  \displaystyle {-i\over \tau(n,k,j)+k} \sqrt{\beta n\over 2} 
\psi_{n-1}(x_1)+{i\over \tau(n,k,j)-k}
\sqrt{\beta (n+1)\over 2}\psi_{n+1}(x_1)\\
\psi_n(x_1)
   \\
\displaystyle {i\over \tau(n,k,j)+k} \sqrt{\beta n\over 2} 
\psi_{n-1}(x_1)+{i\over \tau(n,k,j)-k}
\sqrt{\beta (n+1)\over 2}\psi_{n+1}(x_1)
\end{array}
\right).$$
We have $L\Psi_{n,k,j}=C_{n,k,j} e^{ikx_2} V_{n,k,j}$ where $V_{n,k,j}$ denotes
$$ \left(
  \begin{array}{c}
  \displaystyle \psi'_n(x_1)+ik  \sqrt{\frac\beta 2}
  \left(
  \frac{i\sqrt n} {\tau(n,k,j)+k} 
 \psi_{n-1}(x_1)+\frac{i\sqrt{n+1} }{ \tau(n,k,j)-k}
 \psi_{n+1}(x_1) 
\right)\\
  \displaystyle   \sqrt{\frac\beta 2} \left(
  \frac{i\sqrt n} {\tau(n,k,j)+k}  (\beta 
x_1 \psi_{n-1}(x_1)-\psi'_{n-1}(x_1))
+ \frac{i\sqrt{n+1}} {\tau(n,k,j)-k} 
 (\beta x_1\psi_{n+1}(x_1)+\psi'_{n+1}(x_1))\right)
   \\
\displaystyle -\beta x_1 \psi_n(x_1) +ik 
\sqrt{\frac\beta 2}
\left(  \frac{ -i\sqrt n} 
{\tau(n,k,j)+k }  
\psi_{n-1}(x_1)+\frac{ i\sqrt{n+1}} { \tau(n,k,j)-k }
 \psi_{n+1}(x_1)\right)
\end{array}
\right)$$
which can be rewritten using the identities (\ref{psi-identities})
$$ \left(
  \begin{array}{c}
  \displaystyle {\tau(n,k,j)\over \tau(n,k,j)+k} \sqrt{\beta n\over 2} 
\psi_{n-1}(x_1)-{\tau(n,k,j)\over
\tau(n,k,j)-k}
\sqrt{\beta (n+1)\over 2}\psi_{n+1}(x_1)\\
  \displaystyle  {i\over \tau(n,k,j)+k} \beta n \psi_n(x_1)
+{i\over
\tau(n,k,j)-k}
\beta (n+1) \psi_n (x_1)
   \\
\displaystyle  -{\tau(n,k,j)\over \tau(n,k,j)+k} \sqrt{\beta n\over 2}
\psi_{n-1}(x_1)-{\tau(n,k,j)\over \tau(n,k,j)-k}
\sqrt{\beta (n+1)\over 2}\psi_{n+1}(x_1)
\end{array}
\right).$$
As $\tau(n,k,j)$ satisfies (\ref{tau-def1}), we have
$${i\over \tau(n,k,j)+k} \beta n+{i\over
\tau(n,k,j)-k}
\beta (n+1)= i {(2n+1)\beta \tau(n,k,j) +\beta k \over 
\tau(n,k,j)^2-k^2}=i\tau(n,k,j)$$
from which we deduce that
$$L\Psi_{n,k,j} =i\tau(n,k,j) \Psi_{n,k,j} \hbox{ for all } 
(n,k,j)\in \N\times \Z \times \{-1,1\} \cup
\N^*\times \Z^* \times \{0\}.$$

\medskip
 For the Kelvin modes we start from  formula (\ref{psi-def3})
$$\Psi_{0,k,0} ={1\over \sqrt{4\pi}} e^{ikx_2}\left(
  \begin{array}{c}
  \displaystyle \psi_0\\
0
   \\
\displaystyle \psi_0
\end{array}
\right).$$
We have
$$L\Psi_{0,k,0} ={1\over \sqrt{4\pi}} e^{ikx_2}\left(
  \begin{array}{c}
  \displaystyle ik\psi_0\\
\beta x_1 \psi_0 +\psi_0'
   \\
\displaystyle ik \psi_0
\end{array}
\right)= {ik\over \sqrt{4\pi}} e^{ikx_2}\left(
  \begin{array}{c}
  \displaystyle \psi_0\\
0
   \\
\displaystyle  \psi_0
\end{array}
\right),$$
or equivalently
$$L\Psi_{0,k,0}=ik \Psi_{0,k,0} \hbox{ for all } k\in \Z.$$

\medskip
 For the non-oscillating modes we start from 
formula (\ref{psi-def2})
$$\Psi_{n,0,0} = \frac{1}{\sqrt{2\pi (2n+1)}}\left(
  \begin{array}{c}
  \displaystyle - \sqrt{  (n+1)\over 2} \psi_{n-1} -
\sqrt{ n\over 2}\psi_{n+1} \\
0
   \\
\displaystyle  \sqrt{  n+1 \over 2} \psi_{n-1}-
\sqrt{  n\over 2}\psi_{n+1} 
\end{array}
\right).$$
An easy computation shows that
$$L\Psi_{n,0,0} =\frac{1}{\sqrt{2\pi (2n+1)}}\left(
  \begin{array}{c}
  0\\
\displaystyle   \sqrt{  (n+1)\over 2} (   x_1 
\psi_{n-1} -\psi'_{n-1} ) -
\sqrt{  n\over 2}(x_1 \psi_{n+1} +\psi'_{n+1} )
   \\
0
\end{array}
\right)$$
which is zero by (\ref{psi-identities}). Thus,
$$L\Psi_{n,0,0}=0 \hbox{ for all } n\in \N^*.$$

\bigskip
$\bullet$ {$(\Psi_{n,k,j})$ is an orthonormal family}$ $

By identity (\ref{eigenvectors}) and the fact that
$((2\pi)^{-1/2}e^{ikx_2})_{k \in \Z}$ and $( \psi_n(x_1))_{n \in \N}$ are 
respectively Hilbertian basis of
$L^2(\T)$ and $L^2(\R)$ we are going to deduce that
$(\Psi_{n,k,j})$ is an orthonormal family.

\medskip
 In   formula (\ref{psi-def1}), we choose
\begin{equation}\label{Cnkj}
C_{n,k,j}   =(2\pi)^{-1/2}\left({\beta n\over (\tau(n,k,j)+k)^2}
+{\beta (n+1 ) \over (\tau(n,k,j)-k)^2} +1
\right)^{-1/2}
\end{equation}
so that
$$ \| \Psi_{n,k,j}\|^2_{L^2(\R\times \T)} =1,$$
for  all $(n,k,j)\in \N\times \Z \times \{-1,1\} \cup
\N^*\times \Z^* \times \{0\}.$ 

In the same way,  it is immediate to check that 
$$\|\Psi_{n,0,0}\|^2_{L^2(\R\times \T)}=1,$$
for all $n\in \N^*$, and that
$$\|\Psi_{0,k,0}\|^2_{L^2(\R\times \T)}=1,$$
for all $k\in \Z$.

\medskip
 In order to establish the orthogonality property
we proceed in two steps. 

If $\tau(n,k,j)\neq \tau(n',k',j') $, as $L$ is a 
skew-symmetric operator, we have
$$\begin{aligned}
i\tau(n,k,j) \big( \Psi_{n,k,j} |\Psi_{n',k',j'}\big)
&=- \big( L\Psi_{n,k,j} |\Psi_{n',k',j'}\big)\\
&= \big( \Psi_{n,k,j} |L\Psi_{n',k',j'}\big)\\
&=i\tau(n',k',j') \big( \Psi_{n,k,j} |\Psi_{n',k',j'}\big)
\end{aligned}$$
from which we deduce that
$$\big( \Psi_{n,k,j} |\Psi_{n',k',j'}\big)=0.$$

If  $\tau(n,k,j)=\tau(n',k',j') $, we first note that 
$$\hbox{ for } k\neq k',\quad \big( \Psi_{n,k,j}
|\Psi_{n',k',j'}\big)=0$$ using the orthogonality of $
e^{ikx_2}$ and $e^{ik'x_2} $.
So we are left with the case when~$k = k'$. First, if~$\tau(n,k,j)=\tau(n',k ,j') =\tau\neq 0$ and~$n \neq 0$, then 
Lemma~\ref{polynomial} implies that~$n=n'$ and~$j = j'$. Then in
 the case when~$\tau(0,k,j)=\tau(0,k ,j') =\tau\neq 0$, with~$j \neq j'$,  we just have to consider the explicit definition 
 of~$\Psi_{0,k,j}$ and~$\Psi_{0,k,j'}$ given in Paragraph~\ref{precisedescription} to find that
 $$
 \big( \Psi_{0,k,j} |\Psi_{0,k ,j'}\big)=0.
 $$
Finally, if~$\tau(n,k,j)=\tau(n',k',j') =0$, we have~$k=k'=0$ and~$j=j'=0$ and we deduce
 from formula~(\ref{psi-def2}) that
$$\hbox{ for } n\neq n',\quad\big( \Psi_{n,0,0}
|\Psi_{n',0,0}\big)=0.$$

We thus conclude that
$$\big( \Psi_{n,k,j} |\Psi_{n',k',j'}\big)=0,$$
as soon as $(n,k,j)\neq (n',k',j')$.

\bigskip
$\bullet$ {$(\Psi_{n,k,j})$ spans $L^2(\R\times \T)$}$ $

It remains      therefore to see that any vector $\Phi$ of
$L^2(\R\times \T)$ can  be
decomposed                             on the family
$(\Psi_{n,k,j})$.

We first decompose each component on the Hermite-Fourier
basis
$$\Phi (x_{1},x_{2})= \frac1{\sqrt{2\pi}} \sum_{n,k} e^{ikx_2} \left(
  \begin{array}{c}
  \displaystyle \hat \Phi_0(k,n) \psi_n(x_1)\\
\displaystyle\hat \Phi_1(k,n) \psi_n(x_1)
   \\
\displaystyle\hat \Phi_2(k,n) \psi_n(x_1)
\end{array}
\right),$$
which can be rewritten
$$\begin{aligned}
 &\frac1{2\sqrt{2\pi}}  \sumetage{n>0}{k} e^{ikx_2} \left(
  \begin{array}{c}
  \displaystyle   (\hat \Phi_0(k,n+1)+\hat
\Phi_2(k,n+1))
\psi_{n+1}(x_1)+   (\hat \Phi_0(k,n-1)-\hat
\Phi_2(k,n-1))\psi_{n-1}(x_1)\\
\displaystyle 2\hat \Phi_1(k,n) \psi_n(x_1)
   \\
\displaystyle  (\hat \Phi_0(k,n+1)+\hat
\Phi_2(k,n+1))
\psi_{n+1}(x_1)-   (\hat \Phi_0(k,n-1)-\hat
\Phi_2(k,n-1))\psi_{n-1}(x_1)
\end{array}
\right)\\
&+\frac1{2\sqrt{2\pi}} 
\sum_{k} e^{ikx_2} \left(
  \begin{array}{c}
  \displaystyle   (\hat \Phi_0(k,0)+\hat
\Phi_2(k,0))
\psi_0(x_1)+  (\hat \Phi_0(k,1)+\hat
\Phi_2(k,1))\psi_{1}(x_1)\\
\displaystyle2\hat \Phi_1(k,0) \psi_0(x_1)
   \\
\displaystyle   (\hat \Phi_0(k,0)+\hat
\Phi_2(k,0))
\psi_{0}(x_1)+   (\hat \Phi_0(k,1)+\hat
\Phi_2(k,1))\psi_{1}(x_1)
\end{array}
\right).
\end{aligned}$$

\medskip
We then introduce for all $(n,k)\in \N\times \Z$ the matrix
$M_{n,k}             \in  M_3(\R) $ defined by
\begin{equation}
\label{m-def1}
M_{n,k}=\left(\begin{array}{ccc}
\displaystyle {-iC_{n,k,-1} \sqrt{\beta n/2}\over
\tau(n,k,-1)+k}    &\displaystyle{-iC_{n,k,0} \sqrt{\beta
n/2}\over
\tau(n,k,0)+k} &\displaystyle {-iC_{n,k,1} \sqrt{\beta
n/2}\over
\tau(n,k,1)+k} \\ C_{n,k,-1}&C_{n,k,0}&C_{n,k,1}\\
\displaystyle
{iC_{n,k,-1} \sqrt{\beta (n+1)/2}\over
\tau(n,k,-1)-k}    &\displaystyle {iC_{n,k,0} \sqrt{\beta
(n+1)/2}\over
\tau(n,k,0)-k} &\displaystyle {iC_{n,k,1} \sqrt{\beta
(n+1)/2}\over \tau(n,k,1)-k}
\end{array}\right)
\end{equation}if $n\neq 0$ and $k\neq 0$, by
\begin{equation}
\label{m-def2}
M_{n,0}=\left(\begin{array}{ccc}
\displaystyle {-iC_{n,0,-1} \sqrt{\beta n/2}\over
\tau(n,0,-1)}    &\displaystyle -C_{n,0,0} \sqrt{\beta
(n+1)/2} &\displaystyle {-iC_{n,0,1} \sqrt{\beta
n/2}\over
\tau(n,0,1)} \\ C_{n,0,-1}&0&C_{n,0,1}\\
\displaystyle
{iC_{n,0,-1} \sqrt{\beta (n+1)/2}\over
\tau(n,0,-1)}    &\displaystyle - C_{n,0,0}\sqrt{\beta
n/2} &\displaystyle {iC_{n,0,1} \sqrt{\beta
(n+1)/2}\over \tau(n,0,1)}
\end{array}\right)
\end{equation}
if $n\neq 0$ and by
\begin{equation}
\label{m-def3}
M_{0,k}=\left(\begin{array}{ccc}
\displaystyle 0  &\displaystyle \sqrt{1/4\pi} &0 \\
C_{0,k,-1}&0&C_{0,k,1}\\
\displaystyle
{iC_{0,k,-1} \sqrt{\beta/2 }\over
\tau(0,k,-1)-k}    &\displaystyle 0 &\displaystyle
{iC_{0,k,1} \sqrt{\beta/2 }\over \tau(0,k,1)-k}
\end{array}\right).
\end{equation}
As the eigenvectors $\Psi_{n,k,-1}$, $\Psi_{n,k,0}$ and $\Psi_{n,k,1}$ are orthogonal in
$L^2(\R\times \T)$, these matrices are necessarily invertible.

\medskip
We conclude by checking that one can write
$$\Phi =\sum_{n,k,j} \varphi_{n,k,j} \Psi_{n,k,j}$$
where $\varphi_{n,k,j}$ is defined by\label{varphi-not}
$$\left(\begin{array}{c}\displaystyle
\varphi_{n,k,-1}\\ \displaystyle \varphi_{n,k,0}\\
\displaystyle
\varphi_{n,k,1} \end{array}\right) =\frac1{\sqrt{2\pi}} 
M_{n,k}^{-1}              \left(\begin{array}{c}\displaystyle
\frac12 (\hat \Phi_0(k,n-1)  -\hat\Phi_2(k,n-1))\\
\displaystyle
\hat \Phi_1(k,n)\\
\displaystyle
\frac12 (\hat \Phi_0(k,n+1)  +\hat\Phi_2(k,n+1))
\end{array}\right) $$
for $n\neq 0$, and by
               $$\left(\begin{array}{c}\displaystyle
\varphi_{0,k,-1}\\ \displaystyle \varphi_{0,k,0}\\
\displaystyle
\varphi_{0,k,1} \end{array}\right) =\frac1{\sqrt{2\pi}} 
M_{0,k}^{-1}              \left(\begin{array}{c}\displaystyle
\frac12 (\hat \Phi_0(k,0)  +\hat\Phi_2(k,0))\\
\displaystyle
\hat \Phi_1(k,0)\\
\displaystyle
\frac12 (\hat \Phi_0(k,1)  +\hat\Phi_2(k,1))
\end{array}\right) .$$

\bigskip
$\bullet$ The regularity estimates are obtained using the explicit formulas (\ref{psi-def1}),
(\ref{psi-def2}) and (\ref{psi-def3}), as well as the following bounds on the elementary
Fourier and Hermite functions~:
$$\|e^{ikx_2} \|_{\dot H^s(\T)}=|k|^s \|e^{ikx_2} \|_{ L^2(\T)},\quad \|e^{ikx_2} \|_{\dot
 W^{s,\infty}(\T)}=|k|^s \|e^{ikx_2} \|_{ L^\infty(\T)}$$
and
$$\| \psi_n\|_{H^s(\R)} \sim  (1+n)^{s/2}\sup_n \| \psi_n\|_{L^2(\R)} ,\quad \|
\psi_n\|_{  W^{s,\infty}(\R)}\leq C_s  (1+n)^{s/2}
\sup_n \|\psi_n\|_{L^\infty(\R)}$$
coming from identities (\ref{psi-identities}) by a simple recurrence.
The crucial point is therefore to have a uniform $L^\infty$-bound on the Hermite functions,
which is stated for instance in~\cite{lebedev}:
\begin{equation}
\label{Linfty-bounds-psi}
\forall n\in \N, \quad \|\psi_n\|_{L^\infty(\R)} \leq C_\infty \hbox{ with } C_\infty\sim
1.086435.
\end{equation}

\medskip

Finally let us prove that the eigenspace associated with a nonzero
 eigenvalue is of finite dimension. Suppose by contradiction that
 there is~$\lambda \neq 0$ and
  a sequence~$(n_M,k_M,j_M)_{M \in\N}$ in~$\N \times \Z \times
 \{-1,0,1\}$ such that
$$
\tau (n_M,k_M,j_M) = \lambda \quad \mbox{and}\quad n_M + |k_M|
\rightarrow \infty,\:   \mbox{as} \:  M \rightarrow \infty.
$$
By  Lemma~\ref{polynomial},  as~$n$ or~$|k|$ goes
to infinity, the eigenvalue~$\tau (n,k,j)$ goes to zero or
to~$\pm\infty$, which contradicts the assumption that~$\tau (n_M,k_M,j_M) = \lambda $.
 
This concludes the proof of Proposition \ref{diag-prop}.  
\end{proof}

\bigskip
As the behaviour of the eigenmodes are expected to depend strongly of 
their type, i.e. of the class of the corresponding
eigenvalue, we split $L^2(\R\times \T)$ into five supplementary 
subsets, namely the Poincar{\'e} modes, the Rossby
modes, the mixed Poincar{\'e}-Rossby modes, the Kelvin modes and the non-oscillating modes.

\begin{Def} \label{PRKN}
{ 
With the  above notation, let us define
$$
\begin{aligned}
P&=Vect\Bigl\{\Psi_{n,k,j}\,/\, (n,k,j)\in \N^*\times \Z \times \{-1,1\} \cup \{0\} \times
\{(k,-\mbox{sign}(k)) _{ k\in \Z^*}\} \cup \{0\} \times\{0\} \times\{-1,1\}
\Bigr\},\\
R&=Vect\{\Psi_{n,k,0}\,/\, (n,k)\in \N^*\times \Z^* \},\\
 M&=Vect\{\Psi_{0,k,j}\,/\, k\in \Z^*, j= \mbox{sign}(k)\},\\
K&=Vect\{\Psi_{0,k,0}\,/\, k\in \Z^* \},
\end{aligned}
$$
so that
$L^2(\R\times \T)=P\oplus R\oplus M\oplus K \oplus \Ker L.$
Then we denote by $\Pi_P$ (resp. $\Pi_R,\Pi_M, \Pi_K$ and $\Pi_0$) the $L^2$ 
orthogonal projection on $P$ (resp. on $R$, $M$, $K$ and
$\Ker L$).

Moreover we define~${\mathfrak S}$\label{sig-not} the set of all eigenvalues of~$L$, as well as the following
subsets of~$\N \times \Z \times \{-1,0,1\}$:
$$
\begin{aligned}
{\mathfrak S}_{P} & = \left\{\tau (n,k,j) \: \Big / \: (n,k,j) \in   \N^*\times \Z \times \{-1,1\}
\right\}
\cup
 \left\{\tau (0,k, -\mbox{sign}(k)) \: \Big / \:  k\in \Z^* \right\}\cup \{ \pm \sqrt{\beta} \}, \\
{\mathfrak S}_{R}& =\left\{\tau (n,k,0) \: \Big / \: (n,k,j) \in \N^*\times \Z^*  \right\} , \quad \mbox{and} \quad {\mathfrak S}_{K}  =  \Z^* .
\end{aligned}
$$
}
\end{Def}
 Finally it can  be useful for the rest of the study to sum up the previous notation in the following picture.
 \renewcommand{\arraystretch}{1.5}
\begin{table}[htbp]
\begin{center}
\begin{tabular}{|c|c|c|c|c|c|}
\hline
wave & $n$ & $k$ & $j$ & definition of  $\Psi_{n,k,j} $& definition of  $\tau( n,k,j ) $\\
\hline
Poincar{\'e} & $\N^{*}$ & $\Z$& \{-1,1\} &(\ref{psi-def1}) page \pageref{pagepoincarenknotzero} & 
$\displaystyle\tau (n,k,\pm 1) \sim \pm\sqrt{k^2+\beta(2n+1)}  $ \\
& $\{0\}$ & $\Z^{*}$& $- \mbox{sign}(k) $ & (\ref{psi-def1}) page \pageref{pagepoincarenzero} &$
 \displaystyle\tau (0,k,-\mbox{sign}( k)) \sim - k$\\
 & $\{0\}$ & $\{0\}$& $ \{-1,1\}$ & (\ref{psi-def11}) page \pageref{pagepoincarenzero} &$
 \displaystyle\tau (0,0,\pm 1) = \pm \sqrt \beta $\\
\hline
Mixed &$\{0\}$ & $\Z^{*}$& $ \mbox{sign}(k) $ &(\ref{psi-def1}) page \pageref{pagemixed} & $\displaystyle\tau (0,k, 
\mbox{sign} ( k) )\sim
 \frac \beta k$\\
\hline
Kelvin &$\{0\}$ & $\Z^{*}$&$\{0\}$&(\ref{psi-def3}) page \pageref{pagekelvin} &$ \tau (0,k,0) = k$\\
\hline
Rossby & $\N^{*}$& $\Z^{*}$& $\{0\}$&(\ref{psi-def1}) page
\pageref{pagerossby} & $ \displaystyle \tau (n,k,0)  \sim\frac{\beta k}{
k^{2} + \beta (2n+1)}$\\
\hline
non oscillating &$\N $&$\{0\}$ &$\{0\}$ &(\ref{psi-def2}) page \pageref{pagenonoscillating} &$\tau (n,0,0) = 0$\\
\hline
\end{tabular} 
\end{center}
 \caption{Description of the waves}
\end{table}

\subsection{Orthogonality properties of the eigenvectors}
In this section we are going to give some additional properties on the~$\Psi_{n,k,j}$ defined above, which will be useful in the 
next chapters. We will
 write~$\Pi_{n,k,j}$ for  the projection  on the
eigenmode $\Psi_{n,k,j}$ of $L$,
 and~$\Pi_\lambda$  for the projection on the
eigenspace associated with the eigenvalue~$i\lambda$ of $L$\label{pinkj-not}.
The main result is the following, which states an orthogonality property for the ageostrophic modes (meaning
the eigenvectors in~$(Ker L)^\perp$). Note that there is no analogue of that result for geostrophic modes.
\begin{Prop}\label{orthogonality}
{ 
Let~$s \geq 0$ be a given real number. There is a constant~$C> 0$ such that for any non
zero eigenvalue~$i \lambda$  of~$L$ and for any three component vector field~$\Phi$ in~$(\Ker L)^{\perp}$, we have
\begin{equation}\label{quasiortho}
C_s^{-1} \sum 
_{\tau(n,k,j) =\lambda} \|\Pi_{n,k,j}  \Phi \|^2_{\dot H^s(\R\times 
\T)} \leq \| \Pi_\lambda   \Phi\|_{\dot H^s(\R\times \T)}^2 \leq C_s\sum 
_{\tau(n,k,j) =\lambda} \|\Pi_{n,k,j}  \Phi \|^2_{\dot H^s(\R\times 
\T)}. 
\end{equation}
}
\end{Prop}
\begin{proof}
Let~$\Phi$ in~$(\Ker L)^{\perp}$ be given and let~$s$ be any integer (the result for all~$s \geq 0$ will follow
by interpolation). We have 
$$
\begin{aligned}
\left\|\partial^{s} (\Pi _\lambda\Phi)
\right\|_{L^2(\R\times \T)}^2=& \sum_{\tau(n,k,j)=\lambda }\| \partial^{s}
(\Pi_{n,k,j}\Phi)
\|_{L^2(\R\times \T)}^2\\
&+\sum _{\tau(n,k,j)=\tau(n^*,k^*,j^*)=\lambda,\atop
(n,k,j)\neq (n^*,k^*,j^*)} \left(\partial^{s} (\Pi_{n,k,j}\Phi)|\partial^{s}
(\Pi_{n^*,k^*,j^*}\Phi)\right)_{L^2(\R\times \T)}.
\end{aligned}$$
 Of course,
\begin{equation}
\label{k-ortho}
\left( \partial^{s}\Psi_{n,k,j}|\partial^{s}
\Psi_{n^*,k^*,j^*}\right)_{L^2(\R\times \T)}=0 \hbox{ if } k\neq
k^*.
\end{equation}
Moreover we know by Proposition~\ref{diag-prop}, page~\pageref{diag-prop}
 that  if $\tau(n,k,j)=\tau(n^*,k,j^*)=
\lambda\neq 0$ and~$n \neq 0$, then necessarily~$n=n^*$ and~$j=j^*$. Therefore one has in fact
$$\longformule{
\sum _{\tau(n,k,j)=\tau(n^*,k^*,j^*)=\lambda,\atop
(n,k,j)\neq (n^*,k^*,j^*)} \left( \partial^{s} (\Pi_{n,k,j}\Phi)|\partial^{s}
(\Pi_{n^*,k^*,j^*}\Phi)\right)_{L^2(\R\times \T)}} {= \sum _{\tau(0,k,j)=\tau(0,k ,j^*)=\lambda,\atop
 j \neq j^* } \left( \partial^{s} (\Pi_{0,k,j}\Phi)|\partial^{s}(\Pi_{0,k ,j^*}\Phi)\right)_{L^2(\R\times \T)}.}
$$
But according to Remark~\ref{racinesconfondues} page~\pageref{racinesconfondues}, such 
a  situation occurs only if~$2k^{2} = \beta$, in which
case~$\tau(0,k,j) $ is equal to~$ k$. So    there
is at most one possible value  for~$k$ ($k = \lambda$) which occurs only in the 
case when~$\lambda = \pm \sqrt{\beta /2}$. In this last case, we have obviously
$$
\begin{aligned}\| \Pi_\lambda   \Phi\|_{\dot H^s(\R\times \T)}^2&
\sim \| \Pi_\lambda   \Phi\|_{L^2(\R\times \T)}^2\\
&=\sum _{\tau(n,k,j)
=\lambda}
 \|\Pi_{n,k,j}  \Phi \|^2_{L^2(\R\times 
\T)}\\
&\sim \sum_{\tau(n,k,j) =\lambda} \|\Pi_{n,k,j}  \Phi \|^2_{\dot H^s(\R\times 
\T)}.
\end{aligned}
$$
 The result follows.
\end{proof}

\begin{Rem}\label{componentseparately}
Note that the same argument allows actually to prove similar estimates for the components separately~:
$$
\begin{aligned}
C_s^{-1} \sum 
_{\tau(n,k,j) =\lambda} \|(\Pi_{n,k,j}  \Phi)' \|^2_{\dot H^s(\R\times 
\T)} \leq \| (\Pi_\lambda   \Phi)'\|_{\dot H^s(\R\times \T)}^2 \leq C_s\sum 
_{\tau(n,k,j) =\lambda} \|(\Pi_{n,k,j}  \Phi)' \|^2_{\dot H^s(\R\times 
\T)},\\
 C_s^{-1} \sum 
_{\tau(n,k,j) =\lambda} \|(\Pi_{n,k,j}  \Phi)_0 \|^2_{\dot H^s(\R\times 
\T)} \leq \| (\Pi_\lambda   \Phi)_0\|_{\dot H^s(\R\times \T)}^2 \leq C_s\sum 
_{\tau(n,k,j) =\lambda} \|(\Pi_{n,k,j}  \Phi)_0 \|^2_{\dot H^s(\R\times 
\T)},
\end{aligned}
$$
denoting by $\Phi_0$ the first coordinate and by $\Phi'$ the two other coordinates of $\Phi$\label{comp-not}.
\end{Rem}

\section{The filtering operator and the formal limit system} \label{sctschochet}
\setcounter{equation}{0}
In the previous paragraph we have presented a Hilbertian basis
of $L^2(\R\times \T) $ consisting in eigenvectors of the
singular penalization $L$. We are then able to define, in the
spirit of S. Schochet
\cite{schochet}, the ``filtering operator" associated with
the system.

Let $\LL$ be the semi-group generated by $L$~: we write
$\LL(t)=\exp \left(-tL\right)$\label{LL-not}. Then, for any
three component  vector field $\Phi\in L^2(\R\times \T) $,
we have
\begin{equation}
\label{LL-def}\LL(t)\Phi =\sum_{i\lambda \in {\mathfrak S}}
e^{-it\lambda}
\Pi_\lambda \Phi,
\end{equation}
where $\Pi_\lambda$ denotes the $L^2$ orthogonal projection
on the eigenspace of $L$ corresponding to the eigenvalue~$i\lambda$, and where~${\mathfrak S}$ denotes the set of all the eigenvalues of~$L$.

\bigskip
Now let us consider $(\eta_\eps, u_\eps)$ a weak solution to
(\ref{SW-eps}), which is formally equivalent to
\begin{equation}
\label{SW-epsu}
\begin{aligned}
\d_t \eta_\eps +\frac1\eps \DIV \Bigl((1+\eps \eta_\eps) u_\eps\Bigr) =0, \nonumber\\
\d_t  u_\eps + u_\eps\cdot \nabla u_\eps+  
\frac{ \beta x_1}{ \eps}  u_\eps^\perp+
{1\over \eps} 
\nabla
\eta_\eps
-\frac\nu {1+\e\eta_\eps}\Delta u_\eps =0, \nonumber & \\
\eta_{ \eps|t=0}=\eta_\eps^0,\quad u_{ \eps|t=0} =u_\eps^0 &,
\end{aligned}
\end{equation}
and let us define\label{phieps-not} 
\begin{equation}
\label{phi-def}
\Phi_\eps  =\LL\left(-{t\over \eps}  \right) (\eta_\eps, u_\eps).
\end{equation}
 Conjugating
formally equation (\ref{SW-epsu}) by the semi-group leads to
\begin{equation}
\label{filtered}
\d_t \Phi_\eps  + \LL\left(-{t\over \eps}  \right)
Q\left(\LL\left({t\over \eps}  \right)\Phi_\eps,
\LL\left({t\over \eps}  \right)\Phi_\eps\right)  -\nu \LL\left(-{t\over
\eps}  \right)\Delta' \LL\left({t\over \eps}  \right)\Phi_\eps=R_\eps,
\end{equation}
where $\Delta'$ and $Q$\label{op-not} are  the linear and symmetric bilinear operator
defined by
\begin{equation}
\label{Q-lap}
\Delta ' \Phi =(0,\Delta \Phi')                 \hbox{ and }
Q(\Phi,\Phi)=(\nabla \cdot(\Phi_0 \Phi'), (\Phi'\cdot \nabla) \Phi')
\end{equation}
 and
$$R_\eps =\LL\left(-{t\over \eps}\right)(0, -\nu {\eps \eta_\eps \over 1+\eps \eta_\eps}
\Delta u_\eps).$$

  We therefore expect to get a bound on the
time derivative of 
$ \Phi_\eps$ in some space of distributions.
A formal passage to the limit in~(\ref{filtered}) as $\eps $ goes to zero     (based on formula (\ref{LL-def}) and on a
nonstationary phase argument)           leads then to 
\begin{equation}\label{lim-filtered}
\d_t \Phi +Q_L (\Phi,\Phi) -\nu \Delta'_L \Phi =0,
\end{equation}
where $\Delta'_L$ and $Q_L$ denote  the linear and symmetric bilinear operator
defined by\label{opL-not}
\begin{equation}
\label{Q-lap-L}
\Delta _L'\Phi = \sum_{i\lambda \in {\mathfrak S}} \Pi_\lambda \Delta' \Pi_\lambda  
\Phi            
\hbox{ and } Q_L(\Phi,\Phi)=\sum_{i\lambda, i\mu,i\tilde \mu \in {\mathfrak S} \atop
\lambda=\mu+\tilde \mu} 
\Pi_\lambda Q(\Pi_\mu \Phi, \Pi_{\tilde \mu}\Phi).
\end{equation}
The study of~(\ref{lim-filtered}) is the object of
Chapter~\ref{envelope}. The proof that~(\ref{lim-filtered}) is indeed
the limit system to~(\ref{filtered}) is the object of
Chapter~\ref{convergence}.

In the next section we study the resonances associated with the
operator~$ L$: more precisely we describe in what cases the equality
$$
\tau(n,k,j)+\tau(n^*,k^*,j^*)=\tau(m,k+k^*,\ell)
$$
can hold. That will be very important in the rest of the study, to
understand the structure of the nonlinear terms in~(\ref{lim-filtered}).

\section{Interactions between equatorial waves}\label{Resonance}
In this section we will study the nonlinear term
in~(\ref{lim-filtered}). We will first study the resonances of~$ L$,
and then prove that the projection of~(\ref{lim-filtered}) onto the
kernel of~$ L$ is a linear equation.

\subsection{Study of the resonances}
Let us prove the following result.
\begin{Prop}\label{resonances-prop}{ 
Except for a countable number of $\beta$ and with the notation  of  Section 
\ref{precisedescription},  the following condition of 
non resonance
holds for all $n,n^*,m \in \N$, all~$k,k^*
\in \Z$ and all~$j,j^*,\ell \in \{-1,0,1\}$:
$$
\tau(n,k,j)+\tau(n^*,k^*,j^*)=\tau(m,k+k^*,\ell)
$$
implies
$$\begin{aligned}
\hbox{ either } \tau(n,k,j)=0 \hbox{ or }\tau(n^*,k^*,j^*)=0 \hbox{ or
}\tau(m,k+k^*,\ell)=0 ,
\\
\hbox{ or } \tau(n,k,j),\tau(n^*,k^*,j^*),\tau(m,k+k^*,\ell) \in {\mathfrak S}_K, 
\end{aligned}
$$
meaning that, among the ageostrophic modes, only  three Kelvin waves may interact.}
\end{Prop}
\begin{proof}
Let us start by noticing that by definition of Kelvin waves, Kelvin resonances necessarily take place simply
because they correspond to   convolution in Fourier space. 

Before starting with technical results, let us describe the main ideas of the proof.
The crucial argument is that the eigenvalues of the penalization 
operator $L$ are
defined as the roots of a countable number of polynomials whose 
coefficients depend
(linearly) on the ratio $\beta$. In particular, for fixed $n,n^*,m \in \N$ and
$k,k^*
\in \Z$, the occurence of a resonant triad
$$\tau(n,k,j)+\tau(n^*,k^*,j^*)=\tau(m,k+k^*,\ell)$$
is controlled by the cancellation of some polynomial~$P_{n,n^*,m,k,k^*} (\beta)$. 
Therefore, either
this polynomial has a finite number of zeros, or it is identically zero. The
difficulty here is that we are not able to eliminate the second 
possibility using
only the asymptotics $\beta \to \infty$. We therefore also
study the asymptotics~$\beta \to 0$, and in the case when~$n= 0$ or~$n^{*} = 0$, we have to refine 
the previous
argument introducing an auxiliary polynomial.

\bigskip
$\bullet$ {\underline{\it Definition of the polynomial $P_{n,n^*,m,k,k^*} (\beta)$} }

For fixed $n,n^*,m \in \N$ and
$k,k^*
\in \Z$, it is   natural to consider the following quantity
$$
P_{n,n^*,m,k,k^*} (\beta)= \prod_{j,j^*,\ell\in\{-1,0,1\}}
\big(\tau(n,k,j)+\tau(n^*,k^*,j^*)-\tau(m,k+k^*,\ell)\big)\,.
$$
Considerations of symmetry show that this quantity can be rewritten 
as a polynomial of the
symmetric functions of $(\tau(n,k,j))_{j\in \{-1,0,1\}}$, the 
symmetric functions of
$(\tau(n^*,k^*,j^*))_{j^*\in \{-1,0,1\}}$ and the symmetric functions of
$(\tau(m,k+k^*,\ell))_{\ell\in \{-1,0,1\}}
$.

Therefore, as the eigenvalues $(\tau(n,k,j))_{j\in \{-1,0,1\}}$ of the linear
penalization $L$ are defined as the three roots of a polynomial (\ref{tau-def1}) whose 
coefficients
depend (linearly) on $\beta$
$$
\tau^3 -(k^2+\beta(2n+1))\tau+\beta k=0,
$$
the symmetric functions of $(\tau(n,k,j))_{j\in \{-1,0,1\}}$ satisfy
\begin{equation}
\label{sym-functions}
\begin{aligned}
\sum_{j\in\{-1,0,1\} }\tau(n,k,j)=0 ,\\
\sum_{j\in\{-1,0,1\}} \prod_{j'\neq 
j}\tau(n,k,j')=-(k^2+(2n+1)\beta) ,\\
\prod_{j\in\{-1,0,1\}} \tau(n,k,j)=-\beta k,
\end{aligned}
\end{equation}
from which we deduce that $P_{n,n^*,m,k,k^*} (\beta)$ is a polynomial 
(of degree at
most 27) with respect to $\beta$.

In particular, for fixed $n,n^*,m \in \N$ and
$k,k^*
\in \Z$, either $P_{n,n^*,m,k,k^*} (\beta)$ is identically zero or it has a
finite number of roots. In other words, that means that

(a)
either, for all $\beta \in
\R^*$, there is a resonance of the type
$$\tau(n,k,j)+\tau(n^*,k^*,j^*)=\tau(m,k+k^*,\ell)$$
for some $j,j^*,\ell\in \{-1,0,1\}$,

(b) or, except for a finite number of $\beta$, such
resonances do not occur.

\bigskip
$\bullet$ {\underline{\it Asymptotic behaviour of $P_{n,n^*,m,k,k^*} (\beta)$ as $\beta \to
\infty$}}$ $

In order to discard one of these alternatives, we are   interested in the
asymptotic behaviour of the polynomial~$P_{n,n^*,m,k,k^*} (\beta)$ as $\beta \to
\infty$.

We start by describing the asymptotic behaviour of each root 
$(\tau(n,k,j))_{j\in
\{-1,0,1\}}$ as $\beta \to \infty$.

\begin{Lem}\label{DL-tau} { 
With the notation  of Paragraph \ref{precisedescription}, for all $k\in 
\Z$ and all
$n\in \N$, the following expansions hold as $\beta\to \infty$~:
\begin{equation}
\label{asympt-expansions}
\begin{aligned}
\tau(n,k,1)= \sqrt{(2n+1)\beta} -{k\over 2(2n+1)} +o(1),\\
\tau(n,k,-1)= -\sqrt{(2n+1)\beta} -{k\over 2(2n+1)} +o(1),\\
\tau(n,k,0)={k\over 2n+1}-{4n(n+1) k^3\over (2n+1)^4 \beta} +o\left(\frac1\beta\right) \cdotp
\end{aligned}
\end{equation}}
\end{Lem}

\begin{proof} 
We start with the most general case, namely the case when $k\neq 0$.  We proceed by successive
approximations. As the product of the roots $-\beta k$ tends to infinity as $\beta \to \infty$, there is at least
one root which tends to infinity. Therefore, we get at leading order
$$\tau^3-\beta(2n+1) \tau = 0,$$
which implies that the Poincar{\'e} and mixed Poincar{\'e}-Rossby modes are approximately given by
$$\tau(n,k,\pm 1) \sim\pm \sqrt{(2n+1)\beta}.$$
Plugging this Ansatz in the formula
$$\tau^2=(2n+1)\beta +k^2 -{\beta k\over \tau}=(2n+1)\beta \left( 1- {k\over (2n+1)\tau} +{k^2\over \beta(2n+1)}
\right)
$$
provides the next order approximation of the Poincar{\'e} modes, namely
$$\tau(n,k,j) \sim j \sqrt{(2n+1)\beta}-{k\over 2(2n+1)}\cdotp$$
Then, as the sum of the roots is zero (see (\ref{sym-functions})), we deduce that the third mode, i.e. the Kelvin or Rossby  mode,
satisfies
$$\tau(k,n,0)={k\over 2n+1}+o(1).$$
Plugging this Ansatz in the formula
$$ \tau ={\beta k+\tau^3 \over (2n+1)\beta +k^2}$$
leads then to
$$\tau(n,k,0) = {k\over 2n+1} \left( 1+{1\over \beta k} {k^3\over (2n+1)^3} -{k^2\over (2n+1)\beta}  \right) +o\left(
\frac1\beta\right).$$

\bigskip\bigskip
The other case (when $k= 0$) is dealt with in a very simple way. The Poincar{\'e} modes are
exactly
$$\tau(n,0,\pm 1) =\pm \sqrt{(2n+1)\beta},$$
 whereas the third mode is   zero
$$\tau(n,0,0)=0,$$
and thus they satisfy the general identities (\ref{asympt-expansions}).

The result is proved.
\end{proof}

\bigskip
Equipped with this technical lemma, we are now able to characterize the asymptotic behaviour of most of the factors
$$\tau(n,k,j)+\tau(n^*,k^*,j^*)-\tau(m,k+k^*,\ell)$$
in $P_{n,n^*,m,k,k^*}(\beta)$ as $\beta \to \infty$.

\begin{Lem}\label{triad} { 
With the notations of Paragraph \ref{precisedescription},  any triad of non zero modes
$$ \big(\tau(n,k,j),\tau(n^*,k^*,j^*),\tau(m,k+k^*,\ell)\big)$$
with  $k,k^*\in 
\Z$ and 
$n,n^*,m\in \N$, which is not constituted of three Kelvin or three Rossby modes, is asymptotically non
resonant as $\beta \to \infty$.

More precisely the following expansions hold as $\beta \to \infty$~:

(i) for three Poincar{\'e} or mixed Poincar{\'e}-Rossby modes ($j\neq 0$ and $j^*\neq 0$ and $\ell\neq 0$)
$$ \tau(n,k,j)+\tau(n^*,k^*,j^*)-\tau(m,k+k^*,\ell)\sim \sqrt{\beta} \left(
j\sqrt{2n+1}+j^*\sqrt{2n^*+1}-\ell\sqrt{2m+1}\right)\,;$$

(ii) for one Poincar{\'e}  or mixed Poincar{\'e}-Rossby mode and two Rossby or Kelvin or zero modes
$$\tau(n,k,j)+\tau(n^*,k^*,j^*)-\tau(m,k+k^*,\ell) \sim \sqrt{\beta} \left(
j\sqrt{2n+1}+j^*\sqrt{2n^*+1}-\ell\sqrt{2m+1}\right)\,;$$

(iii) for two Poincar{\'e} or mixed Poincar{\'e}-Rossby modes and one Rossby or Kelvin mode
$$\exists C\equiv C(n,n^*,m,k,k^*) >0,\quad \big | \tau(n,k,j)+\tau(n^*,k^*,j^*)-\tau(m,k+k^*,\ell) \big | \geq
C \,;$$

(iv) for two Kelvin modes and one Rossby mode
$$\exists C\equiv C(n,n^*,m,k,k^*) >0,\quad \big | \tau(n,k,j)+\tau(n^*,k^*,j^*)-\tau(m,k+k^*,\ell) \big | \geq
\frac C\beta \,;$$

(v) for two Rossby modes and one Kelvin mode
$$\exists C\equiv C(n,n^*,m,k,k^*) >0,\quad \big | \tau(n,k,j)+\tau(n^*,k^*,j^*)-\tau(m,k+k^*,\ell) \big | \geq
\frac C\beta  \cdotp$$
}\end{Lem}

\begin{proof} 
The proof of these results is based on Lemma \ref{DL-tau}.

(i)
In the case of three Poincar{\'e}  or mixed Poincar{\'e}-Rossby modes, Lemma
\ref{DL-tau} provides
$$\tau(n,k,j)+\tau(n^*,k^*,j^*)-\tau(m,k+k^*,\ell) =\sqrt{\beta} \left(
j\sqrt{2n+1}+j^*\sqrt{2n^*+1}-\ell\sqrt{2m+1}\right)+o(\sqrt{\beta}),$$
and it is easy to check, using considerations of parity, that the constant $$\left(
j\sqrt{2n+1}+j^*\sqrt{2n^*+1}-\ell\sqrt{2m+1}\right)$$ cannot be zero.

(ii)
In the case of one Poincar{\'e}  or mixed Poincar{\'e}-Rossby mode, we have one term which is exactly of order $\sqrt{\beta}$ whereas
the others are negligible compared with $\sqrt{\beta}$, thus the sum is equivalent to the Poincar{\'e} mode, and  the
same formula holds
$$\tau(n,k,j)+\tau(n^*,k^*,j^*)-\tau(m,k+k^*,\ell) =\sqrt{\beta} \left(
j\sqrt{2n+1}+j^*\sqrt{2n^*+1}-\ell\sqrt{2m+1}\right)+o(\sqrt{\beta}).
$$

\bigskip

(iii) The third case is a bit more difficult to deal with, since the leading order terms can cancel each other out.
Without loss of generality, we can assume that $\ell=0$ and $j,j^*\neq 0$ (the other cases being obtained by
exchanging $j,j^*$ and $-\ell$).

 If $j=j^*$, or if $j+j^*=0$ and $n\neq n^*$, the same arguments as previously show that
the same formula holds
$$\tau(n,k,j)+\tau(n^*,k^*,j^*)-\tau(m,k+k^*,\ell) \sim\sqrt{\beta} \left(
j\sqrt{2n+1}+j^*\sqrt{2n^*+1}-\ell\sqrt{2m+1}\right),
$$
since the factor of $\sqrt{\beta}$ is not zero.

If $j+j^*=0$ and $n=n^*$, the factor of $\sqrt{\beta}$ cancels and we have to determine the next term in the
asymptotic expansion~:
$$\tau(n,k,j)+\tau(n^*,k^*,j^*)-\tau(m,k+k^*,\ell)= -{k+k^*\over 2(2n+1)}-{k+k^* \over 2m+1}+o(1) .$$
 Considerations of parity show therefore that the limit cannot be
zero if
$k+k^*\neq 0$, or equivalently if
$\tau(m,k+k^*,0)\neq 0$.

\bigskip
(iv) In the case of one Rossby and two Kelvin modes, we are not able in general to prove that the leading order term, i.e. the
limit as $\beta \to \infty$ of $\tau(n,k,j)+\tau(n^*,k^*,j^*)-\tau(m,k+k^*,\ell)$ is not zero. But we can look directly at the
second term of the expansion, i.e. the factor of $\beta^{-1}$~: 
$$ \omega_1
=-{4k^3n(n+1)\over (2n+1)^4}-{4(k^*)^3n^*(n^*+1)\over (2n^*+1)^4}+ {4(k+k^*)^3m(m+1)\over (2m+1)^4}\cdotp
$$
Considering one Rossby and two Kelvin modes means that $k,k^*$ and $k+k^* $ are not zero, and that exactly two indices among
$n$, 
$n^*$ and $m$ are zero. Thus  $\omega_1\neq 0$ and
$$ \big | \tau(n,k,j)+\tau(n^*,k^*,j^*)-\tau(m,k+k^*,\ell) \big | \geq
\frac{ |\omega_1|}{2\beta }$$
for $\beta$ large enough.

\bigskip
(v) The last situation is the most difficult to deal with, since the only thing we will be able to prove is that the two first
terms of the asymptotic expansion of $\tau(n,k,j)+\tau(n^*,k^*,j^*)-\tau(m,k+k^*,\ell)$ with respect to $\beta$ cannot cancel
together.
By Lemma \ref{DL-tau}, we deduce that for one Kelvin and two Rossby modes
$$\tau(n,k,j)+\tau(n^*,k^*,j^*)-\tau(m,k+k^*,\ell)=\omega_0 +{\omega_1 \over \beta} +o\left({1\over \beta}\right)$$
with 
$$\omega_0={k\over 2n+1}+{k^*\over 2n^*+1}-{k+k^* \over 2m+1}\virgp$$
and 
$$\omega_1= -{4k^3n(n+1)\over (2n+1)^4}-{4(k^*)^3n^*(n^*+1)\over (2n^*+1)^4}+ {4(k+k^*)^3m(m+1)\over (2m+1)^4}
\cdotp$$
Recall moreover that $k,k^*$ and $k+k^* $ are not zero, and that exactly one index among
$n$, 
$n^*$ and $m$ is zero.
Without loss of generality, we can assume that $m=0$ and $n,n^*\neq 0$ (the other cases being obtained by
exchanging $n,n^*$ and $m$).

Then, if $\omega_0=\omega_1=0$,
$$
\begin{aligned}
\,&{kn\over 2n+1 }+{k^* n^* \over 2n^*+1}=0,\\
&{k^3n(n+1)\over (2n+1)^4}+{(k^*)^3n^*(n^*+1)\over (2n^*+1)^4}=0,
\end{aligned}$$
from which we deduce that
$$ {(n+1) \over n^2 (2n+1)}- {(n^*+1) \over (n^*)^2 (2n^*+1)}=0.$$
Therefore, as the function
$$x\mapsto {x+1\over x^2(2x+1)}$$
decreases strictly  on $\R^+$, we get 
$n=n^*$ and thus $k=-k^*$, which contradicts the fact that $k+k^*\neq 0$.

We conclude that either $\omega_0\neq 0$ or $\omega_1\neq 0$, so that
$$ \big | \tau(n,k,j)+\tau(n^*,k^*,j^*)-\tau(m,k+k^*,\ell) \big | \geq
\frac{ |\omega_1|}{2\beta }$$
for $\beta$ large enough.
Lemma~\ref{triad} is proved.
\end{proof}

\bigskip
Let us go back to the proof of Proposition \ref{resonances-prop}, and  first consider the case when  $k\neq 0$, $k^*\neq 0$ and $k+k^*\neq 0$.
In view of Lemma \ref{triad},   the asymptotic behaviour of $P_{n,n^*,m,k,k^*}
(\beta)$ as $\beta \to \infty$ is completely determined by the behaviour of the factor 
$$\sigma_{n,n^*,m,k,k^*}(\beta)=\tau(n,k,0)+\tau(n^*,k^*,0)-\tau(m,k+k^*,0).$$

Indeed, $P_{n,n^*,m,k,k^*} (\beta)$ is defined as a product, eight factors of which involve triads of type (i),
six of which involve triads of type (ii), twelve of which involve triads of type (iii)   and the last factor of which
is $\sigma_{k,k^*,n,n^*,m} (\beta)$. By Lemma \ref{triad} we then deduce that there exists a nonnegative
constant~$C$ (depending on $k,k^*,n,n^*,m$) such that
$$\big|P_{n,n^*,m,k,k^*} (\beta)\big|\geq C \beta^7 \big|\sigma_{n,n^*,m,k,k^*} (\beta)\big|.$$

If one or two among $n$, $n^*$ and $m$ are zero, properties (iv) and (v) in Lemma \ref{triad} allow to conclude that for $\beta$
large enough
$$\big|P_{n,n^*,m,k,k^*} (\beta)\big|\geq C \beta^{6},$$
and thus $P_{,n,n^*,m,k,k^*}$ has a finite number of roots.

If $n,n^*,m$ are all equal to zero or $n,n^*,m $ are all non zero, we cannot conclude as no estimate 
on $\sigma_{n,n^*,m,k,k^*} (\beta)$ at
infinity is available. Therefore, either $\sigma_{n,n^*,m,k,k^*} (\beta)$ is identically zero for $\beta$ large enough, or
$P_{n,n^*,m,k,k^*} (\beta)$ has a finite number of roots.

\bigskip
Thus at this stage, in order to prove Proposition  \ref{resonances-prop}, it remains

$(1) \quad $ to consider the case when $k,k^*,k+k^*\neq 0$ and $\sigma_{n,n^*,m,k,k^*} (\beta)$
is identically zero for~$\beta $ large enough, with~$n,n^*,m$  all zero
or all non zero;
 
$(2) \quad $ to study the case when $k$ or $k^*$ or $k+k^*$ is zero (in order to establish that only the triads
involving a zero mode may be resonant).

\bigskip
$\bullet$ {\underline{\it Conclusion in the case of~(1)}}$ $

In the case when~$n,n^*,m $ are all   zero, then the resonances corresponding to~$\sigma_{0,0,0,k,k^*}
 (\beta) = 0$
are precisely Kelvin resonances, which cannot be
removed. 

In the case when~$n,n^*,m $ are all non zero, then
  $\sigma_{n,n^*,m,k,k^*} (\beta)$ is an
analytic function of $\beta$ (the roots of (\ref{tau-def1})  -- defined explicitely with Cardan's formula -- do not
cross each other according to Lemma~\ref{polynomial}, 
and thus depend analytically on~$\beta$)~: in particular,
 if $\sigma_{n,n^*,m,k,k^*} (\beta)$ 
cancels for
$\beta$ large enough, then it is identically zero. 
 Let us describe the asymptotics of the roots as~$\beta$ goes to zero.
\begin{Lem}\label{DL-taubetazero} 
With the notation  of Paragraph \ref{precisedescription}, for all $k\in 
\Z$ and all
$n\in \N^*$, the following expansions hold as $\beta\to 0$~:
$$\tau(n,k,0)=  \frac\beta k +o(1) .$$
\end{Lem}
\begin{proof}
Since the product of the roots goes to zero as~$\beta$ goes to  zero, we infer that at least one root goes to  
zero with~$\beta$. Let us consider that root.
Since~$\beta (2n+1)$ is negligible with respect to~$k^{2}$ and~$\tau^{3}$ is negligible with respect
 to~$k^{2} \tau$, we find that 
$$
 k^{2}\tau  -   k \beta \sim 0,
$$
so that  one root is equivalent to~$\frac \beta k$ as~$\beta$ goes to zero.
 It is easy to see that  the two other
roots are then equivalent to~$\pm k$, so that we do have~$\tau (n,k,0) \sim \frac \beta k $ (we recall that 
for~$n \neq 0$, the
roots are numbered in increasing order).  The lemma is proved.
\end{proof}
Now going back to the study of case~$(1)$, in view of Lemma~\ref{DL-taubetazero}  it is obvious 
that~$\sigma_{n,n^*,m,k,k^*} (\beta)$ cannot vanish indentically.

\bigskip
$\bullet$ {\underline{\it Conclusion in the case of~(2)}}$ $

In this situation, we need to   refine the previous analysis by introducing an auxiliary polynomial.  We thus define
$$I_{k,k^*}=\left\{(j,j^*,\ell)\in \{-1,0,1\}^3 \,/\,\ell\neq 0 \hbox{ if }k+k^*=0,\quad  j\neq 0 \hbox{ if }k=0\hbox{ and } j^*\neq
0\hbox{ if } k^*=0\right\}$$
and
$$Q_{n,n^*,m,k,k^*}
 (\beta) = 
\prod_{(j,j^*,\ell)\in I_{k,k^*}}
\big(\tau(n,k,j)+\tau(n^*,k^*,j^*)-\tau(m,k+k^*,\ell')\big).
$$
That corresponds to the remaining possible resonances, where we have omitted the trivial case when 
one wave is geostrophic ($\tau = 0$).

As previously, considerations of symmetry show that this quantity can be rewritten 
in terms of the
symmetric functions of $(\tau(n,k,j))_{j\in \{-1,0,1\}}$ (or $(\tau(n,k,j))_{j\in \{-1,1\}}$ if $k=0$), the 
symmetric functions of
$(\tau(n^*,k^*,j^*))_{j^*\in \{-1,0,1\}}$ (or $(\tau(n^*,k^*,j^*))_{j^*\in \{-1,1\}}$ if $k^*=0$) and the
symmetric functions of
$(\tau(m,k+k^*,\ell))_{\ell\in \{-1,1\}}
$ (or $(\tau(m,k+k^*,\ell))_{\ell\in \{-1,0,1\}}
$ if $k+k^*=0$). Noticing that the symmetric functions of~$(\tau(n, 0,j))_{j\in \{-1,1\}}$ are affine
 in~$\beta$, we conclude that~$Q_{n,n^*,m,k,k^*}(\beta)$  is a polynomial in $\beta$.

The asymptotic analysis of the various factors as $\beta \to \infty$ shows that
$$\left|Q_{n,n^*,m,k,k^*}(\beta)\right| \geq C \beta^3$$
for $\beta$ large enough. Therefore, $Q_{n,n^*,m,k,k^*} (\beta)$ has a finite number of roots, meaning that there exist a finite number of
$\beta$ such that resonant triads with $k=0$ or $k^*=0$ or $ k+k^*=0$ (other than the triads involving a
non-oscillating mode) can occur.

\bigskip
We have therefore proved that

$(1) \quad $ in the case when $k,k^*,k+k^*\neq 0$ and $\sigma_{n,n^*,m,k,k^*} (\beta)$
is identically zero, only  the triads involving three Kelvin
modes are resonant for an infinite  number of $\beta$;
 
$(2) \quad $  when $k$ or $k^*$ or $k+k^*$ is zero, 
only triads involving zero modes are resonant for an infinite  number of $\beta$.

Combining this result with the conclusion of the previous paragraph achieves the proof
 of Proposition~\ref{resonances-prop}.
\end{proof}

\subsection{The special case of~$ \Ker L$}
\label{algebraic section}
In this short section we are going to write an algebraic computation
which in particular allows to derive the following proposition. 
\begin{Prop}\label{limitgeolinear}
{ 
Let~$ \Phi$ and~$ \Phi^*$ be two smooth vector fields. Then for every~$n \in \N$, we have
$$
\left(\Psi_{n,0,0} | Q_L ( \Phi, \Phi^*) \right)_{L^2(\R \times \T)}= 0.
$$
}
\end{Prop}
\begin{Rem}
\label{remgeolinear}
That proposition implies that the projection of the limit
system~(\ref{lim-filtered}) onto~$ \Ker L$ can be formally written
$$
\partial_t \Pi_0 \Phi - \nu \Delta'_L \Pi_0 \Phi = 0.
$$
\end{Rem}
\begin{proof} We are going to prove a more general result, computing
  the quantity
$$
\Bigl(  \Phi_\lambda  | Q_L(\Phi_\mu  ,\Phi_{\tilde\mu}) \Bigr)_{L^2(\R \times \T)}
$$
where~$  \Phi_\lambda$, $  \Phi_\mu$ and $  \Phi_{\tilde\mu}$ are three
eigenmodes of $L$ associated respectively
with the eigenvalues $i\lambda$, $i\mu$ and $i\tilde\mu$ where
$\lambda=\mu+\tilde\mu$. 
The proposition corresponds of course to the case when~$ \lambda = 0$.

We have
$$
\Bigl(  \Phi_\lambda  | Q_L(\Phi_\mu  ,\Phi_{\tilde\mu}) \Bigr)_{L^2(\R \times \T)}
=\Bigl(  \Phi_\lambda  | Q(\Phi_\mu  ,\Phi_{\tilde\mu}) \Bigr)_{L^2(\R \times \T)}
$$
hence denoting by~$\bar  \Phi_\lambda$ \label{ccPhi} the complex conjugate of~$ \Phi_\lambda$, we get
$$
\begin{aligned}
\Bigl(  \Phi_\lambda  & | Q_L(\Phi_\mu  ,\Phi_{\tilde\mu}) \Bigr)_{L^2(\R \times \T)}
\\
&=\frac12\int \left(\bar \Phi_{\lambda,0} \nabla \cdot (\Phi_{\mu,0} \Phi_{\tilde\mu}'+ \Phi_{\tilde\mu, 0} \Phi_{\mu}')+ \bar
\Phi_{\lambda}'\cdot (\Phi_{\mu}' \cdot \nabla
\Phi_{\tilde\mu}' +\Phi_{\tilde\mu}' \cdot \nabla \Phi_{\mu}') \right) dx\\
&=\frac12 \int \left(-\nabla\bar \Phi_{\lambda,0}  \cdot (\Phi_{\mu,0} \Phi_{\tilde\mu}'+\Phi_{\tilde\mu, 0}  \Phi_{\mu}')+
\bar \Phi_{\lambda}'\cdot (\nabla(\Phi_{\mu}'\cdot
\Phi_{\tilde\mu}') +\Phi_{\mu}'^\perp\nabla^\perp \cdot \Phi_{ \tilde \mu}'+\Phi_{\tilde\mu}'^\perp \nabla^\perp \cdot \Phi_{ \mu}')
\right) dx\\
&=-\frac12\int \left((\beta x_1 \bar \Phi_{\lambda}'^\perp+\nabla\bar \Phi_{\lambda,0} ) \cdot
(\Phi_{\mu,0} \Phi_{\tilde\mu}'+ \Phi_{\tilde\mu,0} \Phi_{\mu}')+(\nabla
\cdot \bar \Phi_{\lambda}') \Phi_{\mu}' \cdot \Phi_{\tilde\mu}' \right) dx\\
& \quad\quad\quad \quad\quad\quad \quad+\frac12\int \bar \Phi_{\lambda}' \cdot (\Phi_{\mu}'^\perp (
\nabla^\perp \cdot \Phi_{\tilde \mu}'-\beta
x_1\Phi_{\tilde\mu,0})+\Phi_{\tilde\mu}'^\perp(\nabla^\perp \cdot \Phi_{ \mu}' -\beta x_1 \Phi_{\mu,0})) dx.
\end{aligned}
$$
Using the identities
$$
\begin{aligned}
 \nabla \cdot \Phi_{\lambda}' &=i\lambda \Phi_{\lambda,0},\\
\beta x_1 \Phi_{\lambda}'^\perp +\nabla \Phi_{\lambda,0} &=i\lambda \Phi_{\lambda}',
\end{aligned}
$$
as well as their combination
$$
\beta \Phi_{\lambda,1} =i\lambda ( \nabla^\perp \cdot \Phi_\lambda' -\beta x_1
\Phi_{\lambda,0})
$$
and  similar formulas for $\Phi_{\mu}$ and $\Phi_{\tilde \mu}$, we get
$$\begin{aligned}
&\Bigl(  \Phi_\lambda   | Q_L(\Phi_\mu  ,\Phi_{\tilde\mu}) \Bigr)_{L^2(\R \times \T)}
\\
&=\frac12\int \left(i\lambda \bar \Phi_{\lambda}' \cdot (\Phi_{\mu,0} \Phi_{\tilde\mu}'+ \Phi_{\tilde \mu, 0}
\Phi_{\mu}')+i\lambda \bar \Phi_{\lambda,0} \Phi_{\mu}' \cdot
\Phi_{\lambda}'
\right) dx\\
&-\frac1{2\beta}\int \bar \Phi_{\lambda,2} (i\mu + i \tilde \mu) (\nabla^\perp \cdot \Phi_{  \mu}' -\beta x_1 \Phi_{\mu,0})
(\nabla^\perp \cdot \Phi_{\tilde \mu}'-\beta
x_1\Phi_{\tilde\mu,0})  dx\\
&-\frac1{2\beta}\int i\lambda (\nabla^\perp \cdot \bar \Phi_{\lambda}'-\beta x_1 \bar\Phi_{\lambda,0})
(\Phi_{\mu,2} (\nabla^\perp \cdot \Phi_{\tilde \mu}'-\beta
x_1\Phi_{\tilde\mu,0})+\Phi_{\tilde \mu,2}(\nabla^\perp \cdot \Phi_{  \mu}' -\beta x_1 \Phi_{\mu,0})) dx
\end{aligned}
$$
from which we deduce
\begin{equation}
\label{quasigeostrophic}
\begin{aligned}
\Bigl(  \Phi_\lambda  & | Q_L(\Phi_\mu  ,\Phi_{\tilde\mu}) \Bigr)_{L^2(\R \times \T)}\\
 &  =\frac{i\lambda}2\int \left( \bar \Phi_{\lambda}' \cdot (\Phi_{\mu,0} \Phi_{\tilde\mu}'+\Phi_{\tilde\mu,0} \Phi_{\mu}')+\bar
\Phi_{\lambda,0} \Phi_{\mu}' \cdot
\Phi_{\lambda}'
\right) dx\\
&-\frac{i\lambda}{2\beta}\int \bar \Phi_{\lambda,2}  (\nabla^\perp \cdot \Phi_{  \mu}' -\beta x_1 \Phi_{\mu,0})
(\nabla^\perp \cdot \Phi_{\tilde \mu}'-\beta
x_1\Phi_{\tilde\mu,0})dx\\
&-\frac{i\lambda}{2\beta}\int  (\nabla^\perp \cdot \bar \Phi_{\lambda}' -\beta x_1 \bar\Phi_{\lambda,0})
(\Phi_{\mu,2} (\nabla^\perp \cdot \Phi_{\tilde \mu}'-\beta
x_1\Phi_{\tilde\mu,0})+\Phi_{\tilde \mu,2}(\nabla^\perp \cdot \Phi_{  \mu}'-\beta x_1 \Phi_{\mu,0})) dx.
\end{aligned}
\end{equation}
In particular for $\lambda=0$ this quantity is always zero, which proves Proposition~\ref{limitgeolinear}.
\end{proof}


\chapter{The envelope equations}\label{envelope}
\setcounter{equation}{0}

The aim of this  chapter is to study the system $(SW_0)$
obtained formally
page~\pageref{lim-filtered} as the limit of
the filtered system (\ref{filtered}) as $\eps \to 0$. 
Let us recall the system:
$$
\begin{aligned}
& (SW_0) \quad \quad 
\left\{
\begin{array}{r}
 \d_t \Phi +Q_L (\Phi,\Phi) -\nu \Delta'_L \Phi
=0\\
\displaystyle  \Phi_{|t=0}=(\eta^0,u^0),
\end{array}
\right.
\end{aligned}
$$
where $\Delta'_L$ and $Q_L$ denote  the linear and symmetric bilinear operator
defined by~(\ref{Q-lap-L}) page~\pageref{Q-lap-L}. 

Two different types of wellposedness results will be proved on~$(SW_{0})$: first we  will
  prove the existence of  weak solutions in~$L^{2} $ and of a unique, strong solution  if
  the data is smooth enough (on a short time interval, which becomes infinite
  for small data). Then 
  we will show that except for a countable number of~$\beta$, the strong solutions exists globally in 
  time as soon as the initial data is only in~$L^{2}$, of arbitrary norm.
  
  The statements of both theorems can be found in Paragraph~\ref{statementenveloppe}, and their
proofs are respectively the object of Paragraphs~\ref{theoremallbeta} and~\ref{theoremgeneric}.
    In order to
establish those results we will  need to define, in Paragraph~\ref{sectionsuitable}, suitable function spaces,
compatible with the penalization  operator~$L$ as well as the
diffusion operator. Some  technical preliminaries devoted to those spaces are proved in Paragraph~\ref{preliminaries}:
in particular
in Paragraph~\ref{sectiontrilinear} we prove the  continuity
of the bilinear operator~$Q_{L}$ in those function spaces.  Finally the last part of this chapter is devoted
to an additional smoothing property on the divergence.

\section{Definition of suitable functional  spaces}\label{sectionsuitable}
\setcounter{equation}{0}
By construction the operators $\Delta'_L$ and $Q_L$ appearing in the
limiting filtered system $(SW_0)$
are defined in terms of the projections
$(\Pi_\lambda)_{i\lambda \in {\mathfrak S}}$ on the eigenspaces of
$L$. In particular, they are not expected to satisfy ``good"
commutation
properties with the usual derivation $\nabla$.
Therefore in order to establish a priori estimates on the solutions
to $(SW_0)$ we have to introduce some weighted Sobolev
spaces
 associated with some derivation-like operator
which acts separately on each eigenmode of $L$.
 
 Let us therefore introduce the following norms. We will
 write as previously~$\Pi_{n,k,j}$ for  the projection  on the
eigenmode $\Psi_{n,k,j}$ of $L$
 and~$\Pi_\lambda$  for the projection on the
eigenspace associated with the eigenvalue~$i\lambda$ of $L$. Finally we
define\label{S-not}
$$
S = \N \times \Z \times \{-1,0,1\}.
$$
\begin{Def} \label{HsL}
{
Let~$s  \geq 0$ be a given real number. We define the space~$H^{s}_{L}$
 as the subspace of~$(L^{2}(\R \times \T))^{3} $ given by the
 following norm:
 $$
\|\Phi\|_{H^{s}_{L}} \eqdefa  \left(
\sum_{(n,k,j) \in S}(1 + n + k^2)^{s} \|\Pi_{n,k,j} \Phi\|_{L^{2}(\R  \times \T)}^{2}
\right)^{\frac12}.
 $$
}
\end{Def}
Due to the definition of the eigenvectors of~$L$ seen in the previous
chapter, one can prove the following proposition.
\begin{Prop}\label{equivnorm}
{
Let~$s \geq 0$ be given. Then one has
the following property:
$$
\forall \Phi \in  H^{s}_{L}, \quad \|\Phi\|_{ H^{s}_{L} }   \sim
   \|(\mbox{Id}-\Delta + \beta^{2} x_{1}^{2}) ^{s/2} \Phi\|_{L^{2}(\R  \times \T)}.
$$
In particular, $H^{s}_{L}$ is continuously embedded in~$H^{s}(\R \times  \T)$, and for
all compact subsets $\Omega$ of $\R \times \T$,  $H^s_0(\Omega)$ is continuously embedded
in~$H^{s}_L$.

Moreover for all~$\Phi \in H^s_L \cap (\Ker
  L)^\perp$, we have
$$
\|\Phi\|_{ H^s_L} \sim   \left(
\sum_{i \lambda \in {\mathfrak S} \setminus \{0\}}
 \|\Pi_{\lambda} \Phi\|_{H^{s}(\R \times \T)}^{2}
\right)^{\frac12} .
$$
Finally if~$\Phi$ belongs to~$K \cup P$, as defined in Definition~\ref{PRKN} page~\pageref{PRKN}, then
$$
\|\Phi\|_{ H^s_L} \sim   \left(
\sum_{i \lambda \in {\mathfrak S}\setminus \{0\}}
 (1 + \lambda^2)^s\|\Pi_{\lambda} \Phi\|_{L^2(\R \times \T)}^{2}
\right)^{\frac12} .
$$
 }
\end{Prop}
\begin{proof}
Let us first prove the first equivalence: let~$\Phi \in H^{s}_{L}$ be  given. Then we have
$$
\| (\mbox{Id}-\Delta + \beta^{2} x_{1}^{2}) ^{s/2}\Phi\|_{L^2}^2 =  \Bigl\| \sum_{(n,k,j)
\in S}
(\mbox{Id}-\Delta + \beta^{2} x_{1}^{2}) ^{s/2} \Pi_{n,k,j}  \Phi\Bigr\|_{L^{2}}^2.
$$ By the identity
$$
-\psi_n''+\beta x_1^2 \psi_n=\beta(2n+1) \psi_n
$$
the orthogonality of the family $(\psi_n)_{n\in \N}$ and the explicit
formulas (\ref{psi-def1}), (\ref{psi-def2}) and (\ref{psi-def3})  for~$\Psi_{n,k,j}$, we
  infer that, for all integers~$\sigma$,
  $$ \| (\mbox{Id}-\Delta + \beta^{2} x_{1}^{2}) ^{\sigma}  \Psi_{n,k,j}-(1+n+k^2)^\sigma
\Psi_{n,k,j}\|_{L^2} \leq C
(1+n+k^2)^{\sigma-1} ,$$
which implies in particular that
 $$
\Bigl|\Bigl( (\mbox{Id} \: - \Delta + \beta
x_1^2)^\sigma \Psi_{n,k,j} |  (\mbox{Id} \: - \Delta + \beta
x_1^2)^\sigma \Psi_{n,k,j^*}  \Bigr)_{L^2(\R\times \T)}\Bigr| \leq C  (1+n+k^2)^{\sigma-1}.
$$
On the other hand,
  $$
\Bigl( (\mbox{Id} \: - \Delta + \beta
x_1^2)^\sigma \Psi_{n,k,j} |  (\mbox{Id} \: - \Delta + \beta
x_1^2)^\sigma \Psi_{n^*,k^*,j^*}  \Bigr)_{L^2(\R\times \T)} = 0\hbox{  if }
n \neq n^*\hbox{ or } k\neq k^*,
$$
so we find   that
$$\begin{aligned}
\|(\mbox{Id} \: - \Delta + \beta x_1^2)^{\sigma} \Phi\|_{L^2(\R\times  \T)}^2&=\sum_{n,k,j}
\sum_{ j^*} \left((\mbox{Id}
\: - \Delta + \beta x_1^2)^{\sigma}
 \Pi_{n,k,j} \Phi |  (\mbox{Id} \: - \Delta + \beta  x_1^2)^{\sigma}\Pi_{n,k,j^*} 
\Phi\right)_{ L^2}\\
&\sim \sum_{n,k,j} (1+n+k^2)^{2\sigma} \| \Pi_{n,k,j}   \Phi
\|_{ L^2(\R\times \T)}^2.
\end{aligned}
$$
We then obtain the first equivalence for all $s\geq 0$ by interpolation.

Then, from the   inequality
$$\forall \Phi \in H^s_L,\quad \| \Phi\|_{H^s(\R\times \T)} \leq C\|  (\mbox{Id} \: -
\Delta + \beta x_1^2)^{s/2}
\Phi\|_{L^2(\R\times \T)}$$
along with the fact that
and for all~$\Phi \in C^\infty(\R\times \T)$ supported in~$[-R,R]\times\T $,
$$  \|
(\mbox{Id} \: - \Delta +
\beta x_1^2)^{s/2}
\Phi\|_{L^2(\R\times \T)}\leq C(1+R^2)^{s/2} \|
\Phi\|_{H^s(\R\times \T)}$$
we get the embeddings $H^s_0(\Omega) \subset H^s_L \subset H^s(\R\times  \T)$ for all
$\Omega \subset \subset
\R\times \T$.

\bigskip
The second result of the proposition is easy,  using
  Proposition~\ref{orthogonality} page~\pageref{orthogonality}:
$$\forall i\lambda\in {\mathfrak S}\setminus \{0\},\quad {1\over C_s} 
\sum_{\tau(n,k,j)=\lambda} \| \Pi_{n,k,j}
\Phi \|^2_{H^s(\R\times \T)}\leq
 \| \Pi_\lambda
\Phi \|^2_{H^s(\R\times \T)} \leq C_s  \sum_{\tau(n,k,j)=\lambda} \|  \Pi_{n,k,j}
\Phi \|^2_{H^s(\R\times \T)},$$
and recalling that by  Proposition \ref{diag-prop}
 page~\pageref{diag-prop}, we have
$${1\over C_s} (1+ n+k^2)^{s/2} \leq  \| \Psi_{n,k,j}
\|_{H^s(\R\times \T)} \leq C_s (1+ n+k^2)^{s/2}.$$

Finally the last result, concerning Kelvin  and Poincar\'e   modes is simply due to Lemma~\ref{polynomial} and Proposition~\ref{diag-prop} .

The proposition is proved.
\end{proof}

\begin{Rem}\label{decroissance}

The $H^s_L$ estimates are both regularity and decay  estimates. In particular,
the embedding
$H^s_L\subset L^2(\R\times \T)$ is compact, and we have the following  equality
$$\bigcap_{s\geq 0} H^s_L={\mathcal S}(\R\times \T).$$

Note that these spaces are also used by Dutrifoy and Majda \cite{dutrifoymajda} to study the uniform wellposedness of a non viscous version of $(SW_\e)$.
\end{Rem}

\section{Statement of the   wellposedness  result}\setcounter{equation}{0}\label{statementenveloppe}
The main results of this chapter are the following two theorems. We have
written~$\Pi_{\perp} \Phi$ for  the projection
of~$\Phi$ onto~$(\Ker L)^{\perp}$\label{pip-not}.
In the next theorem, we state the global existence of weak solutions and the local in time existence (and uniqueness) of
strong solutions.
\begin{Thm}[Wellposedness results for all~$\beta$]\label{limit-systemallbeta} 
{There is a constant~$C$ such that the following results hold.
Let $\Phi^0\in L^2(\R\times \T;\R^3)$ be given. Then

$\bullet$ there exists a  global
weak solution $\Phi\in L^\infty (\R^+;L^2(\R\times \T))$
to $(SW_0)$, such that~$\Pi_{\perp} \Phi$ belongs to the  space~$L^2 
(\R^{+};H^{1}_{L})$, and  which  satisfies for   every~$t\geq 0$ the  energy
estimate
$$
\frac12\| \Phi(t)\|_{L^2(\R\times \T)}^2 +\nu\int_0^t \|\nabla (\Pi_0  \Phi)'
(t')\|_{L^2(\R\times \T)}^2dt' +{\nu\over C}\int_0^t \|\nabla  (\Pi_\perp \Phi)
(t')\|_{L^2(\R\times \T) }^2 dt'
 \leq \frac12 \|\Phi^0\|_{L^2(\R\times \T)}^2.
$$

$\bullet$ if we further assume that~$\Pi_{0} \Phi^{0} $ belongs to~$H^{s}_{L}$ for~$s\geq 0$, then~$\Pi_{0}
 \Phi$ (which is unique)
belongs to~$L^{\infty}_{loc}(\R^{+};H^{s}_{L})$.

$\bullet$ if~$\Pi_{0} \Phi^{0} $ belongs to~$L^{2}(\R\times \T)$ and~$\Pi_{\perp} \Phi^0$ belongs to~$  H^{1/2}_{L}$,
then there exists a maximal time interval~$[0,T^*[$, with~$T^*=+\infty$ under the smallness  assumption
$$
\| \Pi_0 \Phi^0\|_{L^2(\R\times  \T)} +\|\Pi_\perp\Phi^0\|_{   H^{1/2}_L
}\leq  C^{-1}\nu,
$$
 such that
$\Phi$ is the unique   solution to $(SW_0)$,
and~$ \Pi_\perp \Phi$ belongs to~$
L^\infty_{loc}([0,T^*[,H^{1/2}_{L}) \cap  L^2_{loc}([0,T^*[,H^{3/2}_{L})$.

$\bullet$  if~$\Pi_{\perp} \Phi^0$ belongs to~$ H^{s}_{L}$    for
some~$ 1/2 \leq  s \leq 1$,
then~$\Pi_{\perp} \Phi$ belongs to~$
L^\infty_{loc}([0,T^*[,H^{s}_{L}) \cap  L^2_{loc}([0,T^*[,H^{s+1}_{L})$.
 }
 \end{Thm}
 The previous theorem is much improved if a countable set of values for~$\beta$ is removed.
 
\begin{Thm}[Wellposedness results for generic~$\beta$]\label{limit-systemgeneric}
{There is a constant~$C$ and a countable subset~${\mathcal N}$ of~$\R^{+}$ such that for any~$\beta \in \R^{+} \setminus {\mathcal N}$, 
the following result holds. 
Let $\Phi^0\in L^2(\R\times \T;\R^3)$ be given. Then~$(SW_0)$ is globally wellposed, in the 
sense that there is a unique, global solution~$\Phi$ in~$L^\infty (\R^+;L^2(\R\times \T))$ such that~$\Pi_{\perp} \Phi$ belongs to the  space~$L^2 
(\R^{+};H^{1}_{L})$, and
which satisfies the energy inequality of Theorem~\ref{limit-systemallbeta}.

$\bullet$ if we further assume that~$\Pi_{\perp} \Phi^0$ belongs to~$  H^{s}_{L}$, for~$0 \leq s \leq 1$,  
then~$ \Pi_\perp \Phi$ belongs to~$
L^\infty_{loc}(\R^{+},H^{s}_{L}) \cap  L^2_{loc}(\R^{+},H^{s+1}_{L})$. }
 \end{Thm}

\begin{Rem} {
 These results are based on a precise study of the structure of
$(SW_0)$, and in particular of the
ageostrophic  part of that equation,  meaning its projection
onto~$(\Ker L)^{\perp}$. One can prove in particular that the  ageostrophic part
of~$(SW_0)$
is in fact fully parabolic. That should be compared to the case of the
incompressible limit of the compressible Navier-Stokes equations,
where again the limit system is
parabolic, contrary to the original compressible system
(see~\cite{danchinperiodic},
\cite{gallagherkyoto}, \cite{masmoudi}). Note however
that~$(SW_0)$   actually satisfies the  same type of
trilinear estimates as the three-dimensional
incompressible Navier-Stokes system, which accounts for the fact that in Theorem~\ref{limit-systemallbeta}
unique solutions are only obtained  for a short life span
(despite the fact that the space variable  runs in the two dimensional  domain~$\R \times
\T$). In the case of Theorem~\ref{limit-systemgeneric}, we use the study of resonances of the previous chapter which 
shows that the limit system is linear, except for its projection onto Kelvin modes; but Kelvin modes are essentially
one-dimensional so energy estimates are much improved compared to the case of Theorem~\ref{limit-systemallbeta}, 
and that is why global wellposedness is true in~$L^{2}$, for arbitrarily large initial data.
}
\end{Rem}

The rest of this chapter is devoted to the proof of those theorems.  
 Some preliminary results are proved in Section~\ref{preliminaries} below, namely
the fact that the ageostrophic part of the limit system is parabolic, along with trilinear estimates.
 In Section~\ref{theoremallbeta} we prove Theorem~\ref{limit-systemallbeta}, whereas  the proof 
 of Theorem~\ref{limit-systemgeneric} 
 can be found
 in Section~\ref{theoremgeneric}. The last section will be devoted to an additional regularity result, giving an estimate of the divergence of $\Phi'$ in both cases, which will be useful in the next chapter.
 
\section{Preliminary results}\label{preliminaries}\setcounter{equation}{0}
Let us prove some results that will be used throughout this chapter: in Section~\ref{parabolicity} below, we 
prove that the limit system, projected onto~$(\Ker L)^{\perp}$  is parabolic. In 
Section~\ref{sectiontrilinear} we prove crucial trilinear estimates.
\subsection{Parabolicity of the ageostrophic limit  equation}\label{parabolicity}
In this section we are going to prove that the projection of the limit  system
onto~$(\Ker L)^{\perp}$
is parabolic. To obtain that result, the
    important remark is that, for each eigenmode of
$L$, the
first and third components of the eigenvectors (corresponding to
$\eta$ and $u_2$) have very similar behaviours, and thus controlling
the
regularity of the last two components is sufficient to have an estimate 
on~$\Pi_{\perp}\Phi$ in~$H^{1}_{L}$.
  A result of quasi-orthogonality in~$H^s(\R\times
\T)$  of the nonzero eigenmodes of~$L$ leads indeed   to the following  result.
\begin{Lem}\label{parabolic}
Let~$ s
\geq 0$ be given. 
There is a constant~$C_{s}$ such that for any~$\Phi \in (\Ker L)^{\perp}$,  we have
$$
\| \Phi\|_{ H^{s+1}_L }^2\leq C_s  (\Phi| -\Delta'_L 
\Phi)_{H^s_L},
$$
meaning in particular that the projection of the system
$(SW_0)$ onto~$(\Ker L)^{\perp}$ is fully parabolic.
\end{Lem}

\begin{proof}
The proof of that result
consists in using the structure of the eigenmodes to prove
that the diffusion -- acting a priori only on the velocity
field -- has
also a smoothing effect on the pressure, and more precisely that
\begin{equation}
\label{step1}
\forall (n,k,j) \in S, \quad
\|  \Pi _{n,k,j} \Phi \|_{ H^s(\R\times \T)} \leq C'\| (\Pi _{n,k,j}  \Phi)'
\|_{H^s(\R\times
\T)},
\end{equation}
where $\Pi_{n,k,j}$  denotes the projection on the
eigenmode $\Psi_{n,k,j}$ of $L$ (with the notation  of
the previous chapter) and $C'$ is a nonnegative constant (independent
of $n$, $k$ and $j$).
 By formulas (\ref{psi-def1}), (\ref{psi-def2}) and (\ref{psi-def3})
we deduce that for any integer~$s$, we have
$$\|\d_2^s( \Psi_{n,k,j})_0\|_{L^2(\R\times \T)}=\|\d_2^s(
\Psi_{n,k,j})_2\|_{L^2(\R\times \T)}$$
using the orthogonality of $\psi_{n-1}$ and
$\psi_{n+1}$, and that
$${1\over C }\|\d_1^s( \Psi_{n,k,j})_2\|_{L^2(\R\times \T)}\leq
\|\d_1^s( \Psi_{n,k,j})_0\|_{L^2(\R\times \T)}\leq C  \|\d_1^s(
\Psi_{n,k,j})_2\|_{L^2(\R\times \T)}.$$
This implies in particular (\ref{step1}).

By Remark \ref{componentseparately}, page  \pageref{componentseparately} we then deduce
that for all $\Phi \in
(\Ker L)^\perp$
$$
\begin{aligned}
  (\Phi| -\Delta'_L
\Phi)_{L^2(\R\times \T)} &= \sum_{i \lambda \in {\mathfrak S}\setminus  \{0\}}
  (\Pi_{\lambda }\Phi| -\Delta'_L   (\Pi_{\lambda }
\Phi))_{L^2(\R\times \T)} =\sum_{i \lambda \in {\mathfrak S}\setminus  \{0\}} \|(
\Pi_{\lambda }
\Phi)'\|_{\dot H^1(\R\times \T)}^{2}\\
&\geq {1\over C} \sum_{(n,k,j)\in S^*}\|( \Pi_{n,k,j}
\Phi)'\|_{\dot H^1(\R\times \T)}^{2}\\
&\geq {1\over CC'}\sum_{(n,k,j)\in S^*}\| \Pi_{n,k,j}
\Phi\|_{\dot H^1(\R\times \T)}^{2}\\
&\geq {1\over CC'C_1'}\| \Phi\|_{H^1_L}^2
\end{aligned}
$$
recalling that
$${1\over C_s' } (1+ n+k^2)^{s/2} \leq  \| \Psi_{n,k,j}
\|_{\dot H^s(\R\times \T)} \leq C_s'   (1+ n+k^2)^{s/2}.$$
We therefore obtain the first inequality using Proposition  \ref{equivnorm}.

In a similar way, by Proposition~\ref{diag-prop},  page~\pageref{diag-prop}, we have
 for all $\Phi \in
(\Ker L)^\perp$
$$
\begin{aligned}
  (\Phi| -\Delta'_L
\Phi)_{H^s_L} &= \sum_{(n,k,j)\in S^*}(1+n+k^2)^s
  (\Pi_{n,k,j }\Phi| -\Delta'_L
\Phi)_{L^2(\R\times \T)} \\
&=\sum_{(n,k,j)\in S^*}(1+n+k^2)^s\sum_{\tau(n^*,k^*,j^*)=\tau(n,k,j)}
  (\Pi_{n,k,j }\Phi| -\Delta'
(\Pi_{n^*,k^*,j^*}\Phi))_{L^2(\R\times \T)}\\
&=\sumetage{(n,k,j)\in S^*}{\tau (n,k,j) \neq \pm \sqrt{\beta/2}}
(1+n+k^2)^s
  \|(\Pi_{n,k,j }\Phi)'\|_{\dot H^1(\R\times \T)}^2 +R_\beta
\end{aligned}
$$
where $R_\beta$ is the contribution of the modes~$\pm \sqrt{\beta/2}$, defined by
$$
\begin{aligned}
R_\beta &= \sum_{k = \pm \sqrt{\beta/2}} (1 + \frac\beta2)^s \sum_{(j,j^*) \in (0,\mbox{sign} k)^2}
\left(
\Pi_{0,k,j} \Phi | - \Delta' \Pi_{0,k,j^*} \Phi
\right) \\
&= \sum_{k = \pm \sqrt{\beta/2}} (1 + \frac\beta2)^s  
\left\| \sum_{(j,j^*) \in (0,\mbox{sign} k)^2} (\Pi_{0,k,j} \Phi )'\right\|_{\dot H^1}^2.
\end{aligned}
$$
 Using (\ref{step1}) and Proposition~\ref{orthogonality}
 page~\pageref{orthogonality}
 leads then to the expected estimate.
\end{proof}

\begin{Rem}
Note that these inequalities indicate in particular that the notion of  homogeneous or
inhomogeneous  spaces
does not make sense for these weighted Sobolev spaces.

Moreover we recall that there is no analogue of
Proposition~\ref {orthogonality}  in the case of geostrophic modes, so that the geostrophic equation does
not have that ellipticity property.

\end{Rem}

\subsection{Derivation of the   trilinear estimates}  \label{sectiontrilinear}
An important step in the proof of Theorems~\ref{limit-systemallbeta} and~\ref{limit-systemgeneric} consists   in
establishing some control on  the
nonlinear term arising in  $(SW_0)$ in terms of the Sobolev
norms introduced in Section~\ref{sectionsuitable}. Such estimates are
obtained using classical  para-differential methods. Obviously a more  general statement
could be written, at the price of more technicalities. In order to
keep the proof as simple as  possible we choose to
state only those  estimates that will be used in the
following.

\begin{Prop}\label{trilinear}
{
Denote by $Q_L$ the limit  nonlinear operator defined by  (\ref{Q-lap-L}), and let~$\alpha$ be
any
real number greater than~$3/2$.
Then the following trilinear estimates hold~:
$$
\begin{aligned}
\left|\left( \Phi_{*} |  Q_L(\Phi,\Phi^*) \right)_{L^2(\R\times \T)}
\right|\leq C &\| \Pi_{\perp} \Phi_* \|_{  H^1_L } ^{3/4}
\|   \Pi_{\perp} \Phi_* \|_{
L^2(\R\times \T)} ^{1/4} \|  \Pi_{\perp}   \Phi
\|_{ H^1_L } ^{3/4} \|  \Pi_{\perp}  \Phi^* \|_{ H^1_L }^{3/4} \\
&\quad \quad \quad   \times \left(\|  \Pi_{\perp}   \Phi^* \|_{   H^1_L } ^{1/4} \|
 \Pi_{\perp}   \Phi \|_{ L^2(\R\times \T)} ^{1/4}+\|  \Pi_{\perp}    \Phi \|_{  H^1_L }
^{1/4} \|
  \Pi_{\perp}  \Phi ^*\|_{ L^2(\R\times \T)} ^{1/4}\right) \\
     &\!\!\!\!\!\!\!\!\!\!\!\!\! \!\!\!\!\!\!\! \!\!\!\!\!\!+C \| \Pi_{\perp}  \Phi_*\|_{L^2(\R\times \T)}\left( \| \Pi_0 
\Phi\|_{L^2(\R\times \T)}\| \Pi_{\perp}
\Phi^*\|_{ H^{1}_L } + \| \Pi_0 \Phi^*\|_{L^2(\R\times \T)}\|  \Pi_{\perp} \Phi\|_{
H^{1}_L } \right)\\
\left|\left( \Phi_{*} |  Q_L(\Phi, \Phi^*) \right)_{L^2(\R\times \T)}
\right|\leq C &\|\Pi_{\perp}  \Phi \|_{   H^{3/2}_L } \| \Pi_{\perp}  \Phi_*\|_{
 H^1_L }^{1/2}
\|\Pi_{\perp}
\Phi_*\|_{L^2(\R\times \T)}^{1/2}\|\Pi_{\perp}   \Phi^*\|_{   H^1_L }  ^{1/2}\|
\Pi_{\perp} \Phi^*\|_{L^2(\R\times \T)} ^{1/2}\\
& \quad \quad \quad \quad   + C\|\Pi_{\perp}  \Phi\|_{  
  H^1_L }\| \Pi_{\perp} \Phi^*\|_{
 H^1_L }
\|\Pi_{\perp}  \Phi_*\|_{  H^1_L }^{1/2}\|\Pi_{\perp}
\Phi_*\|_{L^2(\R\times \T)}^{1/2}\\
   &\!\!\!\!\!\!\!\!\!\!\!\!\! \!\!\!\!\!\!\!\!\!\!\!\!\!   + C    \| \Phi_*\|_{L^2(\R\times \T)}
   \left( \| \Pi_0  \Phi\|_{L^2(\R\times \T)}\|
\Pi_{\perp}
\Phi^*\|_{ H^{1}_L } + \| \Pi_0 \Phi^*\|_{L^2(\R\times \T)}\|  \Pi_{\perp} \Phi\|_{
H^{1}_L } \right) \\
\left|\left( \Phi_{*} |  Q_L(\Phi,\Phi^*) \right)_{L^2(\R\times \T)}
\right|\leq C_\alpha &
 \|  \Pi_{\perp}  \Phi_*
\|_{ L^{2}(\R\times \T)}
\| \Pi_{\perp}  \Phi \|_{  H^\alpha_L }  \| \Pi_{\perp}  \Phi^{*} \|_{   H^\alpha_L }   \\
  &\!\!\!\!\!\!\!\!\!\!\!\!\!  \!\!\!\!\!\!\!\!\!\!\!\!\! \!\!\!\!\!\!\! \!\!\!\!\!\!+C _\alpha   \| \Pi_{\perp}  \Phi_*\|_{L^2(\R\times \T)}\left( \| \Pi_0 
\Phi\|_{L^2(\R\times \T)}\| \Pi_{\perp}
\Phi^*\|_{ H^{1}_L } + \| \Pi_0 \Phi^*\|_{L^2(\R\times \T)}\|  \Pi_{\perp} \Phi\|_{
H^{1}_L } \right),
\end{aligned}
$$ and for all $s\leq 1$
 $$
\begin{aligned}
\left|\left( \Phi_{*} |  Q_L(\Phi,\Phi^*) \right)_{
H^{s}_L  } \right|\leq  C & 
\left(\|\Pi_{\perp}
\Phi
\|_{ H^{3/2}_L }  \|  \Pi_{\perp}  \Phi^* \|_{
H^{1}_L }
  +
\|\Pi_{\perp}
\Phi^*
\|_{ H^{3/2}_L }   \|  \Pi_{\perp}  \Phi \|_{
H^{1}_L }\right) \| \Pi_{\perp}   \Phi_* \|_{ H^{2s}_L } \\
+
C  &\| \Pi_{\perp}\Phi_*\|_{ H^{s+1}_L }\left( \| \Pi_0  \Phi\|_{L^2(\R\times \T)}\|
\Pi_{\perp}\Phi^*\|_{ H^{s}_L} + \| \Pi_0 \Phi^*\|_{L^2(\R\times \T)}\|
\Pi_{\perp}  \Phi\|_{ H^{s}_L } \right)\\
\leq  C & 
\left(\|\Pi_{\perp}
\Phi
\|_{ H^{3/2}_L }  \|  \Pi_{\perp}  \Phi^* \|_{
H^{s}_L }^{s} \|  \Pi_{\perp}  \Phi^* \|_{
H^{s+1}_L }^{1-s}
  +
\|\Pi_{\perp}
\Phi^*
\|_{ H^{3/2}_L }   \|  \Pi_{\perp}  \Phi \|_{
H^{s}_L }^{s} \|  \Pi_{\perp}  \Phi \|_{
H^{s+1}_L }^{1-s}\right)\\
& \quad \quad \quad \quad \times  \| \Pi_{\perp}   \Phi_* \|_{ H^{s}_L }^{1-s} \| \Pi_{\perp}   \Phi_* \|_{ H^{s+1}_L }^{s} \\
+
C  &\| \Pi_{\perp}\Phi_*\|_{ H^{s+1}_L }\left( \| \Pi_0  \Phi\|_{L^2(\R\times \T)}\|
\Pi_{\perp}\Phi^*\|_{ H^{s}_L} + \| \Pi_0 \Phi^*\|_{L^2(\R\times \T)}\|
\Pi_{\perp}  \Phi\|_{ H^{s}_L } \right)
\,.
\end{aligned}
$$
}
\end{Prop}
\begin{Rem}\label{restriction}
{
1.The estimates presented in that proposition are   exactly the analogue  of the
 usual trilinear estimate for the three-dimensional Navier-Stokes  equations. For
instance in three space dimensions one has
$$
\left|\left( \Phi_{*} |   \mbox{div} \: (\Phi \otimes  \Phi^*)   \right)_{L^2(\R^{3})}
\right| \leq C  \|  \Phi_* \|_{\dot  H^{\frac34}(\R^{3})}\Bigl (\|   \Phi  \|_{\dot
H^{\frac34}(\R^{3})} \|\nabla \Phi^*\|_{L^{2}(\R^{3})}
+\|  \Phi^* \|_{\dot H^{\frac34}(\R^{3})}
 \|\nabla \Phi \|_{L^{2}(\R^{3})}\Bigr)
$$
whereas in two space dimensions one would expect
$$
\left|\left( \Phi_{*} |   \mbox{div} \: (\Phi \otimes  \Phi^*)  \right)_{L^2(\R^2)}
\right| \leq C  \|  \Phi_* \|_{\dot  H^{\frac12}(\R^{2})}\Bigl (
\|  \Phi  \|_{\dot H^{\frac12}(\R^{2})} \|\nabla  \Phi^*\|_{L^{2}(\R^{2})}
+\|  \Phi^* \|_{\dot H^{\frac12}(\R^{2})}
 \|\nabla \Phi \|_{L^{2}(\R^{2})}\Bigr).
$$
The reason for the
loss of one half derivative compared to the usual two dimensional case  is linked to the
fact that differentiation with respect
to~$x_{1}$ corresponds to a multiplication by~$\sqrt n$ instead of~$n$.

2.  The restriction $s\leq 1$ is due do the particular
 structure of the nonlinear  term, in particular to the coupling
 between Rossby modes. For the Rossby mode
 associated to the eigenvalue~$i\lambda$ the regularity is indeed measured by  $1/ \lambda$.
Therefore the condition of resonance $\lambda =\mu+\tilde \mu$ (which
 is of course not equivalent  to $\lambda ^{-1}=\mu^{-1}+\tilde
 \mu^{-1}$) does not allow to  distribute the derivatives as in the
 usual paradifferential  calculus. 
 Note nevertheless that the  computation (\ref{quasigeostrophic}) page \pageref{quasigeostrophic}
 allows actually to distribute one derivative and to obtain a trilinear estimate of the form
 $$\left|\left( \Phi |  Q_L(\Phi,\Phi) \right)_{
H^{s}_L  } \right|\leq  C \| \Pi_\perp \Phi \|_{H^{s+1}_L}\| \Pi_\perp \Phi \|_{H^s_L} \left( \| \Pi_\perp \Phi \|_{H^{3/2}_L} +\|\Pi_0 \Phi\|_{L^2(\R\times \T)}\right),$$
 for all $s\leq 2$.

}
\end{Rem}

\begin{proof}
The method used to establish these estimates is rather standard~: we
decompose each vector on the eigenmodes of $L$, then compute each
elementary
trilinear term, and finally determine summability conditions. The
fundamental result we will use to estimate the sums is the following (see~\cite{liebloss})
\begin{equation}
\label{convolutionfaible}
\begin{aligned}
\forall v\in \ell^p (\N\times \Z\times \{-1,0,1\}), \: \forall w\in
\ell^{q,\infty}(\N\times \Z\times \{-1,0,1\}),\:  v* w\in  \ell^r(\N\times
\Z\times \{-1,0,1\})\\
\hbox{ with }p,q,r \in ]1,+\infty[\hbox{ and }\frac1r=\frac1p+\frac1q-1,
\end{aligned}
\end{equation}\label{ellqfaible}
where the convolution is to be understood in~$k$ and~$n$,
coupled with the classical result
\begin{equation}
\label{l3/2}
((1+n+k^2)^{-1}) \in \ell^{3/2,\infty}(\N\times
\Z\times \{-1,0,1\}).
\end{equation}
In   the sequel we will use the following notation \label{S*-not}
$$S=\N\times \Z\times \{-1,0,1\} \: \hbox{ and } \:  S^*=\N\times  \Z\times
\{-1,0,1\}\setminus \N\times \{0\}\times \{0\}.$$

 We have by definition of the space~$H^s_{L}$,
 $$
 \left( \Phi_{*} |  Q_L(\Phi,\Phi^*) \right)_{H^s_{L}}  =  \sum_{(n_*,k_*,j_*) \in S} (1
+ n_* + k_*^{2})^{s}
\left( \Pi_{n_*,k_*,j_*}  \Phi_{*} \Big| \Pi_{n_*,k_*,j_*}  Q_L(\Phi,\Phi^*)
\right)_{L^2}.
 $$
 We can then write
 $$
  \left( \Phi_{*} |  Q_L(\Phi,\Phi^*) \right)_{H^s_{L}} =  \!\!\!\!\sum_{(n_*,k_*,j_*) \in S}
  \!\!\!\!\!\sumetage{i\lambda,i\lambda^* \in {\mathfrak S}}{ \lambda +
  \lambda^* = \tau (n_*,k_*,j_*) }
  (1 + n_* + k_*^{2})^{s}
   \left( \Pi_{n_*,k_*,j_*}  \Phi_{*} \Big|   \Pi_{n_*,k_*,j_*} Q  (\Pi_{\lambda}  \Phi,
\Pi_{\lambda^*} \Phi^*)\right)_{L^2} .
 $$
  
\medskip
$\bullet$ Let us start by estimating
   the purely ageostrophic part of $Q_L$, denoted~$ \tilde Q_{L}$ and  defined by
\begin{equation}
\label{tQL-def}
   \tilde Q_L (\Phi,\Phi^*)=\sum_{i\lambda, i\mu,i\mu^* \in {\mathfrak
S}\setminus \{0\} \atop \lambda =\mu+\mu^*} \Pi_\lambda Q
(\Pi_\mu
\Phi,
\Pi_{\mu^*}\Phi^*).
\end{equation}
We have
$$
 \left( \Phi_{*} | \tilde  Q_L(\Phi,\Phi^*) \right)_{H^s_{L}} = \!\!\!\!\!\!\!\!\!\!\!\!
  \sumetage{ k + k^* = k_*}{ \tau ( n, k, j)  + \tau ( n^*, k^*, j^*) =  \tau
(n_*,k_*,j_* )} \!\!\!\!\!\!\!\!\!\!\!\!\!\!\!\!\!\!\!\!\!\!\!\!\!
   (1 + n_* + k_*^{2})^{s}
  \left(  \Pi_{n_*,k_*,j_*}  \Phi_{*}  \Big|  \Pi_{n_*,k_*,j_*} Q  (\Pi_{  n, k, j} \Phi,
\Pi_{ n^*, k^*, j^*} \Phi^*)\right)_{L^2} ,
$$
where the eigenvalues $i\tau(n,k,j),i\tau(n_*,k_*,j_*)$ and
$i\tau(n^*,k^*,j^*)$ run over ${\mathfrak S}\setminus\{0\}$.

Thus using the regularity estimates on the eigenvectors  $(\Psi_{n,k,j})$ of $L$
stated in Proposition \ref{diag-prop}, page~\pageref{diag-prop}, we
get (writing to simplify~$ \tau$ for~$ \tau (n,k,j)$, and similarly~$
\tau_* = \tau (n_*,k_*, j_* )$ and~$ \tau^* = \tau (n^* ,k^* ,j^* )$)
$$\begin{aligned}
&\left|  \sum_{{\tau_*=\tau +\tau^*\atop
k_*=k+k^* }\atop
i\tau,i\tau_*, i\tau^* \in {\mathfrak S}\setminus \{0\}}
   (1 + n_* + k_*^{2})^{s}
  \left(  \Pi_{n_*,k_*,j_*}  \Phi_{*} \Big| \Pi_{n_*,k_*,j_*} Q (\Pi_{   n, k, j} \Phi,
\Pi_{ n^*, k^*, j^*} \Phi^*)\right)_{L^2} 
\right|\\
&\leq C\sum_{{\tau_*=\tau +\tau^*\atop
k_*=k+k^* }\atop
i\tau,i\tau_*, i\tau^* \in {\mathfrak S}\setminus \{0\}}
 (1 + n_* + k_*^{2})^{s}
|{(\varphi_*)}_{n_*,k_*,j_*}| |\varphi_{n,k,j}|
|\varphi^*_{n^*,k^*,j^*}|  ( (n+k^2)^{1/2} +(n^*+(k^*)^2)^{1/2}),
\end{aligned}
$$
where $\varphi_{n,k,j}$ is defined as in the proof of Proposition
 \ref{diag-prop} by
$$\varphi_{n,k,j}= \left( \Psi_{n,k,j}|\Phi \right)_{L^2(\R\times  \T)}.$$

For the sake of clearness, we will simplify (abusively) the notations
as follows~: we will denote  respectively
by $\varphi$, $\varphi_*$ and $\varphi^*$ the coefficients
$\varphi_{n,k,j}$, $(\varphi_*)_{n_*, k_*,
j_*}$ and
$\varphi^*_{n^*,k^*,j^*}$, and by $\sum$ the sum over $(n,k,j)(n_*,
k_*, j_*) (n^*,k^*,j^*)  \in (S^*)
^3$ satisfying the following constraints
$$ \tau(n_*,k_*, j_*)=\tau(n,k,j)+\tau(n^*,k^*,j^*)\hbox{ and }
k_*=k+k^*.$$

 It is fundamental for the following estimates to notice that those
  constraints in fact imply  that  when~$k $ and~$k^*$ as well
  as~$j,j^*$ and~$j_*$ are
 fixed, then the condition~$\tau(n_*,  k + k^*,  j_*)=\tau(n,k,j)+\tau(n^*,k^*,j^*)$
implies that a given~$n$ and~$n^*$
  constrain the value of~$n_*$.  Indeed we
 recall that
according to Lemma~\ref{polynomial}, page~\pageref{polynomial}, there
  is only one value of~$n_*$  associated
with one value of~$\tau(n_*,k_*, j_*) \neq 0$.  Note however that  contrary to the usual
case when there is an actual convolution (as is  the case
for the Fourier variable~$k$ here), we have no ovbious estimate on~$n_*$, as a   function of~$n$
and~$n^*$. So the
usual methods of distribution of derivatives  cannot be fully used here,
as derivatives in the~$x_{1}$ direction, acting on~$\Phi_*$ cannot be  traded for
derivatives on~$\Phi$ or~$\Phi^*$.

\medskip
(i) \underbar{If $s = 0$},
by the Cauchy-Schwarz inequality, we obtain
$$
\begin{aligned}
&\sum  |\varphi_*| |\varphi| |\varphi^*| (n+k^2)^{1/2}\\
 &\leq \left( \sum   (n+k^2)^{1/2} |\varphi| (1+  n_*+k_*^2)^{3/4}|\varphi_*| ^{3/2}
(1+n^*+(k^*)^2)^{-3/4}|\varphi^*|^{1/2}\right)^{1/2}
\\
&\times \left( \sum   (n+k^2)^{1/2} |\varphi| (1+  n_*+k_*^2)^{-3/4}|\varphi_*| ^{1/2}
(1+n^*+(k^*)^2)^{3/4}|\varphi^*|^{3/2}\right)^{1/2}.
\end{aligned}
$$

By (\ref{convolutionfaible}) and (\ref{l3/2}), we therefore get
$$
\begin{aligned}
&\sum  |\varphi_*| |\varphi| |\varphi^*| (n+k^2)^{1/2}\\
 &\leq \left(\| (n+k^2)^{1/2} \varphi\|_{\ell^2(S^*)}\| (1+  n_*+k_*^2)^{1/2}\varphi_*
\|_{\ell^2(S^*)}^{3/2}\| (1+n^*+(k^*)^2)^{-3/4}\|_{\ell^{2,\infty}(S)}
 \|\varphi^*\|_{\ell^2(S^*)}^{1/2} \right)^{1/2}
\\
&\times \left(\| (n+k^2)^{1/2} \varphi\|_{\ell^2(S^*)}\| (1+n_*+
 k_*^2)^{-3/4}\|_{\ell^{2,\infty}(S)}  \|
\varphi_*
\|_{\ell^2(S^*)}^{1/2}\|(1+n^*+(k^*)^2)^{1/
2}\varphi^*\|_{\ell^2(S^*)}^{3/2}\right)^{1/2}\\
\end{aligned}
$$
and a similar estimate for the symmetric term.

By definition, 
\begin{equation}
\label{equivalence-sum}
\| (1+n+k^2)^{1/2} \varphi\|_{\ell^2(S^*)} \leq C \|\Pi_{\perp} \Phi  \|_{
 H^1_L }.
\end{equation}
 Plugging this estimate in the previous inequality leads to
\begin{equation}
\label{tQ-est}
\begin{aligned}
\,&\left|\left( \Phi_{*} | \tilde Q_L(\Phi,\Phi^*)
\right)_{L^2(\R\times \T)}\right| \leq C\| \Pi_{\perp} \Phi_* \|_{   H^1_L } ^{3/4}
\|   \Pi_{\perp} \Phi_* \|_{
L^2(\R\times \T)} ^{1/4} \|  \Pi_{\perp}   \Phi
\|_{ H^1_L } ^{3/4} \|  \Pi_{\perp}  \Phi^* \|_{ H^1_L }^{3/4} \\
&\quad \quad\quad\quad\quad \quad  \times \left(\|  \Pi_{\perp}    \Phi^* \|_{  H^1_L }
^{1/4} \|
 \Pi_{\perp}   \Phi \|_{ L^2(\R\times \T)} ^{1/4}+\|  \Pi_{\perp}    \Phi \|_{  H^1_L }
^{1/4} \|
  \Pi_{\perp}  \Phi ^*\|_{ L^2(\R\times \T)} ^{1/4}\right).
\end{aligned}
\end{equation}

\medskip
(ii) \underbar{If $s = 0$}, another way to estimate the $L^2$ scalar  product is as
follows
$$
\begin{aligned}
&\sum  |\varphi_*| |\varphi| |\varphi^*| (n+k^2)^{1/2}\\
&\leq \left( \sum \left( (1+n+k^2)^{3/4} |\varphi|\right)^{3/2} 
\left((1+n_*+k_*^2)^{1/4} |\varphi_*|\right)^{3/2}  (1+n^*+(k^*)^2)^{-3/4} \right)^{1/3}\\
&\times \left( \sum \left( (1+n+k^2)^{3/4} |\varphi|\right)^{3/2} 
\left((1+n^*+(k^*)^2)^{1/4} |
\varphi^*|\right)^{3/2} (1+n_*+k_*^2)^{-3/4} \right)^{1/3}\\
&\times \left(\sum  (1+n+k^2)^{-3/4} \left( (1+n_*+k_*^2)^{1/4}  |\varphi_*|\right)^{3/2}
\left((1+n^*+(k^*)^2)^{1/4} |
\varphi^*|\right)^{3/2} \right)^{1/3}
\end{aligned}
$$
and
$$
\begin{aligned}
&\sum  |\varphi_*| |\varphi| |\varphi^*| (n^*+(k^*)^2)^{1/2}\\
&\leq \left( \sum  (n^*+(k^*)^2)^{1/2} |\varphi^*| \left(  (1+n+k^2)^{1/2}
|\varphi|\right)^{7/6}\left((1+n_*+k_*^2)^{1/4}  |\varphi_*|\right) ^{7/6}(1+n_*
+k_*^2)^{-1/2}
\right)^{1/2}\\
&\times  \left(\sum  (n^*+(k^*)^2)^{1/2} |\varphi^*| \left(  (1+n+k^2)^{1/2}
|\varphi|\right)^{5/6}\left((1+n_*+k_*^2)^{1/4}  |\varphi_*|\right)^{5/6} (1+ n +
k^2)^{-1/4}
\right)^{1/2}.
\end{aligned}
$$

By (\ref{convolutionfaible}) and (\ref{l3/2}), we therefore get
$$
\longformule{ \sum  |\varphi_*| |\varphi| |\varphi^*| (n+k^2)^{1/2}\\
 \leq \| (1+n+k^2)^{3/4} \varphi\|_{\ell^2(S^*)} \|(1+n_*+k_*^2)^{1/4} 
\varphi_*\|_{\ell^2(S^*)}}{ \times\| (1+n^*+(k^*)^2)^{1/4}
\varphi^*\|_{\ell^2(S^*)}\|(1+n+k^2)^{-3/4}\|_{\ell^{2,\infty}(S^*)}}
$$
and
$$
\longformule {\sum  |\varphi_*| |\varphi| |\varphi^*|  (n^*+(k^*)^2)^{1/2} \leq \| 
(n^*+(k^*)^2)^{1/2}
 \varphi^*\|_{\ell^2(S^*)} \| (1+n+k^2)^{1/2} 
\varphi\|_{\ell^2(S^*)}}{\times\|(1+n_*+k_*^2)^{1/4} 
\varphi_*\|_{\ell^2(S^*)}\|(1+n+k^2)^{-3/4}\|_{\ell^{2,\infty}(S^*)}   .}
$$

So we find
\begin{equation}
\label{tQ-est*}
\begin{aligned}
\left|\left( \Phi_{*} |  \tilde Q_L(\Phi, \Phi^*) \right)_{L^2(\R\times  \T)}
\right|&\leq C \|\Pi_{\perp}  \Phi \|_{   H^{3/2}_L } \| \Pi_{\perp}  \Phi_*\|_{
 H^1_L }^{1/2}
\|\Pi_{\perp}
\Phi_*\|_{L^2(\R\times \T)}^{1/2}\|\Pi_{\perp}   \Phi^*\|_{   H^1_L }  ^{1/2}\|
\Pi_{\perp} \Phi^*\|_{L^2(\R\times \T)} ^{1/2}\\
& \quad \quad \quad \quad \quad + C\|\Pi_{\perp}  \Phi\|_{  H^1_L }\|  \Pi_{\perp}
\Phi^*\|_{  H^1_L }
\|\Pi_{\perp}  \Phi_*\|_{  H^1_L }^{1/2}\|\Pi_{\perp}
\Phi_*\|_{L^2(\R\times \T)}^{1/2}.
\end{aligned}
\end{equation}

\medskip
(iii) \underbar{If $s = 0$} and we have additional regularity on~$\Phi$  and~$\Phi^*$,
then again
one can write a different estimate. We can write indeed
\begin{eqnarray}\label{additionalestimate}
 &&\sum  |\varphi_*| |\varphi| |\varphi^*| (n^*+(k^*)^2)^{1/2}
 \leq  
 \| \varphi_*\|_{\ell^2(S^*)} \|  (1+n^*+(k^*)^2)^{\alpha/2}
 \varphi^*\|_{\ell^2(S^*)} \\
 & \times&  \|  (1+n^*+(k^*)^2)^{1/2 -  \alpha/2}\|_{\ell^{\infty}(S^*)} 
  \|  (1+n +k^2)^{s/2}
 \varphi \|_{\ell^2(S^*)}  \|  (1+n +k^2)^{-\alpha/2}\|_{\ell^{2}(S^*)},   \nonumber
\end{eqnarray}
which, coupled with the similar estimate for  the symmetric term, gives  the expected
result.

\medskip
(iv) \underline{If $s \leq 1$},
we have  by H{\"o}lder's inequality
$$
\begin{aligned}
&\sum   (1 +   n_* +  k_*^{2})^{s}  |  \varphi_*| |\varphi|  |\varphi^*|
(n+k^2)^{1/2}\\
&   \! \! \! \! \!\leq \left( \sum   (1 +   n_* +  k_*^{2})^{s}
 |  \varphi_*|\left(  (1+n+k^2)^{3/4}|\varphi|\right)^{7/6}\left(  (1+n^*+(k^*)^2)^{1/2}
|\varphi^*|
\right)^{7/6} (1+n+k^2)^{-1/2}
\right)^{1/2}\\
&  \! \! \! \! \!\times \left( \sum  (1 +   n_* +  k_*^{2})^{s}
 |  \varphi_*|\left(  (1+n+k^2)^{3/4}|\varphi|\right)^{5/6}\left(  (1+n^*+(k^*)^2)^{1/2}
|\varphi^*|
\right)^{5/6} (1+n^*+(k^*)^2)^{-1}
\right)^{1/2}\\
\end{aligned}
$$
from which we deduce
$$
\longformule {
 \sum (1 +   n_* +  k_*^{2})^{s}  |  \varphi_*| |\varphi|  |\varphi^*|
(n+k^2)^{1/2}\\
 \leq C\| (1 +   n_* +  k_*^{2})^{s}    \varphi_*\|_{\ell^2(S^*)}  \|
(1+n+k^2)^{3/4}\varphi\|_{\ell^2(S^*)}}
{ \times \|(1+n^*+(k^*)^2)^{1/2} \varphi^*\|_{\ell^2(S^*)}
\|(1+n+k^2)^{-3/4}\|_{\ell^{2,\infty}(S^*)}}
 $$
and a similar estimate for the symmetric term.

Therefore, we get
$$\left|\left( \Phi_{*} |  \tilde Q_L(\Phi,\Phi^*) \right)_{
H^{s}_L } \right|\leq C  \| \Pi_{\perp}   \Phi_* \|_{ H^{2s}_L }
(\|\Pi_{\perp}
\Phi
\|_{ H^{3/2}_L }  \|  \Pi_{\perp}  \Phi^* \|_{
H^{1}_L }
  + \|\Pi_{\perp}
\Phi^*
\|_{ H^{3/2}_L }  \| \Pi_{\perp}   \Phi \|_{
H^{1}_L })\,,
$$
and we conclude by interpolation
\begin{equation}
\label{tQ-est2*}
\begin{aligned}
 \left|\left( \Phi_{*} |  \tilde Q_L(\Phi,\Phi^*) \right)_{
H^{s}_L } \right|&\leq C \left(\|\Pi_{\perp}
\Phi
\|_{ H^{3/2}_L }  \|  \Pi_{\perp}  \Phi^* \|_{
H^{s}_L }^{s} \|  \Pi_{\perp}  \Phi^* \|_{
H^{s+1}_L }^{1-s} +
\|\Pi_{\perp}
\Phi^*
\|_{ H^{3/2}_L }   \|  \Pi_{\perp}  \Phi \|_{
H^{s}_L }^{s} \|  \Pi_{\perp}  \Phi \|_{
H^{s+1}_L }^{1-s}\right)\\
&\quad \quad \times  \| \Pi_{\perp}   \Phi_* \|_{ H^{s}_L }^{1-s} \| \Pi_{\perp}   \Phi_* \|_{ H^{s+1}_L }^{s} 
\,.
\end{aligned}
\end{equation}

\medskip
$\bullet$ Proposition~\ref{limitgeolinear}  shows that
$$Q_L (\Phi,\Phi^*) = Q_L(\Pi_0 \Phi,   \Pi_{\perp} \Phi^*) + Q_L(
  \Pi_{\perp} \Phi, \Pi_0 \Phi^*) + \tilde Q_L (\Phi, \Phi^*).$$

In order to end the proof of the proposition, it remains therefore to  estimate the
terms coupling the geostrophic and ageostrophic parts. We start
by noticing that the constraint~$\tau (n,k,j) = \tau (n',k,j')$ implies  by
Proposition~\ref{polynomial}
page~\pageref{polynomial} that
necessarily~$n = n'$ and~$j = j'$, except if~$n = 0$. But according to 
Remark~\ref{racinesconfondues} that
case  corresponds to two different values of~$j$ for the same  eigenvalue only if~$2k^{2}
= \beta$, in which
case the multiple root is~$ k = \pm \sqrt {\beta/2}$. That means that  one can write
$$
\begin{aligned}
 |\left(\Phi_* | Q_L(\Pi_0 \Phi,   \Pi_{\perp}  \Phi^*)\right)_{H^{s}_{L}}|\leq &
 \left|\sum_{(n,k,j)\in S^*}(1 + n + k^{2})^{s}\left( \Pi_{ n, k,
j} \Phi_* | Q( \Pi_0 \Phi,\Pi_{n,k,j} \Phi^*)  \right)_{L^2(\R\times\T)}\right|\\
& +\sum_{j=\pm 1}\left|(1 +{\beta\over 2})^{s}\left( \Pi_{ 0,  j\sqrt{\beta/2},
0} \Phi_* | Q( \Pi_0 \Phi, \Pi_{ 0, \sqrt{\beta/2},
j}  \Phi^*) \right)_{L^2(\R\times\T)}\right|\\
&+\sum_{j=\pm 1}\left|(1 +{\beta\over 2})^{s}\left( \Pi_{ 0,  j\sqrt{\beta/2},
j} \Phi_* | Q( \Pi_0 \Phi, \Pi_{ 0, \sqrt{\beta/2},
0}  \Phi^*) \right)_{L^2(\R\times\T)}\right|.
\end{aligned}
$$
Integrating by parts when the derivative acts on $\Pi_0 \Phi$, we get
$$\  \left|\left( \Pi_{ n, k,
j} \Phi_* | Q( \Pi_0 \Phi,\Pi_{n,k,j^*} \Phi^*) \right)_{L^2(\R\times
\T)}\right| \leq  C |(\varphi_*)_{n,k,j}|
|(\varphi^*)_{n,k,j^*}|
(n+k^2)^{1/2}\| \Pi_0
\Phi\|_{L^2(\R\times \T)} .
$$
By the Cauchy-Schwarz inequality, we then get
$$\begin{aligned}
\,&|\left(\Phi_* | Q_L(\Pi_0 \Phi,   \Pi_{\perp}  \Phi^*)\right)_{H^{s}_{L}}| \\
&\leq   C\| \Pi_0 \Phi\|_{L^2(\R\times \T)} \|(1 + n_* + k_*^{2})^{s/2} 
\varphi_*\|_{\ell^2(S^*)} \|
 (1+n^*+(k^*)^2)^{(s+1)/2}  \varphi^*\|_{\ell^2(S^*)}
\end{aligned}
$$
Remark that  the derivatives can be distributed either on $\Phi_*$ or  on  $\Phi^*$.

We finally deduce that
\begin{equation}
\label{Q0-est}
\begin{aligned}
\,&\left|\left( \Phi_{*} | \tilde Q_L(\Pi_0\Phi,\Phi^*)
\right)_{ H^{s}_L }+\left( \Phi_{*} | \tilde Q_L(\Phi,\Pi_0\Phi^*)
\right)_{ H^{s}_L}\right|\\
&\leq  C  \| \Pi_{\perp}\Phi_*\|_{ H^{s}_L }\left( \| \Pi_0  \Phi\|_{L^2(\R\times \T)}\|
\Pi_{\perp}\Phi^*\|_{ H^{s+1}_L} + \| \Pi_0 \Phi^*\|_{L^2(\R\times  \T)}\|
\Pi_{\perp}  \Phi\|_{ H^{s+1}_L } \right)
\end{aligned}
\end{equation}
Note that this term is not zero as in the case of the usual Sobolev  spaces, because the
spectrum of $L$ is not symmetric with
respect to $0$.

One should also remark that in~(\ref{Q0-est}), no derivative acts on  any vector field
in~$\Ker L$. This
can seem somewhat surprising, but is due to the very strong constraint  induced by the
resonance: instead of
a summation over three types of indexes (namely~$(n,k)$, $(n^*,k^*)$  and~$(n_*, k_*)$),
one
only sums over~$n$.

\medskip
Combining (\ref{tQ-est}),    (\ref{tQ-est*})  and~(\ref{additionalestimate}) with 
(\ref{Q0-est})
(with~$ s=0$) provides the first estimate of the proposition, while (\ref{tQ-est2*}) and (\ref{Q0-est}) (with~$ s \leq 1$)  give the second
one.

The proposition is proved.
\end{proof}

\section{Proof of Theorem 2}\label{theoremallbeta}\setcounter{equation}{0}
The proof of Theorem~\ref{limit-systemallbeta} is divided into four steps. In Paragraph~\ref{sectionweak} is proved
the existence of weak solutions, and the propagation of regularity of the geostrophic part is proved in 
Paragraph~\ref{regularity}. The construction of strong solutions is performed in Paragraph~\ref{sectionstrong}, 
while the   propagation of regularity of the ageostrophic part is proved in 
Paragraph~\ref{regularitystrong}.

\subsection{Weak solutions} \label{sectionweak}
In this section we are going to prove the existence of weak solutions  to the limit
filtered system~$(SW_0)$. We follow the lines of the classical proof  of the Leray theorem, stating  the existence of weak solutions to the Navier-Stokes equations.

\medskip

$\bullet$ {\it Definition of the approximation scheme}
\index{K_N@$K_N$}

Denote by $K_N$ the truncation operator\label{KN-not} defined by
\begin{equation}\label{definitionKNenvelope}
K_N =\sumetage{(n,k,j) \in S}{(n + k^2)^{1/2} \leq N} \Pi_{n,k,j}.
\end{equation}
Clearly  the operator~$
 K_N Q_L(K_N\Phi,K_N\Phi)-\nu K_N\Delta'_L K_N\Phi$
is continuous on $L^2(\R\times \T)$ (with a norm depending on $N$).
Therefore, we deduce from
the Cauchy-Lipschitz theorem that there exists a unique maximal\label{Phi(N))}
solution $\Phi_{(N)} \in
C([0,T_N[, L^2(\R\times \T))$ to
\begin{equation}
\label{Friedrichs}
\begin{aligned}
\d_t \Phi_{(N)}+K_N Q_L(K_N\Phi_{(N)},K_N\Phi_{(N)})-\nu K_N\Delta'_L  K_N\Phi_{(N)}=0,\\
\Phi_{(N)}(t=0)=K_N \Phi^0\,.
\end{aligned}
\end{equation}
Note that the uniqueness implies in particular that  $\Phi_{(N)}=K_N\Phi_{(N)}$.

Now let us write an energy estimate on~(\ref{Friedrichs}). The
quadratic form being  skew-symmetric in~$ L^2$ (this can easily be   seen by
its definition as the limit of the filtered quadratic form
 in~(\ref{filtered}) page~\pageref{filtered}) we find that
$$
\frac12 \frac d{dt} \|\Phi_{(N)}(t)\|_{L^2(\R\times \T)}^2 - \nu
(\Delta'_L \Phi_{(N)} | \Phi_{(N)} )_{L^2(\R\times \T)} = 0.
$$
Applying  Lemma~\ref{parabolic} implies that
\label{deltaprimterm}
$$
 -
(\Delta'_L \Pi_\perp\Phi_{(N)} | \Pi_\perp \Phi_{(N)} )_{L^2(\R\times  \T)}  \geq C_0^{-1} \|
 \Pi_\perp \Phi_{(N)} \|_{H^1(\R\times \T)} ^{2}
$$
so we infer by Gronwall's lemma that
\begin{eqnarray}
\label{approximate-energy}
\|\Phi_{(N)}(t)\|_{L^2(\R\times \T)}^2 &+& 2\nu \int_0^t \|
\nabla (\Pi_0 \Phi_{(N)})'(t')\|_{L^2(\R\times
\T)}^2 dt' +
 2{\nu \over C_0} \int_0^t \| \Pi_\perp\Phi_{(N)} (t') \|_{  H^1(\R\times \T)}^2 dt'
\nonumber \\
&\leq&  \|K_N \Phi^0\|_{L^2(\R\times \T)}^2 \leq  \|
\Phi^0\|_{L^2(\R\times \T)}^2,
\end{eqnarray}
so that
 the  approximate solution is defined globally, i.e.,  $T_N=+\infty$.
Moreover the proof of Lemma~\ref{parabolic} also implies that
$$
\Pi_\perp\Phi_{(N)} \: \mbox{is bounded in} \:  L^{2}(\R^{+};H^1_{L}).
$$
\medskip
$\bullet$ {\it Existence of a weak solution}

We will only sketch the proof of the existence of a weak solution, as
it is very similar to the case of the 3D incompressible Navier-Stokes equations.
 By (\ref{approximate-energy}) we deduce that
$$
\begin{aligned}
((\Phi_{(N)})_0) \hbox{ is
uniformly bounded in }L^\infty (\R^+,L^2(\R\times \T))  \: \: \:    \\
(\Phi'_{(N)}) \hbox{ is
uniformly bounded in }L^\infty (\R^+,L^2(\R\times \T)) \cap
L^2_{loc}(\R^+, H^1(\R\times \T)) \: \: \: \mbox{and}  \\
 \Pi_\perp \Phi_{(N)} \hbox{ is
uniformly bounded in }L^\infty (\R^+,L^2(\R\times \T)) \cap
L^2(\R^+, H^1_L).
\end{aligned}
$$

For any~$h>0$, denote by $\delta_h \Phi_{(N)}(t,x)=\Phi_{(N)}(t+h,x)-\Phi_{(N)}(t,x)$.  Then
$$
\begin{aligned}
\| &\delta_h \Phi_{(N)} (t)\|_{L^2([0,T]\times \R\times \T)}^2 \\
&=\int_0^T \left( \delta_h \Phi_{(N)} (t)\Big|\int_t^{t+h} \d_t  \Phi_{(N)}(s) ds
\right)_{L^2(\R\times
\T)}dt\\
&=-\int_0^T \int_t^{t+h}\left( \delta _h \Phi_{(N)} (t) | Q_L
(\Phi_{(N)}(s),\Phi_{(N)}(s))\right)_{L^2(\R\times \T)}
dsdt\\
&+\nu \int_0^T \int_t^{t+h}\left( \delta _h \Phi_{(N)} (t)| \Delta'_
L\Phi_{(N)}(s)\right)_{L^2(\R\times \T)}
dsdt
\end{aligned}
$$
  By Proposition~\ref{trilinear} and the positivity of $-\Delta'_L$, we
  deduce that
$$
\begin{aligned}
\| &\delta_h \Phi_{(N)} (t)\|_{L^2([0,T]\times \R\times \T)}^2 \\
&\leq C\int_0^T \int_t^{t+h}
\| \Pi_\perp \delta_h \Phi_{(N)} (t)\|_{  H^1_L}^{3/4}
\| \Pi_\perp \delta_h \Phi_{(N)} (t)\|_{ L^2(\R\times\T)}^{1/4}
\| \Pi_\perp  \Phi_{(N)}(s)\|_{   H^1_L } ^{7/4} \\
& \quad \quad \quad \quad \quad \quad \quad \quad \quad \quad \quad  \quad \quad \quad
\quad \quad \quad \quad \quad \quad
\times \| \Pi_\perp\Phi_{(N)}(s)\|_{L^2(\R\times \T)}^{1/4} dsdt\\
&+  C\int_0^T \int_t^{t+h}
\| \delta_h \Pi_\perp\Phi_{(N)} (t)\|_{L^{2}}\|\Pi_0
\Phi_{(N)}(s)\|_{L^2(\R\times\T)}
\| \Pi_\perp  \Phi_{(N)}(s)\|_{ H^{1}_{L}} dsdt \\
&+\nu \int_0^T \int_t^{t+h}
\left(\Phi_{(N)} (s)| \Delta'_
L\Phi_{(N)}(s)\right)_{L^2(\R\times \T)}^{1/2}  \left( \delta _h \Phi_{(N)} (t)| \Delta'_
L\delta_h \Phi_{(N)}(t)\right)^{1/2}_{L^2(\R\times \T)} dsdt\\
\end{aligned}
$$
Therefore, using the uniform $L^\infty(\R^+,L^2(\R\times \T))$ bounds on $\Phi_{(N)}$ and $\delta_h \Phi_{(N)}$,  and the uniform $L^2(\R^+, H^1_L)$ bounds on $\Pi_\perp \Phi_{(N)}$ and $\Pi_\perp \delta_h \Phi_{(N)}$ coming from the energy estimate, we get by H\" older's inequality
$$
\begin{aligned}
\| &\delta_h \Phi_{(N)} (t)\|_{L^2([0,T]\times \R\times \T)}^2 \\
& \leq  C T ^{5/8} \| \Pi_\perp \delta_h \Phi_{(N)} \|_{ L^2([0,T], H^1_L)} ^{3/4} \| \Pi_\perp\Phi_{(N)} \|_{ L^2([0,T], H^1_L)} ^{7/4} h^{1/8} + C T  \| \Pi_\perp \delta_h \Phi_{(N)} \|_{ L^2([0,T], H^1_L)} h^{1/2} \\
&+C\nu T ^{1/2} \left(\Phi_{(N)} | \Delta'_
L\Phi_{(N)}\right)_{L^2(\R\times \R\times \T)}^2h^{1/2}\,,
\end{aligned}
$$
and thus
$$
\| \delta_h \Phi_{(N)} (t)\|_{L^2([0,T]\times \R\times \T)}^2 \leq C_T  h^{1/8}.
$$
By interpolation, on gets therefore  that (up to extraction)
\begin{eqnarray*}
\label{strong-cv}
\Pi_0\Phi_{(N)}&\rightharpoonup &\Pi_0\Phi\hbox{ weakly in 
}L^2(\R^+,L^2_{loc}(\R\times
\T))\nonumber \\
\Phi_{(N)}'& \to& \Phi' \hbox{ strongly in  }L^2_{loc}(\R^+,L^2_{loc}(\R\times
\T))\nonumber \\
\Pi_\perp \Phi_{(N)} & \to& \Pi_\perp \Phi  \hbox{ strongly in 
}L^2_{loc}(\R^+,L^2(\R\times
\T)).
\end{eqnarray*}
Note that, because of Remark \ref{decroissance}, the last convergence  is actually global
in space.
We are then able, as in the usual case of the 3D Navier-Stokes
 equations,
 to take limits in the weak formulation of~(\ref{Friedrichs}), which proves
that $\Phi$ is a weak solution to $(SW_0)$.

\bigskip
$\bullet$ {\it Strong-weak uniqueness}

In general such a weak solution is not unique and the Cauchy problem
is not well-posed in  $L^2(\R\times \T)$.
Nevertheless we have the following strong-weak uniqueness principle.

\begin{Prop}
There is a positive constant~$C$ and a nondecreasing, positive  function~$ C(t)$ such
that the
following holds.
Let $\Phi$ and $ \Phi_* $ be two weak solutions to
$(SW_0)$ with respective  initial data
$\Phi^0$ and $ \Phi^0_*$, satisfying the energy estimate.
 Assume that there exists some $T>0$ such that~$\Pi_\perp \Phi$ belongs to~$ L^\infty([0,T],H^{1/2}_L)\cap L^2([0,T],  H^{3/2}_L  )$. Then  for all $t\in  [0,T]$, the function $\delta \Phi= \Phi_*-\Phi$
 satisfies
$$
\|\delta \Phi (t)\|_{L^2} ^2\leq \|
\delta \Phi (0)\|_{L^2}^2\exp\Bigl( C  (t) (1+ \|\Phi^0\|_{L^2}^2)
+ C_{\nu} (1+\|\Pi_\perp\Phi \|^2_{L^{\infty}( [0,t];H^{1/2}_L)}) \int_0^t
\|\Pi_\perp\Phi(t')\|^2_{ H^{3/2}_L} dt' \Bigr).
$$
In particular, $ \Phi_*=\Phi$ on $[0,T]\times \R\times \T$ if $  \Phi^0_*=\Phi^0$.
\end{Prop}
\begin{proof}
In order to establish the stability inequality we start by writing  (formally) the
equation on $\delta \Phi= \Phi_*-\Phi$
\begin{equation}
\label{delta-phi}
\d _t \delta \Phi + Q_L (\delta \Phi,\delta \Phi) +2Q_L (\delta \Phi,
\Phi)-\nu \Delta'_L \delta  \Phi=0 .
\end{equation}
Proposition~\ref{trilinear} implies that
$$
\begin{aligned}
|\left(\delta \Phi |  Q_L (\delta \Phi,\Phi)\right)_{L^2}| \leq
C \|\Pi_\perp  \Phi\|_{H^{3/2}_L}  \|\Pi_\perp \delta
\Phi\|_{H^{1}_L} \|\Pi_\perp \delta   \Phi\|_{L^2}
& \\
 {} + C  \|\Pi_\perp \Phi\|_{H^{1}_L} \|\Pi_\perp
\delta\Phi\|_{H^{1}_L}^{3/2}   \|\Pi_\perp \delta   \Phi\|_{L^2}^{1/2}
& \\
 {} + C \|\Pi_\perp\delta   \Phi\|_{L^2} (
\|\Pi_0 \Phi\|_{L^2}  \|\Pi_\perp \delta\Phi\|_{H^{1}_L} +
\|\Pi_0\delta \Phi\|_{L^2}
\|\Pi_\perp \Phi\|_{H^{1}_L}
).&
\end{aligned}
$$
We then deduce (using the same argument as in the construction of a
weak solution page~\pageref{deltaprimterm} for the~$ \Delta'_{L}$  term), that
$$
\begin{aligned}
 \| \delta \Phi(t)\|_{L^2(\R\times \T)}^2 -\|\delta  \Phi^0\|_{L^2(\R\times \T)}^2
+  2\nu \int_0^t  \| \Pi_0(\delta \Phi)'(t')\|_{L^2(\R\times \T)}  dt'+{\nu\over C}
\int_0^t \| \Pi_\perp \delta \Phi (t')\|_{H^1_L }^2dt'  &\\
 \leq  C_\nu \int_0^t \left( \|\Pi_\perp \Phi (t')\|_{H^{3/2}_L}^2 +
\|\Pi_\perp \Phi (t')\|_{H^{1}_L}^4 \right) \|\Pi_\perp
\delta \Phi (t')\|_{L^2}^2   dt'&\\
 {} +  C_\nu \int_0^t
\left( \|\Pi_0 \Phi (t')\|_{L^2}^2  \|\Pi_\perp
\delta \Phi (t')\|_{L^2}^2 +  \|\Pi_\perp \Phi (t')\|_{H^{1}_L}^2  \|\Pi_0\delta
\Phi\|^2_{L^2}\right)
  dt' ,&
\end{aligned}
$$
using the embedding $H^1_L\subset L^2$.
Gronwall's lemma yields
$$
\| \delta \Phi(t)\|_{L^2} ^2\leq \| \delta \Phi^0\|_{L^2}^2\exp\Bigl(  C_\nu   \int_0^t
\bigl(1+
\|\Pi_0 \Phi(\tau)\|^2_{ L^{2} (\R\times \T)} +
\|\Pi_\perp\Phi(\tau)\|_{ H^{1}_L}^4  +
\|\Pi_\perp\Phi(\tau)\|^2_{ H^{3/2}_L} \bigr)d\tau \Bigr),
$$
and the conclusion comes from the fact that~$\Pi_\perp \Phi$ belongs
to~$ L^4([0,T],H^1_L)$ by interpolation
between~$L^\infty([0,T],H^{1/2}_L) $ and~$
L^2([0,T],H^{3/2}_L)$, along with the energy estimate on~$\Phi $. The  proposition is
proved.
\end{proof}

\subsection{Propagation of the geostrophic regularity}\label{regularity}
The following regularity result for the geostrophic equation is inspired by the
Weyl-H{\"o}rmander symbolic calculus, even if that theory does not seem to
be appliable directly due to the possible singularity at $x_1=0$.

Using the formula giving $\Pi_0$ in Proposition~\ref{explicitPi0}, we first see that the geostrophic equation
$$\d_t  \Phi-\nu \Pi_0\Delta'\Pi_0  \Phi=0$$
can be brought back (at least formally) to the scalar equation
$$\d_t u_2 -\nu  D(DD^T+Id)^{-1} D^T\Delta u_2=0,$$
where we recall that~$D$ is the differential operator defined by $\displaystyle D\cdot=\d_1\left({\cdot 
\over
\beta x_1}\right)$. 
Then by a simple change of variables this scalar equation becomes
$$\d_t \varphi -\nu A\varphi =0,$$
where $A$ is some self-adjoint scalar pseudo-differential operator (possibly singular at
$x_1=0$), the principal symbol of which is given by
$$ a(x_1,\xi_1)={\xi_1^4 \over \beta^2 x_1^2+\xi_1^2}$$
neglecting the possible singularity at $x_1=0$. Then, in order to propagate Sobolev regularity on $\varphi$, we
 should have to control some commutator of the type $[A,\nabla^s]$, which is not so easy because $A$ cannot
  be written simply in terms of the usual derivatives $\nabla$.

In order to find a convenient way to measure the regularity, we therefore use formally the results of symbolic calculus.
 Note that the Weyl-H\"ormander theory is used here just  to guide intuition,
 the result of propagation   being actually proved by explicit computations. 
 In order to determine the class of operators $A$ should belong to, we have first to characterize the metric.
Computing the partial derivatives of $a$ with respect to $x_1$ and $\xi_1$
$${\d_{x_1}a\over a}(x_1,\xi_1)=-{2\beta^2 x_1 \over \beta^2x_1^2+\xi_1^2},\quad
{\d_{\xi_1}a\over a}(x_1,\xi_1)={4\over \xi_1} -{2 \xi_1
\over
\beta^2x_1^2+\xi_1^2}$$
shows that the H{\"o}rmander metric to be considered is the one associated to the harmonic
oscillator
$$g(dx_1,d\xi_1)=\beta^2\frac{ dx_1^2}{1+x_1^2}+ \frac{d\xi_1^2}{1+\xi_1^2}\cdotp$$
Then it is natural to measure the regularity by powers of the harmonic oscillator, and
therefore to study the propagation equation
$$\d_t (-\d^2_{x_1x_1} +\beta^2 x_1^2)^s \varphi - \nu A (-\d^2_{x_1x_1} +\beta^2
x_1^2)^s\varphi=\nu [(-\d^2_{x_1x_1} +\beta^2 x_1^2)^s,A]
\varphi
$$
The fundamental result of the Weyl-H{\"o}rmander theory states the following: if $A$ is
a pseudo-differential operator (meaning in particular
that there is no singularity at $x_1=0$), the commutator  occuring in the right-hand side
of the previous equation is a pseudo-differential
operator of lower order (for the metric
$g$), meaning that we expect
$\| [(-\d^2_{x_1x_1} +\beta^2 x_1^2)^s,A] \varphi\|_{L^2(\R)}$ to be controlled by
$\|(-\d^2_{x_1x_1} +\beta^2 x_1^2)^s \varphi\|_{L^2(\R)}$ and
$-\big((-\d^2_{x_1x_1} +\beta^2 x_1^2)^s \varphi| A(-\d^2_{x_1x_1} +\beta^2 x_1^2)^s
\varphi\big)$.

Nevertheless, as we are not able to prove in a simple way that there is no singularity at
$x_1=0$, we shall not use the general theory of
pseudo-differential operators and will proceed instead using explicit computations.
We have seen in the previous chapter that the family $(\Psi_{n,0,0})_{n\in \N}$ defined by
$$ \Psi_{n,0,0} (x_1)=\left(
  \begin{array}{c}
  \displaystyle - \sqrt{ (n+1)\over 2(2n+1)} \psi_{n-1}(x_1)-
\sqrt{ n\over 2(2n+1)}\psi_{n+1}(x_1)\\
0
   \\
\displaystyle  \sqrt{ (n+1)\over 2(2n+1)} \psi_{n-1}(x_1)-
\sqrt{ n\over 2(2n+1)}\psi_{n+1}(x_1)
\end{array}
\right)$$
constitutes an Hermitian basis of $\Ker L$, and that, due to the properties of the
Hermite functions,
$$\forall n\in \N,\quad \left\| (-\d^2_{x_1x_1} +\beta^2 x_1^2)^s
\Psi_{n,0,0}\right\|_{L^2(\R)} \sim (1+n)^s \left\|
\Psi_{n,0,0}\right\|_{L^2(\R)}.$$
Therefore
it is natural to study the propagation of the $H^s_L$
norm of~$\Pi_0\Phi$, recalling that\label{bvarphi-not}
$$
\|\Pi_0\Phi\|_{H^s_L } ^2 =\sum_{n \in \N} (1+n)^s |\underline\varphi_n|^2,
$$
where we have defined
$$\underline\varphi_n=(\Psi_{n,0,0}| \Pi_0\Phi ).$$
In the following we will denote by $N_s$ the operator defined by\label{Ns-not}
$$\forall n\in \N,\quad 
N_s \Pi_0\Phi = \sum (1+n)^s \underline\varphi_n  \Psi_{n,0,0},$$
so that
$$\|N_s  \Pi_0\Phi\|_{L^2(\R)}=\|  \Pi_0\Phi\|_{H^{2s}_L}.$$
We have
$$
\d_t N_s \Pi_0\Phi- \nu\Pi_0 \Delta' N_s \Pi_0\Phi= \nu [N_s,\Pi_0
\Delta'] \bar\Phi.
$$
\medskip
$\bullet$
We start by computing $\Pi_0 \Delta' \Psi_{n,0,0}$. From (\ref{psi-identities}) page \pageref{psi-identities} we deduce
that
$$\d_{x_1} \Psi_{n,0,0}'(x_1)=\left(
  \begin{array}{c}
  \displaystyle 0\\

   \\
\displaystyle  \sqrt{\beta (n+1)(n-1)\over 4(2n+1)}\psi_{n-2}(x_1)-\sqrt{\beta
n(n+1)\over 2n+1 } \psi_n(x_1)+\sqrt{\beta n(n+2)\over 4(2n+1)}
\psi_{n+2}(x_1)
\end{array}
\right)$$
and
$$\d^2_{x_1 x_1}\Psi_{n,0,0}'(x_1)={\beta\over 2\sqrt{2(2n+1)}} \left(
  \begin{array}{c}
  \displaystyle 0\\

   \\
\displaystyle  \sqrt{(n+1)(n-1)(n-2)} \psi_{n-3}(x_1)-(3n-1)\sqrt{n+1} \psi_{n-1}(x_1)
\\+(3n+4) \sqrt{n} \psi_{n+1}(x_1)-\sqrt{n(n+2)(n+3)}
\psi_{n+3}(x_1)
\end{array}
\right)$$
with the usual convention that $\psi_n\equiv 0$ for $n<0$.
Therefore, using the orthogonality of the Hermite functions in $L^2(\R)$, we get
\begin{equation*}
\label{psi-n-image}
\Pi_0\Delta' \Psi_{n,0,0}=\alpha_{n}^{(-4)} \Psi_{n-4,0,0} +\alpha_{n}^{(-2)}
\Psi_{n-2,0,0}+\alpha_{n}^{(0)} \Psi_{n,0,0}+\alpha_{n}^{(2)}
\Psi_{n+2,0,0}+\alpha_{n}^{(4)} \Psi_{n+4,0,0}
\end{equation*}
with
\begin{equation}
\label{alpha-def}
\begin{aligned}
\alpha_{n}^{(-4)} &= -{\beta \over 4} \sqrt{(n-4)(n-2)(n-1)(n+1)\over (2n+1)(2n-7)}\virgp\\
\alpha_{n}^{(-2)}&={\beta\over 4} (4n-2) \sqrt{(n-2)(n+1)\over (2n+1)(2n-3)}\virgp\\
\alpha_n^{(0)} &=-{\beta\over 4} {6n^2+6n-1 \over (2n+1)}\virgp\\
\alpha_n^{(2)}&={\beta\over 4} (4n+6) \sqrt{n(n+3)\over (2n+1)(2n+5)}\virgp\\
\alpha_n^{(4)}&=-{\beta\over 4} \sqrt{n(n+2)(n+3)(n+5)\over (2n+1)(2n+9)}\cdotp
\end{aligned}
\end{equation}

\medskip
$\bullet$
From the previous computation we deduce that
$$\begin{array}{rl}
[N_s, \Pi_0 \Delta'] \Psi_{n,0,0}&=((n-3)^s-(n+1)^s)\alpha_{n}^{(-4)} \Psi_{n-4,0,0}
+((n-1)^s-(n+1)^s)\alpha_{n}^{(-2)}
\Psi_{n-2,0,0}\\
&+((n+2)^s-(n+1)^s)\alpha_{n}^{(2)}
\Psi_{n+2,0,0}+((n+5)^s-(n+1)^s)\alpha_{n}^{(4)} \Psi_{n+4,0,0}
\end{array}$$
Thus, using the definition (\ref{alpha-def}) of the coefficients $\alpha$, we get
$$|\alpha_n|\leq C (n+1)$$
and
$$
\left\| [N_s, \Pi_0 \Delta'] \Psi_{n,0,0} \right\|_{L^2(\R)} \leq C_s (n+1)^s\leq C_s \|N_s
\Psi_{n,0,0} \|_{L^2(\R)}.$$
Because of the quasi-orthogonality of $([N_s, \Pi_0 \Delta'] \Psi_{n,0,0})_{n\in \N}$, we
have actually the more general commutator estimate
$$
\begin{aligned}
\| &[N_s, \Pi_0 \Delta'] \Pi_0\Phi \|^2_{L^2(\R)}\\
&=\sum _n\Big|((n+1)^s-(n-3)^s) \alpha_{n-4}^{(-4)} \underline\varphi_{n+4}+((n+1)^s-(n-1)^s)
\alpha_{n-2}^{(-2)} \underline\varphi_{n+2}\\
&\qquad\quad  +((n+1)^s-(n+3)^s) \alpha_{n+2}^{(2)} \underline\varphi_{n-2}+((n+1)^s-(n+5)^s)
\alpha_{n+4}^{(4)} \underline\varphi_{n-4}\Big|^2\\
&\leq C_s \sum _n (n+1)^{2s} \left( |\underline\varphi_{n+4}|^2+|\underline\varphi_{n+2}|^2+|\underline\varphi_{n-2}|^2+|\underline\varphi_{n-4}|^2\right)\\
& \leq C_s \sum_n (n+1)^{2s} |\varphi_n|^2
\end{aligned}
$$
which can be rewritten
\begin{equation}
\label{commutator}
\left\| 
[N_s, \Pi_0 \Delta'] \Pi_0\Phi
 \right\|_{L^2(\R)}
\leq C_s \| N_s \Pi_0\Phi\|_{L^2(\R)}.
\end{equation}
Note that, due to the particular choice of the operator $N_s$, there is some additional
cancellation, meaning that the commutator $[N_s, \Pi_0
\Delta']$ which is expected to be a pseudodifferential operator of order $(2s+1)$ is
actually of order $2s$.

\medskip
$\bullet$
It is now very easy to propagate regularity using Gronwall's lemma. We
recall that
$$\d_t N_s \Pi_0\Phi- \nu\Pi_0 \Delta' N_s \Pi_0\Phi = \nu [N_s,\Pi_0 \Delta'] \Pi_0\Phi,$$
from which we deduce that
$$\| N_s \Pi_0\Phi (t)\|_{L^2(\R)}^2 +\nu \int_0^t \| \nabla (N_s \Pi_0\Phi)'(\tau)\|^2_{L^2(\R)} d\tau \leq \| N_s \Pi_0 \Phi^0
\|_{L^2(\R)}^2+C_s
\int_0^t \| N_s \Pi_0\Phi (\tau)\|_{L^2(\R)}^2d\tau $$
and finally
\begin{equation}
\label{prop-Ns}
\| N_s \Pi_0\Phi (t)\|_{L^2(\R)}^2 +\nu \int_0^t \| \nabla (N_s \Pi_0\Phi)'(\tau)\|^2_{L^2(\R)} d\tau \leq \| N_s \Pi_0 \Phi^0 \|_{L^2(\R)}^2\exp
(C_s t).
\end{equation}
This concludes the proof. 
\subsection{Local strong solutions}
\label{sectionstrong}
In this section we are going to prove the existence of unique, strong
solutions for smooth  enough initial data. As in the case of weak
solutions discussed in Section~\ref{sectionweak} above, we will not write the full proof, but
detail the estimates enabling one to use the usual Fujita-Kato  theory of strong
solutions to the 3D Navier-Stokes equations (see~\cite{cdggbook} for  instance).

\medskip
$\bullet$ {\it Global existence of strong solutions for small data}

We prove here that under a suitable smallness assumption there exists a  (unique) global
strong
solution to $(SW_0)$ such that~$ \Pi_\perp \Phi$ belongs 
to~$L^\infty(\R^+,
H^{1/2}_L)
\cap L^2(\R^+, H^{3/2}_L)$.

As previously we start from the solutions $\Phi_{(N)}$ of the  approximation scheme
(\ref{Friedrichs}). We have
of course
$$ \frac d{dt}
\|\Pi_0 \Phi_{(N)}(t) \|_{L^2(\R\times \T)}^2 +2\nu \|\nabla (\Pi_0  \Phi_{(N)})'(t)
\|_{L^2(\R\times
\T)}^2 \leq 0 $$
and, by Proposition \ref{trilinear},
$$\begin{aligned}
  \frac d{dt}&\| \Pi_\perp \Phi_{(N)}(t)\|_{  H^{1/2}_L  }^2
- 2\nu \left ( \Pi_\perp  \Phi_{(N)} | \Delta'_L  \Pi_\perp   \Phi_{(N)}\right)_{ 
H^{1/2}_L } (t)
\\
&\leq -  2\left(  \Pi_\perp  \Phi_{(N)} | \Pi_\perp Q_L ( \Phi_{(N)},  \Phi_{(N)})
\right)_{  H^{1/2}_L }(t)\\
&\leq C \Bigl(\| \Pi_\perp \Phi_{(N)}\|_{   H^{3/2}_L }^2\| \Pi_\perp  \Phi_{(N)}\|_{ 
H^{1/2}_L }
+  \|\Pi_0 \Phi_{(N)}\|_{L^2(\R\times \T)}\|  \Pi_\perp
\Phi_{(N)}\|_{   H^{3/2}_L }
\| \Pi_\perp   \Phi_{(N)}\|_{  H^{1/2}_L(\R\times
\T)}\Bigr)(t).
\end{aligned}
$$
By Lemma~\ref{parabolic} and the obvious embedding $H^{3/2}_L\subset  H^{1/2}_L$,
that inequality  can be written
\begin{equation}
\label{eneplusplus}
\begin{aligned}
  \frac d{dt} \| \Pi_\perp \Phi_{(N)}(t)\|_{  H^{1/2}_L  }^2
+ 2\frac \nu {C_{1/2}} \| \Pi_\perp \Phi_{(N)} (t)\|_{  H^{3/2}_L  }^2    \\
\leq C \|
\Pi_\perp \Phi_{(N)}\|_{   H^{3/2}_L }^2 \Bigl(\| \Pi_\perp  \Phi_{(N)}\|_{  H^{1/2}_L }
+  \|\Pi_0 \Phi_{(N)}\|_{L^2(\R\times \T)}\Bigr) (t) .
\end{aligned}
\end{equation}
As usual we notice that this inequality is useful only  if~$\|
 \Pi_\perp \Phi_{(N)}(t)\|_{  H^{1/2}_L}$ and~$ \|\Pi_0  \Phi_{(N)}\|_{L^2(\R\times \T)}
$ are small, which is a typical
 phenomena of global  results under a smallness condition.

Define
$$D_N =\left\{ t\in \R^+\,/\, \forall t'\leq t,\, \| \Pi_0
  \Phi_{(N)}(t')\|_{L^2(\R\times \T)}^2
+\| \Pi_\perp \Phi_{(N)}(t')\|_{  H^{1/2}_L}^2\leq \left({\nu \over 
2CC_{1/2}}\right)^2\right\},$$
where  $C$ is the constant appearing in~(\ref{eneplusplus}), and
let us impose the following
smallness assumption on the initial data:
$$
\| \Pi_0 \Phi_{(N)}^0\|_{L^2(\R\times  \T)}^2+\|\Pi_\perp\Phi_{(N)}^0\|^2_{   H^{1/2}_L
}\leq  \left({\nu \over 4CC_{1/2}}\right)^2.
$$
Then clearly~$ D_N$ is not empty.
 By construction, $\Pi_0\Phi_{(N)} $ belongs to~$ C(\R^+, L^2(\R\times  \T))$ and
$\Pi_\perp \Phi_{(N)}$ belongs to~$C(\R^+,  H^{1/2}_L  )$ thus $D_N$ is  a closed set.

Denote by~$ T_N = \max D_N$.
If $T_N<+\infty$, then
$$
\| \Pi_0
  \Phi_{(N)}(T_N)\|^2_{L^2(\R\times \T)} +  \|\Pi_\perp    \Phi_{(N)}(T_N)\|^2_{ 
H^{1/2}_L}\leq  \left({\nu \over 2CC_{1/2}}\right)^2
$$
and we deduce from (\ref{eneplusplus}) that
$$
{d\over dt}\left( \| \Pi_0 \Phi_{(N)}\|_{L^2(\R\times  \T)}^2+\|\Pi_\perp
\Phi_{(N)}\|^2_{  H^{1/2}_L(\R\times
\T)}\right)(T_N)< 0,
$$
which is in contradiction with the definition of~$T_N=\max D_N$.  Therefore $T_N=+\infty$.

Then we deduce immediately that for all $t\in \R^+$
$$
\|\Pi_\perp \Phi_{(N)}(t)\|^2_{  H^{1/2}_L}+ \frac \nu {C_{1/2}} \int_0^t
\|\Pi_\perp \Phi_{(N)}(t')\|^2_{  H^{3/2}_L} dt'
\leq  \|\Pi_\perp \Phi_{(N)}^0\|^2_{  H^{1/2}_L}.
$$
Up to the extraction of a subsequence that converges to a Leray  solution $\Phi$ of
$(SW_0)$, the
previous estimate implies that
$$
\Pi_\perp \Phi\in L^\infty(\R^+,  H^{1/2}_L)\cap L^2(\R^+,
H^{3/2}_L).
$$
The strong-weak stability principle established in the previous
section provides then the  uniqueness of such a
solution.

\bigskip
$\bullet$ {\it Local existence of strong solutions}

Let us now consider the case of large data. The idea (see for
 instance~\cite{cdggbook})
 is to split $\Phi_{(N)}$ in two parts as follows
$$\Phi_{(N)}= \Phi_{(N)}^< +\Phi_{(N)}^>$$
where  $\Phi_{(N)}^<$ is the unique solution to
\begin{equation}
\label{semigroup}
\begin{aligned}
\,&\d_t \Phi_{(N)}^< -\nu \Delta'_L \Phi_{(N)}^< =0\\
&\Phi_{(N)}^<(0)=\sumetage{(n,k,j) \in S}{(n+k^2)^{1/2} \leq A}  \Pi_{n,k,j} \Pi_\perp
\Phi_{(N)}^0+\Pi_0
\Phi_{(N)},
\end{aligned}
\end{equation}
and $A>0$ is a truncation parameter to be determined (independent of~$
N$).
Using Proposition \ref{parabolic}, it is easy to check that
$$\begin{aligned}
\| \Pi_\perp \Phi_{(N)}^<\|_{L^\infty(\R^+,L^2(\R\times \T))} \leq
C
\|\Phi^0\|_{L^2(\R\times \T)},\\
\| \Pi_\perp \Phi_{(N)}^<\|_{L^\infty(\R^+,H^s_L)} \leq C A^s
\|\Phi^0\|_{L^2(\R\times \T)}.
\end{aligned}
$$

 By (\ref{Friedrichs}) and (\ref{semigroup}) we deduce the equation  satisfied by
$\Phi_{(N)}^>$~:
\begin{equation}
\label{lin-Friedrichs}
\begin{aligned}
\,&\d_t \Phi_{(N)}^>-\nu \Delta'_L \Phi_{(N)}^> +J_NQ_L(   \Phi_{(N)}^>,\Phi_{(N)}^>)+2
J_N Q_L (  \Phi_{(N)}^>,\Phi_{(N)}^<)=-J_N
Q_L(  \Phi_{(N)}^<,\Phi_{(N)}^<)\\
&\Phi_{(N)}^> (0)=\Phi_{(N)}^0 - \left(\sumetage{(n,k,j) \in  S}{(n+|k|)^{1/2} \leq A}
\Pi_{n,k,j} \Pi_\perp \Phi_{(N)}^0+\Pi_0
\Phi_{(N)} \right).
\end{aligned}
\end{equation}
We are going to show that~$\Phi_{(N)}^>$  remains small in~$  H^{1/2}_L  $ on a time
interval which does
not depend on~$N$.

Let us write an energy inequality in~$ H^{1/2}_L $ on  (\ref{lin-Friedrichs}): we have
$$
\longformule{
\frac12 {d\over dt} \|\Phi_{(N)}^>(t)\|^2_{  H^{1/2}_L } +  \frac \nu {C_{1/2}}
\|\Phi_{(N)}^>(t)\|^2_{  H^{3/2}_L }
=-(J_N Q_L(  \Phi_{(N)}^>,\Phi_{(N)}^>)|
\Phi_{(N)}^>)_{  H^{1/2}_L }(t)  }{ + 2 (J_N Q_L (   \Phi_{(N)}^>,\Phi_{(N)}^<)|
\Phi_{(N)}^>)_{  H^{1/2}_L }(t)- (J_N
Q_L(  \Phi_{(N)}^<,\Phi_{(N)}^<) | \Phi_{(N)}^>)_{  H^{1/2}_L }(t),
}
$$
and Proposition~\ref{trilinear} yields
$$
\begin{aligned}
&
{d\over dt} \|\Phi_{(N)}^>(t)\|^2_{  H^{1/2}_L } +  2\frac \nu {C_{1/2}}
\|\Phi_{(N)}^>(t)\|^2_{  H^{3/2}_L } \leq  C  \|\Phi_{(N)}^>(t)\|_{
  H^{1/2}_L } \|\Phi_{(N)}^>(t)\|^2_{  H^{3/2}_L }\\
 &+  C  \|\Phi_{(N)}^>(t)\|_{  H^{1}_L } \bigl(
 \|\Phi_{(N)}^>(t)\|_{  H^{3/2}_L } \|\Pi_\perp \Phi_{(N)}^<(t)\|_{   H^{1}_L }
+  \|\Phi_{(N)}^>(t)\|_{  H^{1}_L } \|\Pi_\perp \Phi_{(N)}^<(t)\|_{
  H^{3/2}_L }\bigr)\\
& + C  \|\Phi_{(N)}^>(t)\|_{  H^{1}_L } \bigl(
 \|\Phi_{(N)}^>(t)\|_{  H^{1}_L } \|\Pi_0   \Phi_{(N)}(t)\|_{L^2(\R\times\T)}
+  \|\Pi_\perp \Phi_{(N)}^<(t)\|_{
   H^{1}_L }  \|\Pi_\perp \Phi_{(N)}^<(t)\|_{  H^{3/2}_L }\bigr)\\
& + C
 \|\Phi_{(N)}^>(t)\|_{  H^{1/2}_L } \|\Pi_0   \Phi_{(N)}(t)\|_{L^2(\R\times\T)}
 \|\Pi_\perp \Phi_{(N)}^<(t)\|_{  H^{3/2}_L }.
\end{aligned}
$$
Now consider the set
$$
D_N =\left\{ t\in \R^+\,/\, \forall t'\leq t,\, \|
 \Phi_{(N)}^>(t')\|_{  H^{1/2}_L}\leq {\nu
 \over CC_{1/2}}\right\},
$$
(where $C$ is the constant appearing in the right-hand side of the previous inequality)
and~$ T_N=\sup D_N $. We are going to prove that there exists $T>0$  such that
$$\forall N\in \N^*, \quad T_N \geq T.$$
We notice that if~$ A$ is chosen large enough (independently of~$ N$),
then~$ D_N$ is not empty: we can indeed choose~$ A$ so that
\begin{equation}\label{Alargeenough}
 \|\Phi_{(N)}^>(0)\|_{  H^{1/2}_L} \leq {\nu \over 4CC_{1/2}}, \quad  \forall N
 \in \N.
\end{equation}
As long as~$ t \leq T_N$, we can write the previous inequality in the
following way:
$$
\begin{aligned}
{d\over dt} \|\Phi_{(N)}^>(t)\|^2_{  H^{1/2}_L } +  \frac \nu {2C_{1/2}}
\|\Phi_{(N)}^>(t)\|^2_{  H^{3/2}_L } &\leq \frac C{\nu^3}   \|\Phi_{(N)}^>(t)\|^2_{
  H^{1/2}_L }  \|\Pi_\perp \Phi_{(N)}^<(t)\|_{
  H^{1}_L }^4 \\
&\!\!\!\!\!\!+\frac C\nu  \|\Phi_{(N)}^>(t)\|^2_{
  H^{1/2}_L } \!\left(  \|\Pi_\perp \Phi_{(N)}^<(t)\|_{
  H^{3/2}_L }^2\|+\|\Pi_0  \Phi_{(N)}(t)\|_{L^2(\R\times\T)}^2\right)  \\
  &\!\!\!\!\!\!+\frac
C{\nu^{1/3}}  \|\Phi_{(N)}^>(t)\|_{
  H^{1/2}_L }^{2/3} \|\Pi_\perp \Phi_{(N)}^<(t)\|_{
  H^{1}_L }^{4/3}  \|\Pi_\perp \Phi_{(N)}^<(t)\|_{
  H^{3/2}_L }^{4/3}\\
&\!\!\!\!\!\!+ C
 \|\Phi_{(N)}^>(t)\|_{  H^{1/2}_L } \|\Pi_0   \Phi_{(N)}(t)\|_{L^2(\R\times\T)}
 \|\Pi_\perp \Phi_{(N)}^<(t)\|_{  H^{3/2}_L }
\end{aligned}
$$
So we get
$$
{d\over dt} \|\Phi_{(N)}^>(t)\|^2_{  H^{1/2}_L } +  \frac \nu {2C_{1/2}}
\|\Phi_{(N)}^>(t)\|^2_{  H^{3/2}_L } \leq  (1+\|\Phi_{(N)}^>(t)\|^2_{   H^{1/2}_L})
F(A,\Phi^0,\nu)  ,
$$
where
$$
F(A,\Phi^0,\nu) = C\left( 1+\frac 1{\nu^3}\right) (1+A^4)  (1+\|\Phi^0\|_{L^2(\R\times\T)}^4) .
$$
Gronwall's lemma enables us to infer that for all~$t \leq T_N$,
$$
 \|\Phi_{(N)}^>(t)\|^2_{  H^{1/2}_L }
\leq
(1+\|\Phi_{(N)}^>(0)\|^2_{  H^{1/2}_L } )\exp \left(t F(A,\Phi^0,\nu)  \right)
-1
$$
Since~$A$ is chosen so that~(\ref{Alargeenough}) is satisfied, it
suffices now to choose~$ T$ in such a way that
$$
\exp \left(T F(A,\Phi^0,\nu) \right) \leq 2,\quad \exp \left(T  F(A,\Phi^0,\nu) \right)
-1\leq \frac12 \left({\nu \over C C_{1/2}}\right)^2
$$
so that
 for any~$ t \leq T$,
$$
\|\Phi_{(N)}^>(t)\|^2_{  H^{1/2}_L } \leq   \left({\nu \over  CC_{1/2}}\right)^2
$$
hence necessarily
$$\forall N\in \N^*,\quad T_N \geq T.$$

\medskip
$\bullet$
Gathering those results, we infer that any limiting point of  $\Phi_{(N)}^<+\Phi_{(N)}^>$
(in particular
any Leray solution to $(SW_0)$) satisfies
$$\begin{aligned}
\Pi_0\Phi \in L^\infty(\R^+,L^2(\R\times \T)) , \quad (\Pi_0\Phi)'
\in L^2(\R^+,H^1(\R\times \T))\\
\Pi_\perp \Phi \in L^\infty([0,T],  H^{1/2}_L ) \cap L^2([0,T],   H^{3/2}_L  ).
\end{aligned}
$$
The weak-strong stability principle gives then the uniqueness of such a  solution on
$[0,T]$.

\subsection{Propagation of the ageostrophic regularity}\label{regularitystrong}
In order to construct a strong approximation of the filtered solution to the Saint-Venant system in the
next chapter, we will actually need further regularity on the solution of the limit filtered system, which
 is obtained in a very standard way from the trilinear estimates stated in Proposition \ref{trilinear}. So suppose 
 that~$\Pi_\perp \Phi^0 $ belongs to~$ H^s_L$, with~$1/2 \leq s \leq 1$, and consider as previously the sequence $(\Phi_{(N)})$ of approximate solutions to $(SW_0)$ 
defined by~(\ref{Friedrichs}). From Lemma~\ref{parabolic} and the last estimate in 
Proposition~\ref{trilinear},  we deduce that  for all~$s\leq 1$,
$$ 
\begin{aligned}
 {d\over dt} &\| \Pi_\perp \Phi_{(N)}(t)\|_{H^s_L}^2 +{ 2\nu \over C_s}\| \Pi_\perp\Phi_{(N)}(t)\|_{H^{s+1}_L}^2\\
& \leq 
C \| \Pi_\perp \Phi_{(N)}(t)\|_{H^s_L} \| \Pi_\perp \Phi_{(N)}(t)\|_{H^{s+1}_L}\left(\| \Pi_\perp \Phi_{(N)}(t)\|_{H^{3/2}_L}+\|\Pi_0 \Phi_{(N)}(t)\|_{L^2(\R\times \T)} \right)\\
& \leq {\nu \over C_s} \| \Pi_\perp\Phi_{(N)}(t)\|_{H^{s+1}_L}^2 +{C\over \nu} \| \Pi_\perp \Phi_{(N)}(t)\|_{H^s_L} ^2\left(\| \Pi_\perp \Phi_{(N)}(t)\|_{H^{3/2}_L}+\|\Pi_0 \Phi_{(N)}(t)\|_{L^2(\R\times \T)} \right)^2\,,
\end{aligned}
$$
and thus by Gronwall's lemma
$$ 
\begin{aligned}
\| \Pi_\perp \Phi_{(N)}(t)\|_{H^s_L}^2 &+{ \nu \over C_s}\int_0^t \| \Pi_\perp\Phi_{(N)}(t')\|_{H^{s+1}_L}^2dt'\\
 &\leq  \| \Pi_\perp \Phi_{(N)}^0\|_{H^s_L}^2\exp\left( \frac C\nu \int_0^t
  \left(\| \Pi_\perp \Phi_{(N)}(t')\|_{H^{3/2}_L}+\|\Pi_0 \Phi_{(N)}(t')\|_{L^2(\R\times \T)} \right)^2 dt'\right)\,.
\end{aligned}
$$
Taking limits as $N\to \infty$ in the previous inequality shows that
$$ \Pi_\perp \Phi \in L^\infty_{loc}([0,T^*[, H^s_L) \cap L^2_{loc}([0,T^*[, H^{s+1}_L),$$
and that
$$ 
\begin{aligned}
\| \Pi_\perp \Phi(t)\|_{H^s_L}^2 &+{ \nu \over C_s}\int_0^t \| \Pi_\perp\Phi(t')\|_{H^{s+1}_L}^2dt' \\
&\leq  \| \Pi_\perp \Phi^0\|_{H^s_L}^2\exp\left( \frac C\nu \int_0^t \left(\| \Pi_\perp \Phi(t')\|_{H^{3/2}_L}+\|\Pi_0 \Phi(t')\|_{L^2(\R\times \T)} \right)^2 dt'\right)\,,
\end{aligned}
$$
which proves the  propagation of regularity result, and completes the  proof of Theorem~\ref{limit-systemallbeta}. 
\hfill $\Box$

\section{Proof of Theorem 3}\label{theoremgeneric}\setcounter{equation}{0}
In this section we shall prove Theorem~\ref{limit-systemgeneric}:   Paragraph~\ref{globalgeneric} is devoted
to the global wellposedness result in~$L^{2}$,  while the propagation of regularity results are given in
 Paragraphs~\ref{sectionpropaghs}
and~\ref{sectionpropag}.

\subsection{Global wellposedness}\label{globalgeneric}
In this section we shall prove the first part of Theorem~\ref{limit-systemgeneric}, namely the fact that except
for a countable number of~$\beta$, the limit system is globally wellposed in~$L^{2}$.  This turns out
to be an easy matter in view of the resonance results obtained in Section~\ref{Resonance} in the previous chapter.

Indeed Proposition~\ref{resonances-prop}, page~\pageref{resonances-prop} indicates that except for
a countable set of values for~$\beta$, the limit system reduces to the following:
$$
\begin{aligned}
\d _t \Pi_0 \Phi -\nu \Pi_0 \Delta'_L  \Phi=0,\\
\d _t \Pi_R \Phi +2Q_L( \Pi_0 \Phi,\Pi_R \Phi) -\nu \Pi_R
\Delta'_L
\Phi=0,\\
\d _t \Pi_M \Phi +2Q_L( \Pi_0 \Phi,\Pi_M \Phi)-\nu \Pi_M \Delta'_L
\Phi=0,\\
\d _t \Pi_P \Phi +2Q_L( \Pi_0 \Phi,\Pi_P \Phi)-\nu \Pi_P \Delta'_L
\Phi=0,\\
\d _t \Pi_K \Phi +2Q_L( \Pi_0 \Phi,\Pi_K \Phi)+Q_L(\Pi_K \Phi,\Pi_K  \Phi)-\nu \Pi_K
\Delta'_L
\Phi=0.\\
\end{aligned}
$$
So the limit system is a linear  equation on all modes but Kelvin modes; Kelvin modes being essentially one-dimensional,
 it will be easy to prove the wellposedness of the system. In fact the only point to be proved is the uniqueness of the
 solution, since existence was proved in the previous section. Uniqueness will be immediate in the case of all non-Kelvin
 modes so let us concentrate on the equation on~$\Pi_{K} \Phi$. Let~$\Phi$ and~$\Phi_*$ be two solutions, and
 define~$\delta \Phi = \Phi_* - \Phi$. Then
 $$
 \longformule{
 \d _t \Pi_K\delta \Phi +2Q_L( \Pi_0\delta \Phi,\Pi_K \delta \Phi)
 + 2 Q_L( \Pi_0  \Phi,\Pi_K \delta \Phi) +  2Q_L( \Pi_0 \delta \Phi,\Pi_K \Phi)
 +Q_L(\Pi_K \delta \Phi,\Pi_K  \delta \Phi)}{
  + 2Q_L(\Pi_K \Phi,\Pi_K  \delta \Phi)  
 -\nu \Pi_K
\Delta'_L
\delta \Phi=0.}
 $$
 Now let us write an energy estimate in~$L^{2}$, in the spirit of the computations of Section~\ref{sectionweak}
 above. We note that in the case of~$\Pi_{K}$, the decomposition on eigenmodes of~$L$ simply corresponds
 to the Fourier decomposition, so that 
$$\left( \Pi_K \delta \Phi | Q_L (\Pi_K \Phi,\Pi_K\delta \Phi)\right)_{L^{2}(\R\times  \T)} =\left( \Pi_K
\delta \Phi | Q (\Pi_K
\Phi,\Pi_K\delta \Phi)\right)_{L^{2}( \R\times \T)} .$$
Moreover in that case, the usual Sobolev spaces~$H^{s}$ coincide with the~$H^{s}_{L} $ spaces, since~$n = 0$.
Using the fact that~$\Pi_{0}$ projects onto~$x_{2}$-independent functions and that the dependence in~$x_{1}$ 
of~$\Pi_{K} \delta \Phi$ is that of the Gaussian~$\psi_{0}$, it is easy to see that
$$
\left| \left(\Pi_K \delta \Phi | Q_L( \Pi_0 \Phi,\Pi_K \delta \Phi) \right)_{L^{2}(\R\times  \T)}\right|
\leq  \|\Pi_{K} \delta \Phi\|_{H^{1} (\R \times \T)} \|\Pi_0 \Phi\|_{L^{2}(\R)}\|\Pi_K \delta \Phi\|_{L^{\infty}(\R;L^{2}(\T))},
$$
and that
\begin{eqnarray*}
\left| \left( \Pi_K
\delta \Phi | Q (\Pi_K
\Phi,\Pi_K\delta \Phi)\right)_{L^{2}( \R\times \T)}\right| &\leq & 
 \|\Pi_{K} \delta \Phi\|_{H^{1} (\R \times \T) } \|\Pi_K
\Phi\|_{L^{\infty}( \R\times \T)}  \|\Pi_{K} \delta \Phi\|_{L^{2}( \R\times \T) } \\
 &\leq & \frac \nu {2C_{1}}  \|\Pi_K
\delta \Phi\|_{H^{1} ( \R\times \T)}^2+ \frac C \nu  \|\Pi_K
\Phi\|_{L^{\infty}( \R\times \T)}^2  \|\Pi_{K} \delta \Phi\|_{L^{2}( \R\times \T) }^2.
\end{eqnarray*}
We recall indeed that
$$\Pi_K \: \cdot= \sum_{k\in \Z^*} (\Psi_{0,k,0 \:}| \:\cdot \:)_{L^2(\R\times \T)} \Psi_{0,k,0} \hbox{ where }  \Psi_{0,k,0}(x_1,x_2)={1\over \sqrt{4\pi}} e^{ikx_2}\left(
  \begin{array}{c}
  \displaystyle \psi_0(x_1)\\
0
   \\
\displaystyle \psi_0(x_1)
\end{array}
\right).$$ 
One sees easily that~$  \|\Pi_K
\Phi\|_{L^{\infty}( \R\times \T)} \leq C \|\Pi_K
\Phi\|_{  H^{1} (\R\times \T) }$, 
so finally an energy estimate (coupled with a Gronwall lemma) gives
 $$
 \longformule{
 \| \Pi_K\delta \Phi(t)\|_{L^{2}( \R\times \T) }^2 + \frac\nu C \int_{0}^{t} \| \Pi_K\delta \Phi(t')\|_{ H^{1}(\R\times \T)  }^2
 \: dt'  \leq  \| \Pi_K\delta \Phi(0)\|_{L^{2}( \R\times \T) }^2 }{\times \exp \left(
  \frac C \nu \int_{0}^{t} \|\Pi_K
\Phi(t')\|_{H^{1}}^2 \: dt' \right).}
 $$
 Since~$\Pi_K\delta \Phi(0) = 0$, uniqueness follows from the energy bound on~$\Pi_K\Phi$. Note that we have recovered
 here the usual, two-dimensional Navier-Stokes type estimates since in the case of purely Kelvin modes, the quadratic
 form and the spaces involved are the same as in the    Navier-Stokes case. In fact Kelvin modes are even one-dimensional (up to a multiplication by~$\psi_0(x_1$)) so it is possible to improve those estimates (that will
 be done in the last section of this chapter).

\subsection{Propagation of regularity} \label{sectionpropaghs} 
Let us now prove that if the initial data~$\Pi_{\perp} \Phi^{0}$ is in~$H^{s}_{L}$ with~$s \in [0,1]$, then that
regularity is propagated to~$\Pi_{\perp} \Phi$. Recalling the special form of the limit system for almost
all~$\beta$, the only equation we need to study is the one on~$\Pi_{K} \Phi$. Indeed denoting
$$
\Psi =( \Pi_{R} + \Pi_{P} + \Pi_{M}) \Phi,
$$
we have
$$
\d _t \Psi   +2Q_L( \Pi_0 \Phi,\Psi) -\nu  
\Delta'_L
\Psi=0,
$$
and   Proposition~\ref{trilinear} gives directly
$$
\|\Psi(t)\|_{H^{s}_{L}}^{2} + \frac\nu{C } \int_{0}^{t} \|\Psi(t')\|_{H^{s+1}_{L}}^{2} \leq
\|\Psi(0)\|_{H^{s}_{L}}^{2} \exp \left(\frac C \nu \|\Pi_0 \Phi\|^{2}_{L^{2}([0,T];L^{2})} \right).
$$
Now let us turn to the Kelvin modes.  In that case we simply use again  the fact that~$H^{s} $ and~$H^{s}_{L}$ spaces
coincide in the case of~$\Pi_{K}$. We have, using
  Proposition~\ref{trilinear} again,
$$
\longformule
{\|\Pi_{K} \Phi (t)\|_{H^{s} ( \R\times \T)}^{2} + \frac\nu {C_{1}}  \int_{0}^{t}  \|\Pi_{K}
 \Phi(t')\|_{H^{s+1}( \R\times \T) }^{2} \leq\left|
\int_{0}^{t} \left( \Pi_K
\Phi | Q (\Pi_K
\Phi,\Pi_K\Phi)\right)_{H^{s}( \R\times \T)}
(t')\: dt'\right| }{
+ \: C \int_{0}^{t}\|\Pi_0 \Phi (t')\|_{L^{2}} \|\Pi_{K} \Phi (t')\|_{H^{s} ( \R\times \T)} 
\|\Pi_{K} \Phi (t')\|_{H^{s+1}
 ( \R\times \T)} \: dt'.}
$$
Two-dimensional product rules give
\begin{eqnarray*}
\left|\left( \Pi_K
\Phi | Q (\Pi_K
\Phi,\Pi_K\Phi)\right)_{H^{s}( \R\times \T)} \right|
&\leq &
\|\Pi_{K} \Phi\|_{H^{s+1} ( \R\times \T)}
   \|\Pi_{K} \Phi \otimes \Pi_{K} \Phi\|_{H^{s } ( \R\times \T)}  \\
  &\leq &\|\Pi_{K} \Phi\|_{H^{s+1} ( \R\times \T)}    \|\Pi_{K} \Phi\|_{H^{s +\frac12} ( \R\times \T)}
    \|\Pi_{K} \Phi\|_{H^{ \frac12} ( \R\times \T)}\\
      &\leq &\|\Pi_{K} \Phi\|_{H^{s+1} ( \R\times \T)}^{\frac32}\|\Pi_{K} \Phi\|_{H^{s }
       ( \R\times \T)}^{\frac12}  \|\Pi_{K} \Phi\|_{H^{ \frac12} ( \R\times \T)}\\
   &\leq & \frac\nu{2C_{1}}     \|\Pi_{K} \Phi\|_{H^{s+1} ( \R\times \T)}^{2}
   + \frac C\nu  \|\Pi_{K} \Phi\|_{H^{s }
       ( \R\times \T)}^{2}\|\Pi_{K} \Phi\|_{H^{ \frac12} ( \R\times \T)}^{4}
\end{eqnarray*}
so that finally
$$
\longformule{
\|\Pi_{K} \Phi (t)\|_{H^{s} ( \R\times \T)}^{2} + \frac\nu {4C_{1}}  \int_{0^{t}} \|\Pi_{K}
 \Phi(t')\|_{H^{s+1}( \R\times \T) }^{2} \leq \frac C{\nu} \int_{0}^{t}  \|\Pi_0 \Phi(t')\|_{L^{2}}^{2} \|\Pi_{K} \Phi(t')\|_{H^{s} ( \R\times \T)}^{2}
  \: dt' }
 {+ \frac C{\nu} \int_{0}^{t} 
   \|\Pi_{K} \Phi(t')\|_{H^{s }
       ( \R\times \T)}^{2}\|\Pi_{K} \Phi(t')\|_{H^{ \frac12} ( \R\times \T)}^{4} \: dt' ,}
$$
and the result follows from Gronwall's lemma and the energy estimate on~$\Pi_{K} \Phi$. 

That concludes the proof of Theorem~\ref{limit-systemgeneric}. \hfill $\Box $

\section{A regularity result for the divergence} \label{sectionpropag}
In this section we are going to prove an additional regularity result for the system $(SW_0)$, which will be useful to study the strong asymptotics of the rotating shallow-water system in the next chapter.

 \begin{Prop}\label{divergence-prop}
 Let $\Phi^0$ belong to $L^2(\R\times \T)$.

 For all $\beta \in \R^+_*$, if $\Pi_\perp \Phi^0$ belongs to $H^{1/2}_L$ and $(\Pi_{P} + \Pi_{K}) \Phi^{0}$ belongs to~$ H^{\alpha}_{L}$ for~$\alpha > 3/2$,  
then the solution $\Phi$   of $(SW_0)$ with initial data $\Phi^0$ defined on $[0,T¬*[$ satisfies for all~$t
\in[0,T^*[$
$$
\int_0^t \|\nabla \cdot \Phi'(t')\|_{L^\infty(\R\times \T)}dt' <+\infty,
$$
where we recall that $\Phi'$ denotes the two last components of $\Phi$.

  Furthermore the regularity assumption on the initial data can be relaxed for all  but a countable number of $\beta$. Indeed, for all $\beta \in \R^+\setminus {\mathcal N}$ where  $\mathcal N$ is the countable subset of $\R^+$ defined in Theorem \ref{limit-systemgeneric},  if $(\Pi_{P} + \Pi_{K}) \Phi^{0}$ belongs to~$ H^{\alpha}_{L}$ for~$\alpha > 1/2$,  
then  the solution $\Phi$   of $(SW_0)$ with initial data $\Phi^0$  satisfies for all~$t
\in \R^+$
$$
\int_0^t \|\nabla \cdot \Phi'(t')\|_{L^\infty(\R\times \T)}dt' <+\infty.$$
\end{Prop}

\begin{proof}
Such a result is established by decoupling the equations on the
various parts of $\Phi$, and proving that  the regularity is
propagated for the Poincar{\'e} and Kelvin modes, while a smoothing
property on the divergence holds  for the nonoscillating,
Rossby and mixed Rossby-Poincar{\'e} modes.

Let us decompose~$\Phi$ on the supplementary subsets  $\Ker L$, $ R$, $M$, $P$ and $K$
$$\Phi=\Pi_0 \Phi+\Pi_R \Phi +\Pi_M \Phi +\Pi_P \Phi+\Pi_K \Phi,$$
  and  estimate each projection
separately.

\medskip
$\bullet$  For nonoscillating, Rossby and mixed modes
  the smoothness
 of the divergence is not due to a propagation result but to a stationnary property of the eigenvectors.

 Clearly, by definition of~$ \Ker L$, $\nabla \cdot  (\Pi_{0} \Phi)' = 0$.

\medskip
  Let us consider the Rossby and mixed Rossby-Poincar{\'e}  modes. By 
definition of the  Rossby modes we deduce the following
relation
$$\forall \Phi \in  R,\quad \nabla \cdot \Phi' =\sum_{i\lambda \in
{\mathfrak S}_R } \nabla \cdot \Phi'_\lambda=\sum_{i\lambda \in  {\mathfrak
S}_R } i\lambda (\Phi_\lambda)_0$$
with the notation $\Phi_\lambda=\Pi_\lambda \Phi$.
It follows that
\begin{eqnarray*}
 \|\nabla \cdot \Phi'\|_{H^{2}}^{2}
&= &\Bigl\|\sum_{i\lambda \in
{\mathfrak S}_R} i\lambda (\Phi_\lambda)_0 \Bigr \|_{H^{2}}^{2} \\
&\leq & C \sum_{i\lambda \in
{\mathfrak S}_R} \left \| \lambda
(\Phi_\lambda)_0 \right\| _{H^2(\R\times \T)} ^2,
\end{eqnarray*}
by Proposition \ref{equivnorm}. But, as Rossby waves  correspond to~$ j =
0$, we have
$$
 \| \lambda (\Phi_\lambda)_0 \| _{H^2(\R\times \T)} ^2 \leq C
 |\lambda|^2 \sum_{\tau (n,k,0) = \lambda} \| (\Pi_{n,k,0}\Phi)_0\|  _{H^2(\R\times \T)}
^2,
$$
using  Remark~\ref{componentseparately}. Recalling the explicit
form of~$  (\Psi_{n,k,j})_0$, we see that
$$
 \| (\Pi_{n,k,0}\Phi)_0\| _{H^2(\R\times \T)} ^2 \leq (1 + n + k^2)  \| 
(\Pi_{n,k,0}\Phi)_0\| _{H^1(\R\times \T)}
^2.
$$
But  for Rossby modes, the following asymptotics hold
as~$ |k| $ or~$n  $ goes to infinity:
$$
\lambda = \tau (n,k,0)  \sim\frac{\beta k}{
k^{2} + \beta (2n+1)}\cdotp
$$
So we infer that as~$ |k| $ or~$n  $ goes to infinity,
$$
 |\lambda|^2 \| (\Pi_{n,k,0}\Phi)_0\| _{H^2(\R\times \T)} ^2 \leq C  \| 
(\Pi_{n,k,0}\Phi)_0\| _{H^1(\R\times \T)}
^2.
$$
Finally we infer that
\begin{eqnarray*}
 \|\nabla \cdot \Phi'\|_{H^{2}}^{2}  &\leq & C \sum_{i\lambda \in
{\mathfrak S}_R} \left \| \lambda
(\Phi_\lambda)_0 \right\| _{H^2(\R\times \T)} ^2 \\
 &\leq & C \sum_{(n,k,0) \in {\mathfrak S}_R }  \|  (\Pi_{n,k,0}\Phi)_0\|_{H^1(\R\times
\T)} ^2 \\
 &\leq & C \|\Phi\|_{H^1_L} ^2.
\end{eqnarray*}
By the embedding of~$H^{2}(\R\times \T)$ into~$L^{\infty}(\R\times \T)$  we conclude
that~$\nabla\cdot  (\Pi_{R} \Phi)'$ belongs to the space~$L^2([0,T];  L^\infty(\T \times
\R))$.
The same result can easily be extended to the mixed Poincar{\'e}-Rossby
modes (it is in fact easier since~$ n = 0$ in that case) and we obtain
$$
\| \Pi_M \Phi \|_{L^2([0,T],   H^1_L )} \leq C_T,
\quad \|\nabla \cdot ( \Pi_M \Phi)'
\|_{L^2([0,T], L^{\infty}(\R\times \T))}\leq C_T\,.
$$
Finally  we deduce that
$$
\|\nabla \cdot (( \Pi_{0} +\Pi_R+ \Pi_M) \Phi)'
\|_{L^2([0,T], L^\infty(\R\times \T))}\leq C_T\,.
$$

\bigskip
$\bullet$ In order to establish a similar estimate for the Poincar\'e and Kelvin modes, we prove that the equation governing these modes propagates the $H^s_L$  regularity without restriction on $s$. Indeed we have seen in Remark \ref{restriction} that the restriction to $s\leq 1$ for the propagation of regularity stated in
Theorem~\ref{limit-systemallbeta}  is due to the coupling between Rossby modes.

The idea here is to study the propagation of the following norm on $(\Ker L)^\perp$
$$\| \Phi\|^2_s= \sum_{\lambda\in {\mathfrak S}}(1+\lambda^2)^s\| \Pi_\lambda \Phi\|^2_{L^2(\R\times \T)},$$
which controls the $H^s_L$  regularity of the Poincar\'e and Kelvin modes only, due to the following easy
estimate (see Proposition~\ref{equivnorm}):
\begin{equation}
\label{estimatesHsL}C^{-1} \| (\Pi_K+\Pi_P)\Phi\|_{H_L^s} \leq \| \Pi_\perp \Phi\|_s\leq C (\| (\Pi_K+\Pi_P)\Phi\|_{H_L^s}+ \| \Pi_\perp \Phi \|_{L^2(\R\times \T)}) .
\end{equation}This norm  is convenient to deal with the condition of resonance occuring in the nonlinear term of $(SW_0)$.

\medskip
Similar arguments as in   Proposition \ref{trilinear}  allow to write the following trilinear estimate:
\begin{equation}
\label{trilinear-modified}
|(\Phi |\tilde Q_L(\Phi,\Phi))_s | \leq C_s \| \Pi_\perp \Phi\|_s (\|\Pi_\perp \Phi\|_{s+1}+\|\Pi_\perp \Phi\|_{H^1_L}) (\| \Pi_\perp \Phi\|_{H^{3/2}_L}+ \| \Pi_0 \Phi\|_{L^2(\R\times \T)} ).
\end{equation}
With the same notation as in the proof of Proposition  \ref{trilinear}, we have indeed
$$\begin{aligned}
&(\Phi_*|\tilde Q_L(\Phi,\Phi^*))_s \eqdefa \left|  \sum_{{\tau_*=\tau +\tau^*\atop
k_*=k+k^* }\atop
i\tau,i\tau_*, i\tau^* \in {\mathfrak S}\setminus \{0\}}
   (1 + \tau_*^{2})^{s}
  \left(  \Pi_{n_*,k_*,j_*}  \Phi_{*} \Big| \Pi_{n_*,k_*,j_*} Q (\Pi_{   n, k, j} \Phi,
\Pi_{ n^*, k^*, j^*} \Phi^*)\right)_{L^2} 
\right|\\
&\leq C_s \sum (1 + \tau_*^{2})^{s/2} \Big( (1+\tau^2)^{s/2} + (1+(\tau^*)^2)^{s/2}\Big)
|\varphi_*| |\varphi |
|\varphi^*|  ( (n+k^2)^{1/2} +(n^*+(k^*)^2)^{1/2}),
\end{aligned}
$$
from which we deduce  that
$$
|(\Phi |\tilde Q_L(\Phi,\Phi))_s \leq C_s \| \Pi_\perp \Phi\|_s (\|\Pi_\perp \Phi\|_{s+1}+\|\Pi_\perp \Phi\|_{H^1_L}) \| \Pi_\perp \Phi\|_{H^{3/2}_L} . $$
In the same way, for the geostrophic part, we have
$$
(\Phi|\tilde Q_L(\Pi_0\Phi,\Pi_\perp \Phi))_s\leq 
C_s \| \Pi_\perp \Phi\|_s (\|\Pi_\perp \Phi\|_{s+1}+\|\Pi_\perp \Phi\|_{H^1_L}) \| \Pi_0 \Phi\|_{L^2(\R\times \T)} .$$
 Note that the~$H^1$ norm appearing above is  used to control the gradient of the Rossby and mixed
 modes.

\medskip
Using (\ref{trilinear-modified}) and Gronwall's lemma, we get  the following propagation result
$$ 
\begin{aligned}
\| \Pi_\perp \Phi(t)\|_s^2 &+{ \nu \over C_s}\int_0^t\left(  \| \Pi_\perp\Phi(t')\|_{{s+1}}^2+\| \Pi_\perp \Phi(t')\|_{H^1_L}^2 dt'\right)\\
 &\leq  \| \Pi_\perp \Phi^0\|_{s}^2\exp\left( \frac {C_s} \nu \int_0^t
  \left(\| \Pi_\perp \Phi(t')\|_{H^{3/2}_L}+\|\Pi_0 \Phi(t')\|_{L^2(\R\times \T)} \right)^2 dt'\right)\,.
\end{aligned}
$$
(see Paragraph \ref{regularitystrong} for the detailed proof), and therefore by~(\ref{estimatesHsL}) we get
$$
\begin{aligned}
&\| (\Pi_K+\Pi_P) \Phi(t)\|_{H^s_L}^2+{ \nu \over C_s}\int_0^t  \| (\Pi_K+\Pi_P)\Phi(t')\|_{H^{s+1}_L}^2 dt'\\
 &\leq  (\| (\Pi_K+\Pi_P) \Phi^0\|_{H^s_L}^2+ \| \Phi^0\|_{L^2(\R\times \T)} ^2)\exp\left( \frac {C_s} \nu \int_0^t
  \left(\| \Pi_\perp \Phi (t')\|_{H^{3/2}_L}+\|\Pi_0 \Phi(t')\|_{L^2(\R\times \T)} \right)^2 dt'\right)\,.
\end{aligned}
$$

\medskip
By Sobolev embeddings we then deduce the expected control on the divergence
$$
\|\nabla \cdot (( \Pi_{K} +\Pi_P) \Phi)'
\|_{L^2([0,T], L^\infty(\R\times \T))}\leq C_T
$$
for all $T \leq T^*$ where $T^*$ is the lifespan of the strong solution $\Phi$.

That concludes the proof of the proposition in the case of general~Ê$\beta$.

\bigskip
$\bullet$ By Proposition~\ref{resonances-prop} page~\pageref{resonances-prop},  we recall that,
for all  $\beta \in \R^+ \setminus {\mathcal N}$, the  only possible resonances
are Kelvin resonances.  

Let us now consider the equation governing the Poincar{\'e}  modes which can be
seen as a
linear parabolic equation whose coefficients depend on $\Pi_0 \Phi$.  Therefore it is
very easy to
propagate regularity once one has noticed that the multiplication by  $|\lambda|$ for the
Poincar{\'e}
mode $\Pi_\lambda \Phi$ is ``equivalent" to a derivation.

Introduce as previously  the notation
$$\Phi=\sum_{(n,k,j)\in S} \varphi_{n,k,j} \Psi_{n,k,j},$$
so that
$$\Pi_P \Phi= \sum_{(n,k,j) \in S_P} \varphi_{n,k,j} \Psi_{n,k,j},$$
where
$$S_P=\N^*\times \Z\times \{-1,1\} \cup \{0\} \times \Z^+_* \times  \{1\}\cup \{0\}
\times \Z^-_* \times \{-1\} \cup {0} \times {0} \times \{-1,1\}.$$

We can use Proposition~\ref{diag-prop}  page \pageref{diag-prop} to deduce that
for each $(n,k,j)$ in~$ S_P$ the equation  governing~$\varphi_{n,k,j}$  can be decoupled
(recall  that~$\Pi_{0} \Phi  $ only
depends on~$x_{1}$):
$$
\d_t \varphi_{n,k,j} -\nu\varphi_{n,k,j}  ( \Psi_{n,k,j}|\Delta'  \Psi_{n,k,j}
)_{L^2(\R\times \T)}
=- 2\varphi_{n,k,j}
 (\Psi_{n,k,j} | Q(\Psi_{n,k,j},\Pi_0 \Phi))_{L^2(\R\times \T)}
$$
which can be rewritten
$$
\longformule{\d_t \left( \varphi_{n,k,j}\exp\left( -\nu t(  \Psi_{n,k,j}|\Delta'
\Psi_{n,k,j} )_{L^2(\R\times \T)}
 \right)\right)}
{=-2\varphi_{n,k,j}(\Psi_{n,k,j} | Q(\Psi_{n,k,j},\Pi_0  \Phi))_{L^2(\R\times \T)}
\exp\left( -\nu t( \Psi_{n,k,j}|\Delta' \Psi_{n,k,j} )_{L^2(\R\times  \T)} \right).
}$$

By Gronwall's lemma and the estimates
$$\begin{aligned}
\left| (\Psi_{n,k,j} | Q(\Psi_{n,k,j},\Pi_0 \Phi))_{L^2(\R\times \T)}  \right| \leq
C_1(n+k^2)^{1/2},\\
-( \Psi_{n,k,j}|\Delta' \Psi_{n,k,j} )_{L^2(\R\times \T)}  \geq  C_2(n+k^2),
\end{aligned}
$$
we then deduce that there exists a nonnegative constant $C_\nu$  (depending only on
$\nu$) such that,
\begin{equation}
\label{heat-estimate}
\forall (n,k,j)\in S_P,\quad |\varphi_{n,k,j}(t)|\leq  |\varphi_{n,k,j}(0)| \exp(-C_\nu
(n+k^2)t).
\end{equation}
We have
$$
\begin{aligned}
\left\| \nabla \cdot (\Pi_P \Phi)' (t) \right\|_{L^\infty(\R\times  \T)}& \leq
\sum_{(n,k,j)\in S_P}|\varphi_{n,k,j}(t)|\left\| \nabla  \cdot (\Psi_{n,k,j})' (t)
\right\|_{L^\infty(\R\times \T)}\\
&\leq C \sum_{(n,k,j)\in S_P}|\varphi_{n,k,j}(t)| (n+k^2)^{1/2}
\end{aligned}
$$
since $(\Psi_{n,k,j})$ is uniformly bounded in $L^\infty(\R\times \T)$.  Thus, by
(\ref{heat-estimate}),
$$\left\| \nabla \cdot (\Pi_P \Phi)' (t) \right\|_{L^\infty(\R\times  \T)}\leq C
\sum_{(n,k,j)\in S_P}|\varphi_{n,k,j}(0)| \exp(-C_\nu  (n+k^2)t)(n+k^2)^{1/2}.$$
Integrating with respect to time leads then to
$$
\longformule
{\left\| \nabla \cdot (\Pi_P \Phi)' \right\|_{L^1([0,T];L^\infty(\R  \times \T))}
\leq C'_\nu \sum_{(n,k,j)\in S_P}|\varphi_{n,k,j}(0)|(n+k^2)^{-1/2} }
{\leq C'_\nu \left(\sum_{(n,k,j)\in  S_P}|\varphi_{n,k,j}(0)|^2(n+k^2)^{\alpha} \right)
^{1/2} \left(
\sum_{(n,k,j)\in S_P}(n+k^2)^{-1-\alpha}\right)^{1/2},}
$$
from which we deduce that for $\alpha>1/2$,
$$\left\| \nabla \cdot (\Pi_P \Phi)'  \right\|_{L^1([0,T],L^\infty(\R\times \T))}\leq C
\|\Pi_P \Phi^0\|_{   H^\alpha_L }$$
where $C$ depends only on $\nu$ and $\alpha$.

\medskip
 It remains then to establish the propagation of regularity for
the Kelvin part of the equation, which is nonlinear and has no  smoothing effect for the
divergence
as the Rossby part.

The crucial point here is to recall as above that this equation is actually  one-dimensional
(modulo a smooth function with respect to
$x_1$). The propagation of regularity result proved in Paragraph~\ref {sectionpropaghs} 
implies that as soon as the initial data is in~$H^{\alpha}_{L}$ with~$0 \leq \alpha \leq 1$, then the
solution lies in~$L^{2}(\R^{+};H^{\alpha+1}_{L})$. In particular~$ \nabla \cdot (\Pi_K \Phi)' $
lies in~$L^{2}(\R^{+};H^{\alpha }_{L})$. So if~$\alpha>1/2$, using the fact that~$H^{\alpha}(\T)$ is 
embedded in~$L^{\infty}(\T)$, the result follows directly. 

Proposition~\ref{divergence-prop} is proved.
\end{proof}

\chapter{Convergence results}\label{convergence}
The aim of this chapter is to study the asymptotics of the rotating
shallow-water system~(\ref{SW-eps}) presented
page~\pageref{SW-eps}. In particular we will see that  the
system~(\ref{lim-filtered}) obtained formally
in Section~\ref{sctschochet}, page~\pageref{lim-filtered}, is indeed  
the limit system, after application of the filtering operator~$ \exp
(-tL/\e) $. In order to simplify the
presentation let us recall here the two main systems we will be
considering in this chapter, namely the shallow-water system\label{SWeps}
$$
\begin{aligned}
& (SW_\eps) \quad \quad 
\left\{
\begin{array}{r}
\displaystyle \d_t \eta +\frac1\eps \DIV \Bigl((1+\eps \eta) u\Bigr) =0, \nonumber\\
\displaystyle \d_t \Bigl((1+\eps \eta)u\Bigr) + \nabla \cdot \Bigl((1+\eps \eta)u\otimes u\Bigr) +  
\frac{ \beta x_1}{ \eps} (1+\eps \eta)u^\perp+
{1\over \eps}(1+\eps \eta)
\nabla
\eta
-\nu \Delta u =0, \nonumber \\
\displaystyle  \eta_{|t=0}=\eta^0,\quad u_{|t=0} =u^0,
\end{array}
\right. \\
\\
&\!\!\!  \mbox{and the limit system} \\
\\\label{SW0}
& (SW_0) \quad \quad 
\left\{
\begin{array}{r}
 \d_t \Phi +Q_L (\Phi,\Phi) -\nu \Delta'_L \Phi
=0\\
\displaystyle  \Phi_{|t=0}=(\eta^0,u^0),
\end{array}
\right.
\end{aligned}
$$
where $\Delta'_L$ and $Q_L$ denote  the linear and symmetric bilinear operator
defined by~(\ref{Q-lap-L}) page~\pageref{Q-lap-L}. 

We also recall the formal equivalent form of~$(SW_\eps)$,
$$
\begin{aligned}
&\ 
\left\{
\begin{array}{r}
\displaystyle \d_t (\eta, u )+\frac1\eps L(\eta,u)  
+ Q\left((\eta, u ),(\eta, u )\right) - \nu \Delta' (\eta,u)  = R
\nonumber \\
\displaystyle ( \eta, u)_{|t=0}=(\eta^0,u^0),
\end{array}
\right. 
\end{aligned}
$$
where~$Q$ and~$\Delta'$ are defined by~(\ref{Q-lap}) page~\pageref{Q-lap} 
 and
$$R  = (0, -\nu {\eps \eta  \over 1+\eps \eta }
\Delta u ).$$

The study of the asymptotics of~$(SW_\eps)$ will be achieved through three different methods, which provide
three different types of results. In Section~\ref{weakweakcv} we
describe the weak limit of the weak solutions to~$(SW_\eps)$ as~$
\e$ goes to zero, which is proved to satisfy the geostrophic equation
studied in the previous chapter, i.e., the projection
of~$(SW_0)$ onto~$ \Ker L$. The statement is given in
Theorem~\ref{weakcv} below. Then  for smooth enough initial data, 
 we prove in  Section~\ref{strongstrongcv} the  strong convergence of the
filtered sequence of solutions towards the unique solution
of~$(SW_0)$. The precise statement depends on the setting, as in Chapter~\ref{envelope}: 
for all~$\beta>0$ we are only able to prove results locally in time (globally for small data) whereas if a countable
set of values of~$\beta$ is removed, then the convergence is strong for all times, and the smoothness assumptions on 
the initial data are less restrictive (see   Theorems~\ref{strong-convergenceallbeta} and~\ref{strong-convergencegenericbeta}).
 Finally in Section~\ref{hybrid}
we propose an intermediate study between those two asymptotic results, by considering  the asymptotic 
behaviour of the 
filtered sequence~$e^{-tL/\e}(\eta_\e,u_\e) $, where~$ (\eta_\e,u_\e)
$ is a weak solution to~$(SW_\eps)$. We prove a strong convergence result towards a weak solution 
to~(\ref{lim-filtered}), where unfortunately due to the lack of compactness of~$\eta_\e$ in space, 
a defect measure remains 
 (see
Theorem~\ref{thmstrongweakcv}).  In order to circumvent that difficulty we propose an alternate system to the
Saint-Venant equations~$(SW_\eps)$, where capillarity effects are included. Technically the effect of capillarity
is to have a uniform control on~$\e \eta_{\e}$ in strong enough norms so as to obtain an evolution equation for~$u_\e$. A strong
convergence result for~$e^{-tL/\e}(\eta_\e,u_\e) $ is established in that new setting, see
 Theorem~\ref{thmstrongweakcvcapillarity}.

In this chapter, many results and  notation of the
previous chapters will be used. However precise references will be
made each time, so that this chapter can be read independently of the
others (assuming the results of course). 

\section{Weak convergence of weak solutions}\label{weakweakcv}
\setcounter{equation}{0}
The first aim of the chapter is  to describe the weak limit $(\eta,u)$ of $(\eta_\eps,u_\eps)$ 
as $\eps$
goes to zero.
\begin{Thm}[Weak convergence]\label{weakcv}
{\sl
Let $(\eta^0,u^0) \in L^2(\R\times \T)$ and $(\eta^0_\eps,u_\eps^0)$
be such that
\begin{equation}\label{firstcondition}
\begin{aligned}
\frac12 \int \left(
|\eta^0_\eps|^2 +(1 + \e
 \eta^0_\eps)
|u^0_\eps|^2\right)dx\leq \EE^0,\\
(\eta^0_\eps, u^0_\eps) \to (\eta^0,u^0) \hbox{ in } L^2(\R\times \T).
\end{aligned}
\end{equation}
 For all
$\eps >0$, denote by $(\eta_\eps,u_\eps)$ a solution of~$(SW_\e)$ with initial
data~$(\eta_\eps^0,u_\eps^0)$, as constructed in
Corollary~\ref{corexistenceeta} page~\pageref{corexistenceeta}. 
Then  $(\eta_\eps, 
u_\eps)$
converges weakly in~$L^2_{loc}(\R^+\times \R\times \T)$ to the solution
$(\eta,u)\in L^\infty(\R^+,L^2(\R))$, with~$u$ also belonging to~$ L^2(\R^+, \dot 
H^1(\R))$, of the
following linear equation (given in weak formulation)
\begin{equation}\label{constraintsstatement}
u_1=0, \quad
 \beta x_1 u_2 +\d_1 \eta =0,
\end{equation}
and for all~$(\eta^*,u^*)  \in L^2  \times H^1(\R ) $ satisfying 
(\ref{constraintsstatement})
\begin{equation}
\label{heat}
\int (\eta \eta^* + u_{2} u^*_{2}) (t,x) \: dx + \nu \int_{0}^{t} \int \nabla u_{2} \cdot 
\nabla
u^*_{2} (t',x) \: dx \: dt' = \int (\eta^{0} \eta^* + u_{2}^0 u^*_{2}) (x) \: dx .
\end{equation}
}
\end{Thm}
\begin{Rem} { 
$\bullet $ Theorem~\ref{weakcv} shows that the system satisfied by
  the weak limits of~$ \eta_\eps$ and~$ u_\eps$ is linear. There is
  therefore no  convective term in the mean flow:
  system~(\ref{constraintsstatement}, \ref{heat}) actually corresponds
  to the projection of~$ (SW_0)$ onto~$ \Ker L$: as seen in
  Section~\ref{algebraic section}, that projection can indeed be
  formally written
$$
\begin{aligned}
\partial_{t} (\eta,0,u_2) - \nu \Pi_0(0,0,\Delta u_2) = 0,\\
 (\eta,u)(t) = \Pi_0
(\eta,u)(t) \: \forall t \geq 0,\\
(\eta,u)_{|t = 0}=\Pi_0 (\eta^0,u^0).
\end{aligned}
$$

\medskip

$\bullet $ Note that~$(\eta^0,u^0)$ do not necessarily satisfy the 
constraints~(\ref{constraintsstatement}),
so in general~$(\eta,u)_{|t = 0} $ is not equal to~$(\eta^0,u^0)$.

\medskip

$\bullet $  The study of the waves induced by~$L$, in Chapter~\ref{envelope},   revealed the presence of trapped
equatorial waves, which however  do not appear in the mean flow
described by Equation~(\ref{heat}):  no constructive interferences take place in the 
limiting
  process, in other words the fast oscillating modes decouple from the
  mean flow, without creating any additional term in the limit
  system (that feature was already observed in~\cite{gallagher/sr} in
  the case of inhomogeneous rotating fluid equations, modelling the
  ocean or the atmosphere at midlatitudes). This will be obtained by a
  compensated compactness argument in
  Section~\ref{compensatedcompactnesssection}.  
}
\end{Rem} 
\subsection{Constraints on the weak limit}
We recall (see Chapter~\ref{intro}) that the uniform energy bound
on~$(\eta_{\e},u_{\e})$ implies  the existence of a weak limit~$(\eta,u)$. In this
paragraph we are going to prove that the weak limit belongs to~$\Ker L$.

\begin{Prop} 
\label{cor:limit} { 
Let $(\eta^0,u^0)\in L^2(\R\times \T)$. Denote by
  $(\eta_\eps,u_\eps)_{\eps> 0}$ a
family of  solutions of~$(SW_\e)$, and by $(\eta,u)$ any of its limit points. Then, 
$(\eta,u)\in L^\infty(\R^+,L^2(\R))$  belongs to~$\Ker L$, and in particular satisfies the constraints
\begin{equation}
\label{constraintlimit}
u_1=0,\quad \beta x_1 u_2+\d_1 \eta =0.
\end{equation}
}
\end{Prop}
\begin{proof} Let $\chi,\psi\in \DD(\R^+\times \R\times \T)$ be any  test functions. 
Multiplying
  the conservation of mass in~$(SW_\e)$ by $\eps \chi$
 and integrating with respect to all variables leads to 
$$
\iint \left( \eps \eta_\eps  \d _t \chi +(1+\eps \eta_\eps) u_\eps
\cdot
\nabla
\chi\right)dxdt =0 .
$$
Because of the bounds coming from the energy estimate~(\ref{uniform}), we can take limits 
in the previous
identity as $\eps$ goes to $0$ to get
$$\iint u\cdot \nabla \chi dxdt=0.$$
Similarly, multiplying the conservation of momentum by $\eps \psi$ and integrating with 
respect to all variables leads to
$$
\iint
\Bigl(
\eps (1+\eps \eta_\eps)u_\e \partial_t \psi +
 \eps(1+\eps \eta_\eps)u_\e\cdot ( u_\e \cdot \nabla) \psi  + \beta x_1(1+\eps \eta_\eps) 
u_\eps \cdot
\psi^\perp + (1+\frac\eps{2} \eta_\eps) \eta_\e \DIV \psi+\nu u_\eps \cdot \Delta 
\psi\Bigr) dxdt=0.
$$
Once again the bounds coming from the energy estimate~(\ref{uniform})
  enable us to take the limit as $\eps$ goes to $0$, and find that
$$ \iint (\eta \nabla \cdot \psi+\beta x_1 u\cdot
\psi^\perp)dxdt=0.$$

It follows that  $(\eta(t),u(t))$ belongs to $\Ker L$ for almost all~$t\in \R^+$, and we 
conclude
by Proposition~\ref{kernel} page~\pageref{kernel} that~$(\eta,u)$ does
not depend on~$x_2$  and satisfies the
constraints~(\ref{constraintlimit}).
\end{proof}

To go further in the description of the weak limit $(\eta,u)$, we have to
isolate the fast oscillations generated by the singular perturbation $L$,
which produce ``big" terms in $(SW_\e)$, but converge weakly to 0.

Therefore, a natural idea consists in introducing the following decomposition
$$(\eta_\eps,u_\eps) =\Pi_0( \eta_\eps,  u_\eps)
+\Pi_\perp(\eta_\eps,u_\eps) ,$$
where $\Pi_0$ is the
$L^2$ orthogonal projection onto $\Ker L$ and~$\Pi_\perp $ the
$L^2$ orthogonal projection onto $(\Ker L)^\perp$.

The idea to get the mean motion is then to apply $\Pi_0$ to $(SW_\e)$~: since $L$ is
a skew-symmetric operator,  we have~$
\Pi_0 L=0$ and we expect $\d_t \Pi_0(\eta_\eps,  u_\eps) $ to be uniformly bounded in some
distribution space. The difficulty comes from the fact that one has no uniform spatial regularity on~$\eta_{\e}$. 
That is why we
will actually consider the weak form of the evolution equations; the point is then to 
take limits in
the nonlinear terms.

\subsection{Rough description of the oscillations}
The analysis of the nonlinear terms lies essentially on the structure
of the oscillations. A rough  description of those
fast oscillations will be enough to prove that they do not produce any
constructive interference, and  therefore do not
appear in the equation governing the mean (geostrophic) motion. The
much more  precise description given in 
Chapter~\ref{equatorialwaves} will not be used in this section, but
will be necessary 
 to discuss the strong asymptotic behaviour of the solutions in the
 next sections.

\noindent In the following statement we have considered a regularization kernel defined 
as follows: let~$\kappa$ be\label{kappa}
a function
of $ C^\infty_c(\R^2 ,\R^+)$ such that $\kappa(x)=0$ if $|x|\geq 1$ and $\int \kappa dx 
=1$. Then for any~$\delta> 0$ we
define~$\kappa_\delta$   by
$$
\kappa_\delta (x) =\delta^{-2} \kappa (\delta^{-1} x).
$$

\begin{Prop}\label{descriptionoscillations}
 { 
Let $(\eta^0,u^0) \in L^2(\R\times \T)$ and $(\eta^0_\eps,u_\eps^0)$
satisfy assumptions~(\ref{firstcondition}), and denote by
$((\eta_\eps,u_\eps))_{\eps>0}$  a family of solutions
of~$(SW_\e)$ with respective initial data~$(\eta_\eps^0,u_\eps^0)$.

Then $\eta_\e^{\delta} =\kappa_\delta \star \eta_\e$ and $m_\e^{\delta} =\kappa_\delta 
\star ((1+\eps
\eta_\eps)u_\e)=u_\eps^\delta +\eps (\eta_\eps u_\eps)^\delta$
 satisfy, for all $T>0$, the uniform convergences for all $\Omega \subset \subset 
\R\times \T$
\begin{equation}
\label{delta-cv}
\begin{aligned}
\| \eta_\eps-\eta_\eps^\delta\|_{L^\infty(\R^+, H^s(\Omega))} \to 0 \hbox{ as } \delta 
\to 0
\hbox{ uniformly in }\eps>0,  \hbox{ for all } s<0,\\
\| u_\eps-u_\eps^\delta\|_{L^2([0,T]; H^s(\Omega))} \to 0 \hbox{ as } \delta \to 0
\hbox{ uniformly in }\eps>0,  \hbox{ for all } s<1,\\
\| \eta_\eps u_\eps -(\eta_\eps u_\eps)^\delta\|_{L^2([0,T]; H^s(\Omega))} \to 0 \hbox{ 
as } \delta \to 0
\hbox{ uniformly in }\eps>0,  \hbox{ for all } s<0,\\
\end{aligned}
\end{equation}
as well as the approximate wave equations
\begin{equation}\label{oscillation}
\begin{aligned}
\eps\d_t \eta_\eps^\delta +\DIV m_\eps^\delta = 0,\\
\eps \d_t m_\e^{\delta} + \beta x_1 (m_{\eps}^\delta)^\perp + \nabla \eta_{\eps}^\delta = 
\eps s_\e^\delta + \delta
\sigma_{\e}^{\delta} ,
\end{aligned}
\end{equation}
denoting  by $s_\eps^\delta $ and  $\sigma_{\e}^{\delta} $  some
quantities satisfying, for all~$T > 0$,
\begin{equation}
\label{wave-remainder}
\begin{aligned}
 \sup_{\delta>0} \sup_{\eps>0} \|\sigma_{\e}^{\delta}\|_{L^{2}([0,T];H^{1}(\R\times \T))} 
< +\infty,
\\
\forall \delta > 0, \quad \sup_{\eps>0}
\|
 s_\eps^\delta \|_{L^1 ([0,T];H^{1}( \R\times \T))}<\infty   .
\end{aligned}
\end{equation}
In particular the approximate vorticity~$\omega_{\e}^\delta  =\nabla^\perp \cdot 
m_\eps^\delta  $
 satisfies
\begin{equation}\label{oscillation2}
\e \partial_{t} (\omega_{\e}^\delta  -\beta  x_{1} \eta_{\e}^\delta ) + \beta 
m_{\eps,1}^\delta = \e q_{\e}^\delta +
\delta p_{\e}^{\delta},
\end{equation}
 with,  for all~$T > 0$,
\begin{equation}
\label{remainder2}
\supetage{\e > 0}{\delta > 0} \|p_{\e}^{\delta}\|_{L^{2}([0,T];L^{2}(\R\times \T))} < 
+\infty \hbox{ and }\forall
\delta >0,\quad
\sup_{\e > 0}\| q_{\e}^\delta \|_{L^{1}([0,T];L^{2}(\R\times \T))} < +\infty .
\end{equation}
}\end{Prop}

\begin{proof} We proceed in two steps, first stating the wave equations for 
$(\eta_\eps,m_\eps)$, then
introducing the regularization $(\eta_\eps^\delta,m_\eps^\delta)$.

$\bullet$
The first step consists in establishing some bounds for
$$\eps \d_t (\eta_\eps,m_\eps) +L (\eta_\eps,m_\eps).$$
We have
$$
 \e \partial_t \eta_\e +  \DIV m_\e = 0,
$$
and
$$
\e\partial_t m_\e +  \beta x_1 m_\e^\perp + \nabla \eta_\e = \e r_\e,
$$
with
$$
 r_\e = - \nabla \cdot (m_\eps \otimes u_\eps)- \eta_\eps \nabla \eta_\eps +  \nu
\Delta u_\e   .
$$
Let us now find a bound for~$ r_\e$. It is made of three
contributions. The easiest to handle is $\Delta
u_\e $. Indeed~$ u_\e$ is bounded in~$ L^2(\R^+,\dot H^1(\R\times \T))$, so
 $
\Delta u_\e $ is bounded in $ L^2(\R^+,H^{-1}(\R\times \T))$.

 Next let us consider
the nonlinear terms
$$ - \nabla \cdot (m_\eps \otimes u_\eps) - \eta_\eps \nabla \eta_\e=
- \nabla \cdot (m_\eps \otimes u_\eps) 
-\frac12\nabla  \eta_\eps^2.   $$ By
the energy bound (\ref{uniform}) we infer that they are bounded 
in~$L^\infty(\R^+,W^{-1,1} (\R\times \T)) $.

Therefore in particular
\begin{equation}
\label{r-control}
\| r_\eps\|_{L^2([0,T]; H^{-5/2}(\R\times \T))} \leq C_T.
\end{equation}

\medskip
$\bullet$
Now let us proceed to the regularization. We recall that the energy inequality 
(\ref{uniform}) provides the following uniform
bounds
$$\begin{aligned}
\,& \| \sqrt{1+\eps \eta_\eps} u_\eps\|_{L^\infty(\R^+, L^2(\R\times \T))} \leq C,\\
& \| \eta_\eps\|_{L^\infty(\R^+, L^2(\R\times \T))} \leq C,\\
& \| u_\eps \|_{L^2(\R^+, \dot H^1(\R\times \T))} \leq C.
\end{aligned}
$$
In particular we have
$$
\|u_\eps\|_{L^2(\R\times \T)}^2 \leq \| \sqrt{1+\eps \eta_\eps} u_\eps\|_{ L^2(\R\times 
\T)}^2 +C\eps \|\eta_\eps\|_{
L^2(\R\times \T)} \| u_\eps\|_{L^2(\R\times \T)} \| u_\eps\|_{\dot H^1(\R\times \T)} ,$$
thus by the Cauchy-Schwarz inequality
$$
\frac12\|u_\eps\|_{L^2(\R\times \T)}^2 \leq \| \sqrt{1+\eps \eta_\eps} u_\eps\|_{ 
L^2(\R\times \T)}^2+ 8C^2 \eps^2\|\eta_\eps\|_{
L^2(\R\times \T)}^2 \|\nabla  u_\eps\|_{L^2(\R\times \T)}^2.$$
and finally
\begin{equation}
\label{H1bound}
 \| u_\eps \|_{L^2([0,T];  H^1(\R\times \T))} \leq C_{T}.
\end{equation}
We also deduce from the usual product laws that
\begin{equation}
\label{Hsboundetau}
\|\eta_\eps u_\eps \|_{L^2([0,T];H^s(\R\times \T))} \leq C_s \hbox{ for all } s<0. 
\end{equation}

The convergences (\ref{delta-cv}) are then obtained by the Rellich Kondrachov theorem.

\medskip
$\bullet$ By convolution we get (with obvious notation)
$$
 \e \partial_t \eta_\e^{\delta} +  \DIV u_\e^{\delta} = 0,
$$
and
$$
\e\partial_t m_\e^{\delta} + \beta  x_1 (m_\e^{\delta}) ^\perp + \nabla \eta_\e^{\delta}  
= \e r_\e^{\delta}  +
\beta x_1 (m_\e^{\delta}) ^\perp - (\beta x_1 m_\e^\delta ) ^{\perp}.
$$
 Notice that
$$
\begin{aligned}
x_1 m_\e^{\delta} (x) - ( x_1 m_\e)^{\delta} (x) &=  \int \kappa_{\delta}(y) m_\e(x-y) 
(x_{1} -
(x_{1}- y_{1} ))
\:  dy \\
&=  \int y_{1}\kappa_{\delta}(y)m_\e(x-y) \: dy.
\end{aligned}
$$
That implies that
\begin{equation}
\label{convolution-error}
 x_1 m_\e^{\delta} (x) - ( x_1 m_\e)^{\delta} (x) =   \delta \kappa^{(1)}_{\delta}(y) 
\star m_{\e},
\end{equation}
where~$ \displaystyle \kappa^{(1)}_{\delta}(y)  =  \delta^{-2}   \kappa^{(1)}
 ( \delta^{-1} x)$ and~$\kappa^{(1)} (x)
= x_{1} \kappa (x)$. We therefore infer that
$$
\beta x_1 (m_\e^{\delta})^{\perp} (x) - \beta (x_1 m_\e^{\perp})^{\delta} (x) = \delta
\sigma_{\e}^{\delta}(x) + \e \beta \delta \kappa^{(1)}_{\delta} * (\eta_{\e} u_{\e}),
$$
where for all~$T>0$,
$$
\sup_{\delta > 0} \sup_{\e > 0} \|\sigma_{\e}^{\delta}\|_{L^{2}([0,T];H^{1}(\R\times 
\T))} < +\infty.
$$

It remains then to control
$$s_\eps^\delta = r_\e^{\delta} + \beta  
\delta \kappa^{(1)}_{\delta} * (\eta_{\e} u_{\e})
\,.$$
By (\ref{r-control}) and (\ref{Hsboundetau})  
 we get
$$
\forall \delta > 0, \quad \sup_{\eps>0}
\|
 s_\eps^\delta \|_{L^1 ([0,T];H^{1}( \R\times \T))}<\infty   .$$

\medskip
$\bullet$
Taking the vorticity in the second equation of (\ref{oscillation}) leads then to
$$\d_t \omega_\eps^\delta +\nabla^\perp\cdot ( \beta x_1 (m_{\eps}^\delta)^\perp + \nabla 
\eta_{\eps}^\delta) = \eps
\nabla^\perp\cdot s_\e^\delta + \delta
\nabla^\perp\cdot \sigma_{\e}^{\delta},$$
from which we deduce that
$$\d_t \omega_\eps^\delta +\beta x_1 \nabla\cdot  m_\eps ^\delta +\beta m_{\eps,1}^\delta 
 =
\eps
\nabla^\perp\cdot s_\e^\delta + \delta
\nabla^\perp\cdot \sigma_{\e}^{\delta}.$$
Combining this last equation with the first one in (\ref{oscillation}) gives finally
$$\d_t (\omega_\eps^\delta -\beta x_1\eta_\eps ^\delta )+\beta m_{\eps,1}^\delta  =
\eps
\nabla^\perp\cdot s_\e^\delta + \delta
\nabla^\perp\cdot \sigma_{\e}^{\delta},$$
from which we deduce the estimate (\ref{remainder2}) on the remainder.
The proposition is proved.
\end{proof}

\subsection{Proof of Theorem \ref{weakcv}}
 We consider an initial data $(\eta^0,u^0) \in L^2(\R\times \T)$ and a family~$(\eta^0_\eps,u_\eps^0)$ such that
$$\begin{aligned}
\frac12 \int \left(
|\eta^0_\eps|^2 +(1 + \e
 \eta^0_\eps)
|u^0_\eps|^2\right)dx\leq \EE^0,\\
(\eta^0_\eps, u^0_\eps) \to (\eta^0,u^0) \hbox{ in } L^2(\R\times \T).
\end{aligned}
$$
We consider  a family  $((\eta_\eps,u_\eps))_{\eps>0}$ of solutions
of~$(SW_\e)$ with respective initial data~$(\eta_\eps^0,u_\eps^0)$ (given by
Corollary~\ref{corexistenceeta} page~\pageref{corexistenceeta}), and
$(\eta,u)$ any of its limit points.
 Finally we consider~$(\eta^*,u^*)$ in~$(L^{2} \times H^1) (\R)$ such that~$(\eta^*,u^*)$
 belongs to the kernel of~$L$. In particular by Proposition~\ref{kernel} we know that
$$
u_1^*=0 \quad \mbox{and } \quad  \beta  x_1 u_2^{*} +\d_1 \eta^* =0 .
$$
Our aim is to prove that
\begin{equation}
\label{weak-heat}
\int (\eta \eta^* + u_{2} u^*_{2}) (t,x) \: dx + \nu \int_{0}^{t} \int \nabla u_{2} \cdot 
\nabla
u^*_{2} (t',x) \: dx \: dt' = \int (\eta^{0} \eta^* + u_{2 }^0 u^*_{2}) (x) \: dx .
\end{equation}

\medskip
$\bullet$
In order to establish such an identity, the idea is to take limits in the weak form of 
System $(SW_\e)$,
which will require some further regularity on $(\eta^*,u^*)$, and then to extend the 
limiting equality  to
all vector fields~$(\eta^*,u^*)\in L^2\times H^1(\R) \cap \Ker L$ by a density argument.

Note that the classical regularization method cannot be applied here, since the kernel of 
$L$ is not stable by
convolution. In view of the explicit formula (\ref{pi-def}) page~\pageref{pi-def} giving the projector $\Pi_0$ 
(which is written in
terms of the singular pseudo-differential operator $\d_1((\beta x_1)^{-1} \cdot)$), it is 
actually natural to consider the  Hermite functions introduced in the
previous chapter,  and we recall (see~(\ref{psi-def2})
page~\pageref{psi-def2}) that
 any element of $\Ker L$  is a linear 
combination of the
following  
$$(\eta_n,u_n)= 
\frac1{\sqrt{2\pi(2n+1)}}
\left(
  \begin{array}{c}
  \displaystyle - \sqrt{   n+1 \over 2} \psi_{n-1}(x_1)-
\sqrt{ n\over 2}\psi_{n+1}(x_1)\\
0
   \\
\displaystyle  \sqrt{   n+1 \over 2} \psi_{n-1}(x_1)-
\sqrt{  n\over 2}\psi_{n+1}(x_1)
\end{array}
\right) \hbox{ for }n\geq 1,$$
$$\hbox{ and }(\eta_0,u_0)= \left(
  \begin{array}{c}
  \displaystyle \psi_0(x_1)\\
0
   \\
\displaystyle  \psi_0(x_1)
\end{array}
\right).$$
We will therefore restrict our attention to these particular vector fields which are 
smooth and integrable
against any polynomial in $x_1$ (recall that
$$\psi_n(x_1)=\exp\left(-{\beta x_1^2 \over 2} \right) P_n(x_1\sqrt{\beta})$$
where $P_n$ is the $n$-th Hermite polynomial), and then conclude by a density
argument.

\medskip
$\bullet$
 Using the conservations of mass and momentum $(SW_\e)$  it is easy to see that
\begin{eqnarray*}
&&\int \left(\eta_{\e} \eta_n + m_{\e,2} u_{n,2}\right) (t,x) \: dx + \nu \int_{0}^{t} 
\int \nabla u_{\e,2}
 \cdot \nabla
u_{n,2} (t',x) \: dx \: dt'  \\
& =& \int (\eta^{0}_\eps \eta_n + m_{\eps,2}^0 u_{n,2}) (x) \: dx +  \int_{0}^{t} \int 
\left(
m_\eps\cdot ( u_\eps \cdot \nabla  u_n)
\right) (t',x) \: dx   dt' .
\end{eqnarray*}
Now we need to take  limits as~$\e$ goes to zero in all four $\e$-dependent integrals 
appearing in that
expression.

Clearly the three first terms converge to their expected limits, as
$$
\begin{aligned}
\int \left(\eta_{\e} \eta_n + u_{\e,2} u_{n,2}\right) (t,x) \: dx \rightarrow \int 
\left(\eta
 \eta_n
 + u_{2} u_{n,2}\right) (t,x) \: dx\\
\int \left(\eta_{\e}^0 \eta_n + u_{\e,2}^0 u_{n,2}\right) (x) \: dx \rightarrow \int 
\left(\eta^0
 \eta_n
 + u_{2}^0 u_{n,2}\right) (x) \: dx
\end{aligned}
$$
and
$$
 \nu \int_{0}^{t} \int \nabla u_{\e,2}
 \cdot \nabla
u_{n,2} (t',x) \: dx \: dt'  \rightarrow  \nu \int_{0}^{t} \int \nabla u_{2}
 \cdot \nabla
u_{n,2} (t',x) \: dx \: dt'
$$
for all~$t \geq 0$,
as~$\e$ goes to zero.

So the only term we need to worry about is the coupling term
$$
 \int_{0}^{t} \int m_{\eps,2} u_{\eps,1} \d_1 u_{n,2} (t',x) \: dx   dt' .
$$
We will prove in the following lemma that it actually converges to 0~:
$$
\lim_{\e \rightarrow 0}\int_{0}^{t} \int m_{\eps,2} u_{\eps,1} \d_1 u_{n,2} (t',x) \: dx  
 dt' = 0 ,
$$ which is due to the special structure of the oscillations
pointed  out in Proposition~\ref{descriptionoscillations}. This result clearly ends the 
proof of
Theorem~\ref{weakcv}, and is proved in the next paragraph.

\subsection{The compensated compactness argument}\label{compensatedcompactnesssection}
Let us prove the following lemma.
\begin{Lem}\label{weakcp}
{ 
With the previous notation, we have locally uniformly in~$t $
$$
\lim_{\e \rightarrow 0}\int_{0}^{t} \int m_{\eps,2} u_{\eps,1} \d_1 u_{n,2} (t',x) \: dx  
 dt' = 0 .
$$
}
\end{Lem}

\begin{proof}
Let us introduce the same regularization as in Proposition \ref{descriptionoscillations}, 
defining
$$
\eta_{\e}^{\delta} = \eta_{\e} \star \kappa_{\delta},\quad u_{\e}^{\delta} = u_{\e} \star 
\kappa_{\delta} \quad \mbox{and} \quad
m_{\e}^{\delta} = m_{\e} \star \kappa_{\delta}.
$$
Then
\begin{equation}
\label{fourintegrals}
\begin{aligned}
 \int_{0}^{t} \int m_{\eps,2} u_{\eps,1} \d_1 u_{n,2} (t',x) \: dx   dt'
&= \int_{0}^{t} \int m_{\eps,2}^\delta  m_{\eps,1}^\delta \d_1 u_{n,2} (t',x) \: dx   
dt'\\
&+ \int_{0}^{t} \int m_{\eps,2}^\delta (u_{\eps,1}^\delta- m_{\eps,1}^\delta) \d_1 
u_{n,2} (t',x) \: dx   dt'\\
&+ \int_{0}^{t} \int m_{\eps,2}^\delta  (u_{\eps,1} -u_{\eps,1}^\delta )\d_1 u_{n,2} 
(t',x) \: dx   dt'\\
&+ \int_{0}^{t} \int (m_{\eps,2}- m_{\eps,2}^\delta)  u_{\eps,1}\d_1 u_{n,2} (t',x) \: dx 
  dt'.
 \end{aligned}
\end{equation}

\medskip
$\bullet$
By the energy estimates and the bounds on the Hermite functions 
given in Proposition~\ref{diag-prop} page~\pageref{diag-prop}, 
we can prove that the two last integrals converge towards zero 
as $\delta$ goes to
zero uniformly in $\eps$. Indeed for all $\alpha>0$ there exists some bounded subset 
$\Omega_\alpha\times \T$ of $\R\times \T$ such that (recalling that~$n$ is fixed)
$$ \|\d_1 u_{n,2}\|_{W^{1,\infty}(\R\setminus  \Omega_\alpha)}\leq \alpha.$$
Then, for $0<s<1$ and for any~$s'>0$,
$$
\begin{aligned}
\,&\left|
\int_{0}^{t} \int (m_{\eps,2}- m_{\eps,2}^\delta)  u_{\eps,1}\d_1 u_{n,2} (t',x) \: dx   
dt'
\right| \\
&\leq \|m_{\e,2} -  m_{\e,2}^{\delta}\|_{L^2([0,T];H^{-s-s'}(\Omega_\alpha\times \T))}
\|u_{\eps,1}\|_{L^{2}([0,T];H^1(\R\times \T))} \|\d_1 u_{n,2}\|_{W^{1,\infty}(\R)}\\
& + 2 \alpha \|m_{\eps,2}\|_{L^2([0,T];H^{-s-s'}(\R\times \T))}  
\|u_{\eps,1}\|_{L^{2}([0,T];H^1(\R\times \T))}
\end{aligned}
$$
which goes to zero  as~$\alpha$ then $\delta$ go to zero,
uniformly in~$\e$ by (\ref{uniform}), (\ref{H1bound}) and (\ref{delta-cv}).

Similarly,  we get, for $0<s<1$,
$$
\begin{aligned}
\,&\left|
\int_{0}^{t} \int m_{\eps,2}^\delta  (u_{\eps,1} -u_{\eps,1}^\delta )\d_1 u_{n,2} (t',x) 
\: dx   dt'
\right| \\
&\leq \|  m_{\e,2}^{\delta}\|_{L^2([0,T];H^{-s}(\R\times \T))}
\|u_{\eps,1} -u_{\eps,1}^\delta\|_{L^{2}([0,T];H^{s}(\Omega_\alpha\times \T))} \|\d_1 
u_{n,2}\|_{W^{1,\infty}(\R\times \T)}\\
& + 2 \alpha \|m_{\eps,2}\|_{L^2([0,T];H^{-s}(\R\times \T))}  
\|u_{\eps,1}\|_{L^{2}([0,T];H^1(\R\times \T))}
\end{aligned}
$$
which goes to zero  as~$\alpha$ then $\delta$ go to zero,
uniformly in~$\e$ by (\ref{uniform}),  (\ref{H1bound}) and (\ref{delta-cv}).

\medskip
Next we prove that for all $\delta>0$, the second integral in the right-hand side of~(\ref{fourintegrals}) 
goes to zero as $\eps$ goes to~$0$. We have
seen in~(\ref{Hsboundetau}) that $\eta_\eps u_\eps$ and consequently $m_\eps$ are uniformly bounded in 
the space~$L^2([0,T];H^s(\R\times \T))$ for $s<0$. Therefore, for fixed
$\delta>0$,
$(\eta_\eps u_\eps)^\delta$ and  $m_\eps^\delta$ are uniformly bounded in 
$L^2([0,T]\times \R\times \T)$. Then,
$$
\begin{aligned}
\left|
\int_{0}^{t} \int m_{\eps,2}^\delta (u_{\eps,1}^\delta- m_{\eps,1}^\delta) \d_1 u_{n,2} 
(t',x) \: dx   dt'
\right| \\
\leq \eps \|  m_{\e,2}^{\delta}\|_{L^2([0,T]\times \R\times \T)}
\|(\eta_\eps u_{\eps,1})^\delta\|_{L^{2}([0,T]\times \R \times \T)} \|\d_1 
u_{n,2}\|_{L^\infty(\R\times \T)}\\
\end{aligned}
$$
which goes to zero as $\eps \to 0$ for all fixed $\delta>0$.

\medskip
$\bullet$
So finally we need to consider the first term in the right-hand side 
of~(\ref{fourintegrals}).
We are going to prove that the limit of that term is zero using Proposition~\ref{weakcp}. 
Integrating by parts, we have, recalling that~$\omega_{\e}^{\delta} = \nabla^{\perp} m_{\e}^\delta$,
$$
\begin{aligned}
\int_{0}^{t} &\int m_{\eps,2}^\delta  m_{\eps,1}^\delta \d_1 u_{n,2} (t',x) \: dx   dt'\\
&=
 -\int_{0}^{t} \int  \left((\d_1 m_{\eps,2}^\delta ) m_{\eps,1}^\delta+m_{\eps,2}^\delta 
(\d_1 m_{\eps,1}^\delta)\right)
u_{n,2}
 (t',x) \: dx   dt'  \\
& =  -\int_{0}^{t} \int  \left((-\omega_\eps^\delta +\d_2 m_{\eps,1}^\delta ) 
m_{\eps,1}^\delta+m_{\eps,2}^\delta (\nabla \cdot
m_\eps^\delta -\d_2 m_{\eps,2}^\delta)\right) u_{n,2}
 (t',x) \: dx   dt'  \\
& =  -\int_{0}^{t} \int  \left((-\omega_\eps^\delta +\beta x_1 \eta_\eps^\delta ) 
m_{\eps,1}^\delta+\eta_\eps^\delta (-\beta
x_1 m_{\eps,1}^\delta +\d_2 \eta_\eps^\delta) +m_{\eps,2}^\delta (\nabla
\cdot m_\eps^\delta) \right) u_{n,2}
 (t',x) \: dx   dt'\\
&\quad -\frac12 \int_0^t \int \d_2 \left( (m_{\eps,1}^\delta )^2-(m_{\eps,2}^\delta 
)^2-(\eta_\eps^\delta)^2\right) u_{n,2}(t',x)dxdt'
\end{aligned}
$$
and the last term is  zero because $\d_2 u_{n,2}=0.$

Proposition~\ref{weakcp} now implies that
$$
\begin{aligned}
\e \partial_{t} (\omega_{\e}^{\delta}  - \beta x_{1} \eta_{\e}^{\delta} ) + \beta
m_{\e,1}^{\delta} = \e q_{\e}^{\delta} + \delta p_{\e}^{\delta},\\
\eps \d_t m_{\eps, 2}^\delta -\beta
x_1 m_{\eps,1}^\delta +\d_2 \eta_\eps^\delta=\eps s_{\eps,2}^\delta +\delta 
\sigma_{\eps,2}^\delta,\\
\eps \d_t \eta_\eps^\delta+\nabla \cdot m_\eps^\delta=0,
\end{aligned}
$$
where~$ q_{\e}^{\delta}$ and $s_\eps^\delta$  are  bounded respectively 
in~$L^{1}([0,T];L^{2}(\R\times \T))$ and  $L^{1}([0,T];H^1(\R\times
\T))$ for any~$T > 0$ uniformly in~$\e$
 (by
a constant depending on~$\delta$), and where $p_{\e}^{\delta}$ and $\sigma_\eps^\delta$ 
are  uniformly
bounded in~$\e$ and~$\delta$,  respectively in the spaces~$L^{2}([0,T];L^{2}(\R\times 
\T))$ and $L^{2}([0,T];H^1(\R\times \T))$ for any~$T>
0$.  It follows that
$$ \begin{aligned}
\int_{0}^{t} &\int m_{\eps,2}^\delta  m_{\eps,1}^\delta \d_1 u_{n,2} (t',x) \: dx   dt'\\
&=-
 \int_{0}^{t} \int \left( \frac\eps{2\beta} \d_t ( \beta x_{1} \eta_{\e}^{\delta}  - 
\omega_{\e}^{\delta})^2+ \frac\eps \beta ( \beta x_{1}
\eta_{\e}^{\delta}  - \omega_{\e}^{\delta}) q_\eps^\delta+ \frac\delta \beta ( \beta x_{1}
\eta_{\e}^{\delta}  - \omega_{\e}^{\delta}) p_\eps^\delta
\right)
u_{n,2} (t',x) \: dx dt'  \\
&-\int_0^t \int \left(- \eps \d_t(\eta_\eps^\delta m_{\eps,2}^\delta) +\eps 
\eta_\eps^\delta s_{\eps,2}^\delta +\delta
\eta_\eps^\delta \sigma_{\eps,2}^\delta\right) u_{n,2}(t',x) dxdt'
\end{aligned}
$$
Now we notice that
$$
\begin{aligned}
\left|
 \int_{0}^{t} \int  ( \beta x_{1} \eta_{\e}^{\delta}  - \omega_{\e}^{\delta})
 p_{\e}^{\delta} u_{n,2} (t',x) \: dx dt'
\right|
\leq  &C\left( T^{1/2}\| \eta_{\e}^{\delta}\|_{L^{\infty}(\R^+;L^{2}(\R\times \T))} +\|
\omega_{\e}^{\delta}\|_{L^{2}([0,T];L^{2}(\R\times \T))} \right)\\
&\times
\|  (1+x_1^2)^{1/2}u_{n,2}\|_{L^{\infty}(\R\times \T)} \| p_{\e}^{\delta}\|_{
L^{2}([0,T];L^{2}(\R\times \T))} ,
\end{aligned}
$$
and similarly
$$
\left|
 \int_{0}^{t} \int  \eta_\eps^\delta \s_{\eps,2}^\delta u_{n,2} (t',x) \: dx dt'
\right|
\leq  C T^{1/2}\| \eta_{\e}^{\delta}\|_{L^{\infty}(\R^+;L^{2}(\R\times \T))}
\|  u_{n,2}\|_{L^{\infty}(\R\times \T)} \| \s_{\e}^{\delta}\|_{
L^{2}([0,T];L^{2}(\R\times \T))} .
$$
So writing
$$
\| \omega_{\e}^{\delta}\|_{L^{2}([0,T];L^{2}(\R\times \T))}  \leq
 \|\nabla^{\perp} \cdot u_{\e}^{\delta}\|_{L^{2}([0,T];L^{2}(\R\times \T))} + \e \|\nabla^{\perp} \cdot
 (\eta_{\e}^{\delta} u_{\e}^{\delta})\|_{L^{2}([0,T];L^{2}(\R\times \T))},
$$
we infer that
$$
\begin{aligned}
\lim_{\delta \rightarrow 0}\lim_{\e \rightarrow 0} \left(\frac\delta\beta  \left|
 \int_{0}^{t} \int  ( \beta x_{1} \eta_{\e}^{\delta}  - \omega_{\e}^{\delta})
 p_{\e}^{\delta} u_{n,2} (t',x) \: dx dt'
\right|\right) = 0,  \quad \mbox{and} \\
\lim_{\delta \rightarrow 0} \left(\delta  \left|
 \int_{0}^{t} \int \eta_\eps^\delta \s_{\eps,2}^\delta u_{n,2} (t',x) \: dx dt'
\right|\right) = 0, \quad \mbox{uniformly in} \: \e.
\end{aligned}
$$
On the other hand,
$$
\begin{aligned}
&\left| \int_{0}^{t} \int  ( \beta x_{1} \eta_{\e}^{\delta}  - \omega_{\e}^{\delta})
 q_{\e}^{\delta} u_{n,2} (t',x) \: dx dt'
\right|
\leq   C\left( \| \eta_{\e}^{\delta}\|_{L^{\infty}(\R^+;L^{2}(\R\times \T))} +\|
\omega_{\e}^{\delta}\|_{L^{\infty}(\R^+;L^{2}(\R\times \T))}  
\right)\\
&\times
\|  (1+x_1^2)^{1/2}u_{n,2}\|_{L^{\infty}(\R\times \T)} \| q_{\e}^{\delta}\|_{
L^{1}([0,T];L^{2}(\R\times \T))}\\
\leq  &C\left( \| \eta_{\e}^{\delta}\|_{L^{\infty}(\R^+;L^{2}(\R\times \T))} 
+\frac1\delta \|
u_{\e}\|_{L^{\infty}(\R^+;L^{2}(\R\times \T))}
+ \e \|\nabla^{\perp} \cdot (\eta_{\e}^{\delta}u_{\e}^{\delta})\|_{L^{\infty}(\R^+;L^{2}(\R\times \T))}
 \right)\\
&\times
\|  (1+x_1^2)u_{n,2}\|_{L^{\infty}(\R\times \T)} \| q_{\e}^{\delta}\|_{
L^{1}([0,T];L^{2}(\R\times \T))} ,
\end{aligned}
$$
and
$$
\left|
 \int_{0}^{t} \int  \eta_\eps^\delta s_{\eps,2}^\delta u_{n,2} (t',x) \: dx dt'
\right|
\leq  C \| \eta_{\e}^{\delta}\|_{L^{\infty}(\R^+;L^{2}(\R\times \T))}
\| u_{n,2}\|_{L^{\infty}(\R\times \T)} \| s_{\e}^{\delta}\|_{
L^{1}([0,T];L^{2}(\R\times \T))}\\
$$
so
$$
\begin{aligned}
\lim_{\e \rightarrow 0} \left(\frac\e\beta \left|
 \int_{0}^{t} \int  ( \beta x_{1} \eta_{\e}^{\delta}  - \omega_{\e}^{\delta})
 q_{\e}^{\delta} u_{n,2} (t',x) \: dx dt'
\right|\right) = 0, \quad \mbox{for all} \: \delta > 0,\\
\lim_{\e \rightarrow 0} \left(\eps\left|
 \int_{0}^{t} \int \eta_\eps^\delta s_{\eps,2}^\delta u_{n,2} (t',x) \: dx dt'
\right|\right) = 0, \quad \mbox{for all} \: \delta > 0.
\end{aligned}
$$
So  we simply need to let~$\e$ go to zero, then~$\delta$, and the result follows.
\end{proof}

\section{Strong convergence of filtered weak solutions towards a   strong  solution}\label{strongstrongcv}
\setcounter{equation}{0}
In this   paragraph we will prove the following strong convergence
theorems. We recall that~$\Pi_\perp$ denotes the projection onto~$
(\Ker L)^\perp$, and that the spaces~$ H^s_L$ were defined and studied
in Chapter~\ref{envelope} and defined again in this chapter, page~\pageref{redefhsl}. 

The first result we will state concerns the case of smooth enough initial data, and requires no restriction on~$\beta$.
\begin{Thm}[strong convergence  for all~$\beta$]\label{strong-convergenceallbeta}
{  
Let $\Phi^0=(\eta^0,u^0)$ belong to~$ L^2(\R\times \T)$,
and 
 consider a family~$((\eta^0_\eps,u_\eps^0))_{\e > 0}$ such that
$$
\begin{aligned}
\frac12 \int \left(
|\eta^0_\eps|^2 +(1 + \e
 \eta^0_\eps)
|u^0_\eps|^2\right)dx\leq \EE^0 \quad \mbox{and} \\
\frac12 \int \left(
|\eta^0_\eps-\eta^0|^2 +(1 + \e
 \eta^0_\eps)
|u^0_\eps-u^0|^2\right)dx\to 0 \hbox{ as }\eps \to 0.
\end{aligned}
$$
For all
$\eps >0$ denote by $(\eta_\eps,u_\eps)$ a solution of~$(SW_\e)$ with initial
data~$(\eta_\eps^0,u_\eps^0)$.
  Finally suppose that~$\Pi_{\perp} \Phi^0$ belongs to~$ H^{1/2}_{L}$
   and that $ \Pi_{P} \Phi^0$ and~$ \Pi_{K} \Phi^0$ belong  to~$ H^{\alpha}_{L}$
  for some~$\alpha > 3/2$. Then 
  the sequence of filtered solutions $(\Phi_\eps)$ to~$(SW_\e)$ defined by
\begin{equation}
\label{phieps-def}
\Phi_\eps=\LL\left(-{t\over \eps}\right) (\eta_\eps,u_\eps),
\end{equation}
 converges strongly towards~$\Phi$ in $L^2_{loc}
([0,T^*[; L^2(\R\times \T))$, where~$\Phi $ is the unique solution on~$[0,T^*[$
 of~$ (SW_0)$
 constructed in Theorem~\ref{limit-systemallbeta},  page~\pageref{limit-systemallbeta} .}
\end{Thm}

The next theorem requires less assumptions on the initial data; on the other hand one must first remove a countable set of 
values for~$\beta$.
\begin{Thm}[strong convergence  for generic~$\beta$]\label{strong-convergencegenericbeta}
{  
There is a  countable subset~${\mathcal N}$ of~$\R^{+}$ such that for any~$\beta \in \R^{+} \setminus {\mathcal N}$, 
the following result holds. 
Let $\Phi^0\in L^2(\R\times \T;\R^3)$ be given,
and 
 consider a family~$((\eta^0_\eps,u_\eps^0))_{\e > 0}$ such that
\begin{equation}
\label{init-cv}
\begin{aligned}
\frac12 \int \left(
|\eta^0_\eps|^2 +(1 + \e
 \eta^0_\eps)
|u^0_\eps|^2\right)dx\leq \EE^0 \quad \mbox{and} \\
\frac12 \int \left(
|\eta^0_\eps-\eta^0|^2 +(1 + \e
 \eta^0_\eps)
|u^0_\eps-u^0|^2\right)dx\to 0 \hbox{ as }\eps \to 0.
\end{aligned}
\end{equation}
For all
$\eps >0$ denote by $(\eta_\eps,u_\eps)$ a solution of~$(SW_\e)$ with initial
data~$(\eta_\eps^0,u_\eps^0)$.
  Finally suppose that~$ \Pi_{P} \Phi^0$ and~$ \Pi_{K} \Phi^0$ belong  to~$ H^{\alpha}_{L}$
  for some~$\alpha > 1/2$. Then 
  the sequence of filtered solutions $(\Phi_\eps)$ to~$(SW_\e)$ defined by (\ref{phieps-def})
$$\Phi_\eps=\LL\left(-{t\over \eps}\right) (\eta_\eps,u_\eps) $$
 converges strongly towards~$\Phi$ in $L^2_{loc}
(\R^{+}; L^2(\R\times \T))$, where~$\Phi $ is the unique solution
 of~$ (SW_0)$
 constructed in Theorem~\ref{limit-systemgeneric},  page~\pageref{limit-systemgeneric} .}
\end{Thm}

\begin{Rem}\label{onlyoneproof}{ 
$\bullet$ Note that  definition (\ref{phieps-def}) of $\Phi_\eps$ does make sense since as stated
in Corollary~\ref{corexistenceeta}, one has a $L^2$ bound on $u_\eps$.  \medskip

$\bullet$ The strong compactness of $(\Phi_\eps)$ in $L^2_{loc}(\R^+,
L^2(\R\times \T))$  cannot be
obtained directly using some a priori estimates. Indeed we have a
priori no uniform  regularity on
$\eta_\eps$ with respect to the space variable $x$ (besides we expect  the limiting system to
be a mixed hyperbolic-parabolic system). 

\medskip

$\bullet$  The  proof of both  convergence results is   based on a weak-strong stability
property of $(SW_\e)$. It is  therefore crucial to be able to construct a
smooth approximate solution $\Phi_{app}$ to $\Phi_\eps$, writing an
asymptotic  expansion in $\eps$
whose first term is $\Phi$. The regularity assumptions on the initial data stated in   both theorems are precisely that
enabling one to guarantee that the limit system  has a unique, stable solution and propagates regularity. In particular it
should be noted that in both cases, the assumptions on the initial data imply that~$\nabla \cdot \Phi'$ belongs to~$L^{1}([0,T];
L^{\infty}(\R\times \T))$ (see Proposition~\ref{divergence-prop}). Once the setting is posed so that the limit system does satisfy those properties, the proofs are
very much  the same in 
both cases. So in the following we will only prove Theorem~\ref{strong-convergencegenericbeta}, and leave to the
reader the easy adaptations in the case of Theorem~\ref{strong-convergenceallbeta}.

 }
  \end{Rem} 
\begin{proof}
As noted in Remark~\ref{onlyoneproof} above, we will only prove Theorem~\ref{strong-convergencegenericbeta}
here.

The idea is, as usual in filtering methods, to start by 
approximating the solution of the limit system, and then to use a
weak-strong stability method to  conclude.

So let us consider the solution~$\Phi$  of~$ (SW_0)$ constructed in the previous
chapter (see Theorem~\ref{limit-systemgeneric} page~\pageref{limit-systemgeneric}), 
which we truncate in 
the following way:
\begin{equation} 
\label{PhiN}
\Phi_N = J_N \Pi_{\perp} \Phi   + \Pi_0 \Phi_N, 
\end{equation}
where $J_N$ is the spectral truncation defined by
\begin{equation}
\label{JN}
J_N=\sum_{i\lambda \in {\mathfrak S}_N} \Pi_\lambda
\end{equation}
with
\label{sigmaN}
$$
{\mathfrak S}_N = \left\{i\tau (n,k,j) \in {\mathfrak S} \Big/ n \leq N, \:
  |k|\leq N  \right\},
$$
 and  $\Pi_\perp$ denotes as previously  the projection
onto~$(\Ker L)^{\perp}$. Finally~$ \Pi_0 
\Phi_N $ solves
$$
\begin{aligned}
\partial_t \Pi_0 \Phi_N -\nu \Pi_0 \Delta' \Pi_0 \Phi_N = 0  &
\\
\Pi_0  \Phi_{N|t = 0} = \sum_{0 \leq n \leq N} \Pi_{n,0,0} \Phi^0.
\end{aligned}
$$
We recall that~$ \Pi_{n,0,0}$ denotes the projection onto the
eigenvector~$ \Psi_{n,0,0}$ of~$ \Ker L$.
Then  for all fixed
$N\in \N$ we have (see Theorem~\ref{limit-systemallbeta})
\begin{equation}\label{pi0phiregular}
 \Pi_0 \Phi_N  \:  \mbox{belongs to} \: L^\infty(\R^+;H^\sigma_L) ,
 \: \forall \sigma \geq 0.
\end{equation}
Recall that such a result means that~$ \Pi_0 \Phi_N$ is as smooth as needed, and decays as fast as needed when~$x_{1}$
goes to infinity.

Moreover by the stability of the limit system (which is linear) we have of course, for all~$T>0$, 
\begin{equation}\label{convergencepiophin}
\lim_{N \to \infty}\| \Pi_0   \Phi_N  - \Pi_{0} \Phi\|_{L^{\infty}([0,T];L^{2}(\R\times\T))} = 0 .
\end{equation}
 Note also  that for all fixed
$N\in
\N$, using the smoothness and the decay of the eigenvectors of $L$, we
get for  any polynomial~$Q \in \R[X]$ 
\begin{equation} 
\label{phiN-reg}
Q(x_1)  \Phi_{N}\in L^\infty([0,T]; C^{\infty}(\R\times\T))
\end{equation}
 We have moreover,  for all~$T>0$, 
\begin{equation} 
\label{phiN-cv}
 \Bigl( \|\Pi_{\perp}( \Phi -\Phi_N)\|_{L^\infty([0,T]; L^2
(\R\times\T))} + \|(\Pi_{K}+ \Pi_{P})( \Phi -\Phi_N)\|_{L^\infty([0,T]; H^{\alpha}_{L} )} \Bigr)
\to 0 \hbox{ as } N\to
\infty,
\end{equation}
and
\begin{equation}
\label{phiN-cv2}
 \Bigl(
\| \Pi_{\perp}( \Phi -\Phi_N ) \|_{L^2([0,T]; H^{ 1}_L)}
+ \|(\Pi_{K}+ \Pi_{P})( \Phi -\Phi_N)\|_{L^2([0,T]; H^{\alpha+1}_L)}
 \Bigr)
 \to 0 \hbox{ as } N\to
\infty.
\end{equation}
Finally since~$  {J_N}$ commutes with $\Delta'_L$, the vector field~$\Phi_{N}$ satisfies  the
approximate limit  filtered system
\begin{equation}
\label{app-lim-filtered}
\begin{aligned}
\d_t \Phi_N + {J_N} Q_L (\Phi,\Phi)-\nu  \Delta'_L \Phi_N=0,\\
\Phi_{N|t=0} = {J_N} \Phi^0.
\end{aligned}
\end{equation}
Conjugating this equation by the semi-group $\LL$ leads then to
$$\d_t \left(\LL\left({t\over \eps}\right)\Phi_N\right) +\frac1\eps L\left(\LL\left({t\over
\eps}\right)\Phi_N\right)+  {J_N} Q_L \left(\LL\left({t\over
\eps}\right)\Phi,\LL\left({t\over \eps}\right)\Phi\right)-\nu  \Delta'_L \LL\left({t\over
\eps}\right)
\Phi_N=0,$$
using the definitions (\ref{Q-lap-L}) of $Q_L$ and $\Delta'_L$. Let us now rewrite
 this last equation in a convenient way
$$\begin{aligned}
\d_t &\left(\LL\left({t\over \eps}\right)\Phi_N\right) +\frac1\eps L\left(\LL\left({t\over
\eps}\right)\Phi_N\right)+Q \left(\LL\left({t\over
\eps}\right)\Phi_N,\LL\left({t\over \eps}\right)\Phi_N\right)-\nu  \Delta' \LL\left({t\over
\eps}\right)
\Phi_N\\
&=  (Q -Q_L)\left(\LL\left({t\over
\eps}\right)\Phi_N,\LL\left({t\over \eps}\right)\Phi_N\right)-\nu (
\Delta'-\Delta'_L) \LL\left({ t\over
\eps}\right) \Phi_N\\ &+(Id-  {J_N}) Q_{L}\left(\LL\left({t\over
\eps}\right)\Phi ,\LL\left({t\over \eps}\right)\Phi \right)+   Q_L \left(\LL\left({t\over
\eps}\right)(\Phi_N-\Phi),\LL\left({t\over \eps}\right)(\Phi_N+\Phi)\right) .
\end{aligned}
$$
Because of~(\ref{phiN-cv}) and~(\ref{convergencepiophin}), the last   term  in
the right-hand side is  expected to be small
when~$N$ is large, uniformly in~$\e$, and similarly for the third
term, using the stability of  the limit 
system proved in the previous chapter. So we are left with the
first two terms, which as  usual cannot
be dealt with so easily since they do not converge strongly towards zero. However they are 
 fast oscillating terms,  and  will be treated by introducing a small
 quantity~$\eps \phi_N$ (which  will be
small when~$\e$ goes to zero, for each fixed~$N$), 
so that
$$  (\d_t+\frac1\eps L) \left(\LL\left({t\over \eps}\right) \eps \phi_N\right) \sim -  (Q
-Q_L)\left(\LL\left({t\over
\eps}\right)\Phi_N,\LL\left({t\over \eps}\right)\Phi_N\right)+\nu
(\Delta'-\Delta'_L)\LL\left({t\over
\eps}\right)
\Phi_N.$$
Let us now define
\begin{equation}\label{defpetitPhiN}
\phi_N= -\sum_{\lambda \neq \mu+\tilde \mu \atop i\lambda \in {\mathfrak S}, i\mu, i\tilde \mu 
\in  {{\mathfrak S}_N}}
{ e^{i\frac t\eps (\lambda -\mu-\tilde \mu)}\over i(\lambda -\mu-\tilde \mu)} 
\Pi_\lambda   Q( \Pi_\mu \Phi_N,
\Pi_{\tilde \mu }\Phi_N) 
+\nu \sum_{\lambda \neq \mu,\atop
i\lambda \in {\mathfrak S}, i\mu \in {\mathfrak S}_N  }
{ e^{i\frac t\eps (\lambda -\mu)}
\over i(\lambda
-\mu)} \Pi_\lambda \Delta' \Pi_\mu \Phi_N, 
\end{equation}
 and consider 
$$
\Phi_{\eps,N}=\Phi_N+\eps \phi_N.
$$
Let us prove the following result.
\begin{Prop}\label{approximate-sol}
{ 
For all but a countable number of~$\beta$, the following result
holds. 
Consider a vector field~$\Phi^0=(\eta_0, u_0) \in L^2(\R \times \T)$, with~$ (\Pi_P + \Pi_{K})
\Phi^0 $ in~$
    H^{\alpha}_L$
for some $\alpha>1/2$.  Denote by $\Phi$ the associate solution of 
$(SW_0)$. Then 
 there exists  a
family 
$(\eta_{\eps,N}, u_{\eps,N})=\LL     \left( {t\over
\eps}\right)\Phi_{\eps,N} 
$, bounded  in the space~$L^\infty_{loc}( \R^{+}, L^{2}) 
\cap L^2_{loc}( \R^{+},   H^1  ) $, 
 such that~$ (\Pi_P + \Pi_{K} )(\eta_{\eps,N}, u_{\eps,N})$ is
uniformly bounded in the space~$L^\infty_{loc}(\R^{+},     H^{\alpha}_L) 
\cap L^2_{loc}(\R^{+},   H^{\alpha+1}_L ) $,  and satisfying
the following properties:

$\bullet$  $\Phi_{\eps,N} $  behaves asymptotically as $ \Phi$ as $\eps \to 0$ and $N\to \infty$~:
\begin{equation}\label{asympt-app}
\forall T>0,\quad 
\lim_{N\to \infty} \lim_{\eps \to 0} \left \| \Phi_{\eps,N}-\Phi\right\| _{L^\infty([0,T];
 L^2(\R\times\T))} = 0;
\end{equation}

$\bullet$    for all  $N\in \N$, $(\eta_{\eps,N}, u_{\eps,N}) $ is smooth:
 for all~$ T>0$ and all~$ Q\in \R[X]$,
\begin{equation}\label{reg-app}
    Q(x_1)(\eta_{\eps,N}, u_{\eps,N})  \:  \mbox{is bounded in} \: 
L^\infty([0,T];C^{\infty}(
\R\times\T)) , \:  \mbox{uniformly in} \: \e;
\end{equation}

$\bullet$ $(\eta_{\eps,N}, u_{\eps,N})$ satisfies the uniform regularity estimate
\begin{equation}\label{unifregularity}
\forall T>0,\quad \sup_{N \in \N} \lim_{\eps \to 0} \|\nabla \cdot  u_{\eps,N} \|_{
 L^1([0,T];L^\infty(\R\times\T))}\leq C_T;
\end{equation}

$\bullet$ $(\eta_{\eps,N}, u_{\eps,N}) $     satisfies approximatively the viscous Saint-Venant
system $(SW_\e)$~:
\begin{equation}\label{SW-app}
\d_t (\eta_{\eps,N}, u_{\eps,N}) +\frac1\eps L(\eta_{\eps,N}, u_{\eps,N})+Q
\left((\eta_{\eps,N}, u_{\eps,N}),(\eta_{\eps,N}, u_{\eps,N})\right)-\nu  \Delta'
(\eta_{\eps,N}, u_{\eps,N})=R_{\eps,N}
\end{equation}
where $R_{\eps,N}$  goes to 0 as $\eps \to 0$ then $N\to \infty$:
\begin{equation}\label{repsNgoestozero}
\lim_{N\to \infty} \lim_{\eps \to 0}  \left(\| R_{\eps, N}\|_{L^1([0,T];L^2(\R\times\T))} 
+\eps \|R_{\eps,N}\|_{L^\infty([0,T]\times\R\times\T)}  \right) =0  .
\end{equation}
}
\end{Prop}
\begin{proof}
Let us define $\Phi_N$ as in (\ref{PhiN}) and~$\phi_N $ as
in~(\ref{defpetitPhiN}). We can write\label{phiN1}
$$
\phi_N = \phi_N^{(1)} +\phi_N^{(2)}, \quad \mbox{with}
$$
$$
\begin{aligned}
\phi_N^{(1)}&= -\sum_{\lambda \neq \mu+\tilde \mu \atop i\lambda 
\in {\mathfrak S}, i\mu, i\tilde \mu \in  {{\mathfrak S}}_N}
{ e^{i\frac t\eps (\lambda -\mu-\tilde \mu)}\over i(\lambda
  -\mu-\tilde \mu)}  \Pi_\lambda    Q( \Pi_\mu \Phi_N,
\Pi_{\tilde \mu }\Phi_N) \\
\phi_N^{(2)}&= \nu \sum_{\lambda \neq \mu,\atop
i\lambda \in {\mathfrak S}, i\mu \in  {{\mathfrak S}}_N}{
e^{i\frac t\eps (\lambda  -\mu)}\over i(\lambda
-\mu)} \Pi_\lambda \Delta' \Pi_\mu \Phi_N.
\end{aligned}
$$
We will check that the approximate solution $\Phi_{\eps,N}$ defined by
\begin{equation}
\label{phi-epsN}
\Phi_{\eps,N}=\Phi_N+\eps \phi_N
\end{equation}
satisfies the required properties.

It will be useful to notice that there are two positive functions~$
\Lambda_{\pm} (N)$ such that if~$ \mu$ belongs to~${\mathfrak S}_N $,
then either~$ \mu = 0$ or
$$
0 < \Lambda_{-} (N) \leq \mu \leq \Lambda_{+} (N) < +\infty.
$$
It will also  be useful to recall that, considering
 the asymptotics of~$\lambda = \tau (n,k,j)$ as~$k$ or~$n$ go to infinity, 
  the operator~$\displaystyle\sum_{\lambda \neq 0} \frac{\Pi_{\lambda} }\lambda $ is continuous
 from~$H^{\sigma}_{L}$
to~$H^{\sigma-1}_{L}$ for any given~$\sigma \in \R$. 

Finally we recall that the spectrum of~$L$ admits only~$0$ and~$\infty$ as accumulation points.

$\bullet$ The correction $\phi_N$ is defined as the sum of two terms.

Let us consider the first one, $\phi_N^{(1)}$. It can 
in turn be written
$$
\phi_N^{(1)} =  \Pi_{0}\phi_N^{(1)} + \Pi_{\perp}\phi_N^{(1)} .
$$
The first part, $\Pi_{0}\phi_N^{(1)}$, is easy  to handle since~$\Pi_{0} $ is of course
continuous from~$H^\sigma_{L}$ to~$H^\sigma_{L}$ for any~$\sigma$. 

Let us now study~$ \Pi_{\perp}\phi_N^{(1)} $.
Clearly $   Q(\Pi_{\mu} \Phi_{N},
\Pi_{\tilde\mu}\Phi_{N})$ 
is in~$H^{\sigma}_{L}$ for any~$\sigma \geq 0$.  We infer that~$ \Pi_{\perp}\phi_N^{(1)}$ 
belongs to~$H^{\sigma}_{L}$ for any~$\sigma \geq 0$. 
Indeed if~$\mu + \tilde \mu \neq 0$ then~$  | \lambda -\mu-\tilde \mu |$ is bounded from below
since~Ê$\mu$ and~$\tilde \mu$ are in~${\mathfrak S}_N$ and the   accumulation points of~$ \lambda$
are 0 and~$\infty$.  On the other hand if~$\mu + \tilde \mu = 0$ then we use    the continuity 
property of~$\displaystyle\sum_{\lambda \neq 0} \frac{\Pi_{\lambda} }\lambda$ recalled
above. 

So we find that for all polynomials $Q\in \R[X]$,  
$$Q(x_1)  \phi_N^{(1)} \in
L^\infty([0,T];C^\infty(
\R
\times \T)),$$
as well as
$$Q(x_1)\LL\left({t\over \eps}\right)  \phi_N^{(1)} \in
L^\infty([0,T];C^\infty(
\R
\times \T)).$$

This obviously implies that for all $k\in \N$
$$\forall N\in \N, \quad \lim_{\eps \to 0}\| \eps
Q(x_1)\LL\left({t\over \eps}\right)  \phi_N^{(1)} \|_{L^\infty( [0,T]; C^k(\R\times \T))}=0.$$

The second term $\phi_N^{(2)}$ is dealt with similarly, splitting it into two terms:
$$ \phi_N^{(2)}= \nu \sum_{\lambda \neq \mu,\atop
i\lambda \in {\mathfrak S}, i\mu \in   {\mathfrak S}_N\setminus \{0\}}{
e^{i\frac t\eps (\lambda -\mu)}\over i( \lambda
-\mu)} \Pi_\lambda \Delta' \Pi_\mu \Phi_N +\nu \sum_{\lambda \neq 0 ,\atop
i\lambda \in {\mathfrak S}}{ e^{i\frac t\eps (\lambda )}\over i\lambda
} \Pi_\lambda \Delta' \Pi_0 \Phi_N.
$$
Because of the relations (\ref{psi-identities}) satisfied by the
Hermite functions, it is easy to see that the  first
contribution can be rewritten as a finite combination of some
eigenvectors of $L$ (which are smooth functions rapidly decaying in
$x_1$), and the   second contribution is dealt with again by using the fact  that the
 operator~$\displaystyle\sum_{\lambda \neq 0} \frac{\Pi_{\lambda} }\lambda $ is continuous
 from~$H^{\sigma}_{L}$
to~$H^{\sigma-1}_{L}$.

We conclude that for all $Q\in \R[X]$
$$Q(x_1)\LL\left({t\over \eps}\right) \phi_N^{(2)} \in
L^\infty([0,T];C^\infty(
\R
\times \T)).$$
and thus
$$\forall N\in \N, \quad \lim_{\eps \to 0}\| \eps Q(x_1)
\LL\left({t\over \eps}\right) \phi_N^{(2)} \|_{ L^\infty([0,T]; C^{k}(\R\times \T))}=0, \: \forall
k \in \N.$$

Combining these results with (\ref{pi0phiregular}), (\ref{convergencepiophin}),  (\ref{phiN-reg}) and (\ref{phiN-cv})
leads to (\ref{asympt-app})  and (\ref{reg-app}).

$\bullet$ The uniform regularity estimate~(\ref{unifregularity}) is obtained in a very
similar way. Of course the regularity of the  correction
established previously  shows that 
its contribution to $\nabla \cdot u_{\eps,N}$ converges to zero as~$\eps$ goes to~$0$ in the sense of smooth  functions
rapidly decaying with respect to $x_1$. Therefore the only point to be checked is that
$$\nabla \cdot \left( \LL\left({t\over \eps}\right)\Phi_N\right)'
\hbox{ is uniformly bounded  in
}L^1([0,T];L^\infty(\R\times \T)),$$
which is obtained as the regularity property stated in Proposition \ref{divergence-prop}, remarking that the~$L^\infty$ bound comes from
estimates which are stable by the  truncation $ {J_N}$ and
by conjugation by the  semi-group~$\LL$. 
Indeed for almost all $\beta$, provided that 
$$(\Pi_K+\Pi_P) \Phi^0 \in   H^\alpha_L \hbox{ for some }\alpha>\frac12\virgp$$
any weak solution $\Phi$ to $(SW_0)$  satisfies the
following estimates 
$$
\left\|\nabla \cdot \Phi'\right\|_{ L^1([0,T]; L^\infty(\R\times \T))}\leq C_T 
$$
as well as
$$
\left\|\nabla \cdot \left( \LL\left({t\over
        \eps}\right)\Phi_N\right)'\right\|_{ L^1([0,T];  L^\infty(\R\times \T))}\leq C_T,$$
where $C_T$ depends only on $T\in \R^+$, $\| \Phi^0\|_{L^2(\R\times
  \T)}$ and $\|(\Pi_K+\Pi_P)  \Phi^0\|_{  H^\alpha_L }$ (neither on $N$ nor on~$\eps$.)

$\bullet$ It remains then to establish the equation satisfied by
$(\eta_{\eps,N},u_{\eps,N})$.  A direct computation provides
$$
\begin{aligned}
\eps \d_t \phi_N&= -\sum_{\lambda \neq \mu+\tilde \mu \atop i\lambda \in {\mathfrak S},
 i\mu, i\tilde \mu \in   {{\mathfrak S}}_N}
 e^{i\frac t\eps (\lambda -\mu-\tilde \mu)} \Pi_\lambda 
 Q( \Pi_\mu \Phi_N,
\Pi_{\tilde \mu }\Phi_N)\\
&+\nu \sum_{\lambda \neq \mu,\atop
i\lambda \in {\mathfrak S}, i\mu \in   {{\mathfrak S}}_N}
 e^{i\frac t\eps (\lambda  -\mu)} \Pi_\lambda \Delta' \Pi_\mu
\Phi_N\\
& -2\eps\sum_{\lambda \neq \mu+\tilde \mu \atop i\lambda \in {\mathfrak S}, i\mu,
 i\tilde \mu \in  {{\mathfrak S}}_N}
{ e^{i\frac t\eps (\lambda -\mu-\tilde \mu)}\over i(\lambda
 -\mu-\tilde \mu)} \Pi_\lambda   Q(   \Pi_\mu \d_t \Phi_N, 
\Pi_{\tilde \mu }\Phi_N)\\
&+\eps \nu \sum_{\lambda \neq \mu,\atop
i\lambda \in {\mathfrak S}, i\mu \in   {{\mathfrak S}}_N}{
 e^{i\frac t\eps (\lambda  -\mu)}\over i(\lambda
-\mu)} \Pi_\lambda \Delta' \Pi_\mu \d_t \Phi_N
\end{aligned}
$$
By (\ref{app-lim-filtered}) we infer that $\d_t \Phi_N$ is smooth
and rapidly decaying (recalling in particular that~$\Pi_{0} Q_{L } = 0$
due to Proposition~\ref{limitgeolinear} page~\pageref{limitgeolinear}), and  thus the last two terms
go to zero as $\eps \to 0$ (for all fixed $N$). The previous identity
can therefore be  rewritten
\begin{equation}
\label{varphiN-eq}
\eps \d_t \phi_N= -\LL\left(-{t\over \eps}\right)  
(Q-Q_L) \left( \LL\left({ t\over \eps}\right) \Phi_N,
\LL\left({t\over \eps}\right) \Phi_N \right)
+\nu\LL\left(-{t\over \eps}\right)  (\Delta'-\Delta'_L) \LL\left( {t\over \eps}\right)  \Phi_N+r_{\eps,N}
\end{equation}
where
\begin{equation}
\label{r-est}
\forall k \in \N, \forall N\in \N, \forall Q\in \R[X],\quad \lim_{\eps\to 0}\| Q(x_1)
r_{\eps,N}\|_{L^\infty( [0,T];
C^{k}(\R\times \T))}=0.
\end{equation}
Now let us recall that
\begin{equation}
\label{phiN-eq}\begin{aligned}
\d_t &\left(\LL\left({t\over \eps}\right)\Phi_N\right) +\frac1\eps L\left(\LL\left({t\over
\eps}\right)\Phi_N\right)+Q \left(\LL\left({t\over
\eps}\right)\Phi_N,\LL\left({t\over \eps}\right)\Phi_N\right)-\nu  \Delta' \LL\left({t\over
\eps}\right)
\Phi_N\\
&= (Q -Q_L)\left(\LL\left({t\over
\eps}\right)\Phi_N,\LL\left({t\over \eps}\right)\Phi_N\right)-\nu (
\Delta'-\Delta'_L) \LL \left({t\over
\eps}\right) \Phi_N\\ &
+(Id- {J_N}) Q_{L}\left(\LL\left({t\over
\eps}\right)\Phi ,\LL\left({t\over \eps}\right)\Phi \right)+
  Q_L \left(\LL \left({t\over
\eps}\right)(\Phi_N-\Phi),\LL\left({t\over \eps}\right)(\Phi_N+\Phi)\right).
\end{aligned}
\end{equation}
We thus have, recalling that
$$
 (\eta_{\eps,N},u_{\eps,N}) =  \LL\left({t\over\eps}\right) (\Phi_{N} + \e \phi_{N}), 
$$
\begin{equation}
\label{eq-app}
\begin{aligned}
\d_t (\eta_{\eps,N},u_{\eps,N}) &+\frac1\eps
L(\eta_{\eps,N},u_{\eps,N})+Q \left((\eta_{\eps,N},
  u_{\eps,N}),(\eta_{\eps,N},u_{\eps,N}) \right)-\nu 
\Delta' (\eta_{\eps,N},u_{\eps,N})\\
&=(Id- {J_N}) Q_{L}\left(\LL\left({t\over
\eps}\right)\Phi ,\LL\left({t\over \eps}\right)\Phi \right)+   Q_L \left(\LL\left({t\over
\eps}\right)(\Phi_N-\Phi),\LL\left({t\over \eps}\right)(\Phi_N+\Phi)\right)\\
& +\eps Q \left( \LL\left({t\over \eps}\right) \phi_N,\LL\left({t\over \eps}\right) (2\Phi_N+\eps
\phi_N)\right)-\eps \nu \Delta'\left(\LL\left({t\over \eps}\right)\phi_N\right)+r_{\eps,N}
\end{aligned}
\end{equation}
Note that the regularity estimates on $\Phi_N$ and $\phi_N$ allow to
prove that the two last  explicit terms in the
right-hand side go to zero as $\eps \to 0$ (for all fixed $N$), and
therefore to incorporate them  into the remainder~$r_{\eps,N}$.

The stability of the limiting filtered system $(SW_0)$ allows to prove that the second term in the
right-hand side of (\ref{eq-app}) goes to zero as $N\to \infty$ uniformly in $\eps$. We have indeed 
$$
\left\|  Q_L \left(\LL\left({t\over
\eps}\right)(\Phi_N-\Phi),\LL\left({t\over
\eps}\right)(\Phi_N+\Phi)\right)\right\|_{L^2(\R\times \T)}= \|  Q_L
\left(\Phi_N-\Phi,\Phi_N+\Phi\right)\|_{L^2(\R\times \T)}
$$
and recalling that only the Kelvin waves can have resonances,
$$
\longformule{\| Q_L
\left(\Phi_N-\Phi,\Phi_N+\Phi\right) \|_{L^2(\R\times \T)} \leq  
\|  Q_L
\left(\Pi_K (\Phi_N-\Phi),\Pi_K (\Phi_N+\Phi)\right)\|_{L^2(\R\times \T)} }{ + \| Q_L
\left(\Pi_0 (\Phi_N-\Phi),\Pi_\perp
  (\Phi_N+\Phi)\right)\|_{L^2(\R\times \T)} +   \|  Q_L
\left(\Pi_\perp (\Phi_N-\Phi),\Pi_0 (\Phi_N+\Phi)\right)\|_{L^2(\R\times \T)}  , }
$$
so by Proposition~\ref{trilinear}
page~\pageref{trilinear} and two-dimensional product rules on the Kelvin part (recall as in the previous chapter
that~$H^{s}$ and~$H^{s}_{L}$ spaces coincide in the case of Kelvin modes) we infer that
\begin{eqnarray}
  \|  Q_L
\left(\Phi_N-\Phi,\Phi_N+\Phi\right)\|_{L^2(\R\times \T)}
&\leq &  
 C_\alpha \|\Pi_K (\Phi_N-\Phi)\|_{   H^{\alpha+1}_L} \| \Pi_K (\Phi_N+\Phi)\|_{ 
H^{\alpha }_L} \label{estimateK-K1}\\
&+ &C \|\Pi_0 (\Phi_N-\Phi)\|_{L^2(\R\times\T)}  \| \Pi_\perp (\Phi_N+\Phi)\|_{ 
H^{1}_L}\nonumber\\
& +& C  \|\Pi_0 (\Phi_N+\Phi)\|_{L^2(\R\times\T)}  \| \Pi_\perp (\Phi_N-\Phi)\|_{ 
H^{ 1}_L}\nonumber.
\end{eqnarray} 

So by (\ref{phiN-cv}) and~(\ref{phiN-cv2}) we conclude that
$$ \lim_{N\to \infty} \|  Q_L
\left(\Phi_N-\Phi,\Phi_N+\Phi\right)\|_{L^2([0,T];L^2(\R\times \T))}=0.$$

Let us estimate the first term in the right side of (\ref{eq-app}). We can write as above
(recalling that~$Q(\Pi_{0} \: 
\cdot \: , \Pi_{0}\:  \cdot\:  ) = 0$)
$$
\left\|  Q_{L} \left(\LL\left({t\over
\eps}\right)\Phi ,\LL\left({t\over
\eps}\right)\Phi \right)\right\|_{L^{2}} \leq
 \left\|  Q_{L} \left(\Pi_K  \Phi ,\Pi_K \Phi \right)\right\|_{L^2 } + 2
\left\|  Q_{L} \left(\Pi_0 \Phi ,\Pi_\perp \Phi \right)\right\|_{L^2 }   
$$
so  we find that
\begin{equation} \label{estimateK-K12}
 \left\|  Q_L \left(\LL\left({t\over
\eps}\right)\Phi_N,\LL\left({t\over
\eps}\right)\Phi_N\right)\right\|_{L^2} 
\leq C \left(\| \Pi_K  \LL\left({t\over
\eps}\right)\Phi_N\|^2_{H^{ 1}_L} + \|\Pi_0 \Phi_N\|^2_{L^2} \right)
\end{equation}
and thus
$$\lim_{N\to \infty} \sup_{\eps}\| (Id- {J_N}) Q_{L} \left(\LL\left({t\over
\eps}\right)\Phi ,\LL\left({t\over \eps}\right)\Phi \right)\|_{L^2([0,T];L^2(\R\times \T))}=0.$$

 Note that in the case when $\beta$ belongs to~$\mathcal N$
  (Theorem~\ref{strong-convergenceallbeta}), equations~(\ref{estimateK-K1})
and ~(\ref{estimateK-K12}) must be adapted using the third estimate of Proposition~\ref{trilinear}.

Finally to prove that for all~$N  \in \N$, for all~$ T>0$, the quantity
$$
\e    Q_L \left(\LL\left({t\over
\eps}\right)(\Phi_N-\Phi),\LL\left({t\over
\eps}\right)(\Phi_N+\Phi)\right) +\e   (Id- {J_N}) Q_{L} \left(\LL\left({t\over
\eps}\right)\Phi ,\LL\left({t\over \eps}\right)\Phi \right)
$$
goes to zero as~$\e$ goes to zero, in the space~$L^\infty([0,T] \times \R\times \T)$, we
simply notice that
$$
\longformule{Q_L \left(\LL\left({t\over
\eps}\right)(\Phi_N-\Phi),\LL\left({t\over
\eps}\right)(\Phi_N+\Phi)\right) +  (Id- {J_N}) Q_{L} \left(\LL\left({t\over
\eps}\right)\Phi ,\LL\left({t\over \eps}\right)\Phi \right)}{= J_{N}
Q_L \left(\LL\left({t\over
\eps}\right)(\Phi_N-\Phi),\LL\left({t\over
\eps}\right)(\Phi_N+\Phi)\right) +  (Id- {J_N}) Q_{L} \left(\LL\left({t\over
\eps}\right)\Phi_{N} ,\LL\left({t\over \eps}\right)\Phi_{N} \right)
}
$$
and the convergence result is obvious if one considers the right-hand side, simply because those terms are smooth
for each fixed~$N$. 

Combining all the previous estimates shows that
$(\eta_{\eps,N},u_{\eps,N})$ satisfies the expected approximate
equation~(\ref{SW-app}), where $R_{\eps,N}$ satisfies the expected estimate (\ref{repsNgoestozero}) as well as 
\begin{equation}
\label{epsr-cv}
\lim_{N\to \infty} \lim_{\eps \to 0} \eps \|R_{\eps,N}\|_{L^\infty([0,T]\times \R\times \T)} =0.
\end{equation} 
Proposition~\ref{approximate-sol} is proved.
\end{proof}

 Equipped with that result, we are now ready to prove the strong
 convergence theorem. The method  relies
 on a weak-strong stability method which we shall now detail. We are going to prove that
 \begin{equation}\label{limfinale}
 \lim_{N\to \infty} \lim_{\eps \to 0}   \|(\eta_\eps,u_\eps)-(\eta_{\eps,N},
u_{\eps,N})\|_{L^2([0,T]\times\R\times\T)}
= 0,
\end{equation}
where $(\eta_{\eps,N}, u_{\eps,N})$ is the approximate solution to $(SW_\e)$ defined in
Proposition \ref{approximate-sol}.  Note that
 combining this   estimate with the fact that $(\eta_{\eps,N},u_{\eps,N})$
  is close to $\LL\left({t\over \eps}\right) \Phi$  provides the expected convergence, namely 
  the fact that
 $$\forall T>0,\quad \lim_{\eps \to 0} 
 \left \| (\eta_\eps,u_\eps)-\LL\left({t\over \eps}\right) 
 \Phi \right\|_{L^2([0,T]\times\R\times \T)} =0. $$

The key to the proof of~(\ref{limfinale}) lies in the following proposition.
\begin{Prop}\label{strong-weak}
{ 
There is a constant~$C$ such that the following property holds.
Let~$(\eta^0,u^0)$ and~$(\eta^0_\eps,u_\eps^0)$ satisfy assumption~(\ref{init-cv}), and let~$T>0$
be given.
 For all
$\eps>0$, denote by     $(\eta_\eps,   u_\eps)$ a solution of
$(SW_\e)$        with initial data~$(\eta^0,u^0)$. For any couple of  vector
fields~$(\el,\ul) $ belonging to~$ L^\infty([0,T];C^{\infty}(
\R\times\T))$ and rapidly decaying with  respect to $x_1$,
 define  
$$\EE_\eps(t)= \frac12 \int \left((\eta_\eps-\el)^2  +(1+\eps \eta_\eps)|u_\eps-\ul|^2\right)(t,x) dx
+\nu\int_0^t \int
 | \nabla (u_\eps-\ul )|^2(t',x) dxdt'.
$$ 
Then the following stability inequality
holds for all $t \in [0,T]$:
$$
\begin{aligned}
& \EE_\eps(t)\leq C \EE_\eps(0)\exp\left( \chi (t)\right) +\omega_\eps(t)\\
& +  C\int_0^t e^{  \chi (t)-\chi (t')} \int \left(\d_t
\el+\frac1\eps
\nabla
\cdot
\ul+\nabla
\cdot (\el\ul)
\right)(\el-\eta_\eps )(t',x)dxdt'
\\ & +   C\int_0^t  e^{ \chi (t)-\chi (t') }\int  (1+\eps
\eta_\eps)  \left(\d_t
\ul+\frac1\eps (\beta x_1 \ul^\perp   +\nabla \el) +(\ul\cdot \nabla )\ul-
\nu \Delta \ul\right)\cdot
(\ul-u_\eps)     
  (t',x)dxdt' ,
\end{aligned}
$$
where $\omega_\eps(t)$ depends on~$\ul$ and goes to zero with~$\eps$, uniformly in time,  and where
$$
\chi (t) =C\int_0^t \left( \|
\nabla \cdot \ul\|_{L^\infty(\R\times \T)}+\| \nabla  
\ul \|_{L^2(\R\times\T)}^2  \right)(t') dt'.
$$
}
\end{Prop}
Let us postpone the proof of that result, and end the proof of the strong convergence. We apply that proposition
to~$(\el,\ul) = (\eta_{\eps,N},u_{\eps,N})$, where~$(\eta_{\eps,N},u_{\eps,N})$  is the
 approximate solution 
  given by Proposition \ref{approximate-sol}. We will denote by~$\chi_{\eps,N}$ and~$\EE_{\eps,N}$
  the quantities defined in Proposition~\ref{strong-weak}, where~$(\el,\ul)$ has been 
  replaced by~$ (\eta_{\eps,N},u_{\eps,N})$.
  
 Because of the uniform regularity estimates on~$(\eta_{\eps,N},u_{\eps,N})$, we have
$$
\forall T>0,\quad 
\sup_N \lim_{\eps \to 0} \left(\|\nabla u_{\eps,N} \|_{L^2([0,T]; L^2(\R\times\T))}^2+
 \|\nabla \cdot  u_{\eps,N} \|_{ L^1([0,T];L^\infty(\R\times
\T))}\right)\leq C_T ,$$
so we get a uniform bound on $\chi_{\eps,N}$~:
$$\sup_N  \lim_{\eps \to 0} \|\chi_{\eps,N}\|_{L^\infty([0,T])} \leq C_T.$$
Then, from the initial convergence (\ref{init-cv}) we obtain that 
$$
\forall N\in \N, \quad \EE_{\eps,N}(0)\exp\left( \chi_{\eps,N}(t)\right) \to 0 \hbox{ as }
 \eps \to 0\hbox{ in }L^\infty([0,T]).
 $$
 Moreover by Proposition~\ref{strong-weak} we have
\begin{equation}\label{eqepsN}
\d_t
(\eta_{\eps,N},u_{\eps,N})+{1\over \eps }
L(\eta_{\eps,N},u_{\eps,N}) +Q((\eta_{\eps,N},u_{\eps,N}),(\eta_{\eps,N},u_{\eps,N}))-\nu \Delta'(\eta_{\eps,N},u_{\eps,N})
= R_{\e,N}.
\end{equation}
Let us estimate the contribution of the remainder term. We can write 
$$
\int_0^t e^{ \chi_{\eps,N}(t)-\chi_{\eps,N}(t') }
   \int R_{\eps,N}\cdot \left((\eta_{\eps,N}-\eta_\eps), (1+\eps \eta_\eps)
   (u_{\eps,N}-u_\eps)\right) (t',x)dxdt' = I_{\e,N}^{(1)} (t) +  I_{\e,N}^{(2)}  (t), 
$$
with
\begin{eqnarray*}
I_{\e,N}^{(1)} (t) &\eqdefa &\int_0^t e^{ \chi_{\eps,N}(t)-\chi_{\eps,N}(t') }  \int R_{\eps,N,0} 
(\eta_{\eps,N}-\eta_\eps)  (t',x)dxdt' , \quad \mbox{and}
\\
I_{\e,N}^{(2)} (t) &\eqdefa &\int_0^t e^{ \chi_{\eps,N}(t)-\chi_{\eps,N}(t') }  \int R_{\eps,N}' 
(1+\eps \eta_\eps)
   (u_{\eps,N}-u_\eps)(t',x)dxdt' .
\end{eqnarray*}
The first term can be estimated in the following way:
$$
|I_{\e,N}^{(1)} (t)| \leq C_{T} \| R_{\eps,N}\|_{L^{1}([0,T];L^{2}(\R\times\T))}
 \|\eta_{\eps,N}-\eta_\eps
\|_{L^{\infty}([0,T];L^{2}(\R\times\T))}.
$$
For the second term we can write
$$
|I_{\e,N}^{(2)} (t)|   \leq  C_T \|\sqrt{1 +\e \eta_{\e}} (u_{\eps,N}-u_\eps)\|_{
L^{\infty}([0,T];L^{2}(\R\times \T))} \|\sqrt{1 +\e \eta_{\e}} R_{\eps,N}\|_{L^{1}([0,T];
L^{2}(\R\times\T))}.
$$
Now we can write
$$
 \|\sqrt{1 +\e \eta_{\e}} R_{\eps,N}\|_{ L^{2}(\R\times\T)}^{2} \leq C (
 \|R_{\eps,N}\|_{ L^{2}(\R\times\T)
 }^{2} +
 \e \|\eta_{\e}\|_{ L^{2}(\R\times\T)} \|R_{\eps,N}\|_{ L^{4}(\R\times\T)}^{2}) .
$$
 Since
 $$
 \e \|R_{\eps,N}\|_{ L^{4}(\R\times\T)}^{2} \leq \e\|R_{\eps,N}\|_{ L^{\infty}(\R\times\T)} 
   \|R_{\eps,N}\|_{ L^{2}(\R\times\T)}  ,
 $$
we infer that the quantity~$\e^{\frac12} R_{\eps,N}$ goes to zero as~$\e$ goes to zero and~$N$ goes to infinity,  in the 
space~$L^{2}([0,T];L^{4}(\R\times\T))$,
so in particular
$$
 \lim_{N\to \infty} \lim_{\e \to 0} \e^{\frac12} \|R_{\eps,N}\|_{ L^{1}([0,T];L^{4}(\R\times\T))} = 0.
$$
Finally by the uniform bound on~$\eta_{\e}$ in~$L^{\infty}([0,T];L^{2}(\R\times\T))$ and by the 
smallness assumptions
on~$R_{\e,N}$, we deduce that  
$$
\longformule{
\int_0^t e^{ \chi_{\eps,N}(t)-\chi_{\eps,N}(t') }
   \int R_{\eps,N}\cdot \left((\eta_{\eps,N}-\eta_\eps), (1+\eps \eta_\eps)
   (u_{\eps,N}-u_\eps)\right) (t',x)dxdt'}{
   \leq \frac12 ( \|\eta_{\eps,N}-\eta_\eps
\|^{2}_{L^{\infty}([0,T];L^{2} )} + \|\sqrt{1 +\e \eta_{\e}} (u_{\eps,N}-u_\eps)\|^{2}_{
L^{\infty}([0,T];L^{2} )} ) + \omega_{\e,N}(t) ,
   }
$$
where
$$
 \lim_{N\to \infty} \lim_{\eps \to 0} \|\omega_{\e,N}(t)\|_{L^{\infty}([0,T])} = 0.
$$ 
We now recall that by Proposition~\ref{strong-weak}, using~(\ref{eqepsN}), we have
$$
\longformule{
\EE_{\eps,N}(t) \leq C\EE_{\eps,N}(0)\exp\left( \chi_{\eps,N}(t)\right) + \omega_{\e}(t) }{{}+ C\int_0^t e^{ \chi_{\eps,N}(t)-\chi_{\eps,N}(t') }
   \int R_{\eps,N}\cdot \left((\eta_{\eps,N}-\eta_\eps), (1+\eps \eta_\eps)
   (u_{\eps,N}-u_\eps)\right) (t',x)dxdt' }
$$
where \label{strong-weakN}
$$
\EE_{\eps,N}(t) = \frac12\left( \| (\eta_\eps-
  \eta_{\eps,N})(t)\|_{L^{2}}^2  +\|\sqrt{1+\eps  \eta_\eps}  (u_\eps-
u_{\eps,N})(t)\|_{L^{2}}^2 \right)
+\nu\int_0^t \| \nabla (u_\eps-u_{\eps,N} )(t')\|_{L^2}^2(t' )  dt'.
$$
Putting together the previous results we get that~$\displaystyle
\lim_{N\to \infty} \lim_{\eps  \to 0} \EE_{\eps,N}(t) = 0$
uniformly on~$[0,T]$, hence that
 $$
 \begin{aligned}
 \lim_{N\to \infty} \lim_{\eps \to 0} \|
 \eta_{\eps,N}-\eta_\eps\|_{L^\infty([0,T];L^2( \R\times\T))} =0,\\
 \lim_{N\to \infty} \lim_{\eps \to 0} \| \sqrt{1+\eps \eta_\eps}(u_{\eps,N}-u_\eps)\|_{
 L^\infty([0,T];L^2(\R\times\T))} =0,\\
\lim_{N\to \infty} \lim_{\eps \to 0}
 \|u_{\eps,N}-u_\eps\|_{L^2([0,T];\dot H^1(\R\times\T))} =0.
 \end{aligned}
 $$
By interpolation we therefore find that
 $$ 
 \lim_{N\to \infty} \lim_{\eps \to 0}\left( \| \eta_{\eps,N}-\eta_\eps
 \|_{L^\infty([0,T];L^2(\R\times\T))} +\|u_{\eps,N}-u_\eps\|_{L^2([0,T]; H^1
 (\R\times\T))} \right)=0, 
 $$ 
 hence~(\ref{limfinale}) is proved.  
 \medskip
 
 To conclude the proof of Theorem~\ref{strong-convergencegenericbeta}  it remains
 to prove
 Proposition~\ref{strong-weak}.
 As the energy is a Lyapunov functional for $(SW_\e)$, we have
$$\begin{aligned}
\EE_\eps(t)-\EE_\eps(0) &\leq \int_0^t {d\over dt} \int\left( (\frac12 \el^2-\el \eta_\eps)
+(1+\eps \eta_\eps)(\frac12 |\ul|^2-\ul\cdot u_\eps)\right)(t',x) dxdt'\\ 
& \quad  \quad \quad\quad \quad \quad\quad +\int_0^t
\int
\nu   (\nabla
\ul-2\nabla u_\eps)\cdot \nabla \ul(t',x) dxdt'\\
&\leq \int_0^t \int \left(\d_t \el (\el-\eta_\eps )  +    (1+\eps
  \eta_\eps)  \d_t \ul\cdot  (\ul-u_\eps)     
\right)  (t',x)dxdt'
\\ &\quad    \quad\quad\quad \quad\quad \quad- \int_0^t \int
\left(\d_t \eta_\eps \el              + \d_t ((1+\eps
\eta_\eps )u_\eps) \cdot \ul -\frac{\eps}2 \d _t
\eta_\eps  |\ul|^2)\right)(t',x)dxdt'\\
&\quad  \quad  \quad\quad \quad\quad \quad-\int_0^t
\int
\nu  
\left(\Delta\ul \cdot(\ul-u_\eps)-\Delta u_\eps
\cdot
\ul\right)(t',x) dxdt' .
\end{aligned}$$
Using the conservation of mass and  of momentum 
we get
$$\begin{aligned}
\EE_\eps(t)-\EE_\eps(0) &\leq \int_0^t \int \left(\d_t \el
(\el-\eta_\eps ) +    (1+\eps
\eta_\eps)  (\d_t
\ul-\nu \Delta \ul)\cdot (\ul-u_\eps)     
\right)  (t',x)dxdt' \\
 &+ \int_0^t \int   \frac1\eps
\nabla\cdot  \left((1+\eps
\eta_\eps) u_\eps\right)\left(\el -\frac\eps2   |\ul|^2 \right)(t',x)dxdt'
\\ &+\int_0^t \int           
\left({ (1+\eps \eta_\eps )\over \eps} (\beta x_1         u_\eps^\perp +\nabla \eta_\eps)  + \nabla \cdot
((1+\eps \eta_\eps) u_\eps \otimes u_\eps)\right)
\cdot
\ul (t',x)dxdt'\\
&+\int_0^t \int \eps\nu  \eta_\eps \Delta \ul \cdot (\ul-u_\eps)(t',x) dxdt'.
\end{aligned}$$
Integrating by parts leads then to
\begin{equation}
\label{stab1}
\begin{aligned}
\EE_\eps(t)-\EE_\eps(0)
&\leq \int_0^t \int \left(\d_t \el+\frac1\eps \nabla \cdot \ul+\nabla \cdot (\el\ul)
\right)(\el-\eta_\eps )(t',x)dxdt'
\\ &\!\!+ \!\!  \int_0^t \! \int  (1+\eps \eta_\eps)  \left(\d_t
\ul+\frac1\eps (\beta x_1 \ul^\perp   +\nabla \el) +(\ul\cdot \nabla )
\ul-\nu \Delta \ul\right)\cdot
(\ul-u_\eps)     
  (t',x)dxdt' \\
&-\int_0^t \int         (1+\eps \eta_\eps)D\ul : (\ul
-u_\eps)^{\otimes 2}(t',x)dxdt'\\&-\int_0^t \int\left(\frac12\eta_\eps^2
\nabla
\cdot \ul+(\el-\eta_\eps )
\nabla \cdot (\el \ul)+\eta_\eps \ul\cdot \nabla \el\right) (t',x)dxdt'+R_\eps,
\end{aligned}
\end{equation}
where 
$$
R_\eps(t)=\int_0^t \int \eps\nu 
\eta_\eps
\Delta
\ul
\cdot (\ul-u_\eps)(t',x) dxdt'.
$$
The last term is rewritten in a convenient form by integrating by parts
\begin{equation}
\label{nonlinear-eta}
\begin{aligned}
-\int_0^t \int&\left(\frac12\eta_\eps^2
\nabla
\cdot \ul+(\el-\eta_\eps )\nabla \cdot (\el \ul)+\eta_\eps \ul\cdot \nabla \el\right) (t',x)dxdt'\\
&=-\int_0^t \int\left(\frac12\eta_\eps^2
\nabla
\cdot \ul+(\el-\eta_\eps )(\ul \cdot \nabla \el+\el \nabla \cdot \ul)+\eta_\eps \ul\cdot \nabla \el\right)
(t',x)dxdt'\\
&=-\int_0^t \int\left(\frac12\eta_\eps^2
\nabla
\cdot \ul+(\el-\eta_\eps )\el \nabla \cdot \ul+\frac12 \ul\cdot \nabla \el^2\right) (t',x)dxdt'\\
&=-\int_0^t \int\frac12(\eta_\eps-\el)^2
\nabla
\cdot \ul (t',x)dxdt' .
\end{aligned}
\end{equation}
 Plugging (\ref{nonlinear-eta})   into (\ref{stab1}) leads to
\begin{equation}
\label{stab2}
\begin{aligned}
&\EE_\eps(t)-\EE_\eps(0)
\leq \int_0^t \int \left(\d_t \el+\frac1\eps \nabla \cdot \ul+\nabla \cdot (\el\ul)
\right)(\el-\eta_\eps )(t',x)dxdt'
\\ &+  \int_0^t \int  (1+\eps \eta_\eps)  \left(\d_t
\ul+\frac1\eps (\beta x_1 \ul^\perp   +\nabla \el) +(\ul\cdot \nabla )\ul-\nu \Delta \ul\right)\cdot
(\ul-u_\eps)     
  (t',x)dxdt' \\
&-\int_0^t \int         (1+\eps \eta_\eps) D \ul:(\ul-u_\eps)^{\otimes 2} (t',x)dxdt'
-\int_0^t \int\frac12(\eta_\eps-\el)^2
\nabla
\cdot \ul (t',x)dxdt'+R_\eps(t).
\end{aligned}
\end{equation}

\bigskip
In order to get an inequality of Gronwall type, one has to control the
right hand side in terms of $\EE_\eps$. We start by estimating the flux
term. We have
$$\begin{aligned}
-\int_0^t &\int         (1+\eps \eta_\eps) D \ul:(\ul-u_\eps)^{\otimes
2} (t',x)dxdt' \\
&\leq \int_0^t \left(  \|
\nabla \ul\|_{L^2(\R\times\T)}+\eps 
\|\eta_\eps\|_{L^2(\R\times\T)}\|
\nabla \ul\|_{L^\infty(\R\times\T)}\right) \| \ul-u_\eps\|_{L^4(\R\times\T)}^2(t') dt'\\
&\leq C\int_0^t \left(  \|
\nabla \ul\|_{L^2(\R\times\T)}+\eps 
\|\eta_\eps\|_{L^2(\R\times\T)}\|
\nabla \ul\|_{L^\infty(\R\times\T)}\right) \| \ul-u_\eps\|_{L^2(\R\times\T)} \\
& \quad \quad \quad \quad \quad \quad \quad \quad 
\quad \quad \quad \quad \quad \quad \quad \quad \quad \quad 
\quad \quad \quad \quad \quad \quad \quad \times \| \ul-u_\eps\|_{\dot
H^1(\R\times\T)}(t') dt'
\end{aligned}
$$
and 
$$
\| \ul-u_\eps\|_{L^2(\R\times\T)}^2 \leq  \|\sqrt{1+\eps \eta_\eps} (u_\eps -
\ul)\|^2_{L^2(\R\times\T)}
+ \eps
\|\eta_\eps\|_{L^2(\R\times\T)} \| \ul-u_\eps\|_{L^2(\R\times\T)}
\| \ul-u_\eps\|_{\dot H^1(\R\times\T)}
$$
which implies
$$
\| \ul-u_\eps\|_{L^2(\R\times\T)}^2 \leq  2 \|\sqrt{1+\eps \eta_\eps} (u_\eps -\ul)\|^2_{
L^2(\R\times\T)}
 + 16\eps^2
\|\eta_\eps\|^{2}_{L^2(\R\times\T)} \| \ul-u_\eps\|^2_{\dot H^1(\R\times\T)}.
$$
 Therefore, using the uniform bounds on~$\eta_{\e}$, $\sqrt{1+\e \eta_{\e}} u_{\e}$ and on~$u_{\e}$ given
 by the energy estimate, we gather that
\begin{equation}
\label{rightside1}
\begin{aligned}
-\int_0^t &\int         (1+\eps \eta_\eps) D \ul:(\ul-u_\eps)^{\otimes
2} (t',x)dxdt' \\
&\leq C\int_0^t ( \|\nabla \ul\|_{L^2 } + \e \|\nabla \ul\|_{L^\infty } )
\| \sqrt{1+\eps \eta_\eps}(u_\eps -\ul)\|_{L^2 } \|\ul-u_\eps\|_{\dot H^1 }(t') dt'\\
& \quad \quad \quad \quad \quad \quad \quad \quad \quad 
+C \e \int_0^t  ( \|\nabla \ul\|_{L^2 } + \e \|\nabla \ul\|_{L^\infty } )
  \|\ul-u_\eps\|^{2}_{\dot H^1 }(t') dt'
\\
&\leq \frac\nu 4 \int \|\ul-u_\eps\|_{\dot H^1}^2(t') dt'+ \frac{C}\nu \int \|
\nabla \ul\|_{L^2 }^2   \|\sqrt{1+\eps \eta_\eps}(u_\eps -\ul)
\|_{L^2 }^2(t') dt'   + \omega_{\e}(t).
\end{aligned}
\end{equation}
We also have
\begin{equation}
\label{rightside2}
-\int_0^t \int\frac12(\eta_\eps-\el)^2
\nabla
\cdot \ul (t',x)dxdt' \leq \frac12\int_0^t \| \nabla \cdot \ul\|_{L^\infty(\R\times\T)}
\|\el-\eta_\eps\|_{L^2(\R\times\T)}^2 (t') dt' , 
\end{equation}
so we are left with the study  of the remainder $R_\eps$. We have
$$
R_{\e}(t) \leq \e \nu \|\eta_\eps\|_{L^\infty(\R^+;L^2(\R\times\T))}  \int_0^t
\|\Delta \ul\|_{L^4(\R\times\T)} \|\ul-u_\eps\|_{ L^4(\R\times\T)} 
 (t') dt'.
$$
The above estimate on~$\|\ul-u_\eps\|_{ L^2(\R\times\T )}$ implies in particular that~$\|\ul-u_\eps\|_{ 
L^2(\R\times\T )}$
is bounded in~$L^{2}([0,T])$, hence we get that~$ \|\ul-u_\eps\|_{ L^4(\R\times\T)} $ is also bounded 
in~$L^{2}([0,T])$. So we infer directly that~$R_{\e}(t)$ goes to zero in~$L^{\infty}([0,T])$ as~$\e$ goes to zero. That 
result, joint with~(\ref{rightside1}) and~(\ref{rightside2})  allows to deduce from~(\ref{stab2}) the following estimate:
$$
\begin{aligned}
&\frac12\EE_\eps(t)-\EE_\eps(0)
\leq \int_0^t \int \left(\d_t \el+\frac1\eps \nabla \cdot \ul+\nabla \cdot (\el\ul)
\right)(\el-\eta_\eps )(t',x)dxdt' 
\\ & +  \int_0^t \int  (1+\eps \eta_\eps)  \left(\d_t
\ul+\frac1\eps (\beta x_1 \ul^\perp   +\nabla \el) +(\ul\cdot \nabla )\ul-\nu \Delta \ul\right)\cdot
(\ul-u_\eps)     
  (t',x)dxdt'\\
& +\frac{C}\nu \int \|\nabla \ul\|_{L^2 }^2   \| \sqrt{1+\eps \eta_\eps}(u_\eps -\ul)
\|_{L^2(\R\times \T)}^2(t') dt' +\frac12\int_0^t \| \nabla \cdot \ul\|_{L^\infty }
\|\el-\eta_\eps\|_{L^2 }^2 (t') dt'
+\omega_{\e}(t)
\end{aligned}
$$
thus applying Gronwall's lemma provides the expected stability inequality.  
 Theorem~\ref{strong-convergencegenericbeta} is proved.
\end{proof}

\section{Strong convergence of filtered weak solutions towards a weak solution}
\label{hybrid}
\setcounter{equation}{0}
The aim of this section is to prove an intermediate convergence result, in the sense that we will seek   a strong
convergence result of  the filtered weak solutions, towards a weak solution of the limit system;  thus no 
additional smoothness will be required on the initial data other than~$L^{2}(\R \times \T)$. 
As explained in the introduction
of this chapter, the lack of compactness in the spatial variables of~$\eta_{\e}$ will prevent us from obtaining
at the limit the expected system~$(SW_{0})$: a defect measure remains at the limit, which we are unable
to remove. In order to gain some space compactness and to get rid of that defect measure,
 we propose in the final paragraph of this section 
(Paragraph~\ref{capillarity} below) an alternate model which takes into account capillarity effects, and for 
which one can prove the strong convergence of filtered solutions towards a weak solution of~$(SW_0)$.

  The first result of this paragraph is the following.  
\begin{Thm}[strong convergence towards weak solutions]\label{thmstrongweakcv}
{ 
Let $(\eta^0,u^0) \in L^2(\R\times \T)$ and $(\eta^0_\eps,u_\eps^0)$ satisfy~$(\ref{firstcondition})$. 
 For all
$\eps >0$, denote by~$(\eta_\eps,u_\eps)$ a solution of~$(SW_\e)$ with initial
data~$(\eta_\eps^0,u_\eps^0)$, and by\label{Phieps} $$\Phi_\eps=\LL\left(-{t\over  \eps}\right) 
(\eta_\eps,u_\eps).$$ 
Up  to the extraction of a subsequence, $\Phi_\eps$
converges strongly  in~$L^2_{loc}(\R^+;H^s_{loc}(\R\times \T))$ (for  all~$s<0$) towards some 
  solution~$\Phi$ \label{Phi}
of the following limiting filtered system: for all~$i\lambda \in {\mathfrak S} $, there 
is a bounded measure~$\upsilon_{\lambda} \in {\mathcal M}(\R^{+}\times\R \times \T )$ (which vanishes if~$\lambda = 0$),
 such that  for
all smooth~$\Phi^*_\lambda \in \Ker (L-i\lambda
Id)$,
$$
\longformule{
\int \Phi \cdot \bar\Phi^*_\lambda (x) \: dx -\nu \int_{0}^{t} \int \Delta'_{L} \Phi  \cdot    \bar\Phi^*_\lambda  
(t',x) \: dxdt'}{
+\int_{0}^{t}  \int Q_{L} (\Phi,\Phi) \cdot  \bar\Phi^*_\lambda(t',x) \: dxdt' + \int_{0}^{t}  \int \nabla
\cdot  ( \bar\Phi^*_\lambda )' \upsilon_{\lambda} (dt' dx) = \int \Phi^{0} \cdot\bar \Phi^*_\lambda (x) \: dx
,}
$$
where~$Q_{L} $ and~Ê$\Delta'_{L} $ are defined by~(\ref{Q-lap-L}) page~\pageref{Q-lap-L}, and where~$\Phi^0 = (\eta^0,u^0)$.}
\end{Thm}

\begin{Rem}
{  $\bullet$ Note that, by interpolation with the uniform $L^2_{loc}(\R^+;  H^1(\R \times \T))$ bound 
on $u_\eps$, we get the strong
convergence of $u_\eps$ in $L^2_{loc}(\R^+; L^2(\R\times \T))$~: up to  extraction of a 
subsequence,
$$\left\| u_\eps - \left( \LL\left({t\over \eps } \right)\Phi  
\right)'\right\|_{L^2(\R\times \T)} \to 0 \hbox{ in }
L^2_{loc}(\R^{+}).$$

$\bullet$ As explained above, the presence of the defect measure~$\upsilon_{\lambda}$ at the limit
is due to a possible defect of compactness in space of the sequence~$(\eta_{\e})_{\e >0}$. As the proof of the theorem
will show, that measure is zero if one is able to prove
some equicontinuity in space on~$\eta_{\e}$, or even on~$\e \eta_\e$.  Since we have been unable to prove such a result, we study in the final
paragraph of this section a slightly different model, where   capillarity effects are added   in order to gain
that compactness. Note that the model introduced in Paragraph~\ref{capillarity}   is unfortunately 
not very physical due to the particular form of the capillarity operator
 (see its definition in~(\ref{defcapillarity}) below). 
}
\end{Rem}
 Theorem~\ref{thmstrongweakcv} is proved in Sections~\ref{strongcompactpilambdaphieps} to~\ref{takinglimits}, and 
 the result in the presence of capillarity is stated and proved in Section~\ref{capillarity} .

\subsection{Strong compactness of $\Pi_\lambda \Phi_\eps$}\label{strongcompactpilambdaphieps}
Let us prove the following lemma.

\begin{Lem}\label{strongcomp-pil}
{ 
With the notation of Theorem~\ref{thmstrongweakcv}, the following results hold.

$\bullet$  For all $i\lambda \in {\mathfrak S}\setminus\{0\}$, $\Pi_ \lambda \Phi_\eps$ is 
strongly compact in $L^2([0,T];
H^s(\R\times \T))$ for all $T>0$ and all $s \in \R$;

$\bullet$   $\Pi_0 \Phi_\eps$ is strongly  compact in 
$L^2([0,T];
H^s_{loc}(\R\times
\T))$ for all $T>0$ and all $s<0$.}
\end{Lem}

 \begin{proof} 

$\bullet$
For all $\lambda\neq 0$, we recall that by Proposition~\ref{diag-prop} page~\pageref{diag-prop}, 
 the eigenspace of $L$ associated with the  eigenvalue~$i\lambda$ is a 
  finite dimensional subspace of~$H^{\infty} (\R\times \T)$. Therefore
the only point to be checked is the compactness with respect to time,  which is obtained 
as follows.

Let~$(n,k,j) \in \N \times \Z \times \{-1,0,1\}$ be given, such that~$\lambda = \tau(n,k,j) \neq 0$, and
let~$\Psi_{n,k,j}$ be the corresponding eigenvector. 
Multiplying the system~$(SW_\e)$ by~$\Psi_{n,k,j}$  (which is smooth 
and rapidly decaying as~$|x_1|$ goes to infinity) and integrating with respect to $x$ leads to
$$\begin{aligned}
\d_t \int \left(\eta_\eps (\bar\Psi_{n,k,j })_{0}+m_\eps \cdot \bar\Psi_{n,k,j}'\right) (t,x) \: dx 
+{i\tau(n,k,j)\over \eps} \int
(\eta_\eps
  ( \bar\Psi_{n,k,j })_{0}+m_\eps \cdot \bar \Psi_{n,k,j}')(t,x) \: dx\\
+\nu \int \nabla  u_\eps : \nabla \bar \Psi_{n,k,j}'(t,x) \: dx
-\int m_\eps\cdot (u_\eps \cdot \nabla)  \bar \Psi_{n,k,j}'(t,x) \: dx-\frac12 \int \eta_{\e}^{2}  \nabla \cdot\bar \Psi_{n,k,j}' \: dx
=0
\end{aligned}
$$
where $\bar \Psi_{n,k,j}$ denotes the complex conjugate of $\Psi_{n,k,j}$,
or equivalently
\begin{equation}
\label{pil-eq}
\begin{aligned}
& \d_t \left( \exp\left({it\tau(n,k,j) \over \eps}\right) \int (\eta_ \eps 
 (\bar \Psi_{n,k,j })_{0}+m_\eps \cdot
\bar  \Psi_{n,k,j}')(t,x) \: dx
\right) \\
&
+\nu \int \nabla  \left(  \exp\left({it\tau(n,k,j) \over \eps}\right)  u_\eps \right): 
\nabla \bar  \Psi_{n,k,j}'(t,x) \: dx\\
&-\int  \exp\left({it\tau(n,k,j) \over \eps}\right) 
\Bigl( m_\eps\cdot (u_ \eps \cdot \nabla) \bar \Psi_{n,k,j}' + \frac12\eta_\eps^{2} \nabla \cdot
 \bar \Psi_{n,k,j}'\Bigr) (t,x) \: dx
=0.
\end{aligned}
\end{equation}

By the uniform estimates coming from the energy inequality we then  deduce that
$$
\d_t \left( \exp\left({it\tau(n,k,j) \over \eps}\right) \int (\eta_ \eps
 ( \bar \Psi_{n,k,j })_{0}+m_\eps \cdot
\bar  \Psi_{n,k,j}')(t,x)dx\right) \hbox{ is uniformly bounded in }\eps.
$$
Therefore the family
$$
\left(\exp\left({it\lambda \over \eps}\right) \Pi_\lambda (\eta_\eps,m_ \eps)  \right)_{\e > 0} \hbox{ is compact in
}L^2([0,T];H^{s}(\R\times \T)) \: \mbox{for any} \: s \in \R , 
$$
and  since  $\eps \eta_\eps u_\eps$ converges to 0 in $L^2(\R^ +; H^s(\R\times 
\T))$ for all $s<0$, we deduce that
$$
\exp\left({it\lambda \over \eps}\right) \Pi_\lambda (\eta_\eps,u_ \eps)=\Pi_\lambda 
\Phi_\eps \hbox{ is compact in
}L^2([0,T];H^{s}(\R\times \T)) \:\mbox{for any} \: s \in \R.
$$

$\bullet$
For $\Pi_0 \Phi_\eps=\Pi_0  (\eta_\eps,u_\eps)$ the study is a little more  difficult since the 
compactness with respect to
spatial variables has to be taken into account. By the energy  estimate we have the 
uniform bound
$$
\Phi_\eps \: \hbox{is uniformly bounded in }L^2_{loc}(\R^+, L^2(\R\times \T)).
$$

We recall that we have defined in Section~\ref{sectionsuitable} (Definition~\ref{HsL} 
page~\pageref{HsL}) the space\label{redefhsl}
$$
H^s_L=\left\{ \psi\in L^{2}(\R\times \T)\,/\, \sum_ {n,k,j\in S}(1+n+k^2)^s
\left(\psi|\Psi_{n,k,j}\right)_ {L^2(\R\times \T)}^2<+\infty\right\},
$$
where~$S = \N \times \Z \times \{-1,0,1\}$
or equivalently (see Proposition~\ref{equivnorm} page~\pageref{equivnorm})
$$
H^s_L=\left\{ \psi\in L^{2}(\R\times \T)\,/\, (\mbox{Id} \: -\Delta+ \beta^2x_1^2)^{  s/2 }
\psi\in L^2(\R\times \T)\right\}.
$$

As $(\Psi_{n,0,0})_{n \in \N}$ is a Hilbertian basis of $\Ker L$, we    have for all
$T>0$ and all
$s<0$
$$ \left\| \sum_{n\leq N} \left( \Phi_\eps|\Psi_{n,0,0}\right)_{L^2(\R \times\T)} 
\Psi_{n,0,0}-\Pi_0\Phi_\eps
\right\|_{L^2([0,T]; H^s_L)}\to 0 \hbox{ as }N\to \infty  \hbox{ uniformly 
in }\eps.$$
Let $\Omega$ be any relatively compact open subset of $\R\times \T$.   Proposition~\ref{equivnorm}
page~\pageref{equivnorm} implies
that, for all~$s \geq 0$
$$H^s_0(\Omega) \subset H^s_L\subset H^s(\R\times \T),$$
and conversely for $s\leq 0 $,
\begin{equation}\label{embeddings<0}
H^s(\R\times \T) \subset H^s_L \subset H^s(\Omega).
\end{equation}
Thus  for all~$s<0$ and  all $T>0$, we have
$$ \left\| \sum_{n\leq N} \left( \Phi_\eps|\Psi_{n,0,0}\right)_{L^2(\R \times\T)} 
\Psi_{n,0,0}-\Pi_0\Phi_\eps
\right\|_{L^2([0,T]; H^s (\Omega))}\to 0 \hbox{ as }N\to \infty \hbox { uniformly in 
}\eps.$$

Moreover the same computation as previously shows that for any~$n \in \N$, 
\begin{equation}
\label{pi0-eq}
\begin{aligned}
\d_t  \left(\int (\eta_\eps \bar \eta_{n,0,0}+m_\eps \cdot
\bar u_{n,0,0})(t,x)dx
\right)+\nu \int \nabla    u_\eps : \nabla \bar u_{n,0,0}(t,x)dx\\
-\int   m_\eps\cdot (u_\eps \cdot \nabla)  \bar u_{n,0,0}(t,x)dx=0, 
\end{aligned}
\end{equation}
and, since  $\eps \eta_\eps u_\eps$ converges to 0 in $L^2 (\R^ +;
H^s(\R\times \T))$ for  any~$s < 0$ we get
$$
\sum_{n\leq N} \Pi_{n,0,0} (\eta_\eps,u_\eps) \hbox{ is
compact in }L^2([0,T]\times\R\times \T).
$$

Combining both results shows finally that
$$
\Pi_0 \Phi_\eps \hbox{ is compact in } L^2([0,T]; H^s_{loc}(\R\times \T))
$$
for all $T>0$ and all $s<0$.
 Lemma~\ref{strongcomp-pil} is proved.
\end{proof}
\bigskip\bigskip
As  the spectrum of~$L$, $\mathfrak S$ is countable (see Chapter~\ref{equatorialwaves}), we
are therefore  able to construct (by  diagonal 
extraction) a subsquence of
$\Phi_\eps$, and  some \label{Philambda}$\Phi_\lambda \in \Ker(L-i\lambda Id)$ such  that for all $s<0$ 
and all $T>0$
\begin{equation}\label{limpilambdaphieps}
\forall i\lambda \in {\mathfrak S} ,\quad \Pi_\lambda \Phi_\eps \to  \Phi_\lambda\hbox{ 
in } L^2([0,T];H^s_{loc}(\R\times \T)).
\end{equation}
Note that the $\Phi_\lambda$ defined as the strong limit  of~$\Pi_ \lambda \Phi_\eps$ can 
also be obtained as the
weak limit of~$\exp\left({it\lambda \over \eps}\right)(\eta_\eps, u_ \eps)$. We have 
indeed the following lemma.
\begin{Lem}\label{weak-charlambda}
{ 
With the notation of Theorem \ref{thmstrongweakcv}, consider a  subsequence of
$(\Phi_\eps)_{\e > 0}$, and  some~$\Phi_\lambda $ in~$ \Ker(L-i\lambda Id)$ such  that for all $s<0$ 
and all $T>0$
$$\forall i\lambda \in {\mathfrak S} ,\quad \Pi_\lambda \Phi_\eps \to  \Phi_\lambda\hbox{ 
in } L^2([0,T];H^s_{loc}(\R\times \T)).$$

Then, for all $i\lambda \in {\mathfrak S}$, $\displaystyle e^{ {it\lambda \over  \eps}   }
(\eta_\eps,u_\eps) $ converges
to~$ \Phi_\lambda $  weakly in~$L^2([0,T]\times\R\times \T). $
In particular, for all~$i\lambda \in {\mathfrak S}$, the vector field~$ \Phi_\lambda '$
belongs to~$L^2  ([0,T];H^1(\R\times \T))$. }
\end{Lem}
\begin{proof}
Denote by $(\eta_\lambda, u_\lambda)$ any weak limit point of the  sequence 
$\displaystyle\exp\left({it\lambda \over \eps} \right)
(\eta_\eps,u_\eps) $ (recall that the sequence is bounded in~$L^{2}_{loc}(\R^+\times\R\times \T)$).
Let~$\chi $ and~$\psi$ be any  test function  and vector field in~$ \DD(\R^+\times \R\times \T)$. 
Multiplying
  the conservation of mass in~$(SW_\e)$ by $\eps \chi\exp\left({-
  i  t \lambda \over \eps} 
\right) $
and integrating with respect to all variables leads to
$$
\iint \left( \eps \eta_\eps \exp\left({it\lambda \over \eps} \right)   \left(\d _t 
\chi+{i\lambda \over \eps}\chi\right)
+(1+\eps
\eta_\eps) u_\eps\exp\left({it\lambda \over \eps} \right)
\cdot
\nabla
\chi\right)dxdt =0 .
$$
Because of the bounds coming from the energy estimate~(\ref{uniform}) page~\pageref{uniform}, we can take 
limits in the previous
identity as $\eps$ goes to $0$ to get
$$\iint (u_\lambda \cdot \nabla \chi +i\lambda \eta_\lambda \chi) dxdt=0.$$
Similarly, multiplying the conservation of momentum by $\eps \psi \exp \left({it\lambda 
\over \eps} \right) $ and
integrating with respect to all variables leads to
$$
\begin{aligned}
\iint
\Bigl(
\eps (1+\eps \eta_\eps)u_\e \exp\left({it\lambda \over \eps} \right) \left(\partial_t 
\psi + {it\lambda \over \eps}
\psi\right)+
\eps\exp\left({it\lambda \over \eps} \right)(1+\eps \eta_\eps)u_\e \cdot ( u_\e \cdot 
\nabla) \psi \\ + \beta
x_1\exp\left({it\lambda \over \eps} \right)(1+\eps
\eta_\eps) u_\eps \cdot
\psi^\perp + (1+\frac\eps{2} \eta_\eps) \eta_\e \exp\left({it\lambda  \over \eps} 
\right)\DIV \psi\\
+ \e\nu u_\eps
\exp\left({it\lambda \over \eps} \right) \cdot
\Delta \psi\Bigr) dxdt=0.
\end{aligned}
$$
Once again the bounds coming from the energy estimate~(\ref{uniform})
will enable us to take the limit as $\eps$ goes to $0$, to get
$$ \iint (\eta_\lambda \nabla \cdot \psi+\beta x_1 u_\lambda \cdot
\psi^\perp+i\lambda u_\lambda \psi)dxdt=0.$$

It follows that  $(\eta_\lambda(t),u_\lambda(t))$ belongs to $\Ker (L- i\lambda Id)$ for 
almost all~$t\in \R^+$, and we
conclude by uniqueness of the limit and $L^2$ continuity of $\Pi_ \lambda$ that~$\Phi_\lambda=
(\eta_\lambda, u_\lambda).$
The lemma is proved.   
\end{proof}

\subsection{Strong convergence of $\Phi_\eps$}\label{strongcvphiepssection}
As a corollary of the previous mode by mode convergence results, we  get the following 
convergence for $\Phi_\eps$. 
\begin{Lem}\label{weak-char}
{ 
With the notation of Theorem~\ref{thmstrongweakcv}, the following results hold. 
Consider a subsequence of
$(\Phi_\eps)$, and  some $\Phi_\lambda \in \Ker(L-i\lambda Id)$ such  that  as constructed 
in~(\ref{limpilambdaphieps}), for all $s<0$ 
and all $T>0$
$$\forall i\lambda \in {\mathfrak S} ,\quad \Pi_\lambda \Phi_\eps \to  \Phi_\lambda\hbox{ 
in } L^2([0,T];H^s_{loc}(\R\times \T)).$$

Then,
$$
\begin{aligned}
\Phi_\eps \rightharpoonup \Phi=\sum_{i\lambda \in {\mathfrak S}}\Phi_ \lambda \hbox{ 
weakly in } L^2_{loc}(\R^+;
L^2(\R\times \T)) ,\\
\mbox{and} \quad \Phi_\eps \to \Phi \hbox{ strongly  in }
L^2_{loc}(\R^+; H^s_{loc}(\R\times \T))  \hbox{ for 
all }s<0.
\end{aligned}
$$
Moreover, defining~$K_N$ as in~(\ref{definitionKNenvelope}) page~\pageref{definitionKNenvelope},   we have for any relatively
compact subset~$\Omega$  of~$ \R\times \T$, 
for all~$T>0$  
and for all~$s<0$,
\begin{equation}
\label{approximation-JN}
\|(Id-K_N) \Phi_\eps\|_{L^2([0,T];H^s(\Omega))} + \|(Id-K_N){\mathcal L}  (\frac{t}{\e}  )
  \Phi_\eps\|_{L^2([0,T];H^s(\Omega))}\to 0 \hbox{ as }N\to  \infty , 
\end{equation} 
uniformly in  $\e$.}
\end{Lem}
\begin{proof}
The first convergence statement  comes  directly from the uniform
bound on~$\Phi_\eps$  in the space~$L^2_{loc}(\R^{+};
L^{2}(\R\times \T))$ and
the~$L^2$ continuity of~$\Pi_\lambda$.

In order to establish the strong convergence result, the crucial  argument is to 
approximate (uniformly)~$\Phi_\eps$ by
a finite number of modes, i.e. to prove (\ref{approximation-JN}). The  main idea is the 
same as for the approximation of $\Pi_0 \Phi_\eps$ in Lemma \ref {strongcomp-pil}.
We have for all
$T>0$ and all
$s<0$
$$ \left\| \sum_{n\leq N,|k|\leq N} \left( \Psi_{n,k,j} |\Phi_\eps
\right)_{L^2(\R\times\T)} \Psi_{n,k,j}-\Phi_\eps
\right\|_{L^2([0,T];H^s_L)}\to 0 \hbox{ as }N\to \infty  \hbox{ uniformly 
in }\eps,$$
and similarly  
$$
\left\| \sum_{n\leq N,|k|\leq N} e^{-i \tau (n,k,j)\frac t\e}\left( \Psi_{n,k,j} |\Phi_\eps
\right)_{L^2(\R\times\T)}  \Psi_{n,k,j}- {\mathcal L}(\frac t \e)\Phi_\eps
\right\|_{L^2([0,T];H^s_L)}\to 0 \hbox{ as }N\to \infty,
$$
uniformly in~$\e$. Therefore for all relatively compact subsets~$\Omega $ of~$\R\times \T$, 
the embedding of~$H^s_L 
$ into~$ H^s(\Omega)$ recalled in~(\ref{embeddings<0}) implies that both quantities
$$\displaystyle\sum_{n\leq N,|k|\leq N} \left( \Psi_{n,k,j} |\Phi_\eps
\right)_{L^2(\R\times\T)} \Psi_{n,k,j}-\Phi_\eps
$$ and
$$ \displaystyle\sum_{n\leq N,|k|\leq N} e^{-i \tau (n,k,j)\frac t\e}\left( \Psi_{n,k,j} |\Phi_\eps
\right)_{L^2(\R\times\T)}  \Psi_{n,k,j}- {\mathcal L}(\frac t \e)\Phi_\eps
$$
 converge strongly
towards zero in~$L^2([0,T];H^s(\Omega))$ as~$N$ goes to infinity, uniformly in~$\e$. Finally~(\ref{approximation-JN}) is proved.

The strong convergence is therefore obtained from the following  decomposition:
$$\Phi_\eps-\Phi= (Id-K_N)\Phi_\eps +K_N (\Phi_\eps-\Phi)  -(Id-K_N) \Phi$$
The first term converges to 0 as~$N\to \infty$ uniformly in~$\eps$ in~$L^2_{loc}(\R^+, 
H^s_{loc}(\R\times \T))$ for all~$s<0$ by~(\ref{approximation-JN}). By 
Lemma~\ref{strongcomp-pil}, the  second term (which 
is a finite sum of 
modes) converges to~0 as~$\eps \to 0$ for all fixed~$N$ in~$L^2_{loc}(\R^+;H^s(\R\times \T))$ for 
all~$s<0$. The last term converges to~0 as~$N\to \infty$ in~$L^2_{loc}(\R^ +;
 L^2(\R\times \T))$. Thus taking limits as~$\eps
\to 0$, then as~$N\to \infty$ leads to the expected strong convergence.  
\end{proof}

\subsection{Taking limits in the equation on $\Pi_\lambda \Phi_\eps$}\label{takinglimits}
The next step is then to obtain the evolution equation for each mode $ \Phi_\lambda$, 
taking limits in~(\ref{pil-eq}) and~(\ref{pi0-eq}).  In the following
proposition, we recall that the  first result 
(concerning the geostrophic motion) relies on a compensated
compactness argument,  i.e. on both the 
algebraic structure of the coupling
term and the particular form of the oscillating modes, which implies  that there is no 
contribution of the equatorial
waves to the geostrophic flow. That result was proved in
Section~\ref{weakweakcv} (see also Proposition~\ref{limitgeolinear} page~\pageref{limitgeolinear}). 
Here we will prove the
second part of the statement, concerning the limit ageostrophic motion.

\begin{Prop}\label{takinglimitsprop} { 
With the notation of Theorem \ref{thmstrongweakcv},  there is  a subsequence of
$(\Phi_\eps)$ such that the following result holds. Consider a 
family~$(\Phi_\lambda )_{i\lambda  \in {\mathfrak S}}$ such
 that~$\Phi_\lambda\in \Ker(L-i\lambda Id)$ and such  that for all~$s<0$ 
and all~$T>0$
$$\forall i\lambda \in {\mathfrak S} ,\quad \Pi_\lambda \Phi_\eps \to  \Phi_\lambda\hbox{ 
in } L^2([0,T];H^s_{loc}(\R\times \T)),$$
as constructed in~(\ref{limpilambdaphieps}). 

Then, $\Phi_0=(\eta_{0},u_{0})$ satisfies the geostrophic equation~:
for all~$(\eta^*,u^*) $ belonging  to~$
\Ker L$ and satisfying~$u^*\in H^1(\R\times \T)$,
$$
\int (\eta_{0} \eta^* + u_{0,2} u^*_{2}) (t,x) \: dx + \nu \int_{0}^{t}  \int \nabla u_{0,2} 
\cdot \nabla
u^*_{2} (t',x) \: dx \: dt' = \int (\eta^{0} \eta^* + u_{2}^{0} u^*_ {2}) (x) \: dx .
$$
Moreover for~$\lambda \neq 0$, $\Phi_\lambda=(\Phi_\lambda^0,\Phi_\lambda')$
satisfies the following envelope  equation~: there is a measure~$\upsilon_{\lambda} $ 
in~${\mathcal M}(\R^{+}\times\R \times \T)$,
 such that  for
all smooth~$\Phi^*_\lambda   = 
(   \Phi^*_{\lambda,0},    (\Phi^*_\lambda)') \in \Ker (L-i\lambda
Id)$,
$$
\longformule { \int \Phi_\lambda \cdot    \bar \Phi^*_\lambda  (t,x) \: dx + 
\nu \int_{0}^{t}  \int  \nabla\Phi_\lambda' 
: \nabla     (   \bar \Phi^*_\lambda)' (t',x) \: dx \: dt' 
+ 
 \int_{0}^{t} \int \nabla \cdot   (  \bar \Phi^*_{\lambda})' \upsilon_{\lambda}(dt',dx)  }
 {  + \sum_{i\mu,  i\tilde \mu \in 
{\mathfrak S}\atop
\lambda=\mu+\tilde \mu} 
 \int_{0}^{t}  Q(\Phi_\mu , \Phi_{\tilde \mu}) \cdot     \bar \Phi^*_{\lambda} 
 (t',x) \: dx \: dt'
= \int\Phi^0 \cdot      \bar \Phi^*_\lambda  (x) \: dx 
,}$$
where~$Q$ is defined by  (\ref{Q-lap}) page~\pageref{Q-lap}.
}
\end{Prop}
\begin{proof}
Let us first recall that for $\lambda\neq 0$,  $\Ker(L-i\lambda Id)$ is  constituted of 
smooth, rapidly decaying vector
fields, so that it makes sense to apply $\Pi_\lambda$ to any  distribution.

Starting from  (\ref{pil-eq}) we get  that for all smooth~$    \Phi^*_\lambda  = 
(   \Phi^*_{\lambda,0},    (\Phi^*_\lambda)')\in \Ker (L-i\lambda
Id)$
\begin{equation}
\label{pil-eqsum}
\begin{aligned}
\int \exp\left({it\lambda \over \eps}\right) 
 (\eta_\eps \bar {    \Phi}^*_{\lambda,0}
 +m_\eps \cdot
 (\bar{    \Phi}^*_\lambda)')(t,x)dx-\int (\eta_\eps^0 \bar {    \Phi}^*_{\lambda,0}+m_ \eps^0 \cdot
 (\bar{    \Phi}^*_\lambda)')(x)dx\\
+\nu \int_0^t \int \nabla    \left(\exp\left({it'\lambda \over \eps} \right) 
u_\eps\right): \nabla  
 (\bar{    \Phi}^*_\lambda)'(t',x)dxdt'\\
-\int _0^t\int \exp\left({it'\lambda \over \eps} \right) \Bigl(
m_\eps\cdot (u_\eps \cdot \nabla)  
 (\bar{    \Phi}^*_\lambda)' +  \frac12 
\eta_\eps^{2}\nabla \cdot 
 (\bar{    \Phi}^*_\lambda)' \Bigr)(t',x)dxdt'=0.
\end{aligned}
\end{equation}
Taking limits as $\eps \to 0$ in the three first terms is immediate
using Lemma~\ref{weak-charlambda} and the assumption on the initial data.
The limit as~$\eps \to 0$ in the two nonlinear terms is given in the following proposition.

\begin{Prop}\label{limittwononlinearterms}
With the notation of Proposition~\ref{takinglimitsprop}, we have
$$
\int_0^t\int \exp\left({it'\lambda \over \eps}\right) m_\eps\cdot (u_ \eps \cdot \nabla)  
 (\bar{    \Phi}^*_\lambda)'(t',x)dxdt'
 \to \int_0^t\int \sum_{\mu+\tilde \mu=\lambda \atop   i\mu, 
i\tilde \mu \in
{\mathfrak S}}\Phi_\mu'\cdot (\Phi_{\tilde
\mu}'
\cdot
\nabla)
 (\bar{    \Phi}^*_\lambda)'(t',x)dxdt',
$$
and
$$
\longformule{
\frac12\int _0^t\int \exp\left({it'\lambda \over \eps} \right) 
\eta_\eps^{2}\nabla \cdot 
 (\bar{    \Phi}^*_\lambda)'(t',x)dxdt' \to \frac12
\int_0^t\int
 \sum_{\mu+\tilde \mu=\lambda \atop   i\mu,  i\tilde \mu \in {\mathfrak S}}
 \Phi_{\mu,0} \Phi_{\tilde
\mu,0}
\nabla \cdot   (\bar{    \Phi}^*_\lambda)'  (t',x)dxdt' }{- \int _0^t \nabla \cdot (\bar{    \Phi}^*_\lambda)'
 \upsilon_{\lambda}(dt'dx).}
$$
\end{Prop}
Before proving that result, let us conclude the proof of Proposition~\ref{takinglimitsprop}. It remains to check
that
$$\longformule{
-  \int  \sum_{\mu+\tilde \mu=\lambda \atop  i\lambda, i\mu, 
i\tilde \mu \in
{\mathfrak S}}\Bigl(\Phi_\mu'\cdot (\Phi_{\tilde
\mu}'
\cdot
\nabla)
 (\bar{    \Phi}^*_\lambda)'   +
 \Phi_{\mu,0} \Phi_{\tilde
\mu,0} \nabla \cdot   (\bar{    \Phi}^*_\lambda)'\Bigr)  \: dx}{= \int 
\sum_{\mu+\tilde \mu=\lambda \atop   i\mu, 
i\tilde \mu \in
{\mathfrak S}}\Bigl(
(\Phi_\mu'\cdot \nabla )\Phi_{\tilde \mu}' \cdot  (\bar{    \Phi}^*_\lambda)'   + \nabla \cdot
(\Phi_{\mu,0} \Phi_{\tilde \mu}' )  \bar{    \Phi}^*_{\lambda ,0}
\Bigr) \: dx,}
$$
since
$$
\int 
\sum_{\mu+\tilde \mu=\lambda \atop  i\mu, 
i\tilde \mu \in
{\mathfrak S}}\Bigl(
(\Phi_\mu'\cdot \nabla )\Phi_{\tilde \mu}' \cdot  (\bar{    \Phi}^*_\lambda)'   + \nabla \cdot
(\Phi_{\mu,0} \Phi_{\tilde \mu}' )  \bar{    \Phi}^*_{\lambda ,0}
\Bigr) \: dx = \int 
\sum_{\mu+\tilde \mu=\lambda \atop   i\mu, 
i\tilde \mu \in
{\mathfrak S}} Q(\Phi_\mu,\Phi_{\tilde \mu}) \cdot  \bar{    \Phi}^*_{\lambda }\: dx .
$$
Clearly one has
$$
-  \int \Phi_\mu'\cdot (\Phi_{\tilde \mu}' \cdot \nabla)  (\bar{    \Phi}^*_\lambda)'  \: dx =   \int
  (\Phi_{\tilde \mu}' \cdot \nabla )  \Phi_\mu' \cdot  (\bar{    \Phi}^*_\lambda)'   \: dx
+ \int \Phi_\mu'\cdot (\bar{    \Phi}^*_\lambda)'  \nabla\cdot\Phi_{\tilde \mu}'  \: dx ,
$$
so since~$\mu$ and~$\tilde \mu$ play symmetric roles, we just need to check that
$$\sum_{\mu+\tilde \mu = \lambda \atop  i\mu, 
i\tilde \mu \in
{\mathfrak S}} \Bigl(
  \int \Phi_\mu'\cdot (\bar{    \Phi}^*_\lambda)'  \nabla\cdot\Phi_{\tilde \mu}'  \: dx
- \frac12 \int \:  \Phi_{\mu,0} \Phi_{\tilde
\mu,0} \nabla \cdot   (\bar{    \Phi}^*_\lambda)'  \: dx\Bigr)=\sum_{\mu+\tilde \mu = \lambda \atop  i\lambda, i\mu, 
i\tilde \mu \in
{\mathfrak S}}   \int \nabla \cdot
(\Phi_{\mu,0} \Phi_{\tilde \mu}' )  \bar{    \Phi}^*_{\lambda ,0} \: dx.
$$
Recalling that
$$
\nabla  \bar{    \Phi}^*_{\lambda ,0} = -i \lambda (\bar{    \Phi}^*_\lambda)' - \beta x_{1}
{ (\bar{    \Phi}^*_\lambda)'}^{\perp} = -i (\mu + \tilde \mu)  (\bar{    \Phi}^*_\lambda)' - \beta x_{1}
{ (\bar{    \Phi}^*_\lambda)'}^{\perp} ,
$$
we have
$$
 \int \nabla \cdot
(\Phi_{\mu,0} \Phi_{\tilde \mu}' )  \bar{    \Phi}^*_{\lambda,0} \: dx =  i (\mu + \tilde \mu)
\int 
\Phi_{\mu,0} \Phi_{\tilde \mu}' \cdot  (\bar{    \Phi}^*_\lambda)'  
\: dx+\int 
\Phi_{\mu,0}\beta x_{1}\Phi_{\tilde \mu}' \cdot { (\bar{    \Phi}^*_\lambda)'  }^{\perp}
\: dx.
$$
Then we write  
$$
 i \mu\Phi_{\mu,0} = \nabla \cdot \Phi_{\mu }'
$$
so that
$$
\longformule{
 \int \nabla \cdot
(\Phi_{\mu,0} \Phi_{\tilde \mu}' )  \bar{    \Phi}^*_{\lambda ,0}\: dx =
  \int\nabla \cdot \Phi_{\mu }'  \Phi_{\tilde \mu}' \cdot  (\bar{    \Phi}^*_\lambda)'  \: dx+ i \tilde \mu
\int 
\Phi_{\mu,0} \Phi_{\tilde \mu}' \cdot  (\bar{    \Phi}^*_\lambda)'  
\: dx }{- \int 
\Phi_{\mu,0}\beta x_{1}{\Phi_{\tilde \mu}'}^{\perp} \cdot { (\bar{    \Phi}^*_\lambda)'  }
\: dx.}
$$
Exchanging the roles of~$\mu$ and~$\tilde\mu$ in the first integral we get
$$\longformule{
 \sum_{\mu+\tilde \mu=\lambda \atop   i\mu, 
i\tilde \mu \in
{\mathfrak S}}\int \nabla \cdot
(\Phi_{\mu,0} \Phi_{\tilde \mu}' )   \bar{    \Phi}^*_{\lambda ,0} \: dx =
 \sum_{\mu+\tilde \mu=\lambda \atop  i\mu, 
i\tilde \mu \in
{\mathfrak S}} 
\Bigl(
\int\nabla \cdot \Phi_{\tilde\mu }'  \Phi_{ \mu}' \cdot  (\bar{    \Phi}^*_\lambda)' \: dx}
{+
 i \tilde \mu
\int 
\Phi_{\mu,0} \Phi_{\tilde \mu}' \cdot  (\bar{    \Phi}^*_\lambda)'  
\: dx -\int\Phi_{\mu,0} (i\tilde\mu\Phi_{\tilde\mu}'- \nabla \Phi_{\tilde\mu,0} ) \cdot  (\bar{    \Phi}^*_\lambda)'\: dx 
\Bigr),
}
$$
so
$$
 \sum_{\mu+\tilde \mu=\lambda \atop   i\mu, 
i\tilde \mu \in
{\mathfrak S}}\int \nabla \cdot
(\Phi_{\mu,0} \Phi_{\tilde \mu}' )   \bar{    \Phi}^*_{\lambda ,0} \: dx = 
 \sum_{\mu+\tilde \mu=\lambda \atop   i\mu, 
i\tilde \mu \in
{\mathfrak S}}  \Bigl(\int \nabla \cdot \Phi_{\tilde\mu }'  \Phi_{ \mu}' \cdot  (\bar{    \Phi}^*_\lambda)' \: dx
+ \int \Phi_{\mu,0} \nabla \Phi_{\tilde\mu,0} \cdot  (\bar{    \Phi}^*_\lambda)'\: dx\Bigr) .
$$
The result finally follows from the fact that, by symmetry,
$$
 \sum_{\mu+\tilde \mu = \lambda \atop  i\mu, 
i\tilde \mu \in
{\mathfrak S}}
\int \Phi_{\mu,0} \nabla \Phi_{\tilde\mu,0}  \cdot  (\bar{    \Phi}^*_\lambda)'\: dx
=  \frac12\sum_{\mu+\tilde \mu = \lambda \atop   i\mu, 
i\tilde \mu \in
{\mathfrak S}}
\int  \nabla (\Phi_{\tilde\mu,0} \Phi_{\mu,0} ) \cdot  (\bar{    \Phi}^*_\lambda)'\: dx \\
 $$
which finally implies that
$$
\sum_{\mu+\tilde \mu = \lambda \atop  i\mu, 
i\tilde \mu \in
{\mathfrak S}} \Bigl( \int \Phi_\mu'\cdot (\bar{    \Phi}^*_\lambda)'  \nabla\cdot\Phi_{\tilde \mu}'  \: dx
- \frac12 \int \: \Phi_{\mu,0} \Phi_{\tilde
\mu,0}  \nabla \cdot  (\bar{    \Phi}^*_\lambda)'  \: dx\Bigr)=  \sum_{\mu+\tilde \mu = \lambda \atop    i\mu, 
i\tilde \mu \in
{\mathfrak S}} \int \nabla \cdot
(\Phi_{\mu,0} \Phi_{\tilde \mu}' )  \bar{    \Phi}^*_{\lambda ,0} \: dx.
$$

\medskip

Now let us prove Proposition~\ref{limittwononlinearterms}.
The idea is to  decompose $ \Phi_\eps$ on the  eigenmodes of $L$, by writing
$$(\eta_\eps,u_\eps)(t,x)=\LL\left( {t\over \eps}\right)  \Phi_\eps  (t,x)=\sum_{i\lambda 
\in {\mathfrak S}} e^{-{it\lambda \over \eps}}\Pi_\lambda
 \Phi_\eps (t,x).$$ 
Note in particular that by (\ref{approximation-JN}), for any~$s < 0$, 
$$
(\eta_\eps ,u_\eps)(t)- \LL\left( {t\over \eps}\right) K_N \Phi_\eps  (t) \to 0 \hbox{ 
in }L^2_{loc}(\R^+;
H^s_{loc}(\R\times \T))
$$
as~$N$ goes to infinity, uniformly in $\eps$. 
Let us also introduce the notation\label{PhiepsN}
\begin{eqnarray*} \Phi_{\eps,N} &= & \LL\left(- {t\over \eps} \right) (\eta_{\eps,N} ,u_{\eps,N}) =
K_N\Phi_\eps, \quad \mbox{and} \\
  \Phi  _{\eps,\lambda,N} &= &\Pi_{\lambda} \Phi_{\eps,N}.
\end{eqnarray*}\label{PhiepslambdaN}

\medskip

We will start by considering the first nonlinear term in Proposition~\ref{limittwononlinearterms}, namely
$$
\int _0^t\int \exp\left({it'\lambda \over \eps}\right) m_\eps\cdot (u_ \eps \cdot \nabla)  
 (\bar{    \Phi}^*_\lambda)'(t',x)dxdt' .
$$
We can notice  that
\begin{eqnarray*}
\int _0^t\int \exp\left({it'\lambda \over \eps}\right) m_\eps\cdot (u_ \eps \cdot \nabla)  
(\bar{    \Phi}^*_\lambda)'(t',x)dxdt' \!\!\!\!&=&\!\!\!\!\int _0^t\int 
\exp\left({it'\lambda \over \eps}\right) \eps \eta_\eps  u_\eps\cdot 
(u_\eps \cdot \nabla)  (\bar{    \Phi}^*_\lambda)'(t',x)dxdt' \\
&+& \int _0^t\int \exp\left({it'\lambda \over \eps}\right)   u_\eps\cdot 
(u_\eps \cdot \nabla)  (\bar{    \Phi}^*_\lambda)'(t',x)dxdt' .
\end{eqnarray*}
The uniform bounds coming from the energy estimate  imply clearly  that the first term 
converges to 0 as $\eps\to 0$. 
Then we can  decompose the second contribution in the following way:
\begin{equation}\label{decomposition-nonlinear}
\begin{aligned}
\int _0^t\int \exp\left({it'\lambda \over \eps}\right)   u_\eps\cdot 
(u_\eps \cdot \nabla)  (\bar{    \Phi}^*_\lambda)'(t',x)dxdt'\\
 = \int_0^t \int_{( \R\setminus [-R,R] )\times \T  }
\exp\left({it'\lambda  \over \eps}\right) 
u_\eps\cdot (u_\eps \cdot \nabla)(\bar{    \Phi}^*_\lambda)'(t',x)dxdt'\\
+ \int_0^t \int_{   [-R,R]\times \T}  
\exp\left({it' \lambda  \over \eps} \right) 
 (u_{\eps } - u _{\eps,N})  \cdot (u_\eps \cdot
\nabla) (\bar{    \Phi}^*_\lambda)'(t',x) dxdt'\\
+\int_0^t \int_{[-R,R]\times \T} 
\exp\left({it' \lambda  \over \eps} \right)u _{\eps,N}
\cdot  \left( ( u_{\eps}- u _{\eps,M}  )  \cdot
\nabla\right)
 (\bar{    \Phi}^*_\lambda)'(t',x)dxdt'\\
+\int_0^t \int_{[-R,R]\times \T} 
\exp\left({it' \lambda \over \eps} \right) u _{\eps,N}\cdot 
(u _{\eps,M}\cdot \nabla) (\bar{    \Phi}^*_\lambda)'(t',x)dxdt'.
\end{aligned}
\end{equation}
Let us consider now all the terms in the right-hand side of~(\ref{decomposition-nonlinear}). 
 The uniform bound on~$u_{\e}$  and the decay of~$\Phi^*_\lambda$
 imply  that the first term on the right-hand side 
  converges to 0 as $R\to \infty$ uniformly in $\eps$.

By the inequality
$$\longformule{\left| \int_0^t \int_{[-R,R]\times \T }  
\exp\left({it' \lambda  \over \eps} \right) 
 (u_{\eps } - u _{\eps,N})  \cdot (u_\eps \cdot
\nabla)
 (\bar{    \Phi}^*_\lambda)'(t',x) dxdt'\right|}{
\leq C \|u_{\eps } - u _{\eps,N}\|_{L^2([0,T];H^s([-R,R]\times \T))}
\|u_\eps\|_{L^2([0,T];H^1(\R\times \T))} \| {    \Phi}^*_\lambda 
\|_{W^{2, \infty}(\R\times\T)}, \: \hbox{ with } -1<s < 0, }
 $$
we deduce  that the third
term converges to 0 as $N\to \infty$ uniformly in $\eps$.

Now let us consider the third term on the right-hand side.  Since~$ u_{\eps,N} $ corresponds to the projection
of~$\Phi_{\eps}$ onto a finite number of eigenvectors of~$L$, we deduce that 
$$ 
  \forall N\in \N,\exists C_{N} , \forall \e > 0,  \quad \|  u_{\eps,N}  \|_{L^\infty(\R^+; H^1(\R
   \times \T))} \leq C_N.
$$
Thus
$$\begin{aligned}
\left|\int_0^t \int_{[-R,R]\times \T} 
\exp\left({it' \lambda  \over \eps} \right)u_{\eps,N}
\cdot( (  u_{\eps}- u_{\eps,M}  )  \cdot
\nabla) (\bar{    \Phi}^*_\lambda)'(t',x)dxdt'\right| \\
\leq C_N \| u_{\eps}- u_{\eps,M}\|_{L^2([0,T];H^s([-R,R]\times \T))}\|{    \Phi}^*_\lambda\|_{W^{2,\infty}(\R\times\T)}
\end{aligned}$$
and, for all fixed $N$ and~$R$, the fourth term converges to 0 as $M\to \infty $ uniformly in 
$\eps$.

It remains then to take limits  in the last term of (\ref{decomposition-nonlinear}). It 
can be rewritten
$$\longformule{
\int_0^t \int_{[-R,R]\times \T}\exp\left({it'\lambda \over \eps} \right) u_{\eps,N} \cdot
(u_{\eps,M}
\cdot
\nabla) (\bar{    \Phi}^*_\lambda)'(t',x)dxdt'}{ 
= \int_0^t  \int_{[-R,R]\times \T}\sum_{i\mu ,i \tilde \mu \in 
{\mathfrak S} }
\exp\left({it'(\lambda-\mu-\tilde \mu) \over \eps}\right)  ( \Phi_{\e, \mu,N}) '
\cdot ( \Phi_{\e, \tilde \mu,M} ' \cdot \nabla) (\bar{    \Phi}^*_\lambda)'(t',x)dxdt' .}
$$
This in turn can be written in the following way:
$$\begin{aligned}
\int_0^t  \int_{[-R,R]\times \T}\sum_{i\mu ,i \tilde \mu \in 
{\mathfrak S} }
\exp\left({it'(\lambda-\mu-\tilde \mu) \over \eps}\right)  \Phi_{\e, \mu,N} '
\cdot (  \Phi_{\e, \tilde \mu,M} ' \cdot \nabla)  (\bar{    \Phi}^*_\lambda)'(t',x)dxdt' \\
= \int_0^t  \int_{[-R,R]\times \T}\sum_{i\mu ,i \tilde \mu \in 
{\mathfrak S} } \exp\left({it'(\lambda-\mu-\tilde \mu) \over \eps}\right)  
(\Phi_{\e, \mu,N} ' -  \Phi_{\mu,N} ')
\cdot ( \Phi_{\e, \tilde \mu,M} ' \cdot \nabla) (\bar{    \Phi}^*_\lambda)'(t',x)dxdt' \\
+ \int_0^t  \int_{[-R,R]\times \T}\sum_{i\mu ,i \tilde \mu \in 
{\mathfrak S} } \exp\left({it'(\lambda-\mu-\tilde \mu) \over \eps}\right)   \Phi_{\mu,N} ' 
\cdot ( (\Phi_{\e, \tilde \mu,M} ' -  \Phi_{ \tilde \mu,M}')\cdot \nabla) (\bar{    \Phi}^*_\lambda)'(t',x)dxdt' \\
+ \int_0^t  \int_{[-R,R]\times \T}\sum_{i\mu ,i \tilde \mu \in 
{\mathfrak S} } \exp\left({it'(\lambda-\mu-\tilde \mu) \over \eps}\right)    \Phi_{\mu,N} ' 
\cdot (    \Phi_{ \tilde \mu ,M }' \cdot \nabla)  (\bar{    \Phi}^*_\lambda)'(t',x)dxdt' .
\end{aligned}
$$
We have   denoted
$$
\Phi_{\mu,N} = \Pi_{\mu} \Phi_{N}, \quad \mbox{where} \quad\Phi_{N}=  K_N\Phi .
$$
The first two terms on the right-hand side go to zero as~$\e$ goes to zero, for all given~$N,M$ and~$R$, due to 
the following estimates: for~$-1<s<0$,
$$
\longformule{
\int_0^t  \int_{[-R,R]\times \T}\sum_{i\mu ,i \tilde \mu \in 
{\mathfrak S} }  \left| (\Phi_{\e, \mu,N} ' - \Phi_{\mu,N} ')
\cdot ( \Phi_{\e, \tilde \mu,M} ' \cdot \nabla)  (\bar{    \Phi}^*_\lambda)'(t',x)\right|dxdt' }{\leq C_{N,M}
 \|\Phi_{\e,N} ' - \Phi_{N} '\|_{L^{2}([0,T];H^{s}([-R,R]\times \T))} 
 \|\Phi_{\e, M} \|_{L^{\infty}([0,T];H^{1}(\R \times \T))} 
 \|{   \Phi }^*_\lambda \|_{W^{2,\infty}(\R\times \T)},}
$$
and similarly
$$
\longformule{
\int_0^t  \int_{[-R,R]\times \T}\sum_{i\mu ,i \tilde \mu \in 
{\mathfrak S} }  \left|  \Phi_{\mu,N} ' 
\cdot ( (\Phi_{\e, \tilde \mu,M} '  - \Phi_{ \tilde \mu,M} ' )\cdot \nabla) (\bar{    \Phi}^*_\lambda)'(t',x)\right|dxdt' }
{\leq C_{N,M}
 \|\Phi_{\e,M} ' - \Phi_{M} '\|_{L^{2}([0,T];H^{s}([-R,R]\times \T))} 
 \|\Phi_{\e, N} \|_{L^{\infty}([0,T];H^{1}(\R \times \T))} 
 \|  {   \Phi }^*_\lambda  \|_{W^{2,\infty}(\R\times \T)}.}
$$
Finally let us consider the last term, which can be decomposed in the following way:
$$\begin{aligned}
 \int_0^t  \int_{[-R,R]\times \T}\sum_{i\mu ,i \tilde \mu \in 
{\mathfrak S} }
\exp\left({it'(\lambda-\mu-\tilde
\mu)
\over
\eps}\right) \Phi_{\mu,N}'
\cdot (\Phi_{\tilde \mu,N}'
\cdot
\nabla)
 (\bar{    \Phi}^*_\lambda)'(t',x)dxdt'\\
=\int_0^t  \int_{[-R,R]\times \T}\sum_{i\mu ,i \tilde \mu \in 
{\mathfrak S} \atop \lambda=\mu+\tilde
\mu}
\Phi_{\mu,N}'
\cdot (\Phi_{\tilde \mu,M}'
\cdot
\nabla)
 (\bar{    \Phi}^*_\lambda)'(t',x)dxdt'\\
+\int_0^t  \int_{[-R,R]\times \T}\sum_{i\mu ,i \tilde \mu \in 
{\mathfrak S} \atop \lambda \neq
\mu+\tilde \mu}
\exp\left({it'(\lambda-\mu-\tilde
\mu)
\over
\eps}\right) \Phi_{\mu,N}'
\cdot (\Phi_{\tilde \mu,M}'
\cdot
\nabla)
 (\bar{    \Phi}^*_\lambda)'(t',x)dxdt'.
\end{aligned}
$$
For fixed $N$ and $M$, the nonstationary phase theorem (which corresponds here
to a simple integration by parts in the~$ t'$ variable) shows that the second
term is a finite sum of terms converging to 0 as
$\eps
\to 0$. And the first term (which does not depend on $\eps$)  converges to
$$\int _0^t\int \sum_{\mu+\tilde \mu=\lambda\atop  i\lambda, i\mu, i \tilde \mu \in
{\mathfrak S}}\Phi_\mu'\cdot (\Phi_{\tilde
\mu}'
\cdot
\nabla)
 (\bar{    \Phi}^*_\lambda)'(t',x)dxdt'$$
as $N,M,R\to \infty$, because~$\Phi'_{N}$ converges towards~$\Phi'$ 
strongly in~$L^{2}([0,T];L^{2}(\R \times \T))$ when~$N$ goes to infinity, and then by Lebesgue's theorem   when~$R$
goes to infinity.

Therefore, taking limits as $\eps \to 0$, then $M\to \infty$, then $N\to \infty$, then $R \to \infty$ in
(\ref{decomposition-nonlinear}) leads to
$$
\int _0^t\int \exp\left({it'\lambda \over \eps}\right) m_\eps\cdot (u_ \eps \cdot \nabla)  
 (\bar{    \Phi}^*_\lambda)'(t',x)dx dt'\to \int _0^t\int \sum_{\mu+\tilde \mu=\lambda \atop  i\lambda, i\mu, 
i\tilde \mu \in
{\mathfrak S}}\Phi_\mu'\cdot (\Phi_{\tilde
\mu}'
\cdot
\nabla)
 (\bar{    \Phi}^*_\lambda)'(t',x)dxdt'.
$$
 
 \medskip
 
 Finally let us consider the second term of the proposition, namely
 $$
 \frac12 \int_{0}^t \int  \exp\left({it'\lambda \over \eps}\right)\eta_{\e}^{2} \nabla \cdot (\bar \Phi^*_{\lambda})' (t',x) \: dxdt'.
 $$
 The first step of the above study remains valid, in the sense that one can write
 $$\longformule{
\!\!\!\! \frac12 \int_{0}^t\int  \exp\left({it'\lambda \over \eps}\right) \eta_{\e}^{2} \nabla \cdot (\bar \Phi^*_{\lambda})' (t',x) \: dxdt' 
 = \frac12 \int_{0}^t \int_{\R \setminus [-R,R] \times \T}  \exp\left({it'\lambda \over \eps}\right) \eta_{\e}^{2} \nabla \cdot (\bar \Phi^*_{\lambda})' (t',x) \: dxdt'
} {+\frac12 \int_{0}^t \int_{[-R,R] \times \T}   \exp\left({it'\lambda \over \eps}\right) \eta_{\e}^{2} \nabla \cdot (\bar \Phi^*_{\lambda})' (t',x) 
 \: dxdt',}
 $$
 and the first term converges to zero uniformly in~$\e$ as~$R$ goes to infinity, due to the spatial decay of the eigenvectors
 of~$L$. For such a result, a uniform bound of~$\eta_{\e}$ in~$L^{\infty}(\R^+;L^{2}(\R\times\T))$ is sufficient.
 However the next steps of the above study do not work here, as we have no smoothness on~$\eta_{\e}$ other
 than that energy bound. In order to conclude, let us nevertheless decompose the remaining term as above, for any 
 integers~$N$ and~$M$ to be chosen large enough below:
 \begin{equation}\label{decomposeetasquare}
 \begin{aligned}
 &\frac12 \int_{0}^t \int_{[-R,R] \times \T} \exp\left({it'\lambda \over \eps}\right)\eta_{\e}^{2} \nabla \cdot (\bar \Phi^*_{\lambda})' (t',x) 
 \: dxdt' \\
 &= \frac12 \int_{0}^t \int_{[-R,R] \times \T}  \exp\left({it'\lambda \over \eps}\right) ( \eta_{\e} -  \eta_{\e,N}) \eta_{\e}\nabla \cdot (\bar \Phi^*_{\lambda})' (t',x) 
 \: dxdt'
\\
&+\frac12 \int_{0}^t \int_{[-R,R] \times \T}  \exp\left({it'\lambda \over \eps}\right) \eta_{\e,N}( \eta_{\e} -  \eta_{\e,M}) \nabla \cdot (\bar \Phi^*_{\lambda})' (t',x) 
 \: dxdt'\\
 &+\frac12 \int_{0}^t \int_{[-R,R] \times \T}  \exp\left({it'\lambda \over \eps}\right) \eta_{\e,N} \eta_{\e,M}  \nabla \cdot (\bar \Phi^*_{\lambda})' (t',x) 
 \: dxdt'.
 \end{aligned}
\end{equation}
 The sequence~$\displaystyle
 -\frac12   \exp\left({it'\lambda \over \eps}\right)  ( \eta_{\e} -  \eta_{\e,N}) \eta_{\e}$ is uniformly bounded in~$N \in \N$
and~$\e > 0$ in the space~$L^{1}_{loc}(\R^+\times \R \times \T)$, so up to the extraction of a subsequence 
 it converges weakly, as~$\e$ goes to zero,
 towards a  measure~$\upsilon_{\lambda,N}$, which in turn is uniformly 
 bounded in~${\mathcal M}(\R^+\times \R \times \T)$. Denoting by~$\upsilon_{\lambda}$ its limit
  in~${\mathcal M}(\R^+\times \R \times \T)$ as~$N$ goes to infinity,  we find that
  $$
  \frac12 \int_{0}^t \int_{[-R,R] \times \T} \exp\left({it'\lambda \over \eps}\right) ( \eta_{\e} - 
   \eta_{\e,N}) \eta_{\e}\nabla \cdot (\bar \Phi^*_{\lambda})' (t',x) 
 \: dxdt' \to -   \int_{0}^t \int_{[-R,R] \times \T}\nabla \cdot (\bar \Phi^*_{\lambda})'
   \upsilon_{\lambda} (dt' dx) 
  $$
  as~$\e$ goes to zero and~$N$ goes to infinity,
  which in turn converges to
  $$
  -   \int_{0}^t \int \nabla \cdot (\bar \Phi^*_{\lambda})'
   \upsilon_{\lambda} (dt' dx)
  $$
  as~$R$ goes to infinity, due to the smoothness of~$\nabla \cdot (\bar \Phi^*_{\lambda})' $.
  
  Note that as~${\mathfrak S}$ is countable, one can choose a subsequence such that for all~$i \lambda \in {\mathfrak S}$, 
  the sequence~$\displaystyle
 -\frac12   \exp\left({it'\lambda \over \eps}\right)  ( \eta_{\e} -  \eta_{\e,N}) \eta_{\e}$ converges 
 towards~$\upsilon_{\lambda }$ as~$\e$ goes to zero and~$N$ goes to infinity.

Now let us consider the two last terms in~(\ref{decomposeetasquare}).  We recall that~$\eta_{\e,N}$ corresponds
to the projection of~$\Phi_{\e}$ onto a finite number of eigenvectors of~$L$, so it is smooth for each fixed~$N$. In 
particular we can write, for any~$s<0$,
$$\longformule{
\frac12\left|
\int_0^{t} \int_{[-R,R] \times \T}  \exp\left({it'\lambda \over \eps}\right) \eta_{\e,N}(
 \eta_{\e} -  \eta_{\e,M}) \nabla \cdot (\bar \Phi^*_{\lambda})' (t',x) 
 \: dxdt'
\right| }{\leq C_{N} \| \eta_{\e} -  \eta_{\e,M}\|_{L^{\infty}([0,T];H^{s}([-R,R] \times \T))} 
\|\Phi^*_{\lambda}\|_{W^{2,\infty}(\R\times\T)}.}
$$
  So letting~$M$ go to infinity we find that this term converges to zero uniformly in~$\e$ for each fixed~$N$ and~$R$.
  
  Finally for the last term of~(\ref{decomposeetasquare}) we write similar computations as for the first nonlinear term 
  in Proposition~\ref{limittwononlinearterms}. We have indeed
  $$
  \begin{aligned}
 &\int_{0}^t \int_{[-R,R] \times \T}  \exp\left({it'\lambda \over \eps}\right) \eta_{\e,N} \eta_{\e,M}  
  \nabla \cdot (\bar \Phi^*_{\lambda})' (t',x)  \: dxdt' \\
  &=    \int_{0}^t \int_{[-R,R] \times \T}
  \sum_{i\mu,i\tilde\mu \in {\mathfrak S}}
 \exp\left({it'(\lambda-\mu-\tilde
\mu)
\over
\eps}\right)  (\Phi_{\mu,\e,N})_{0}
  (\Phi_{\tilde\mu,\e,M})_{0} \nabla \cdot (\bar \Phi^*_{\lambda})' (t',x)  \: dxdt'\\
 & =   \int_{0}^t \int_{[-R,R] \times \T}
  \sum_{i\mu,i\tilde\mu \in {\mathfrak S}}
\exp\left({it'(\lambda-\mu-\tilde
\mu)
\over
\eps}\right)  ( \Phi_{\mu,\e,N} -\Phi_{\mu ,N} )_{0}
  (\Phi_{\tilde\mu,\e,M})_{0} \nabla \cdot (\bar \Phi^*_{\lambda})' (t',x)  \: dxdt'\\
&  + \int_{0}^t \int_{[-R,R] \times \T}
  \sum_{i\mu,i\tilde\mu \in {\mathfrak S}}
 \exp\left({it'(\lambda-\mu-\tilde
\mu)
\over
\eps}\right) (\Phi_{\mu ,N})_{0}
  ( \Phi_{\tilde\mu,\e,M} -\Phi_{\tilde\mu ,M} )_{0} \nabla \cdot (\bar \Phi^*_{\lambda})' (t',x)  \: dxdt'\\
   &+ \int_{0}^t \int_{[-R,R] \times \T}
  \sum_{i\mu,i\tilde\mu \in {\mathfrak S}}
 \exp\left({it'(\lambda-\mu-\tilde
\mu)
\over
\eps}\right)  (\Phi_{\mu ,N})_{0}(\Phi_{\tilde\mu ,M} )_{0}
  \nabla \cdot (\bar \Phi^*_{\lambda})' (t',x)  \: dxdt'.
  \end{aligned}
  $$
  The two first terms in the right-hand side are easily shown to converge to zero as~$\e$   goes to zero,
  for each fixed~$N$ and~$M$.
  We have indeed
  $$
  \longformule{
  \left|
  \int_{0}^t \int_{[-R,R] \times \T}
  \sum_{i\mu,i\tilde\mu \in {\mathfrak S}}
\exp\left({it'(\lambda-\mu-\tilde
\mu)
\over
\eps}\right)  ( \Phi_{\mu,\e,N} -\Phi_{\mu ,N} )_{0}
  (\Phi_{\tilde\mu,\e,M})_{0} \nabla \cdot (\bar \Phi^*_{\lambda})' (t',x)  \: dxdt'
  \right|
  }{\leq C_{M}
  \| \Phi_{\mu,\e,N} -\Phi_{\mu ,N}\|_{L^{\infty}([0,T];H^{s}([-R,R] \times \T))}
  \|\Phi^*_{\lambda}\|_{W^{2,\infty}}
,  }
  $$
 and similarly
  $$
  \longformule{
  \left|
 \int_{0}^t \int_{[-R,R] \times \T}
  \sum_{i\mu,i\tilde\mu \in {\mathfrak S}}
 \exp\left({it'(\lambda-\mu-\tilde
\mu)
\over
\eps}\right) (\Phi_{\mu ,N})_{0}
  ( \Phi_{\tilde\mu,\e,M} -\Phi_{\tilde\mu ,M} )_{0} \nabla \cdot (\bar \Phi^*_{\lambda})' (t',x)  \: dxdt'
  \right|
  }{\leq C_{N}
  \| \Phi_{\tilde\mu,\e,M} -\Phi_{\tilde\mu ,M}\|_{L^{\infty}([0,T];H^{s}([-R,R] \times \T))}
  \|\Phi^*_{\lambda}\|_{W^{2,\infty}}
.  }
  $$
 Finally the last term on the right-hand side is dealt with by a nonstationary phase argument, and we have as above, 
 as~$\e$ goes
 to zero and then~$M$, $N$ and~$R$ go successively to infinity, 
 $$\longformule{
\frac12  \int_{0}^t \int_{[-R,R] \times \T}
  \sum_{i\mu,i\tilde\mu \in {\mathfrak S}}
 \exp\left({it'(\lambda-\mu-\tilde
\mu)
\over
\eps}\right)  (\Phi_{\mu ,N})_{0}(\Phi_{\tilde\mu ,M} )_{0}
  \nabla \cdot (\bar \Phi^*_{\lambda})' (t',x)  \: dxdt' }{\to \frac12
  \int_0^t\int
 \sum_{\mu+\tilde \mu=\lambda \atop   i\mu,  i\tilde \mu \in {\mathfrak S}}
 \Phi_{\mu,0} \Phi_{\tilde
\mu,0}  \nabla \cdot(\bar{    \Phi}^*_\lambda)' (t',x)dxdt'.}
 $$
  Proposition~\ref{limittwononlinearterms} is proved, and therefore also
  Theorem~\ref{thmstrongweakcv}.
\end{proof} 

\subsection{The case when capillarity is added}\label{capillarity} 
In this final paragraph we propose an adaptation to the  Saint-Venant model which provides some additional smoothness
on~$\eta_{\e}$, and which enables one to get rid of the defect measure present in the above study.  The model is presented
in the next part, and the convergence result stated and proved below.

\subsubsection{The model}
Let us present an alternative to the Saint-Venant model studied up to now, which presents the advantage of providing the
additional smoothness of~$\e \eta_{\e}$ which is lacking in the original system. Its disadvantage however is that
there is no real physical meaning to the capillarity operator we use in that model. 
With the notation
of Chapter~\ref{intro}, we choose indeed the capillarity operator
\begin{equation}\label{defcapillarity}
K(h)= \kappa (-\Delta)^{2\alpha} h,  
\end{equation}
where~$\kappa > 0$ and~$\alpha> 1/2$ are given constants. 
After rescaling as in  Chapter~\ref{intro}, we find the following system: 
\begin{equation}
\label{SW-epscap}
\begin{aligned}
\d_t \eta +\frac1\eps \DIV \Bigl((1+\eps \eta) u\Bigr) =0,  \\
\d_t  u  + u \cdot \nabla  u +  
\frac{ \beta x_1}{ \eps} u^\perp+
{1\over \eps} 
\nabla
\eta
-\frac\nu{1+\e\eta} \Delta u + \e\kappa\nabla(-\Delta)^{2\alpha} \eta =0,  \\
\eta_{|t=0}=\eta^0,\quad u_{|t=0} =u^0.
\end{aligned}
\end{equation}
 In the next part we discuss the existence of bounded energy solutions to that system of equations (under a smallness assumption), and the following part 
 consists in the proof of the analogue of Theorem~\ref{thmstrongweakcv} in that setting. One should emphasize here that the additional
 capillarity term that is added in the system will not appear in the limit, since it comes as a~$O(\e)$ term. Moreover it
 is a linear term, so it should not change the other  asymptotics proved in this chapter. However its unphysical 
 character (as well as the smallness condition on the initial data) made us  prefer to study the original Saint-Venant system for all  the convergence results of this chapter.

\subsubsection{Existence of solutions} The following theorem is an easy adaptation of the result by D. Bresch
and B. Desjardins in~\cite{breschdesjardins} (see also~\cite{pll} for the
compressible Navier-Stokes system), we give a sketch of the proof below.
 \begin{Thm}[existence of solutions in the case with capillarity]\label{existencecapillary}
There is a constant~$C>0$ such that the following result holds.
 Let~$(\eta_{\e}^{0},u_{\e}^{0})$ be a family of~$H^{2\alpha} \times L^{2}(\R\times\T)$ such that
 for all~$\e>0$,
 $$
 \frac12\int\Bigl( (\eta_{\e}^{0})^{2} + \kappa \e^{2} |(-\Delta)^{\alpha}\eta_{\e}^{0}|^{2} 
 + (1+\e\eta_{\e}^{0}) |u_{\e}^{0}|^{2}\Bigr)(x) \: dx \leq {\mathcal E}^{0}.
 $$
 If~$ {\mathcal E}^{0} \leq C$, then
  there is a family~$(\eta_{\e} ,u_{\e} )$ of weak solutions to~(\ref{SW-epscap}), satisfying the 
 energy estimate
 $$ 
  \frac12\int \Bigl( \eta_{\e}^{2} + \kappa \e^{2} |(-\Delta)^{\alpha}\eta_{\e}|^{2} 
 + (1+\e\eta_{\e} ) |u_{\e} |^{2}\Bigr) (t,x)\: dx   + \nu \int_{0}^{t} \int |\nabla u_{\e} |^{2}(t',x)\: dx  dt'
 \leq {\mathcal E}^{0}.  
 $$
 \end{Thm}
 \begin{proof}
 Weak
solutions can be constructed by a  standard approximation scheme
obtained by regularization : compactness on the approximate
solutions comes from the a priori
bounds derived from the energy  inequality, which
is obtained formally in a classical way by multiplying the
momentum equation by $u$, using the mass conservation and integrating by parts. 
It allows to derive immediately the following a priori bounds
(denoting by~$ \eta$ and~$ u$ approximate solutions) :
$$
\begin{array}{r}
\eta  \in L^\infty (\R^+; L^{2 } (\R \times \T)) \\
\e\eta  \in L^\infty (\R^+; H^{2\alpha} (\R \times \T)) \\
 (1+\e\eta) |u|^2 \in L^\infty (\R^+; L^{1} (\R \times \T)) \\
  |\nabla u|^2 \in L^1 (\R^+; L^{1} (\R \times \T)).
\end{array}
$$
Since~$ \alpha > 1/2$, the first bound implies in particular that
$$
\e\eta  \in L^\infty (\R^+ \times \R \times \T), 
$$
and in particular if~$ \EE_0$ is small enough (compared to the reference height which is~1 here), 
then~$ 1+\e \eta$ is bounded from below. We infer that~$ u$ is bounded in~$
L^\infty (\R^+;L^2 (\R \times \T)) \times L^2 (\R^+,\dot H^1 (\R \times \T))$.

 It is standard
(see~\cite{breschdesjardins}, \cite{pll}) to deduce that~$ 1+\e \eta$ is
compact in~$ L^2_{loc} (\R^+;H^1 (\R \times \T))$, and that~$ u$  is
compact in~$ L^2_{loc} (\R^+;L^2 (\R \times \T))$. 

Taking the limit in the non linear terms is now possible :  we need indeed to deal with the
following nonlinear terms : 
$$
 u \cdot \nabla u, \quad \DIV ((1+\e \eta) u) \quad \mbox{and} \quad \frac{1}{1+\e \eta}
\Delta u .
$$
The compactness of~$1+\e \eta$ and~$ u$ derived above allows to deal with the
two first terms in a standard fashion. For the last one we just have
to recall that~$ 1/(1+\e \eta)$ is bounded in~$  L^\infty (\R^+;H^{2\alpha}  ( \R \times
\T))$. This  completes the proof of Theorem~\ref{existencecapillary}. 
\end{proof}

\subsubsection{Convergence} In this section our aim is to show that the   capillarity term enables us to 
get rid of the defect measure present in the conclusion of Theorem~\ref{thmstrongweakcv} page~\pageref{thmstrongweakcv}. 
As the proof is very similar
to that theorem, up to the compactness of~$\eta_{\e}$, we will not give the full details. The result is the following.
\begin{Thm}[strong convergence  in the case with capillarity]\label{thmstrongweakcvcapillarity}
  Under the assumptions of Theorem~\ref{existencecapillary}, denote by~$(\eta_\eps,u_\eps)$ a solution of~(\ref{SW-epscap}) with initial
data~$(\eta_\eps^0,u_\eps^0)$, and define $$\Phi_\eps=\LL\left(-{t\over  \eps}\right) 
(\eta_\eps,u_\eps).$$ 
Up  to the extraction of a subsequence, $\Phi_\eps$
converges weakly  in~$L^2_{loc}(\R^+;H^s_{loc}(\R\times \T))$ (for  all~$s<0$) toward some 
  solution~$\Phi$
of the following limiting filtered system: for all~$i\lambda  $ in~$ {\mathfrak S} $ and for
all smooth~$\Phi^*_\lambda  $ in~$  \Ker (L-i\lambda
Id)$,
$$
\int \Phi \cdot \bar \Phi^*_\lambda (x) \: dx - \nu \int_{0}^{t} \int \Delta_{L}'
 \Phi  \cdot  \bar \Phi^*_\lambda  
(t',x) \: dxdt'+\int_{0}^{t}  \int Q_{L} (\Phi,\Phi) \cdot  
 \bar \Phi^*_\lambda(t',x) \:dx dt' = \int \Phi^{0} \cdot  \bar \Phi^*_\lambda (x) \: dx,
$$
where~$\Phi^0 = (\eta^0,u^0)$. 
\end{Thm}
\begin{proof}
We will follow the lines of the proof of Theorem~\ref{thmstrongweakcv}. In particular all the results of 
Sections~\ref{strongcompactpilambdaphieps} and~\ref{strongcvphiepssection}
are true in this situation and we will not detail the proofs.Ê So the point, as in Section~\ref{takinglimits}, consists in taking the limit as~$\e$ goes to
zero, of  the 
equation on~$\Pi_\lambda \Phi_\eps$. 

Equation~(\ref{pil-eqsum}) page~\pageref{pil-eqsum} can be written here as follows:
for all smooth~$    \Phi^*_\lambda  = 
(   \Phi^*_{\lambda,0},    (\Phi^*_\lambda)') $ belonging to~$ \Ker (L-i\lambda
Id)$,
\begin{equation}
\label{pil-eqsumcap}
\begin{aligned}
\int \exp\left({it\lambda \over \eps}\right) 
 (\eta_\eps \bar {    \Phi}^*_{\lambda,0}
 +u_\eps \cdot
 (\bar{    \Phi}^*_\lambda)')(t,x)dx-\int (\eta_\eps^0 \bar {    \Phi}^*_{\lambda,0}+u_ \eps^0 \cdot
 (\bar{    \Phi}^*_\lambda)')(x)dx\\
-\e \kappa \int_0^t \int(-\Delta)^{ \alpha} \exp\left({it'\lambda \over \eps}\right)
\eta_{\e} \nabla \cdot (-\Delta)^{ \alpha} (\bar{    \Phi}^*_\lambda)'(t',x) dxdt'\\
- \int_0^t \int \frac\nu{1+\e\eta_{\e}}
 \Delta  \left(\exp\left({it'\lambda \over \eps}\right)u_\eps  \right) \cdot
 (\bar{    \Phi}^*_\lambda)'(t',x)dxdt'\\
+\int _0^t\int \exp\left({it'\lambda \over \eps} \right) ( u_\eps \cdot \nabla )  u_\eps   \cdot
  (\bar{    \Phi}^*_\lambda)'     (t',x)dxdt'
 -\int _0^t\int \exp\left({it'\lambda \over \eps} \right) 
\eta_{\e} u_{\e} \cdot \nabla
  \bar{    \Phi}^*_{\lambda,0} (t',x)dxdt'  =0.
\end{aligned}
\end{equation}

\begin{Rem}
We have chosen to keep the unknowns~$(\eta_{\e},u_{\e})$ and not write the analysis in terms of~$(\eta_{\e},m_{\e})$
as previously (recall that~$m_{\e} = (1+\e\eta_{\e})u_{\e}$): the study of~$m_{\e}$ rather than~$u_{\e}$ 
is indeed unnecessary here
as the factor~$\displaystyle\frac1{1+\e\eta_{\e}}$ which appears in the diffusion term in the equation on~$u_{\e}$
can be controled in this situation, contrary to the previous case. 
The advantage of writing the equations on~$(\eta_{\e},u_{\e})$ is that there is no nonlinear term in~$\eta_{\e}$, contrary to
the previous study, but of course the difficulty is transfered to the study of the diffusion operator; the gain of regularity 
in~$\eta_{\e}$ will appear here.
\end{Rem}

Taking limits as $\eps \to 0$ in the two first terms is immediate.
For the third term, we simply recall that~$\eta_{\e}$ is bounded in~$L^{\infty}(\R^+;L^{2}(\R\times \T))$
and~$\e\eta_{\e}$ is bounded in~$L^{\infty}(\R^+;H^{2\alpha}(\R\times \T))$, so~$\e\eta_{\e}$ goes strongly
to zero in~$L^{\infty}(\R^+;H^{s}(\R\times \T))$ for every~$s<2\alpha$. Since~${    \Phi}^*_\lambda$ is smooth, we infer that
$$
\e \kappa \int_0^t \int (-\Delta)^{ \alpha} \exp\left({it'\lambda \over \eps}\right)
\eta_{\e} \nabla \cdot (-\Delta)^{ \alpha} (\bar{    \Phi}^*_\lambda)' (t',x) dxdt' \to 0, \quad
\mbox{as} \: \e \to 0.
$$
Let us now consider the fourth term,
$$
-\int_0^t \int \frac\nu{1+\e\eta_{\e}}
 \Delta  \left(\exp\left({it'\lambda \over \eps}\right)u_\eps  \right) \cdot
 (\bar{    \Phi}^*_\lambda)'(t',x)dxdt' .
$$
It is here that the presence of capillarity enables us to get a better control. Let us write
$$
\begin{aligned}
&-\int_0^t \int \frac\nu{1+\e\eta_{\e}}
 \Delta  \left(\exp\left({it'\lambda \over \eps}\right)u_\eps  \right) \cdot
 (\bar{    \Phi}^*_\lambda)'(t',x)dxdt' \\
& = \nu\int_0^t \int  
 \nabla \left(\exp\left({it'\lambda \over \eps}\right)u_\eps  \right) : \nabla
 (\bar{    \Phi}^*_\lambda)'(t',x)dxdt'\\
& - \nu \int_0^t \int \nabla \left(\exp\left({it'\lambda \over \eps}\right)u_\eps  \right) : \nabla
 \left( \frac{\e \eta_\eps} {1+\e\eta_{\e}}  (\bar{    \Phi}^*_\lambda)'\right)(t',x)dxdt'.
\end{aligned}
$$
Clearly the first term on the right-hand side converges towards the expected limit: we have
$$
 \nu\int_0^t \int  
 \nabla \left(\exp\left({it'\lambda \over \eps}\right)u_\eps  \right) : \nabla
 (\bar{    \Phi}^*_\lambda)'(t',x)dxdt' \rightarrow  \nu\int_0^t \int  
 \nabla \Phi_{\lambda}'  : \nabla
 (\bar{    \Phi}^*_\lambda)'(t',x)dxdt' ,  \quad \mbox{as} \: \e \to 0.
$$ 
 To study the  second one, we can 
notice that
$$
\nabla
 \left( \frac{\e \eta_\eps} {1+\e\eta_{\e}}  (\bar{    \Phi}^*_\lambda)'\right) = 
 \nabla \left(\frac{\e \eta_\eps} {1+\e\eta_{\e}}   \right)(\bar{    \Phi}^*_\lambda)' + 
 \frac{\e \eta_\eps} {1+\e\eta_{\e}}\nabla (\bar{    \Phi}^*_\lambda)', 
$$
and since the second term on the right-hand side
 is obviously easier to study than the first one,   let us concentrate on the first term. We have
$$
\nabla \frac{\e \eta_{\e} } {1+\e\eta_{\e}} = \frac{\e  \nabla\eta_{\e}} {1+\e\eta_{\e}} - 
\frac{\e^{2} \eta_{\e} \nabla \eta_{\e}} {(1+\e\eta_{\e})^{2}} \cdotp
$$
Since~$\e \eta_{\e}$ is bounded in~$L^{\infty}(\R^+;H^{2\alpha}(\R\times \T))$, we infer easily, 
by product laws in Sobolev spaces, that 
$$
\e^{2} \eta_{\e} \nabla \eta_{\e} \quad \mbox{is bounded in} \quad 
L^{\infty}(\R^+;H^{\sigma}(\R\times \T)), \: \mbox{for some} \: \sigma>0.
$$
But on the other hand~$ \eta_{\e}$ is bounded in~$L^{\infty}(\R^+;L^{2}(\R\times \T))$, so we have
also
$$
\e^{2} \eta_{\e} \nabla \eta_{\e} \to 0\quad \mbox{ in} \quad L^{\infty}(\R^+;H^{2\alpha -2}
(\R\times \T)).
$$
By interpolation we gather that
$$
\e^{2} \eta_{\e} \nabla \eta_{\e} \to 0\quad \mbox{ in} \quad L^{\infty}(\R^+;L^{2}
(\R\times \T)),
$$
and the lower bound on~$1+\e\eta_{\e}$ ensures that
$$
\frac{\e^{2} \eta_{\e} \nabla \eta_{\e}} {(1+\e\eta_{\e})^{2}}   \to 0\quad \mbox{ in} \quad L^{\infty}(\R^+;L^{2}
(\R\times \T)).
$$
The argument is similar (and easier) for the term~$\displaystyle \frac{\e  \nabla\eta_{\e}} {1+\e\eta_{\e}} $, 
so we can conclude that
$$
-\int_0^t \int \frac\nu{1+\e\eta_{\e}}
 \Delta  \left(\exp\left({it'\lambda \over \eps}\right)u_\eps  \right) \cdot
 (\bar{    \Phi}^*_\lambda)'(t',x)dxdt' \to \nu \int_0^t \int \nabla \Phi_{\lambda} : \nabla  (\bar{    \Phi}^*_\lambda)'(t',x)dxdt'.
$$
Finally we are left with the nonlinear terms: let us study the limit of
$$
\int _0^t\int \exp\left({it'\lambda \over \eps} \right) 
(u_\eps\cdot \nabla) u_\eps \cdot  
 (\bar{    \Phi}^*_\lambda)'(t',x)dxdt'
 -\int _0^t\int \exp\left({it'\lambda \over \eps} \right) 
\eta_{\e} u_{\e}\nabla \cdot 
  \bar{    \Phi}^*_{\lambda,0} (t',x)dxdt'.
$$
 The study is very similar to the case studied in Section~\ref{takinglimits} (see 
 Proposition~\ref{limittwononlinearterms}), so we will not give
 all  the details but merely point out the differences.
 First, one can truncate the integral in~$x_{1} \in \R$ to~$x_{1} \in [-R,R]$, where~$R$ is a parameter to be chosen
  large enough in the
 end. As previously that is simply due to the decay of the eigenvectors of~$L$ at infinity. So we are reduced to the study of
 $$
  \int _0^t\int_{[-R,R] \times\T}  \exp\left({it'\lambda \over \eps} \right) 
(u_\eps\cdot \nabla) u_\eps \cdot  
 (\bar{    \Phi}^*_\lambda)'(t',x)dxdt'  \quad \mbox{and}
 $$
 $$
 -\int _0^t\int_{[-R,R] \times\T} \exp\left({it'\lambda \over \eps} \right) 
\eta_{\e} u_{\e} \cdot \nabla
  \bar{    \Phi}^*_{\lambda,0} (t',x)dxdt'.
 $$
 The limit of the first term is obtained in an identical way to Section~\ref{takinglimits}, since~$ u_{\e}  $ satisfies
 the same bounds, so we have
 $$
 \int _0^t\int_{[-R,R] \times\T}  \exp\left({it'\lambda \over \eps} \right) 
(u_\eps\cdot \nabla) u_\eps \cdot  
 (\bar{    \Phi}^*_\lambda)'(t',x)dxdt' \rightarrow \int _0^t\int 
 \sum_{\mu+\tilde \mu=\lambda \atop   i\mu,  i\tilde \mu \in {\mathfrak S}}
( \Phi_{\mu }' \cdot \nabla )\Phi_{\tilde
\mu }'   \cdot(\bar{    \Phi}^*_\lambda)' (t',x)dxdt',
 $$ 
 as  $\e $ goes to  0  and  $R$ goes to infinity.
 
 Now let us concentrate on the last nonlinear term. With the notation of Section~\ref{takinglimits}, we can write
 $$
 \begin{aligned}
 \int _0^t\int_{[-R,R] \times\T} \exp\left({it'\lambda \over \eps} \right) 
\eta_{\e} u_{\e} \cdot \nabla
  \bar{    \Phi}^*_{\lambda,0} (t',x)dxdt'\\
   =  \int _0^t\int_{[-R,R] \times\T} \exp\left({it'\lambda \over \eps} \right) 
(\eta_{\e}- \eta_{\e,N}) u_{\e} \cdot \nabla
  \bar{    \Phi}^*_{\lambda,0} (t',x)dxdt' \\
  +\int _0^t\int_{[-R,R] \times\T} \exp\left({it'\lambda \over \eps} \right) 
 \eta_{\e,N} ( u_{\e}- u_{\e,M} ) \cdot \nabla
  \bar{    \Phi}^*_{\lambda,0} (t',x)dxdt'\\
    +\int _0^t\int_{[-R,R] \times\T} \exp\left({it'\lambda \over \eps} \right) 
 \eta_{\e,N}   u_{\e,M}  \cdot \nabla
  \bar{    \Phi}^*_{\lambda,0} (t',x)dxdt'.
 \end{aligned}
 $$
 The first two terms on the right-hand side converge to zero, due to the following estimates: for some~$-1<s<0$ and for 
 all~$t \in [0,T]$,
 $$
 \longformule{
 \left|
  \int _0^t\int_{[-R,R] \times\T} \exp\left({it'\lambda \over \eps} \right) 
(\eta_{\e}- \eta_{\e,N}) u_{\e} \cdot \nabla
  \bar{    \Phi}^*_{\lambda,0} (t',x)dxdt'  
 \right| }{\leq C_{T}\|\eta_{\e}- \eta_{\e,N}\|_{L^{\infty}([0,T];H^{s}([-R,R] \times\T))}
 \|u_{\e}\|_{L^{2}([0,T];H^{1}([-R,R] \times\T))} \|  \bar{    \Phi}^*_{\lambda } \|_{W^{2,\infty}
 (\R\times\T)},}
 $$
 and similarly
  $$
   \longformule
{ \left|
  \int _0^t\int_{[-R,R] \times\T} \exp\left({it'\lambda \over \eps} \right) 
  \eta_{\e,N}  ( u_{\e}- u_{\e,M} ) \cdot \nabla
  \bar{    \Phi}^*_{\lambda,0} (t',x)dxdt'  
 \right|} {\leq C_{T,N} 
 \|u_{\e}- u_{\e,M}\|_{L^{2}([0,T];H^{s}([-R,R] \times\T))} \|  \bar{    \Phi}^*_{\lambda } \|_{W^{2,\infty}
 (\R\times\T)}.}
 $$
 Finally the limit of the third term is obtained by the (by now) classical nonstationary phase theorem, namely 
 we find, exactly as in the proof of Proposition~\ref{limittwononlinearterms},  that
 $$
 \int _0^t\int_{[-R,R] \times\T} \exp\left({it'\lambda \over \eps} \right) 
 \eta_{\e,N}   u_{\e,M}  \cdot \nabla
  \bar{    \Phi}^*_{\lambda,0} (t',x)dxdt' \rightarrow 
  \int _0^t\int 
 \sum_{\mu+\tilde \mu=\lambda \atop   i\mu,  i\tilde \mu \in {\mathfrak S}}
 \Phi_{\mu,0 }   \Phi_{\tilde
\mu }'  \cdot \nabla  \bar{    \Phi}^*_{\lambda,0}  (t',x)dxdt',
 $$ 
 as  $\e $ goes to  0  and~$M$, $N$ and  $R$ go  to infinity.
 
 That concludes the proof of the theorem.
\end{proof}

\backmatter


\chapter*{Notation Index}

 $\Delta'$, diffusion operator in the Saint-Venant system, p. \pageref{op-not}

 $\Delta'_L$, diffusion operator in the limit filtered system, p. \pageref{opL-not}

\bigskip
$\EE_\eps(t)$,  modulated energy,  p. \pageref{strong-weak}

$\EE_{\eps,N}(t)$,  modulated energy applied to a sequence of approximate solutions,  p. \pageref{strong-weakN}

\bigskip
$\FF$, Fourier transform, p. \pageref{ff-not}

$\FF_2$, Fourier transform with respect to~$x_2$, p. \pageref{defFF2}

 $\widehat f (n,k) $, coefficients of $f$ in the Hermite-Fourier basis, p.
\pageref{widehat-not}

 \bigskip
 $H^s$, inhomogeneous Sobolev space, p. \pageref{Hs-not}

 $\dot H^s$, homogeneous Sobolev space, p. \pageref{Hs-not}

$H^s_0(\Omega)$, for $s>0$, Sobolev space on a bounded set~$\Omega$ with Dirichlet boundary conditions, p. \pageref{hs0}

$H^{-s}(\Omega)$,   for $s>0$, dual space of $H^s_0(\Omega)$, p. \pageref{hs0} 

 $H^s_L$, weighted Sobolev space, p. \pageref{HsL}

 \bigskip
$ J_N$,  spectral truncation,  p. \pageref{JN}

\bigskip
$K_N$, truncation operator in the $(\Psi_{n,k,j})$ basis, p. \pageref{KN-not}

$\kappa$, a regularizing kernel, p.  \pageref{kappa}

 \bigskip
 
 $\ell^{p,\infty}$, weakly convergent series, p. \pageref{ellqfaible}
 
 $L$, singular perturbation, p. \pageref{L-not}

 $\LL$, semi-group generated by $L$, p. \pageref{LL-not}

\bigskip
$N_s$, pseudo-differential operator on $\Ker L$, p. \pageref{Ns-not}

 \bigskip

 $\Pi_{n,k,j}$, projection on $\Psi_{n,k,j}$, p. \pageref{pinkj-not}

 $\Pi_\lambda$, projection on $\Ker (L-i\lambda Id)$, p. \pageref{pinkj-not}

 $\Pi_0$, projection on $\Ker L$, p. \pageref{pi0-not}

 $\Pi_\perp$, projection on $(\Ker L)^\perp$, p. \pageref{pip-not}

 $\Pi_K$, projection on Kelvin modes, p. \pageref{PRKN}

 $\Pi_M$, projection on mixed Rossby-Poincar{\'e} modes, p. \pageref{PRKN}

 $\Pi_P$, projection on Poincar{\'e} modes, p. \pageref{PRKN}

 $\Pi_R$,  projection on Rossby modes, p. \pageref{PRKN}

 \bigskip
 $Q$, quadratic operator in the Saint-Venant system, p. \pageref{op-not}

 $Q_L$, quadratic operator in the limit filtered system, p. \pageref{opL-not}

$\tilde Q_L$, purely ageostrophic quadratic operator in the limit filtered system, p.
\pageref{tQL-def}

 \bigskip
 $S$, set of indices $(n,k,j)$, p. \pageref{S-not}

 $S^*$, set of indices $(n,k,j)$ corresponding to $(\Ker L)^\perp$, p. \pageref{S*-not}

$(SW_\eps)$,  the shallow water system, p. \pageref{SWeps}

$(SW_0)$, the limit system after filtering,  p. \pageref{SW0}

 ${\mathfrak S}$, spectrum of $L$, p. \pageref{sig-not}

 ${\mathfrak S}_K$, subset of $\mathfrak S$ corresponding to Kelvin modes, p.
\pageref{sig-not}

 ${\mathfrak S}_P$,  subset of $\mathfrak S$ corresponding to Poincar{\'e} modes, p.
\pageref{sig-not}

 ${\mathfrak S}_R$,  subset of $\mathfrak S$ corresponding to Rossby modes, p.
\pageref{sig-not}

 ${\mathfrak S}_N$, subset of $\mathfrak S$ defined by a frequency truncation, p.
\pageref{sigmaN}

 \bigskip
 $ i\tau_{n,k,j}$, eigenvalues of $L$, p. \pageref{tau-not}

 \bigskip
$W^{s,\infty}$, Sobolev space, p. \pageref{diag-prop}

 \bigskip
 $\Phi_0$, first coordinate of the three component   vector field $\Phi$, p.
\pageref{comp-not}

 $\Phi' = (\Phi_1,\Phi_2)$, two last coordinates of the three component  vector field $\Phi$, p.
\pageref{comp-not}

${\Phi'} ^\perp = (\Phi_2,-\Phi_1)$, image of~$\Phi'$ by a rotation of angle~$\pi/2$,  p.
\pageref{phiperp}

 $\bar \Phi$, complex conjugate of $\Phi$, p. \pageref{ccPhi}

 $\Phi_\eps=\LL\left(-{t\over  \eps}\right) (\eta_\eps,u_\eps)$,
 where~$(\eta_\eps,u_\eps) $ solves
 the Saint-Venant system, p. \pageref{Phieps}

$\Phi$,  a solution to the limit system~$(SW_0) $,  p. \pageref{Phi}

 $\Phi_\lambda$,  an element of $ \Ker(L-i\lambda Id)$,  p. \pageref{Philambda} 

$\Phi_N$,  an approximation of~$ \Phi$, p. \pageref{PhiN}

$\Phi_{(N)}$,  an approximate Leray solution of~$(SW_{0})$, p. \pageref{Phi(N))}

$\Phi_{\eps,N}$,  an approximation of~$ \Phi_\eps$, p. \pageref{PhiepsN}

$\Phi_{\eps,\lambda, N}$, the projection of~$\Phi_{\eps,N}$ onto~$
 \Ker(L-i\lambda Id)$, p. \pageref{PhiepslambdaN}

 $\underline \Phi$,  a solution to the geostrophic equation,
p.
\pageref{bvarphi-not}

$\underline \varphi_n$, coefficients of $\underline \Phi\in \Ker L$ in the
$(\Psi_{n,0,0})$ basis,
p.
\pageref{bvarphi-not}

$\phi_N$, a corrector to~$\Phi_N $, p. \pageref{defpetitPhiN}

 $\varphi_{n,k,j}$, coefficients of $\Phi$ in the $(\Psi_{n,k,j})$ basis, p.
\pageref{varphi-not}

 \bigskip

 $\psi_n$, Hermite functions, p. \pageref{hermite-not}

 $\Psi_{n,k,j}$, eigenvectors of $L$,
 p. \pageref{psi1-not}--\pageref{psi3-not}


\printindex

\begin{thebibliography}{99}

\bibitem{aubin} J.-P. Aubin, Un th{\'e}or{\`e}me de compacit{\'e}, {\it
    Notes aux Comptes--Rendus de l'Acad{\'e}mie des Sciences de Paris},
    {\bf 309}  (1963), pages 5042--5044.

\bibitem{breschdesjardins} D. Bresch \& B. Desjardins, Existence of
global weak solutions for a 2D viscous shallow water equation and
convergence to the quasi-geostrophic model. {\it Commun. Math. Phys.}
{\bf 238}  (2003), pages  211--223 .

\bibitem{bresch-desjardinslin} D. Bresch, B. Desjardins \& C.K. Lin, On
some compressible fluid models~: Korteweg, lubrication and shallow water
systems. {\it Comm. Partial Diff. Eqs},  {\bf 28} (2003), pages 843--868.

\bibitem{cdggbook} J.-Y. Chemin, B. Desjardins, I. Gallagher and
E. Grenier, {\it Basics of Mathematical Geophysics, 
An introduction to rotating fluids 
and the Navier-Stokes equations}, to appear in  {\it Oxford
University Press}, 2006.

\bibitem{danchinperiodic}  R. Danchin,   Zero Mach number limit for compressible flows with  periodic boundary
conditions. {\it Amer. J. Math.} {\bf 124} (2002), no. 6,  pages 1153--1219.

\bibitem{desjardinsgrenier} B. Desjardins  \&  E Grenier,
 On the homogeneous model of wind-driven ocean circulation. {\it SIAM J. 
Appl. Math.} {\bf 60} (2000), no. 1, 43--60

\bibitem{dutrifoymajda} A. Dutrifoy  \&  A. Majda,  
The dynamics of equatorial long waves: a singular limit with fast variable coefficients, 
{\it Commun. Math. Sci.} {\bf 4} (2006), no. 2, 375--397. 

\bibitem{gallagherkyoto} I. Gallagher,  A Remark on smooth solutions  of the weakly
compressible Navier--Stokes
equations,
{\it  Journal of Mathematics of Kyoto University}, {\bf 40} (2000),  pages 525--540.

\bibitem{gallagher/sr2}  I. Gallagher \& L. Saint-Raymond, On pressureless gases
driven by a strong inhomogeneous magnetic field,  {\it SIAM J. Math. Analysis},
 {\bf 36} (2005), no. 4, pages 1159--1176.

\bibitem{gallagher/sr}  I. Gallagher \& L. Saint-Raymond, Weak convergence
results for inhomogeneous rotating fluid equations,  {\it to appear in Journal
d'Analyse Math{\'e}matique}, 2006.

\bibitem{gallagher/sr3}  I. Gallagher \& L. Saint-Raymond,   {  On the influence of the Earth's 
rotation on geophysical flows}, {\it to appear in Handbook of Mathematical Fluid Dynamics}, Elsevier, 
 Susan Friedlander and Denis Serre editors, 2006.


\bibitem{dgv} D. G{\'e}rard-Varet, Highly rotating fluids in rough
  domains, {\it  Journal de Math{\'e}matiques Pures
et Appliqu{\'e}es}  {\bf 82}  (2003),   
  pages 1453--1498.


\bibitem{GP} J.-F. Gerbeau  \&  B. Perthame, Derivation of viscous Saint-Venant system for laminar shallow water; 
numerical validation, {\it Discrete Contin. Dyn. Syst. } Ser.    B 1  (2001), no. 1, 89--102.

\bibitem{gill} A. E. Gill, {\it Atmosphere-Ocean Dynamics}, { International
Geophysics Series}, {\bf Vol. 30}, 1982.

\bibitem{gill-longuet} A. E. Gill \& M. S. Longuet-Higgins, Resonant
interactions between planetary waves, {\it Proc. Roy. Soc. London}, {\bf
A 299} (1967), pages~120--140.

\bibitem{Greenspan} H.P. Greenspan, {\it The theory of rotating
fluids}, {  Cambridge monographs on mechanics and applied mathematics},
$1969$.

\bibitem{grenier} E. Grenier, Pseudo-differential energy estimates of
singular perturbations. {\it Comm. Pure Appl. Math.}, {\bf 50} (1997),
no. 9, pages 821--865.

\bibitem{hormander} L. H\"ormander, {\it The Analysis of Linear Partial Differential 
Equations} Vol. III, Grundlehren der mathematischen  Wissenschaften {\bf 274}, Springer Verlag, 1985. 


\bibitem{klein/majda} R. Klein \& A. Majda, Systematic multi-scale models for the
tropics. {\it Journal of Atmospheric Sciences},  {\bf 60} (2003), pages 173--196.

\bibitem{lebedev} N. Lebedev, {\it Special functions and their applications}, Dover Publications, 
Inc., New York, 1972. 
 
 
\bibitem{liebloss} E. Lieb \& M. Loss, {\it Analysis}, Graduate Studies in Mathematics {\bf 14}, American 
Mathematical Society, 2001.


\bibitem{lions/temam/wang} J.-L. Lions, R. Temam \& S. Wang, New formulations of the
primitive equations of atmosphere and applications. {\it Nonlinearity} {\bf 5} (1992),
pages 237-288.

\bibitem{pll} P.-L. Lions,  {\it Mathematical Topics in Fluid Mechanics}, Vol.
II, Compressible Models, Oxford Science Publications,  1997.

\bibitem{majda} A. Majda, {\it Introduction to PDEs and Waves for the Atmosphere and Ocean}, Courant Lecture
Notes {\bf 9}, American Mathematical Society, 2003.

\bibitem{masmoudi} N. Masmoudi, Some asymptotic
problems in fluid mechanics.
Evolution equations and their applications in physical and life
sciences (Bad Herrenalb, 1998), pages  395--404, {\it
Lecture Notes in Pure and Appl. Math.}, {\bf 215} (2001), Dekker, New York.

\bibitem{mellet-vasseur} A. Mellet  \&  A. Vasseur, On the isentropic compressible Navier-Stokes equation, 
{\it Preprint}, 2005.  

\bibitem{MOS} W. Magnus, F. Oberhettinger \& R. Soni, {\it Formulas and
Theorems for the Special Functions of Mathematical Physics}, { Springer
Verlag}, 1966.

\bibitem{pedlosky}
J. Pedlosky, {\it Geophysical fluid dynamics}, { Springer},
$1979$.

\bibitem{philander} S. Philander, {\it  El Ni\~no, la Ni\~na, and the Southern
Oscillation}, { Academic Press}, 1990.  


\bibitem{ripa1} P. Ripa, Nonlinear wave-wave interactions in a one-layer reduced-gravity 
model on the equatorial~$\beta$ plane, {\it J. Phys. Oceanogr.}, {\bf 12}(1)  (1982), 97-111. 

\bibitem{ripa2} P. Ripa,   Weak interactions of equatorial waves in a one-layer model. Part 
I: General properties.{\it J. Phys. Oceanogr.}, {\bf 13}(7) (1983), 1208-1226. 


\bibitem{ripa3} P. Ripa,  Weak interactions of equatorial waves in a one-layer model. Part 
II: Applications. {\it J. Phys. Oceanogr.}, {\bf 13}(7), (1983), 1227-1240. 



\bibitem{schochet} S. Schochet, Fast singular limits of hyperbolic PDEs. {\it J. Diff. Equ.} {\bf
114}  (1994), pages 476--512.

\bibitem{TZ} R. Temam \& M. Ziane, Some mathematical Problems in
Geophysical Fluid Dynamics, {\sl Handbook of Mathematical Fluid
Dynamics, vol. III, eds S. Friedlander \& D. Serre}, 2004.



\bibitem{thomson} W. Thomson (Lord Kelvin), On gravitational oscillations
of rotating water, {\it Proc. Roy. Soc. Edinburgh} {\bf 10}  (1879),  pages
92-100.
\end{thebibliography}
\end{document}